\documentclass[12pt, reqno]{amsart}
\usepackage{mathrsfs}
\usepackage{amssymb,amsthm,amsmath}
\usepackage[numbers,sort&compress]{natbib}
\usepackage{amssymb,amsmath}
\usepackage{amsfonts}
\usepackage{mathrsfs}
\usepackage{latexsym}
\usepackage{amssymb}
\usepackage{amsthm}
\usepackage{color}
\usepackage{pdfsync}
\usepackage{indentfirst}
\date{today}

\usepackage{amsmath}
\usepackage{amsthm}

\newtheorem{theorem}{Theorem}
\newcommand{\R}{\mathbb R}

\newtheorem{remark}[theorem]{Remark}
\newtheorem{proposition}[theorem]{Proposition}
\newtheorem{lemma}[theorem]{Lemma}
\newtheorem{definition}[theorem]{Definition}

\newcommand{\beq}{\begin{equation}}
\newcommand{\eeq}{\end{equation}}
\newcommand{\ben}{\begin{eqnarray}}
\newcommand{\een}{\end{eqnarray}}
\newcommand{\beno}{\begin{eqnarray*}}
\newcommand{\eeno}{\end{eqnarray*}}

\numberwithin{equation}{section}

\begin{document}
\title[Suitable weak solution for the 3D chemotaxis-Navier-Stokes equations]{\bf  Global existence of suitable weak solutions to the 3D chemotaxis-Navier-Stokes equations}
\author{Xiaomeng~Chen}
\address[Xiaomeng~Chen]{School of Mathematical Sciences, Dalian University of Technology, Dalian, 116024,  China}
\email{cxm2381033@163.com}

\author{Shuai~Li}
\address[Shuai~Li]{School of Mathematical Sciences, Dalian University of Technology, Dalian, 116024,  China}
\email{leeshy@mail.dlut.edu.cn}

\author{Lili~Wang}
\address[Lili~Wang]{School of Mathematical Sciences, Dalian University of Technology, Dalian, 116024,  China}
\email{wanglili\_@mail.dlut.edu.cn}

\author{Wendong~Wang}
\address[Wendong~Wang]{School of Mathematical Sciences, Dalian University of Technology, Dalian, 116024,  China}
\email{wendong@dlut.edu.cn}
\date{\today}
\maketitle

\begin{abstract}
In 2004, Dombrowski et al. showed that suspensions of aerobic bacteria often develop flows from the interplay of chemotaxis and buoyancy, which is described as the chemotaxis-Navier-Stokes model, and they observed self-concentration occurs as a turbulence by exhibiting transient, reconstituting, high-speed jets, which entrains nearby fluid to produce paired, oppositely signed vortices. In order to investigate the properties  of these vortices (singular points), one approach is to follow the partial regularity theory of Caffarelli-Kohn-Nirenberg to study the singularity properties of suitable weak solutions. 
In this paper, we established the existence of suitable weak solution for the three dimensional  chemotaxis-Navier-Stokes equations,  where the main difficulty is to establish appropriate local energy inequalities.

%The new ingredients  are to establish certain type of local energy inequality and deal with the non-scaling invariant quantity of $n\ln n$, where $n$ represents the cell concentration, which seems to be the first description for the singular set of weak solutions  of the model.
%In this note, we investigate the existence of suitable weak solution for the three dimensional  chemotaxis-Navier-Stokes equations.
%and then achieve three results about internal regularity criterions. We investigate the Hausdorff dimension of these vortices (singular points) by considering partial regularity of weak solutions and obtain the $\frac53$-dimensional Hausdorff measure of the possible singular set is vanishing  at the first blow-up time. The new ingredients are to establish certain type of local energy inequality and deal with the non-scaling invariant quantity of $n\ln n$.
\end{abstract}

\section{Introduction}	
 We consider the following model describing the dynamics of oxygen,
swimming bacteria, and viscous incompressible fluids, which was proposed by Tuval et al.\cite{ILCJR 2005} in $Q_T=\mathbb{R}^3\times(0,T)$:
\begin{eqnarray}\label{eq:GKS}
 \left\{
    \begin{array}{llll}
    \displaystyle \partial_t n+u\cdot \nabla n-\Delta n=-\nabla\cdot(\chi(c)n\nabla c),\\
   % -\Delta u + u \cdot \nabla u + \nabla p = B \cdot \nabla B, ~~ \rm{in} ~~ \mathbb{R}^3, \\
    \displaystyle \partial_t c+ u\cdot \nabla c-  \Delta c=-\kappa(c)n, \\
    \displaystyle \partial_t u+ \mu u\cdot \nabla u-\Delta u+\nabla P =-n\nabla \phi,\\
    \displaystyle \nabla\cdot u=0,
    \end{array}
 \right.
\end{eqnarray}
where  $c(x,t):Q_T\rightarrow{\mathbb R}^{+}$, $n(x,t):Q_T\rightarrow\mathbb{R}^{+}$, $u(x,t):Q_T\rightarrow\mathbb{R}^{3}$ and $P(x,t):Q_T\rightarrow\mathbb{R}$ denote the oxygen concentration, cell concentration, the fluid velocity and the associated pressure, respectively. Moreover, the
gravitational potential $\phi$, the chemotactic sensitivity $\chi(c)\geq 0$ and the per-capita oxygen
consumption rate $\kappa(c)\geq 0$ are supposed to be sufficiently smooth given functions. Here $\mu=1$ or $0$, and the system of (\ref{eq:GKS}) is reduced to the chemotaxis-Stokes model for $\mu=0.$

As Dombrowski et al.  observed in \cite{DCCGK2004} (see also \cite{ILCJR 2005}):
local concentration for  suspensions of aerobic bacteria leads to a jet descending faster
than its surroundings, which entrains nearby fluid to
produce paired, oppositely signed vortices. Thus this system has very strong physical background. It is necessary to study the existence of solutions and the singularity properties of weak solutions.
%In order to investigate the properties  of these vortices (singular points), we follow the partial regularity theory  of Caffarelli-Kohn-Nirenberg in \cite{CKN} to study the singularity properties of suitable weak solutions. 
Recently, Chen-Li-Wang \cite{CLW 2022} studied partial regularity of strong solutions to the simplified 3D chemotaxis-Navier-Stokes equations ($\kappa(c)=c$, $\chi(c)=1$ ) at the first blow-up time. However, it's still unknown whether there exists a suitable weak solution for such system. Here we consider the existence of general  chemotaxis-Navier-Stokes equations, and discuss  partial regularity of these weak solutions in a forthcoming paper. About suitable weak solutions of Navier-Stokes equations, it was started by Scheffer in \cite{SV1,SV2,SV4}, and later Caffarelli-Kohn-Nirenberg \cite{CKN} showed that the set $\mathcal{S}$ of possible interior singular points of a suitable weak solution is one-dimensional parabolic Hausdorff measure zero.
The suitable weak solution is better than Leray-Hopf weak solution introduced by Leray in
\cite{Leray} and  if the local strong solution blows up, then the solution may be continued as a
suitable weak solution (see Proposition 30.1 in \cite{LR}).
More references on simplified proofs and improvements, we refer to Lin \cite{Lin}, Ladyzhenskaya-Seregin \cite{LS}, Tian-Xin \cite{TX}, Seregin \cite{Se}, Gustafson-Kang-Tsai \cite{GKT}, Vasseur \cite{Va} and the references therein.
%Here we consider the partial regularity of the system (\ref{eq:KS}) at the first blow-up time as Dong-Du in \cite{DD}, since the first question of (Q1) is still unknown, which is an open question.

Due to the significance of the biological background (see \cite{DCCGK2004}, \cite{ILCJR 2005}), the model could be used to predict the
large-scale bioconvection affecting clearly the overall oxygen consumption in the
above experiments (see, for example, \cite{AL2010}). Many mathematicians have studied this model and made much progress, 
%such as the existence of weak solutions, the chemotaxis-Navier-Stokes system with a nonlinear diffusion,  blow-up criteria, stability and so on.
and we just mention some related works for the result.
Firstly, for the existence of weak solutions, global classical solutions near constant steady states are constructed for the full
chemotaxis-Navier-Stokes system by Duan-Lorz-Markowich in \cite{DLM2010}. In \cite{AL2010}, for the case of bounded domain of $\mathbb{R}^n$ with $n= 2,3$, the local existence of weak solutions for problem (\ref{eq:GKS}) is obtained by Lorz. Later, Winkler proved the existence of global weak solution in \cite{Winkler2012} by assuming that
\beno
\left(\frac{\kappa}{\chi}\right)'>0,\quad \left(\frac{\kappa}{\chi}\right)''\leq 0,\quad  \left({\kappa}{\chi}\right)'\geq0.
\eeno
By assuming $\chi',\kappa'\geq 0$ and $\kappa(0)=0$, local well-posed results and blow-up criteria were established by Chae-Kang-Lee in \cite{CKL2013}. In \cite{CKL2014}, under the conditions of $\chi,\kappa\in C^m,\kappa(0)=0$, $\chi,\kappa,\chi^{'},\kappa^{'}\geq 0$ and the smallness of $\|c_0\|_{L^\infty}$, Chae-Kang-Lee proved that classical solution of (\ref{eq:GKS}) exists globally. 
Recently, Winkler proved the global existence of weak solutions of the system (\ref{eq:GKS}) in bounded domain with large initial data, and obtained much better a priori estimates such as $\frac{|\nabla c|^4}{c^3}\in L^1$ in \cite{Winkler2016}. Secondly, for the two-dimensional system  of (\ref{eq:GKS}), the system is better understood. Liu and Lorz \cite{LL2011}  proved the global existence of weak solutions to the two-dimensional system  of (\ref{eq:GKS}) for arbitrarily large initial data, under the assumptions on $\chi$ and $f$ made in \cite{DLM2010}. See also \cite{LL2016,DLX2017,WWX2018,WWX2018-2,He2020,WZZ2021} and the references therein.
%For more references  about the existence of solutions, we refer to \cite{CL2016,Winkler2017,KM2019,DL2020,Bl2020,DL2022} and the references therein.
As for the case of the chemotaxis-Navier-Stokes system with a nonlinear diffusion, that means $\Delta n$ is replaced by $\Delta n^m$, there are also many works. one can refer to \cite{FLM2010,TW2012, TW2013,CKK2014,ZK2021} and so on.

As Winkler said in \cite{Winkler2016},
{ ``For the full three-dimensional chemotaxis-Navier-Stokes system, even at the very basic level of global existence in generalized solution frameworks, a satisfactory solution theory is entirely lacking.''} 
%In this paper our aim is to explore partial regularity properties of weak solutions.
For the simplified model of  the case $\kappa(c)=c$ and $\chi(c)=1$, the three dimensional chemotaxis-Navier-Stokes system (\ref{eq:GKS}) is reduced to
\begin{eqnarray}\label{eq:KS}
 \left\{
    \begin{array}{llll}
    \displaystyle \partial_t n+u\cdot \nabla n-\Delta n=-\nabla\cdot(n\nabla c),\\
   % -\Delta u + u \cdot \nabla u + \nabla p = B \cdot \nabla B, ~~ \rm{in} ~~ \mathbb{R}^3, \\
    \displaystyle \partial_t c+ u\cdot \nabla c- \Delta c=-cn, \\
    \displaystyle \partial_t u+u\cdot \nabla u-\Delta u+\nabla p =-n\nabla \phi,~~\nabla\cdot u=0. \\
    %\displaystyle n(x,0)=n_0(x),~~c(x,0)=c_0(x),~~u(x,0)=u_0(x).
    \end{array}
 \right.
\end{eqnarray}

Furthermore, global weak solutions of Leray-Hopf type for this system was obtained in 2D and 3D by Zhang-Zheng \cite{ZZ2014}, He-Zhang  \cite{HZ2017}, and   Kang-Lee-Winkler \cite{KLW2022}, respectively, where they established a priori estimate
\ben\label{eq: a priori-zhang}
\mathcal{U}(t) + \int_0^t \mathcal{V}(t) d\tau \leq C e^{Ct},
\een
where
\beno
\mathcal{U} = \|n\|_{L^1 \cap L \log L} + \|\nabla \sqrt c\|_{L^{2}}^2 + \|u\|_{L^{2}}^2,
\eeno
and
\beno
\mathcal{V} = \|\nabla \sqrt{n+1}\|_{L^{2}}^2 + \|\Delta \sqrt c\|_{L^{2}}^2 + \|\nabla u\|_{L^{2}}^2 + \int_{\mathbb{R}^d} (\sqrt{c})^{-2} |\nabla \sqrt c|^4 dx + \int_{\mathbb{R}^d} n |\nabla \sqrt c|^2 dx,
\eeno
where $d=2,3$ and denote by $ L^{p}(\mathbb{R}^{3})=L^{p} $ for simplicity.
However, up to now more information about these weak solutions is still not known,  especially for the interior singular vortices as described in \cite{DCCGK2004} or the self-organized generation of a persistent hydrodynamic vortex that traps cells near the contact line(see \cite{ILCJR 2005}).

Our main objective of this paper is to present global existence of suitable weak soultion for (\ref{eq:GKS}) in three dimensions.
Assume that
\ben\label{ine:chi}
%\chi(s) \in C^0(\overline{\mathbb{R}^+}) \cap C^2(\mathbb{R^{+}}), ~~~~\chi'(s) \geq  0, ~~~~ \chi''(s) \geq  0,
\chi(s) \in 
%C^1(\overline{\mathbb{R}^+}) \cap
 C^2(\overline{\mathbb{R^{+}}});~~~~~\chi(s)\geq 0;
 % ~~~~ \chi'(0) >  0,~~~~\chi'(s) \geq  0;
% ~~~~ \chi''(s) \geq  0,~~~~ \chi'''(s) \geq  0;
\een
\ben\label{ine:kappa}
\kappa(s) \in C^2(\overline{\mathbb{R^{+}}})
%\cap C^4(\mathbb{R^{+}})
,~~~~\kappa(0)=0,~~~~ \kappa'(s)\geq 0,~~~~\kappa''(s)\geq 0;
%\kappa^{(m)}(s)\geq 0~~{\rm{for}}~~m=1,2,3,4;
%~~~~\kappa'(s) \geq 0, ~~~~\kappa''(s) \geq 0,~~~~\kappa'''(s) \geq 0,~~~~\kappa^{(4)}(s) \geq 0,
\een
and
\ben\label{equ:chi kappa}
\kappa(s)=\Theta_0s\chi(s),
\een
where $ \Theta_0 >0$ is a positive absolute constant.

%Here the definition of suitable weak solutions is as follows.
\begin{definition}\label{sws}
	A triplet $(n, c, u)$ is called a suitable weak solution of the system (\ref{eq:GKS}) if the following holds:\\
	
	(i). For any bounded domain $\Omega \subset \mathbb{R}^3$, \\
	
	$n, n \ln n \in L^\infty_{\rm loc}([0,+\infty);L^1(\Omega))$, $\nabla \sqrt{n} \in L^2_{\rm loc}([0,+\infty);L^2(\Omega))$, \\
	
	$\nabla \sqrt{c} \in L^\infty_{\rm loc}([0,+\infty);L^2(\Omega)) \cap L^2_{\rm loc}([0,+\infty);H^1(\Omega))$, \\
	
	$u\in L_{\rm loc}^\infty([0,+\infty);L^2(\Omega)) \cap L^2_{\rm loc}([0,+\infty);H^1(\Omega))$, $P \in L^\frac32_{\rm loc}([0,+\infty);L^\frac32(\Omega))$; \\
	
	(ii). $(n, c, u)$ solves (\ref{eq:GKS}) %in $\mathbb{R}^3\times\mathbb{R}^{+}$ 
	in the sense of distributions;\\
	
	(iii). $(n, c, u)$ satisfies energy inequality
	\begin{equation*}
		\begin{aligned}
		&||u||_{L^{2}}^{2} + \int_0^t||\nabla u(t)||_{L^{2}}^{2}\\&+~\int_{\mathbb{R}^3} (n+1) \ln (n+1)(\cdot,t) + \int_0^t\int_{\mathbb{R}^3} |\nabla \sqrt{n+1}|^2\\
		&+~ \frac2{\Theta_0}||\nabla \sqrt{c}||_{L^{2}}^{2} +\frac4{3\Theta_0}\int_0^t||\nabla^2 \sqrt{c}||_{L^{2}}^{2} + \frac1{3\Theta_0} \int_0^t \int_{\mathbb{R}^3}(\sqrt{c})^{-2} |\nabla \sqrt{c}|^4  \\
		\leq ~&
		C(\|\nabla \phi\|_{L^\infty},\|n_0\|_{L^1},\|c_{0}\|_{L^{\infty}},\|c_{0}\|_{L^{1}},\|u_{0}\|_{L^{2}},\|(n_{0}+1)\ln (n_{0}+1)\|_{L^{1}},\|\nabla \sqrt{c_0}\|_{L^{2}})(1+t);
		\end{aligned}
	\end{equation*}
	(iv). $(n, c, u)$ satisfies local energy inequality
\begin{equation*}
\begin{aligned}
&\int_{\Omega} (n \ln (n) \psi)(\cdot,t) + 4 \int_{(0,t)\times\Omega} |\nabla \sqrt{n}|^2 \psi+\frac{2}{\Theta_0}  \int_{\Omega} (|\nabla \sqrt{c}|^2 \psi)(\cdot,t)\\&+\frac{4}{3\Theta_0}\int_{(0,t)\times\Omega} |\Delta \sqrt{c}|^2 \psi
+\frac{18}{\Theta_0}\|{c_{0}}\|_{L^{\infty}}\int_{\Omega}(|u|^2)(\cdot,t) \psi\\
& + \frac{18}{\Theta_0}\|{c_{0}}\|_{L^{\infty}}\int_{(0,t)\times\Omega} |\nabla u|^2 \psi+\frac{2}{3\Theta_0}  \int_{(0,t)\times\Omega} (\sqrt{c})^{-2} |\nabla \sqrt{c}|^4 \psi
\\\leq~&\int_{(0,t)\times\Omega} n \ln (n) (\partial_t \psi + \Delta \psi) + \int_{(0,t)\times\Omega} n \ln (n) u \cdot \nabla \psi \\&+\int_{(0,t)\times\Omega}n\chi(c)\nabla c\cdot \nabla \psi +\int_{(0,t)\times \Omega}n\ln n \chi(c)\nabla c\cdot\nabla\psi\\&+ \frac{2}{\Theta_0}\int_{(0,t)\times\Omega} |\nabla \sqrt{c}|^2 (\partial_t \psi + \Delta \psi) + \frac{2}{\Theta_0} \int_{(0,t)\times\Omega} |\nabla \sqrt{c}|^2 u \cdot \nabla \psi\\&
+\frac{18}{\Theta_0}||c_{0}||_{L^{\infty}}\int_{(0,t)\times\Omega} |u|^2 \left(\partial_t \psi +\Delta \psi\right) + \frac{18\mu}{\Theta_0}||c_{0}||_{L^{\infty}}\int_{(0,t)\times\Omega} |u|^2 u\cdot\nabla\psi
\\&+\frac{36}{\Theta_0}||c_{0}||_{L^{\infty}}\int_{(0,t)\times\Omega} (P - \bar{P}) u \cdot \nabla \psi - \frac{36}{\Theta_0}||c_{0}||_{L^{\infty}}\int_{(0,t)\times\Omega} n\nabla\phi \cdot u \psi,
\end{aligned}
\end{equation*}
	where  $\psi\geq 0$ and vanishes in the parabolic boundary of $(0,t)\times\Omega.$
\end{definition}

%Firstly, we give the definition of $L\log L$ norm, which will be used in the following time.
%\begin{definition}\label{def:log}
%	The Zygmund classes with $A(t) = t\log^{+} t$, is defined as the set all functions $f$ such that
%	\beno
%	\int_{\mathbb{R}^3}A(|f(x)|)dx<\infty.
%	\eeno
%	The corresponding Zygmund space $L\log L(\mathbb{R}^3)$ is defined as the linear hull of the Zygmund class,
%	which is equipped with the Luxemburg norm
%	\beno
%	||f||_{L\log L}=\inf\left\{k \big|\int_{\mathbb{R}^3}A(\frac fk)dx\leq 1\right\},
%	\eeno
%	%{\rm in}
%	and
%	\begin{eqnarray}\label{eq:initial}
%	\log^{+}t=\left\{
%	\begin{array}{llll}
%	\displaystyle \log t,~~t\geq 1 ,\\
%	% -\Delta u + u \cdot \nabla u + \nabla p = B \cdot \nabla B, ~~ \rm{in} ~~ \mathbb{R}^3, \\
%	\displaystyle 0,~~{\rm otherwise}.\\
%	%\displaystyle \partial_t u+u\cdot \nabla u-\triangle u+\nabla p =-n\nabla \phi,~~\nabla\cdot u=0 \\
%	\end{array}
%	\right.
%	\end{eqnarray}
%\end{definition}

Our main results are as follows.\\
%{\bf The existence of suitable weak solutions of system $\eqref{eq:GKS}$:}\\
\begin{theorem}\label{suitable weak solution}
	Assume that the initial data  $(n_0,c_0,u_0)$ satisfies
	\ben\label{ine:initial condition assumption 1}
	\left\{
	\begin{array}{llll}
		\displaystyle n_0 \in L^1(\mathbb{R}^3), \quad (n_0+1) \ln (n_0+1)  \in L^1(\mathbb{R}^3),  \quad u_0 \in L^2_\sigma(\mathbb{R}^3);\\
		\displaystyle \nabla \sqrt{c_0} \in L^2(\mathbb{R}^3), \quad c_0 \in L^1 \cap L^\infty(\mathbb{R}^3);\\
		\displaystyle n_0 \geq 0, \quad c_0 \geq 0.\\
	\end{array}
	\right.
	\een
	Moreover, $\nabla \phi\in L^\infty(\mathbb{R}^3)$, $\kappa$ and $\chi$ satisfy $\eqref{ine:chi}$, $\eqref{ine:kappa}$ and $ \eqref{equ:chi kappa}. $
	Then there exists a global suitable weak solution of the system (\ref{eq:GKS}). Moreover, for the chemotaxis-Stokes model of $\mu=0$, assume that $(-\triangle)^{\frac14} u_0 \in L^2(\mathbb{R}^3)$ additionally, and one can get a new a priori estimate of the velocity:
 \ben\label{eq: good esimate of u}
 &&\int_{\mathbb{R}^3}|(-\Delta)^{\frac14} u(\cdot,t)|^2dx+\int_0^t\int_{\mathbb{R}^3}|\nabla(-\Delta)^{\frac14} u|^2dxds\\\nonumber
 &\leq& C(\|\nabla \phi\|_{L^\infty},\|n_0\|_{L^1},\|c_{0}\|_{L^{\infty}\cap L^{1}},\|u_{0}\|_{L^{2}},\|(n_{0}+1)\ln (n_{0}+1)\|_{L^{1}},\|\nabla \sqrt{c_0}\|_{L^{2}})(1+t).
 \een
\end{theorem}

\begin{remark}\label{rem:1}
(i) Note that the right hand side  term of $(\ref{eq:GKS})_1$ is supper-critical, which is cancelled by the right hand term of $c$'s equation.   When $\eqref{equ:chi kappa}$ fails, it is still unknown whether the system of (\ref{eq:GKS}) has a suitable weak solution or even a Leray-Hopf weak solution.

%The existence of suitable weak solution For some equations, it is not easy to study the existence of appropriate solutions or even Leray solutions, such as liquid crystal equations and Keller Segal equations.

(ii) For $\mu=1$, Chae-Kang-Lee in \cite{CKL2013} constructed a global weak solution by assuming $\chi=\mu_0 \kappa$ and using the iterative scheme developed in \cite{DLM2010}, where they obtained the estimates of $\|\nabla c\|_{L^\infty_tL^2_x(\mathbb{R}^3)}$ and $\|\nabla^2 c\|_{L^2((0,t)\times \mathbb{R}^3)}$. Here we introduce a mollifier dealing with the nonlinear terms of $n,u$, and obtain the estimate of $\|({\sqrt{c}})^{-1}|\nabla \sqrt{c}|^2\|_{L^2((0,t)\times \mathbb{R}^3)}$ by assuming \eqref{equ:chi kappa} for the equation of $\sqrt{c}$, which seems to be important in the local energy inequality.
Compared with global weak solutions of the simplified model in \cite{HZ2017} by He-Zhang and \cite{KLW2022} by Kang-Lee-Winkler, we consider the more general model, and obtained local energy inequality by proving the strong convergence of $n\ln n$ in $L^{\frac32}((0,t)\times \mathbb{R}^3))$. Moreover, the uniform a priori energy estimates is linear with respect to  the time $t$, which is slower than the exponential growth (\ref{eq: a priori-zhang})
in \cite{HZ2017} 
or \cite{KLW2022}.

(iii) For $\mu=0$, we get a better a priori estimate of the velocity in (\ref{eq: good esimate of u}), and it's still unknown whether one can get higher regularity estimates for $n$ and $c$. If $\chi(c)$ is replaced by $S(n,c)\geq 0$ and $S(n,c)\leq C(1+n)^{-\alpha}$ with $\alpha>0$, referring to the known results \cite{WC2015}, the solution is globally bounded.

%
%The difference between our result and \cite{Winkler2012} is that we include the $u\cdot \nabla u$ for the fluid equation, while Winkler proved the linearized system. Besides, we have different assumptions for the connection about $\kappa$ and $\chi$. The reason is that we establish the uniform estimate about $\sqrt{c}$, while the author estimate about $c$.

%(ii) The new observation of this theorem is the local energy inequality of (\ref{local energy}) (See Lemma \ref{local energy inequality}), which is indeed a local a priori estimate for weak solutions, which is of independent interest.
%
%(iii) The difficulty mainly lies in dealing with the term including $\ln n$, which is not scaling invariant under the embedding inequality.
% We establish the local a priori estimates of $\int_{B_1} (n \ln n \psi)(\cdot,t)$ firstly by estimating the local energy inequality, then use the equation of $n$ by
%  estimating the term of $\int_{B_1} (n \psi)(\cdot,t)$, which and the embedding inequality in a fixed sphere imply the estimate of $\int_{B_1} n |\ln n| $.
\end{remark}

The proof of Theorem \ref{suitable weak solution} is based on the following global well-posedness result of the following regularized system:
\begin{eqnarray}\label{eq:CNS}
\left\{
\begin{array}{llll}
\displaystyle \partial_t n^{\varepsilon,\tau} + (u^{\varepsilon,\tau} \ast \rho^\varepsilon) \cdot \nabla n^{\varepsilon,\tau} - \Delta n^{\varepsilon,\tau} = - \nabla \cdot \left(\frac{n^{\varepsilon,\tau}}{1+\tau n^{\varepsilon,\tau}} \nabla c^{\varepsilon,\tau} \chi(c^{\varepsilon,\tau})\right),\\
\displaystyle \partial_t c^{\varepsilon,\tau} + u^{\varepsilon,\tau}  \cdot \nabla c^{\varepsilon,\tau} - \Delta c^{\varepsilon,\tau} = - \frac1{\tau}\ln (1+\tau n^{\varepsilon,\tau})\kappa(c^{\varepsilon,\tau}),\\
\displaystyle \partial_t u^{\varepsilon,\tau} + \mu (u^{\varepsilon,\tau} \ast \rho^\varepsilon) \cdot \nabla u^{\varepsilon,\tau} - \Delta u^{\varepsilon,\tau} + \nabla P^{\varepsilon,\tau} = - (n^{\varepsilon,\tau} \nabla \phi) \ast \rho^\varepsilon,\\
\displaystyle \nabla \cdot u^{\varepsilon,\tau} = 0,\\
\displaystyle n^{\varepsilon}(x,0) = n^{\varepsilon}_0(x), \quad c^{\varepsilon}(x,0) = c^{\varepsilon}_0(x), \quad u^{\varepsilon}(x,0) = u^{\varepsilon}_0(x),\\
\end{array}
\right.
\end{eqnarray}
where $\rho^\varepsilon$ is a standard mollifier, $ \tau>0 $ is a sufficiently small constant, and the initial value $n_0^{\varepsilon} $ 
%satisfies
%\ben\label{ine:n L log L}
%n_0^{\varepsilon} \to n_0 \quad {\rm in} \quad L^p \quad {\rm as} \quad \varepsilon \to 0 \quad 1\leq p<\infty,
%\een
,$c_0^{\varepsilon}$ and $u_0^{\varepsilon}$ are defined by $n_0^{\varepsilon}(x) =n_0 \ast \rho^\varepsilon$, $c^{\varepsilon}_0(x) = \left(\sqrt{c_0} \ast \rho^\varepsilon\right)^2$ and $u^{\varepsilon}_0(x) = u_0 \ast \rho^\varepsilon$ respectively. 

The existence of strong solutions is established as follows.
\begin{theorem}\label{lem:2}
	Let $ \nabla \phi\in L^\infty(\mathbb{R}^3) $ and the initial value
	$ (n_{0}^{\varepsilon},c_{0}^{\varepsilon},\nabla c_{0}^{\varepsilon},u_{0}^{\varepsilon})\in H^{2}(\mathbb{R}^{3}),$ $n_{0}^{\varepsilon}, c_{0}^{\varepsilon}\in L^{1}(\mathbb{R}^{3}),$  
	$ \chi, \kappa$ satisfy $\eqref{ine:chi}$ and $\eqref{ine:kappa}$. Assume $ n_{0}^{\varepsilon}(x)\geq 0 $ and $ c_{0}^{\varepsilon}(x)\geq 0 $ for any $ x\in\mathbb{R}^{3}. $ Then there exists a global strong solution $(n^{\varepsilon,\tau},c^{\varepsilon,\tau}, u^{\varepsilon,\tau})$ of system (\ref{eq:CNS}), which satisfies	 
%	\begin{itemize}
%		\item[(i)] 
\\(i)\beno
~~ (n^{\varepsilon,\tau},c^{\varepsilon,\tau},\nabla c^{\varepsilon,\tau},u^{\varepsilon,\tau})\in C([0,T); H^{2}(\mathbb{R}^{3}))\cap L^{2}((0,T); H^{3}(\mathbb{R}^{3}));   
\eeno
and
\beno
&&\|n^{\varepsilon,\tau},c^{\varepsilon,\tau},\nabla c^{\varepsilon,\tau},u^{\varepsilon,\tau}\|_{L^{\infty}_{t}H^{2}_{x}}^{2}+\| n^{\varepsilon,\tau}, c^{\varepsilon,\tau},\nabla c^{\varepsilon,\tau},u^{\varepsilon,\tau} \|_{L^2_{t}H^{2}_{x}}
^{2}\\&\leq&C\left(\varepsilon,\tau,\|\nabla\phi\|_{L^{\infty}}, \|\chi\|_{0},\|\kappa\|_{0},\|n_{0}^{\varepsilon}\|_{L^{1}}, \|E_{0}^{\varepsilon}\|_{H^2}  \right)(T+1)^{54},
%~~{\rm{for}}~~{\rm {any}}~t\in(0,T).
\eeno
for any $ t\in(0,T), $ where
 \beno
\|\chi\|_0=\|\chi\|_{L^\infty(0, \|c_0\|_{L^{\infty}})}+\|\chi'\|_{L^\infty(0, \|c_0\|_{L^{\infty}})}+\|\chi''\|_{L^\infty(0, \|c_0\|_{L^{\infty}})},
\eeno
\beno
\|\kappa\|_0=\|\kappa\|_{L^\infty(0, \|c\|_{L^{\infty}})}+\|\kappa'\|_{L^\infty(0, \|c_0\|_{L^{\infty}})}
+\|\kappa''\|_{L^\infty(0, \|c_0\|_{L^{\infty}})}.
\eeno
		%\item[(ii)] 
(ii) $ \|n^{\varepsilon,\tau}(\cdot,t)\|_{L^{1}}\leq \|n^{\varepsilon}_0\|_{L^{1}},\quad \|c^{\varepsilon,\tau}(\cdot,t)\|_{L^{1}}\leq \|c^{\varepsilon}_0\|_{L^{1}}. $
	%\end{itemize}
\end{theorem}
\begin{remark}
	Indeed, we can consider the following system:
	\begin{eqnarray}\label{eq:CNS n+}
	\left\{
	\begin{array}{llll}
	\displaystyle \partial_t n^{\varepsilon,\tau} + (u^{\varepsilon,\tau} \ast \rho^\varepsilon) \cdot \nabla n^{\varepsilon,\tau} - \Delta n^{\varepsilon,\tau} = - \nabla \cdot \left(\frac{n_{+}^{\varepsilon,\tau}}{1+\tau n_{+}^{\varepsilon,\tau}} \nabla c^{\varepsilon,\tau} \chi(c^{\varepsilon,\tau})\right),\\
	\displaystyle \partial_t c^{\varepsilon,\tau} + u^{\varepsilon,\tau}  \cdot \nabla c^{\varepsilon,\tau} - \Delta c^{\varepsilon,\tau} = - \frac1{\tau}\ln (1+\tau n_{+}^{\varepsilon,\tau})\kappa(c^{\varepsilon,\tau}),\\
	\displaystyle \partial_t u^{\varepsilon,\tau} + (u^{\varepsilon,\tau} \ast \rho^\varepsilon) \cdot \nabla u^{\varepsilon,\tau} - \Delta u^{\varepsilon,\tau} + \nabla P^{\varepsilon,\tau} = - (n^{\varepsilon,\tau} \nabla \phi) \ast \rho^\varepsilon,\\
	\displaystyle \nabla \cdot u^{\varepsilon,\tau} = 0,\\
	\displaystyle n^{\varepsilon}(x,0) = n^{\varepsilon}_0(x), \quad c^{\varepsilon}(x,0) = c^{\varepsilon}_0(x), \quad u^{\varepsilon}(x,0) = u^{\varepsilon}_0(x),\\
	\end{array}
	\right.
	\end{eqnarray}
	where $ n_{+}^{\varepsilon,\tau}=\max\{0, n^{\varepsilon,\tau}\}. $ By Banach fixed theorem there  exists a classical solution of (\ref{eq:CNS n+}) over $ (0,T)$ for a small $T$, and one can prove $ n_{-}^{\varepsilon,\tau}=0 $ over $ (0,T), $ thus the system (\ref{eq:CNS n+}) is reduced to (\ref{eq:CNS}). Hence it suffices to consider the system of (\ref{eq:CNS}).
\end{remark}
\begin{remark}It is worth noting that we did not use a smoothing operator for the nonlinear term in the equation of $c$. One reason is that its regularity can be drived from the equations of $u$ and $n$. Another key reason is that it comes from the elimination estimation of local energy inequalities. That is to say, the bad term  $\int_{(0,t)\times\Omega} |\nabla (u^{\varepsilon,\tau}\ast\rho^\varepsilon)|^2 \psi$ from $(u^{\varepsilon,\tau}\ast\rho^\varepsilon)\cdot\nabla c^{\varepsilon,\tau} $ seems to be not cancelled by the term $\int_{(0,t)\times\Omega} |\nabla u^{\varepsilon,\tau}|^2 \psi$ (see the estimate of $J_2$ and (\ref{ine:energy u})).
\end{remark}
Throughout this article, $C$ denotes an absolute constant independent of $(n,c,u)$ and may be different from line to line.
Denote by $ L^{p}(\mathbb{R}^{3})=L^{p} $ and $ W^{m,p}(\mathbb{R}^{3})=W^{m,p} $ the standard Lebesgue spaces and Sobolev spaces, respectively. In the case $ p=2, $ we denote $ W^{m,2}(\mathbb{R}^{3})=H^{m}(\mathbb{R}^{3})=H^{m}. $ 

The rest of the paper is organized as follows. Section 2 is devoted to obtain the global existence of strong solutions to the regularized system (\ref{eq:CNS}). In Section 3, we prove uniform estimates independent of $\varepsilon$ and $\tau$ for the regularized system (\ref{eq:CNS}). Local energy inequality (\ref{local energy}) is established in Section 4. Besides, we derive convergence of $n^{\varepsilon,\tau}$, $u^{\varepsilon,\tau}$ and $c^{\varepsilon,\tau}$ respectively by Aubin-Lions Lemma in Section 5. Section 6 is to prove 
convergence of local energy inequality when $\tau\rightarrow0$ and $\varepsilon\rightarrow0$, respectively. The proof of Theorem  \ref{suitable weak solution} is in Section 7. Some technical lemmas are shown in the Appendix.

\section{Globally well-posedness theorem to the regularized system (\ref{eq:CNS})}

This section is devoted to proving the global well-posedness of the following regularized system (\ref{eq:CNS}).

%Here, $n^{\varepsilon,\tau}_0(x) = n_0 \ast \rho^\varepsilon$, $c^{\varepsilon,\tau}_0(x) = \left(\sqrt{c_0} \ast \rho^\varepsilon\right)^2$ and $u^{\varepsilon,\tau}_0(x) = u_0 \ast \rho^\varepsilon$.

%The main result is as follows.

%\begin{lemma}\label{lem:3}
%	Letting $E$ is a Banach space, $B$ is bounded bilinear transformation from $E\times E$ to $E$, 
%	\beno
%	\|B(e,f)\|_{E}\leq C_B\|e\|_{E}\|f\|_{E},
%	\eeno 
%	if $0<\delta<\frac1{4C_B}$, $e_0\in E$ and satisfies $\|e_0\|_{E}\leq \delta$, then the equation $e=e_0-B(e,e)$ has a solution with $\|e\|_{E}\leq 2\delta$. The solution is unique in the ball $\bar{B}(0,2\delta)$. Besides, the solution depends on $e_0$ continuously. If $\|f_0\|_{E}\leq \delta$, $f=f_0-B(f,f)$ and $\|f\|_{E}\leq 2\delta$, then 
%	\beno
%	\|e-f\|_{E}\leq\frac1{1-4C_B\delta}\|e_0-f_0\|_{E}.
%	\eeno
%\end{lemma}

%\subsection*{ Global strong solution of (\ref{eq:CNS}).}
%{\bf Step I: Proof of the existence of}

%Assume $ n^{\varepsilon,\tau}\geq 0 $ and $ \tau\in(0,1) $, then clearly
%\begin{equation*}
%	0<\frac{1}{1+\tau n^{\varepsilon,\tau}}\leq 1.
%\end{equation*}
{\bf Proof of Theorem \ref{lem:2}}.
Let  $ R=2\left(\|n_{0}^{\varepsilon}\|_{H^{2}}+\|c_{0}^{\varepsilon}\|_{H^{3}}+\|u_{0}^{\varepsilon}\|_{H^{2}}  \right)<\infty $
and $ T\in (0,1) $ to be decided, and we introduce the Banach space as follows.
$$ X:=L^{\infty}[(0,T);\{L^{\infty}(\mathbb{R}^3)\cap W^{1,2}(\mathbb{R}^{3})\}\times \{\dot{W}^{1,\infty}(\mathbb{R}^{3})\cap W^{2,2}(\mathbb{R}^{3})\}\times W^{2,2}(\mathbb{R}^{3})], $$
along with its closed subset
\ben\label{eq:norm}
S&:=&\{(n^{\varepsilon,\tau},c^{\varepsilon,\tau},u^{\varepsilon,\tau})\in X\big|\sup_{t\in(0, T)}\|n^{\varepsilon,\tau}(\cdot, t)\|_{L^\infty}+\sup_{t\in(0, T)}\|n^{\varepsilon,\tau}(\cdot, t)\|_{W^{1,2}}\nonumber\\&&+\sup_{t\in(0, T)}\|\nabla c^{\varepsilon,\tau}(\cdot,t)\|_{L^{\infty
}}+\sup_{t\in(0, T)}\|c^{\varepsilon,\tau}(\cdot,t)\|_{L^{2}}+\sup_{t\in(0, T)}\|\nabla^{2}c^{\varepsilon,\tau}(\cdot,t)\|_{L^{2}}\nonumber \\&&+\sup_{t\in(0, T)}\|u^{\varepsilon,\tau}(\cdot,t)\|_{L^{2}}+\sup_{t\in(0, T)}\|\nabla^{2}u^{\varepsilon,\tau}(\cdot,t)\|_{L^2}\leq R\}.
\een
For $ (n^{\varepsilon,\tau}, c^{\varepsilon,\tau}, u^{\varepsilon,\tau})\in S $ and $ t\in (0,T), $ we let
\begin{equation*}
\begin{aligned}
&\Phi(n^{\varepsilon,\tau}, c^{\varepsilon,\tau}, u^{\varepsilon,\tau})(\cdot, t)\\=~&\left(\begin{array}{llll}
\Phi_1(n^{\varepsilon,\tau},c^{\varepsilon,\tau},u^{\varepsilon,\tau})(\cdot,t)\\
\Phi_2(n^{\varepsilon,\tau},c^{\varepsilon,\tau},u^{\varepsilon,\tau})(\cdot,t)\\
\Phi_3(n^{\varepsilon,\tau},c^{\varepsilon,\tau},u^{\varepsilon,\tau})(\cdot,t)\\
\end{array}
\right)
\\=~&\left(\begin{array}{llll}
e^{t\Delta}n_0^{\varepsilon,\tau}-\int_0^te^{(t-s)\Delta}\big\{\nabla \cdot \left(\frac{n^{\varepsilon,\tau}}{1+\tau n^{\varepsilon,\tau}} \nabla c^{\varepsilon,\tau} \chi(c^{\varepsilon,\tau})\right)+(u^{\varepsilon,\tau} \ast \rho^\varepsilon) \cdot \nabla n^{\varepsilon,\tau}\big\}(\cdot,s)ds\\
e^{t\Delta}c_0^{\varepsilon,\tau}-\int_0^te^{(t-s)\Delta}\big\{\frac1{\tau}\ln (1+\tau n^{\varepsilon,\tau})\kappa(c^{\varepsilon,\tau})+u^{\varepsilon,\tau}  \cdot \nabla c^{\varepsilon,\tau}\big\}(\cdot,s)ds\\
e^{t\Delta}u_0^{\varepsilon,\tau}-\int_0^te^{(t-s)\Delta}\mathcal{P}\big\{\mu(u^{\varepsilon,\tau} \ast \rho^\varepsilon)\cdot \nabla u^{\varepsilon,\tau}+(n^{\varepsilon,\tau} \nabla \phi) \ast \rho^\varepsilon \big\}(\cdot,s)ds\\
\end{array}
\right),
\end{aligned}
\end{equation*}
where $(e^{t\Delta})_{t\geq 0}$
% $(e^{-tA})_{t\geq 0}$
and $\mathcal{P}$ denote the heat semigroup
%  , the Stokes semigroup,
and
the Leray projection in $L^2$, respectively. The proof is divided into five steps as follows.
%In order to estimate $ \Phi_1,~\Phi_2 $ and $ \Phi_3, $ by Gagliado-Nirenberg inequality, 

%Also, we need the Young's inequality of convolutional form: Let $ f\in L^{r}(\mathbb{R}^{n}), g\in L^{q}(\mathbb{R}^{n}), $ then $ f\ast g\in L^{p}(\mathbb{R}^{n}), $ and
%\begin{equation}\label{4}
%\|f\ast g\|_{L^{p}(\mathbb{R}^{n})}\leq \|f\|_{L^{r}(\mathbb{R}^{n})}\|g\|_{L^{q}(\mathbb{R}^{n})},
%\end{equation}
%where $ 1\leq r, p, q\leq \infty, $ satisfies $ 1+\frac{1}{p}=\frac{1}{r}+\frac{1}{q}. $ 
%In the following estimates, we can take $ \tau\in (0,1) $ as sufficiently small, such that $ \frac12\leq 1+\tau n^{\varepsilon,\tau}\leq 2. $

{\bf 
	%(\uppercase\expandafter{\romannumeral1}) To show
	Step I. $\Phi$ is a mapping from $S$ to $S$.}

{\bf \underline{Estimation of $ \Phi_1 $}.}
%Select $\beta\in(0,1)$ such that $\frac 3{8}<\beta<\frac12$, so that  $D(B^\beta)\hookrightarrow C^0(\mathbb{R}^3)$, where $ B $ stands for the sectorial operator $-\Delta+1$ in $L^4$. Moreover, 
Noting that $ \|e^{t\Delta}n_{0}^{\varepsilon}\|_{L^{\infty}}\leq \|n_{0}^{\varepsilon}\|_{L^{\infty}} $ due to the maximum principle and heat kernel estimates of Lemma \ref{e Delta f } in the appendix, we have 
\begin{equation*}
\begin{aligned}
&\|\Phi_1(n^{\varepsilon,\tau},c^{\varepsilon,\tau},u^{\varepsilon,\tau})(\cdot,t)\|_{L^\infty}\\ \leq~
&\|e^{t\Delta}n_0^{\varepsilon}\|_{L^\infty}+C\int_0^t \left\| e^{-(t-s)\Delta}\nabla\cdot\left(\frac{n^{\varepsilon,\tau}}{1+\tau n^{\varepsilon,\tau}} \nabla c^{\varepsilon,\tau} \chi(c^{\varepsilon,\tau})+(u^{\varepsilon,\tau} \ast \rho^\varepsilon) \cdot n^{\varepsilon,\tau}\right)\right\|_{L^\infty}ds\\\leq~
&\|n_0^{\varepsilon}\|_{L^\infty}+C\int_0^t(t-s)^{-\frac78}\left\|(\frac{n^{\varepsilon,\tau}}{1+\tau n^{\varepsilon,\tau}} \nabla c^{\varepsilon,\tau} \chi(c^{\varepsilon,\tau}))+
(u^{\varepsilon,\tau} \ast \rho^\varepsilon) \cdot n^{\varepsilon,\tau}\right\|_{L^4}ds.
\end{aligned}
\end{equation*}
It follows from (\ref{eq:norm}) that $ \|c^{\varepsilon,\tau}\|_{L^{\infty}}\leq C\|\nabla^{2}c^{\varepsilon,\tau}\|_{L^{2}}^{3/4}\|c^{\varepsilon,\tau}\|_{L^{2}}^{1/4}\leq CR $ due to embedding inequality, then by (\ref{ine:chi}) and (\ref{ine:kappa}) we know
\ben\label{eq:c infty estimate}
&&\|\chi(c^{\varepsilon,\tau})\|_{L^{\infty}}\leq C(R),~~\|\chi'(c^{\varepsilon,\tau})\|_{L^{\infty}}\leq C(R);\nonumber\\
&&\|\kappa(c^{\varepsilon,\tau})\|_{L^{\infty}}\leq C(R),~~\|\kappa'(c^{\varepsilon,\tau})\|_{L^{\infty}}\leq C(R).
\een
Moreover, using \begin{equation}\label{1}
\|\nabla c^{\varepsilon,\tau}\|_{L^{4}}\leq C\|c^{\varepsilon,\tau}\|_{L^{2}}^{\frac18}\|\nabla^{2}c^{\varepsilon,\tau}\|_{L^{2}}^{\frac78},
\end{equation}
\begin{equation}\label{2}
\|u^{\varepsilon,\tau}\|_{L^{4}}\leq C\|u^{\varepsilon,\tau}\|_{L^{2}}^{\frac58}\|\nabla^{2}u^{\varepsilon,\tau}\|_{L^{2}}^{\frac38},
\end{equation}
Young inequality and (\ref{eq:c infty estimate}), we get
\begin{equation*}
\begin{aligned}
&\|\Phi_1(n^{\varepsilon,\tau},c^{\varepsilon,\tau},u^{\varepsilon,\tau})(\cdot,t)\|_{L^\infty}\\
\leq~&\|n_{0}^{\varepsilon}\|_{L^{\infty}}+	C\int_{0}^{t}(t-s)^{-\frac78}\left\{\|n^{\varepsilon,\tau}\|_{L^{\infty}}\|\nabla c^{\varepsilon,\tau}\|_{L^{4}}\|\chi(c^{\varepsilon,\tau})\|_{L^{\infty}}+\|u^{\varepsilon,\tau}\|_{L^{4}}\|\rho^{\varepsilon}\|_{L^{1}}\|n^{\varepsilon,\tau}\|_{L^{\infty}}\right\}ds\\
\leq~
& \|n_0^{\varepsilon}\|_{L^\infty}+C(R)T^{\frac18}\quad\quad~~{\rm{for ~~all}} ~~t\in (0,T).
\end{aligned}
\end{equation*}
Direct calculation and (\ref{eq:norm}) yield that
\begin{equation*}
\begin{aligned}
&\|\Phi_1(n^{\varepsilon,\tau},c^{\varepsilon,\tau},u^{\varepsilon,\tau})(\cdot,t)\|_{L^{2}}\\
%\leq~&\|e^{t\Delta}n_0^{\varepsilon}\|_{L^2}+\left\|\int_0^te^{(t-s)\Delta}\big\{\nabla \cdot \left(\frac{n^{\varepsilon,\tau}}{1+\tau n^{\varepsilon,\tau}} \nabla c^{\varepsilon,\tau}\chi(c^{\varepsilon,\tau})\right)+(u^{\varepsilon,\tau} \ast \rho^\varepsilon) \cdot \nabla n^{\varepsilon,\tau}\big\}(\cdot,s)ds\right\|_{L^2}\\
\leq~
&\|n_{0}^{\varepsilon}\|_{L^{2}}+C \int_{0}^{t}(t-s)^{-\frac12}\left\|\frac{n^{\varepsilon,\tau}}{1+\tau n^{\varepsilon,\tau}}\nabla c^{\varepsilon,\tau}\chi(c^{\varepsilon,\tau})+(u^{\varepsilon,\tau}\ast\rho^{\varepsilon})\cdot\nabla n^{\varepsilon,\tau}\right\|_{L^{2}}ds\\
%\leq~&\|n_{0}^{\varepsilon}\|_{L^{2}}+C\int_{0}^{t}(t-s)^{-\frac12}\left( \|\frac{n^{\varepsilon,\tau}}{1+\tau n^{\varepsilon,\tau}}\|_{L^{2}}\|\nabla c^{\varepsilon,\tau}\|_{L^{\infty}}\|\chi(c^{\varepsilon,\tau})\|_{L^{\infty}}+\|u^{\varepsilon,\tau}\ast\rho^{\varepsilon}\|_{L^{\infty}}\|\nabla n^{\varepsilon,\tau}\|_{L^{2}} \right)ds\\
\leq ~&\|n_{0}^{\varepsilon}\|_{L^{2}}+C\int_{0}^{t}(t-s)^{-\frac12}\left(\|n^{\varepsilon,\tau}\|_{L^{2}}\|\nabla c^{\varepsilon,\tau}\|_{L^{\infty}}\|\chi(c^{\varepsilon,\tau})\|_{L^{\infty}}+\|u^{\varepsilon,\tau}\|_{L^{\infty}}\|\rho^{\varepsilon}\|_{L^{1}}\|\nabla n^{\varepsilon,\tau}\|_{L^{2}}  \right)ds\\\leq~& \|n_{0}^{\varepsilon}\|_{L^{2}}+C(R)T^{\frac12}\quad\quad~~{\rm{for ~~all}} ~~t\in (0,T).
\end{aligned}
\end{equation*}
By (\ref{1}) and 
\begin{equation}\label{3}
\|u^{\varepsilon,\tau}\|_{L^{\infty}}\leq C\|u^{\varepsilon,\tau}\|_{L^{2}}^{\frac14}\|\nabla^{2}u^{\varepsilon,\tau}\|_{L^{2}}^{\frac34},
\end{equation}
noting that $ \nabla\cdot u^{\varepsilon,\tau}=0, $ we have
\begin{equation*}
\begin{aligned}
&\|\nabla\Phi_1(n^{\varepsilon,\tau},c^{\varepsilon,\tau}, u^{\varepsilon,\tau})(\cdot, t)\|_{L^{2}}\\\leq~& \|\nabla n_{0}^{\varepsilon}\|_{L^{2}}+\int_{0}^{t}\left\|\nabla e^{(t-s)\Delta}\nabla\cdot\left(\frac{n^{\varepsilon,\tau}}{1+\tau n^{\varepsilon,\tau}}\cdot\nabla c^{\varepsilon,\tau}\chi(c^{\varepsilon,\tau})+(u^{\varepsilon,\tau}\ast\rho^{\varepsilon})\cdot n^{\varepsilon,\tau}\right)\right\|_{L^{2}}ds\\
%\leq~&\|\nabla n_{0}^{\varepsilon}\|_{L^{2}}+C\int_{0}^{t}(t-s)^{-\frac12}\|\frac{\nabla n^{\varepsilon,\tau}}{(1+\tau n^{\varepsilon,\tau})^{2}}\cdot\nabla c^{\varepsilon,\tau}\chi(c^{\varepsilon,\tau})\\&+\frac{n^{\varepsilon,\tau}}{1+\tau n^{\varepsilon,\tau}}\Delta c^{\varepsilon,\tau}\chi(c^{\varepsilon,\tau})+\frac{n^{\varepsilon,\tau}}{1+\tau n^{\varepsilon,\tau}}|\nabla c^{\varepsilon,\tau}|^{2}\chi'(c^{\varepsilon,\tau})+(u^{\varepsilon,\tau}\ast\rho^{\varepsilon})\cdot\nabla n^{\varepsilon,\tau}\|_{L^{2}}ds\\
\leq~&\|\nabla n_{0}^{\varepsilon}\|_{L^{2}}+C\int_{0}^{t}(t-s)^{-\frac12}( \|\nabla n^{\varepsilon,\tau}\|_{L^{2}}\|\nabla c^{\varepsilon,\tau}\|_{L^{\infty}}\|\chi(c^{\varepsilon,\tau})\|_{L^{\infty}}\\&+~\|n^{\varepsilon,\tau}\|_{\infty}\|\nabla^{2}c^{\varepsilon,\tau}\|_{L^{2}}\|\chi(c^{\varepsilon,\tau})\|_{L^{\infty}}+\|n^{\varepsilon,\tau}\|_{L^{\infty}}\|\nabla c^{\varepsilon,\tau}\|_{L^{4}}^{2}\|\chi'(c^{\varepsilon,\tau})\|_{L^{\infty}}\\&+~\|u^{\varepsilon,\tau}\|_{L^{\infty}}\|\rho^{\varepsilon}\|_{L^{1}}\|\nabla n^{\varepsilon,\tau}\|_{L^{2}} )ds
\\\leq~&\|\nabla n_{0}^{\varepsilon}\|_{L^{2}}+C(R)T^{\frac12}\quad\quad~~{\rm{for ~~all}} ~~t\in (0,T).
\end{aligned}
\end{equation*}

%{\color{blue}
%\begin{equation*}
%	\begin{aligned}
%	&\|\nabla^{2}\Phi_1(n^{\varepsilon}, c^{\varepsilon}, u^{\varepsilon})(\cdot,t)\|_{L^{2}(\mathbb{R}^{3})}\\\leq~&\leq \|\nabla^{2}e^{t\Delta}n_{0}^{\varepsilon}\|_{L^{2}}+\int_{0}^{t}\|\nabla^{2}e^{(t-s)\Delta}\left\{ \nabla\cdot\big(\frac{n^{\varepsilon}}{1+\tau n^{\varepsilon}}\nabla c^{\varepsilon}\chi(c^{\varepsilon})  \big)+(u^{\varepsilon}\ast\rho^{\varepsilon})\cdot\nabla n^{\varepsilon} \right\}(\cdot,s)\|_{L^{2}}ds\\\leq~&C\|\nabla^{2}n_{0}^{\varepsilon}\|_{L^{2}}+C\int_{0}^{t}(t-s)^{-\frac12}\|\nabla\{
%	\frac{\nabla n^{\varepsilon}}{(1+\tau n^{\varepsilon})^{2}}\cdot\nabla c^{\varepsilon}\chi(c^{\varepsilon})+\frac{n^{\varepsilon}}{1+\tau n^{\varepsilon}}\nabla^{2}c^{\varepsilon}\chi(c^{\varepsilon})\\&+\frac{n^{\varepsilon}}{1+\tau n^{\varepsilon}}|\nabla c^{\varepsilon}|^{2}\chi'(c^{\varepsilon})+(u^{\varepsilon}\ast\rho^{\varepsilon})\cdot\nabla n^{\varepsilon}
%	\}\|_{L^{2}}ds
%	\end{aligned}
%\end{equation*}
%
%}

{\bf\underline{Estimation of $ \Phi_2 $}.}
Similar to the estimate of $ \|\Phi_1\|_{L^{\infty}}, $ note that  for any $ n^{\varepsilon,\tau}>0, $ ~ $ \tau>0, $ 
\ben\label{ln 1+f}
\ln(1+\tau n^{\varepsilon,\tau})\leq \tau n^{\varepsilon,\tau}.
\een
Besides, by 
\begin{equation}\label{0}
\|n^{\varepsilon,\tau}\|_{L^{4}}\leq C\|n^{\varepsilon,\tau}\|_{L^{2}}^{\frac14}\|\nabla n^{\varepsilon,\tau}\|_{L^{2}}^{\frac34},
\end{equation}
(\ref{1}) and (\ref{3}), we get
\begin{equation*}
\begin{aligned}
&\|\nabla\Phi_2(n^{\varepsilon,\tau},c^{\varepsilon,\tau},u^{\varepsilon,\tau})(\cdot,t)\|_{L^\infty}\\\leq~& \|\nabla e^{t\Delta}c_0^{\varepsilon}\|_{L^\infty}+\int_{0}^{t}\left\|\nabla e^{(t-s)\Delta}\left\{\frac{1}{\tau}\ln(1+\tau n^{\varepsilon,\tau})\kappa(c^{\varepsilon,\tau})+u^{\varepsilon,\tau}\cdot\nabla c^{\varepsilon,\tau}\right\}(\cdot, s)\right\|_{L^{\infty}}ds
%\\\leq~
%&\|\nabla c_0^{\varepsilon}\|_{L^\infty}+C\int_0^t \left\|B^\beta\nabla e^{(t-s)\Delta}\left(\frac{1}{\tau}\ln(1+\tau n^{\varepsilon,\tau})\kappa(c^{\varepsilon,\tau})+(u^{\varepsilon,\tau}\ast\rho^{\varepsilon})\cdot\nabla c^{\varepsilon,\tau}\right)\right\|_{L^4}ds
\\\leq~
&\|\nabla c_0^{\varepsilon}\|_{L^\infty}+C\int_0^t(t-s)^{-\frac78}\left\|\frac{1}{\tau}\ln(1+\tau n^{\varepsilon,\tau})\kappa(c^{\varepsilon,\tau})+u^{\varepsilon,\tau}\cdot\nabla c^{\varepsilon,\tau}\right\|_{L^4}ds\\\leq ~&
%\|\nabla c_{0}^{\varepsilon}\|_{L^{\infty}}+	C\int_{0}^{t}(t-s)^{-\frac78}\left\{\|\frac{\ln(1+\tau n^{\varepsilon,\tau})}{\tau}\|_{L^{4}}\|\kappa(c^{\varepsilon,\tau})\|_{L^{\infty}}+\|u^{\varepsilon,\tau}\ast\rho^{\varepsilon}\|_{L^{\infty}}\|\nabla c^{\varepsilon,\tau}\|_{L^{4}}\right\}ds\\\leq~&
\|\nabla c_{0}^{\varepsilon}\|_{L^{\infty}}+C\int_{0}^{t}(t-s)^{-\frac78}\left\{\|n^{\varepsilon,\tau}\|_{L^{4}}\|\kappa(c^{\varepsilon,\tau})\|_{L^{\infty}}+\|u^{\varepsilon,\tau}\|_{L^{\infty}}\|\nabla c^{\varepsilon,\tau}\|_{L^{4}}\right\}ds\\\leq~
& \|\nabla c_0^{\varepsilon}\|_{L^\infty}+C(R)T^{\frac18}\quad\quad~~{\rm{for ~~all}} ~~t\in (0,T).
\end{aligned}
\end{equation*}
By (\ref{1}) and (\ref{2}), we have
\begin{equation*}
\begin{aligned}
&\|\Phi_2(n^{\varepsilon,\tau},c^{\varepsilon,\tau},u^{\varepsilon,\tau})(\cdot,t)\|_{L^2}\\\leq~& \| e^{t\Delta}c_0^{\varepsilon}\|_{L^2}+\int_{0}^{t}\left\| e^{(t-s)\Delta}\left\{\frac{1}{\tau}\ln(1+\tau n^{\varepsilon,\tau})\kappa(c^{\varepsilon,\tau})+u^{\varepsilon,\tau}\cdot\nabla c^{\varepsilon,\tau}\right\}(\cdot, s)\right\|_{L^{2}}ds
\\\leq~
%&\|c_0^{\varepsilon}\|_{L^2}+C\int_0^t \left\|\frac{1}{\tau}\ln(1+\tau n^{\varepsilon,\tau})\kappa(c^{\varepsilon,\tau})+(u^{\varepsilon,\tau}\ast\rho^{\varepsilon})\cdot\nabla c^{\varepsilon,\tau}\right\|_{L^2}ds\\\leq~
&\|c_0^{\varepsilon}\|_{L^2}+C\int_0^t\left(\|\frac{\ln(1+\tau n^{\varepsilon,\tau})}{\tau}\|_{L^{2}}\|\kappa(c^{\varepsilon,\tau})\|_{L^{\infty}}+\|u^{\varepsilon,\tau}\|_{L^{4}}\|\nabla c^{\varepsilon,\tau}\|_{L^{4}}  \right)ds
%\\\leq ~
%&\|c_{0}^{\varepsilon}\|_{L^{2}}+	C(\tau)\int_{0}^{t}\left( \|n^{\varepsilon,\tau}\|_{L^{2}}\|\kappa(c^{\varepsilon,\tau})\|_{L^{\infty}}+\|u^{\varepsilon,\tau}\|_{L^{4}}\|\rho^{\varepsilon}\|_{L^{1}}\|\nabla c^{\varepsilon,\tau}\|_{L^{4}} \right)ds
\\\leq~&\|c_0^{\varepsilon}\|_{L^2}+C(R)T\quad\quad~~{\rm{for ~~all}} ~~t\in (0,T).
\end{aligned}
\end{equation*}
By (\ref{3}) and 
\begin{equation}\label{3'}
\|\nabla u^{\varepsilon,\tau}\|_{L^{2}}\leq C\|u^{\varepsilon,\tau}\|_{L^{2}}^{\frac12}\|\nabla^{2}u^{\varepsilon,\tau}\|_{L^{2}}^{\frac12};
\end{equation}
there holds
\begin{equation*}
\begin{aligned}
&\|\nabla^{2}\Phi_2(n^{\varepsilon,\tau},c^{\varepsilon,\tau},u^{\varepsilon,\tau})(\cdot,t)\|_{L^2}\\
%\leq~& \|\nabla^{2} e^{t\Delta}c_0^{\varepsilon}\|_{L^2}+\int_{0}^{t}\left\|\nabla e^{(t-s)\Delta}\nabla\left\{\frac{1}{\tau}\ln(1+\tau n^{\varepsilon,\tau})\kappa(c^{\varepsilon,\tau})+(u^{\varepsilon,\tau}\ast\rho^{\varepsilon})\cdot\nabla c^{\varepsilon,\tau}\right\}(\cdot, s)\right\|_{L^{2}}ds\\
\leq~
&\|\nabla^{2}c_0^{\varepsilon}\|_{L^2}+C\int_0^t (t-s)^{-\frac12}\left\|\frac{\nabla n^{\varepsilon,\tau}}{1+\tau n^{\varepsilon,\tau}}\kappa(c^{\varepsilon,\tau})+\frac{\ln(1+\tau n^{\varepsilon,\tau})}{\tau}\kappa'(c^{\varepsilon,\tau})\nabla c^{\varepsilon,\tau}\right\|_{L^2}\\&+C\int_0^t (t-s)^{-\frac12}\left\|\nabla u^{\varepsilon,\tau}\cdot\nabla c^{\varepsilon,\tau}+u^{\varepsilon,\tau}\cdot\nabla^{2} c^{\varepsilon,\tau}\right\|_{L^{2}} ds\\
%\leq~
%&\|\nabla^{2}c_0^{\varepsilon}\|_{L^2}+C\int_0^t(t-s)^{-\frac12}\{\|\frac{\nabla n^{\varepsilon,\tau}}{1+\tau n^{\varepsilon,\tau}}\|_{L^{2}}\|\kappa(c^{\varepsilon,\tau})\|_{L^{\infty}}\\&+\|\frac{\ln (1+\tau n^{\varepsilon,\tau})}{\tau}\|_{L^{2}}\|\kappa'(c^{\varepsilon,\tau})\|_{L^{\infty}}\|\nabla c^{\varepsilon,\tau}\|_{L^{\infty}}+\|\nabla u^{\varepsilon,\tau}\|_{L^{2}}\|\nabla c^{\varepsilon,\tau}\|_{L^{\infty}}+\|u^{\varepsilon,\tau}\|_{L^{\infty}}\|\nabla^{2}c^{\varepsilon,\tau}\|_{L^{2}}\}ds\\
\leq~&\|\nabla^{2}c_{0}^{\varepsilon}\|_{L^{2}}+C(\tau)\int_{0}^{t}(t-s)^{-\frac12}\left\{\|\nabla n^{\varepsilon,\tau}\|_{L^{2}}\|\kappa(c^{\varepsilon,\tau})\|_{L^{\infty}}+\|n^{\varepsilon,\tau}\|_{L^{2}}\|\kappa'(c^{\varepsilon,\tau})\|_{L^{\infty}}\|\nabla c^{\varepsilon,\tau}\|_{L^{\infty}} \right\}ds \\&+C(\tau)\int_{0}^{t}(t-s)^{-\frac12}\left\{\|\nabla u^{\varepsilon,\tau}\|_{L^{2}}\|\nabla c^{\varepsilon,\tau}\|_{L^{\infty}}+\|u^{\varepsilon,\tau}\|_{L^{\infty}}\|\nabla^{2}c^{\varepsilon,\tau}\|_{L^{2}}\right\}ds\\\leq~
& \|\nabla^{2}c_0^{\varepsilon}\|_{L^2}+C(R)T^{\frac12}\quad\quad~~{\rm{for ~~all}} ~~t\in (0,T).
\end{aligned}
\end{equation*}

{\bf\underline{Estimation of $ \Phi_3 $}.}
By  (\ref{2}) and 
\begin{equation}\label{3''}
\|\nabla u^{\varepsilon,\tau}\|_{L^{4}}\leq C\|u^{\varepsilon,\tau}\|_{L^{2}}^{\frac18}\|\nabla^{2}u^{\varepsilon,\tau}\|_{L^{2}}^{\frac78},
\end{equation} 
we have
\begin{equation*}
\begin{aligned}
&\| \Phi_3(n^{\varepsilon,\tau},c^{\varepsilon,\tau},u^{\varepsilon,\tau})(\cdot,t)\|_{L^2}\\\leq~& 
\| e^{t\Delta}u_0^{\varepsilon}\|_{L^2}+ \int_0^t\left\|e^{(t-s)\Delta}\mathcal{P}\left\{\mu(u^{\varepsilon,\tau} \ast \rho^\varepsilon)\cdot \nabla u^{\varepsilon,\tau}+(n^{\varepsilon,\tau} \nabla \phi) \ast \rho^\varepsilon \right\}(\cdot,s)\right\|_{L^2}ds
%\\
%\leq~& \| u_0^{\varepsilon}\|_{L^2}+C\int_0^t\left\|(u^{\varepsilon,\tau} \ast \rho^\varepsilon)\cdot \nabla u^{\varepsilon,\tau}+(n^{\varepsilon,\tau} \nabla \phi) \ast \rho^\varepsilon\right\|_{L^2} ds
\\
%\leq~& \|u_0^{\varepsilon}\|_{L^2}+C\int_0^t\big(\|u^{\varepsilon,\tau} \ast \rho^\varepsilon\|_{L^{4}}\| \nabla u^{\varepsilon,\tau}\|_{L^4}+\|n^{\varepsilon,\tau}\nabla \phi\|_{L^{2}} \| \rho^\varepsilon\|_{L^1}\big)ds\\
\leq~&  \|u_0^{\varepsilon}\|_{L^2}+C\int_0^t\big(\|u^{\varepsilon,\tau}\|_{L^{4}}\|\rho^{\varepsilon}\|_{L^{1}}\|\nabla u^{\varepsilon,\tau}\|_{L^{4}}+\|n^{\varepsilon,\tau}\|_{L^{2}}\|\nabla\phi\|_{L^{\infty}}\|\rho^{\varepsilon}\|_{L^{1}}\big)ds\\
\leq~& \| u_0^{\varepsilon}\|_{L^2}+C(R) T\quad\quad~~{\rm{for ~~all}} ~~t\in (0,T).
\end{aligned}
\end{equation*}
By (\ref{3''}) and  (\ref{3}), we achieve

\begin{equation*}
\begin{aligned}
&\|\nabla^{2} \Phi_3(n^{\varepsilon,\tau},c^{\varepsilon,\tau},u^{\varepsilon,\tau})(\cdot,t)\|_{L^2}\\\leq~& 
\|\nabla^{2} e^{-tA}u_0^{\varepsilon}\|_{L^2}+ \int_0^t\|\nabla e^{-(t-s)A}\nabla\mathcal{P}\big\{\mu(u^{\varepsilon,\tau} \ast \rho^\varepsilon)\cdot \nabla u^{\varepsilon,\tau}+(n^{\varepsilon,\tau} \nabla \phi) \ast \rho^\varepsilon \big\}(\cdot,s)\|_{L^2}ds
%\\
%\leq~& \|\nabla^{2} u_0^{\varepsilon}\|_{L^2}+C\int_0^t(t-s)^{-\frac12}\|\nabla\{(u^{\varepsilon,\tau} \ast \rho^\varepsilon)\cdot \nabla u^{\varepsilon,\tau}+(n^{\varepsilon,\tau} \nabla \phi) \ast \rho^\varepsilon\}\|_{L^2} ds
\\
\leq~& \|\nabla^{2}u_0^{\varepsilon}\|_{L^2}+C\int_0^t(t-s)^{-\frac12}\big(\|(\nabla u^{\varepsilon,\tau}\ast\rho^{\varepsilon})\cdot\nabla u^{\varepsilon,\tau}+(u^{\varepsilon,\tau}\ast\rho^{\varepsilon})\cdot\nabla^{2}u^{\varepsilon,\tau}\\&+\nabla(n^{\varepsilon,\tau}\nabla\phi)\ast\rho^{\varepsilon}\|_{L^2}\big)ds\\
\leq~&  \|\nabla^{2}u_0^{\varepsilon}\|_{L^2}+C\int_0^t(t-s)^{-\frac12}\big(\|\nabla u^{\varepsilon,\tau}\|_{L^{4}}^{2}\|\rho^{\varepsilon}\|_{L^{1}}+\|u^{\varepsilon,\tau}\|_{L^{\infty}}\|\rho^{\varepsilon}\|_{L^{1}}\|\nabla^{2}u^{\varepsilon,\tau}\|_{L^{2}}\\&+\|n^{\varepsilon,\tau}\|_{L^{2}}\|\nabla\phi\|_{L^{\infty}}\|\nabla\rho^{\varepsilon}\|_{L^{1}}
\big)ds\\
\leq~& \|\nabla^{2} u_0^{\varepsilon}\|_{L^2}+C(R,\varepsilon)T^{\frac12}\quad\quad~~{\rm{for ~~all}} ~~t\in (0,T).
\end{aligned}
\end{equation*}
Combining the estimates of $ \Phi_1 $, $ \Phi_2 $ and $ \Phi_3 $, we get $\Phi$ is a mapping from $S$ to $S$. 

{\bf Step II. To show $\Phi$ acts as a contraction on $ S. $}

Furthermore, using the same idea, we can obtain the following estimates for $ (n_{1}^{\varepsilon,\tau}, c_{1}^{\varepsilon,\tau}, u_{1}^{\varepsilon,\tau}) $, $ (n_{2}^{\varepsilon,\tau}, c_{2}^{\varepsilon,\tau}, u_{2}^{\varepsilon,\tau})\in S $.

{\bf \underline{Estimates of $n^{\varepsilon,\tau}$}.} For $ L^{2} $ norm of $ n^{\varepsilon,\tau}, $ we have 
\begin{equation*}
\begin{aligned} &\|\Phi_1(n_2^{\varepsilon,\tau},c_2^{\varepsilon,\tau},u_2^{\varepsilon,\tau})-\Phi_1(n_1^{\varepsilon,\tau},c_1^{\varepsilon,\tau},u_1^{\varepsilon,\tau})(\cdot,t)\|_{L^2}\\=~&\big{\|}\int_0^te^{(t-s)\Delta}\big\{\nabla \cdot [\frac{n_2^{\varepsilon,\tau}}{1+\tau n_2^{\varepsilon,\tau}} \nabla c_2^{\varepsilon,\tau} \chi(c_2^{\varepsilon,\tau})+(u_2^{\varepsilon,\tau} \ast \rho^\varepsilon) \cdot  n_2^{\varepsilon,\tau}\\&-\frac{n_1^{\varepsilon,\tau}}{1+\tau n_1^{\varepsilon,\tau}} \nabla c_1^{\varepsilon,\tau} \chi(c_1^{\varepsilon,\tau})-(u_1^{\varepsilon,\tau} \ast \rho^\varepsilon) \cdot  n_1^{\varepsilon,\tau}]\big\}ds\big{\|}_{L^2}
\\\leq~&C\int_{0}^{t}(t-s)^{-\frac12}\|
\frac{n_2^{\varepsilon,\tau}}{1+\tau n_2^{\varepsilon,\tau}} \nabla c_2^{\varepsilon,\tau} \chi(c_2^{\varepsilon,\tau})+(u_2^{\varepsilon,\tau} \ast \rho^\varepsilon) \cdot  n_2^{\varepsilon,\tau}\\&-\frac{n_1^{\varepsilon,\tau}}{1+\tau n_1^{\varepsilon,\tau}} \nabla c_1^{\varepsilon,\tau} \chi(c_1^{\varepsilon,\tau})-(u_1^{\varepsilon,\tau} \ast \rho^\varepsilon) \cdot  n_1^{\varepsilon,\tau}\|_{L^{2}}ds\\\leq~& C(R)T^{\frac12}\|(n_2^{\varepsilon,\tau},c_2^{\varepsilon,\tau},u_2^{\varepsilon,\tau})-(n_1^{\varepsilon,\tau},c_1^{\varepsilon,\tau},u_1^{\varepsilon,\tau})\|_{S},
\end{aligned}
\end{equation*}
where we used $ |\chi(c_{1}^{\varepsilon,\tau})-\chi(c_{2}^{\varepsilon,\tau})|\leq \chi'(R)|c_{1}^{\varepsilon,\tau}-c_{2}^{\varepsilon,\tau}|. $
For $L^{\infty}$ norm of  $n^{\varepsilon,\tau}$, we obtain
\begin{equation*}
\begin{aligned}
&\|\Phi_1(n_2^{\varepsilon,\tau},c_2^{\varepsilon,\tau},u_2^{\varepsilon,\tau})(\cdot,t)-\Phi_1(n_1^{\varepsilon,\tau},c_1^{\varepsilon,\tau},u_1^{\varepsilon,\tau})(\cdot,t)\|_{L^{\infty}}\\=~&\|\int_0^te^{(t-s)\Delta}
\big\{\nabla \cdot \left(\frac{n_2^{\varepsilon,\tau}}{1+\tau n_2^{\varepsilon,\tau}} \nabla c_2^{\varepsilon,\tau} \chi(c_2^{\varepsilon,\tau})+(u_2^{\varepsilon,\tau} \ast \rho^\varepsilon) \cdot  n_2^{\varepsilon,\tau}\right)\\&-~\nabla \cdot \left(\frac{n_1^{\varepsilon,\tau}}{1+\tau n_1^{\varepsilon,\tau}} \nabla c_1^{\varepsilon,\tau} \chi(c_1^{\varepsilon,\tau})+(u_1^{\varepsilon,\tau} \ast \rho^\varepsilon) \cdot  n_1^{\varepsilon,\tau}\right)\big\}\|_{L^{\infty}}ds\\\leq~&C(R)T^{\frac18}\|(n_2^{\varepsilon,\tau},c_2^{\varepsilon,\tau},u_2^{\varepsilon,\tau})-(n_1^{\varepsilon,\tau},c_1^{\varepsilon,\tau},u_1^{\varepsilon,\tau})\|_{S}.
\end{aligned}
\end{equation*}
Finally, for ${\dot{W}^{1,2}}$ norm about $n^{\varepsilon,\tau}$, we get
\begin{equation*}
\begin{aligned}
&\|\Phi_1(n_2^{\varepsilon,\tau} ,c_2^{\varepsilon,\tau} ,u_2^{\varepsilon,\tau} )(\cdot,t)-\Phi_1(n_1^{\varepsilon,\tau} ,c_1^{\varepsilon,\tau} ,u_1^{\varepsilon,\tau} )(\cdot,t)\|_{\dot{W}^{1,2}}\\
\leq~& \int_0^t\|\nabla e^{(t-s)\Delta}\big\{\nabla \cdot \left(\frac{n_2^{\varepsilon,\tau} }{1+\tau n_2^{\varepsilon,\tau} } \nabla c_2^{\varepsilon,\tau}  \chi(c_2^{\varepsilon,\tau} )\right)+(u_2^{\varepsilon,\tau}  \ast \rho^\varepsilon) \cdot \nabla n_2^{\varepsilon,\tau}\\
&-~\nabla \cdot \left(\frac{n_1^{\varepsilon,\tau} }{1+\tau n_1^{\varepsilon,\tau} } \nabla c_1^{\varepsilon,\tau}  \chi(c_1^{\varepsilon,\tau} )\right)-(u_1^{\varepsilon,\tau}  \ast \rho^\varepsilon) \cdot \nabla n_1^{\varepsilon,\tau} \big\}(\cdot,s)\|_{L^2}ds\\
\leq~& \int_0^t (t-s)^{-\frac12}\|\nabla \cdot \left(\frac{n_2^{\varepsilon,\tau} }{1+\tau n_2^{\varepsilon,\tau} } \nabla c_2^{\varepsilon,\tau}  \chi(c_2^{\varepsilon,\tau} )\right)+(u_2^{\varepsilon,\tau}  \ast \rho^\varepsilon) \cdot \nabla n_2^{\varepsilon,\tau} \\&-~\nabla \cdot \left(\frac{n_1^{\varepsilon,\tau} }{1+\tau n_1^{\varepsilon,\tau} } \nabla c_1^{\varepsilon,\tau}  \chi(c_1^{\varepsilon,\tau} )\right)-(u_1^{\varepsilon,\tau}  \ast \rho^\varepsilon) \cdot \nabla n_1^{\varepsilon,\tau} (\cdot,s)\|_{L^2}ds\\
\leq~ &C(R)T^{\frac12}\|(n_2^{\varepsilon,\tau} ,c_2^{\varepsilon,\tau} ,u_2^{\varepsilon,\tau} )-(n_1^{\varepsilon,\tau} ,c_1^{\varepsilon,\tau} ,u_1^{\varepsilon,\tau} )\|_{S}.
\end{aligned}
\end{equation*}

{\bf \underline{Estimates of $ c^{\varepsilon,\tau}$}}.
For the estimate of $L^{\infty}$ about $\nabla c^{\varepsilon,\tau}$, there holds
\begin{equation*}
\begin{aligned}
&\left\|\Phi_2(n_2^{\varepsilon,\tau},c_2^{\varepsilon,\tau},u_2^{\varepsilon,\tau})(\cdot,t)-\Phi_2(n_1^{\varepsilon,\tau},c_1^{\varepsilon,\tau},u_1^{\varepsilon,\tau})(\cdot,t)\right\|_{\dot{W}^{1,\infty}}\\
\leq~&\int_0^t\left\|e^{(t-s)\Delta}\left\{\frac1{\tau}\ln (1+\tau n_2^{\varepsilon,\tau})\kappa(c_2^{\varepsilon,\tau})-\frac1{\tau}\ln (1+\tau n_1^{\varepsilon,\tau})\kappa(c_1^{\varepsilon,\tau})\right\}\right\|_{\dot{W}^{1,\infty}}ds\\
&+~
\int_0^t\left\|e^{(t-s)\Delta}\left\{u_2^{\varepsilon,\tau}  \cdot \nabla c_2^{\varepsilon,\tau}-u_1^{\varepsilon,\tau} \cdot \nabla c_1^{\varepsilon,\tau}\right\}\right\|_{\dot{W}^{1,\infty}}ds
\\
\leq~& CT^{\frac18}\|(n_2^{\varepsilon,\tau},c_2^{\varepsilon,\tau},u_2^{\varepsilon,\tau})-(n_1^{\varepsilon,\tau},c_1^{\varepsilon,\tau},u_1^{\varepsilon,\tau})\|_{S},
\end{aligned}
\end{equation*}
where by (\ref{ln 1+f}), we can yield
\begin{equation*}
\begin{aligned}
&\left\|\frac1{\tau}\ln (1+\tau n_2^{\varepsilon,\tau})\kappa(c_2^{\varepsilon,\tau})+u_2^{\varepsilon,\tau}  \cdot \nabla c_2^{\varepsilon,\tau}-\frac1{\tau}\ln (1+\tau n_1^{\varepsilon,\tau})\kappa(c_1^{\varepsilon,\tau})-u_1^{\varepsilon,\tau}  \cdot \nabla c_1^{\varepsilon,\tau}\right\|_{L^4}\\
\leq~& \left\|\frac1{\tau}\ln{\frac{1+\tau n_2^{\varepsilon,\tau}}{1+\tau n_1^{\varepsilon,\tau}}}\kappa(c_2^\varepsilon)\right\|_{L^4}+\left\|\frac1{\tau}\ln (1+\tau n_1^{\varepsilon,\tau})(\kappa(c_2^{\varepsilon,\tau})-\kappa(c_1^{\varepsilon,\tau}))\right\|_{L^4}\\
&+~\left\|(u_2^{\varepsilon,\tau}-u_1^{\varepsilon,\tau})\cdot \nabla c_2^{\varepsilon,\tau}\|_{L^4}+\|u_1^{\varepsilon,\tau} \cdot\nabla(c_2^{\varepsilon,\tau}-c_1^{\varepsilon,\tau})\right\|_{L^4}\\
\leq~&\|n_2^{\varepsilon,\tau}-n_1^{\varepsilon,\tau}\|_{L^4}\|\kappa(c_2^{\varepsilon,\tau})\|_{L^{\infty}}+\|n_1^{\varepsilon,\tau}\|_{L^{\infty}}\|\kappa^{'}(\theta c_2^\varepsilon+(1-\theta)c_1^\varepsilon)(c_2^\varepsilon-c_1^\varepsilon)\|_{L^4}\\
&+~\|u_2^{\varepsilon,\tau}-u_1^{\varepsilon,\tau}\|_{L^4}\|\nabla c_2^{\varepsilon,\tau}\|_{L^{\infty}}+\|u_1^{\varepsilon,\tau}\|_{L^{\infty}}\|\nabla(c_2^{\varepsilon,\tau}-c_1^{\varepsilon,\tau})\|_{L^4}\\
\leq~& C(\kappa(R),\kappa^{'}(R))\|(n_2^{\varepsilon,\tau},c_2^{\varepsilon,\tau},u_2^{\varepsilon,\tau})-(n_1^{\varepsilon,\tau},c_1^{\varepsilon,\tau},u_1^{\varepsilon,\tau})\|_{S}.
\end{aligned}
\end{equation*}
The estimate of $ L^{2} $ about $ c^{\varepsilon,\tau} $ is similar to the estimate of $L^{\infty}$ about $\nabla c^{\varepsilon,\tau}$, we have
\begin{equation*}
\begin{aligned}
&\|\Phi_2(n_2^{\varepsilon,\tau},c_2^{\varepsilon,\tau},u_2^{\varepsilon,\tau})(\cdot,t)-\Phi_2(n_1^{\varepsilon,\tau},c_1^{\varepsilon,\tau},u_1^{\varepsilon,\tau})(\cdot,t)\|_{L^2}\\\leq~&\left\|\int_0^te^{(t-s)\Delta}\left\{\frac1{\tau}\ln (1+\tau n_2^{\varepsilon,\tau})\kappa(c_2^{\varepsilon,\tau})-\frac1{\tau}\ln (1+\tau n_1^{\varepsilon,\tau})\kappa(c_1^{\varepsilon,\tau}) \right\}ds  \right\|_{L^2}\\&+~\left\|\int_0^te^{(t-s)\Delta}\left\{ u_2^{\varepsilon,\tau}  \cdot \nabla c_2^{\varepsilon,\tau}-u_1^{\varepsilon,\tau}\cdot\nabla c_{1}^{\varepsilon,\tau}\right\}ds  \right\|_{L^{2}}\\\leq~&\int_{0}^{t}\left\|e^{(t-s)\Delta}\left(\frac1{\tau}\ln (1+\tau n_2^{\varepsilon,\tau})\kappa(c_2^{\varepsilon,\tau})-\frac1{\tau}\ln (1+\tau n_1^{\varepsilon,\tau})\kappa(c_1^{\varepsilon,\tau})  \right)
\right\|_{L^{2}}ds\\&+~\int_{0}^{t}\left\|e^{(t-s)\Delta}\left(u_2^{\varepsilon,\tau}  \cdot \nabla c_2^{\varepsilon,\tau}-u_1^{\varepsilon,\tau}\cdot\nabla c_{1}^{\varepsilon,\tau}  \right)
\right\|_{L^{2}}ds\\\leq~&C(R)T\|(n_2^{\varepsilon,\tau},c_2^{\varepsilon,\tau},u_2^{\varepsilon,\tau})-(n_1^{\varepsilon,\tau},c_1^{\varepsilon,\tau},u_1^{\varepsilon,\tau})\|_{S}.
\end{aligned}
\end{equation*}
Moreover, we need to estimate the $\dot{W}^{2,2}$ norm for $c^{\varepsilon,\tau}$. We only need to estimate the following norm 
\begin{equation*}
\begin{aligned}
&\left\|\nabla \left\{\frac1{\tau}\ln (1+\tau n_2^{\varepsilon,\tau})\kappa(c_2^{\varepsilon,\tau})+u_2^{\varepsilon,\tau} \cdot \nabla c_2^{\varepsilon,\tau}-\frac1{\tau}\ln (1+\tau n_1^{\varepsilon,\tau})\kappa(c_1^{\varepsilon,\tau})-u_1^{\varepsilon,\tau} \cdot \nabla c_1^{\varepsilon,\tau}\right\}\right\|_{L^{2}}\\
\leq~& \left\|\left(\frac{\nabla n_2^{\varepsilon,\tau}}{1+\tau n_2^{\varepsilon,\tau}}-\frac{\nabla n_1^{\varepsilon,\tau}}{1+\tau n_1^{\varepsilon,\tau}}\right)\kappa(c_2^{\varepsilon,\tau})\right\|_{L^2}+\left\|\frac{\nabla n_1^{\varepsilon,\tau}}{1+\tau n_1^{\varepsilon,\tau}}(\kappa(c_2^{\varepsilon,\tau})-\kappa(c_1^{\varepsilon,\tau}))\right\|_{L^2}\\
&+~\left\|\frac1{\tau}\ln{\frac{1+\tau n_2^{\varepsilon,\tau}}{1+\tau n_1^{\varepsilon,\tau}}}\kappa^{'}(c_2^{\varepsilon,\tau})\nabla c_2^{\varepsilon,\tau}\right\|_{L^2}+\left\|\frac1{\tau}\ln (1+\tau n_1^{\varepsilon,\tau})\nabla c_2^{\varepsilon,\tau}[\kappa^{'}(c_2^{\varepsilon,\tau})-\kappa^{'}(c_1^{\varepsilon,\tau})]\right\|_{L^2}\\
&+~\left\|\frac1{\tau}\ln (1+\tau n_1^{\varepsilon,\tau})(\nabla c_2^{\varepsilon,\tau}-\nabla c_1^{\varepsilon,\tau})\kappa^{'}(c_1^{\varepsilon,\tau})\right\|_{L^2}+\left\|\nabla(u_{2}^{\varepsilon,\tau}-u_{1}^{\varepsilon,\tau})\cdot\nabla c_{2}^{\varepsilon,\tau}\right\|_{L^{2}}\\&+~\left\| \nabla u_{1}^{\varepsilon,\tau}\cdot\nabla (c_{2}^{\varepsilon,\tau}-c_{1}^{\varepsilon,\tau})  \right\|_{L^{2}}+\left\| (u_{2}^{\varepsilon,\tau}-u_{1}^{\varepsilon,\tau})\cdot\nabla^{2}c_{2}^{\varepsilon,\tau} \right\|_{L^{2}}+\left\| u_{1}^{\varepsilon,\tau}\cdot\nabla^{2}(c_{2}^{\varepsilon,\tau}-c_{1}^{\varepsilon,\tau}) \right\|_{L^{2}}\\
\leq~& \left\|\frac{\nabla n_2^{\varepsilon,\tau}}{1+\tau n_2^{\varepsilon,\tau}}-\frac{\nabla n_1^{\varepsilon,\tau}}{1+\tau n_1^{\varepsilon,\tau}}\right\|_{L^2}\|\kappa(c_2^{\varepsilon,\tau})\|_{L^{\infty}}+\left\|\frac{\nabla n_1^{\varepsilon,\tau}}{1+\tau n_1^{\varepsilon,\tau}}\right\|_{L^{\infty}}\|\kappa^{'}\|_{L^{\infty}}\|c_2^{\varepsilon,\tau}-c_1^{\varepsilon,\tau}\|_{L^2}\\
&+~\|\frac1{\tau}\ln{\frac{1+\tau n_2^{\varepsilon,\tau}}{1+\tau n_1^{\varepsilon,\tau}}}\|_{L^2}\|\kappa^{'}(c_2^{\varepsilon,\tau})\|_{L^{\infty}}\|\nabla c_2^{\varepsilon,\tau}\|_{L^{\infty}}\\
&+~\|\frac1{\tau}\ln (1+\tau n_1^{\varepsilon,\tau})\|_{L^{\infty}}\|\nabla c_2^{\varepsilon,\tau}\|_{L^{\infty}}\|\kappa^{''}\|_{L^{\infty}}\|c_2^{\varepsilon,\tau}-c_1^{\varepsilon,\tau}\|_{L^2}\\
&+~\|\frac1{\tau}\ln (1+\tau n_1^{\varepsilon,\tau})\|_{L^{\infty}}\|\nabla c_2^{\varepsilon,\tau}-\nabla c_1^{\varepsilon,\tau}\|_{L^2}\|\kappa^{'}(c_1^{\varepsilon,\tau})\|_{L^{\infty}}\\&+~\|\nabla(u_{2}^{\varepsilon,\tau}-u_{1}^{\varepsilon,\tau})\|_{L^{2}}\|\nabla c^{\varepsilon,\tau}\|_{L^{\infty}}+\|\nabla u_{1}^{\varepsilon,\tau}\|_{L^{2}}\|\nabla(c_{2}^{\varepsilon,\tau}-c_{1}^{\varepsilon,\tau})\|_{L^{\infty}}\\&+~\|u_{2}^{\varepsilon,\tau}-u_{1}^{\varepsilon,\tau}\|_{L^{\infty}}\|\nabla^{2}c_{2}^{\varepsilon,\tau}\|_{L^{2}}+\|u_{1}^{\varepsilon,\tau}\|_{L^{\infty}}\|\nabla^{2}(c_{2}^{\varepsilon,\tau}-c_{1}^{\varepsilon,\tau})\|_{L^{2}}\\
\leq~& C(\kappa(R),\kappa^{'}(R),\kappa^{''}(R))\|(n_2^{\varepsilon,\tau},c_2^{\varepsilon,\tau},u_2^{\varepsilon,\tau})-(n_1^{\varepsilon,\tau},c_1^{\varepsilon,\tau},u_1^{\varepsilon,\tau})\|_{S}.
\end{aligned}
\end{equation*}
%Then, by (\ref{estimate ln}),
%\beno
%&&\|\Phi_2(n_2^{\varepsilon,\tau},c_2^{\varepsilon,\tau},u_2^{\varepsilon,\tau})(\cdot,t)-\Phi_2(n_1^{\varepsilon,\tau},c_1^{\varepsilon,\tau},u_1^{\varepsilon,\tau})(\cdot,t)\|_{\dot{W}^{2,2}}\\
%&&=\|\int_0^t\nabla^2 e^{(t-s)\Delta}\\
%&&\big\{\frac1{\tau}\ln (1+\tau n_2^{\varepsilon,\tau})\kappa(c_2^{\varepsilon,\tau})+(u_2^{\varepsilon,\tau} \ast \rho^\varepsilon) \cdot \nabla c_2^{\varepsilon,\tau}-\frac1{\tau}\ln (1+\tau n_1^{\varepsilon,\tau})\kappa(c_1^{\varepsilon,\tau})-(u_1^{\varepsilon,\tau} \ast \rho^\varepsilon) \cdot \nabla c_1^{\varepsilon,\tau}\big\}(\cdot,s)ds\|_{L^{2}}\\
%&&\leq\int_0^t(t-s)^{-\frac12}\\
%&&\|\nabla \big\{\frac1{\tau}\ln (1+\tau n_2^{\varepsilon,\tau})\kappa(c_2^{\varepsilon,\tau})+(u_2^{\varepsilon,\tau} \ast \rho^\varepsilon) \cdot \nabla c_2^{\varepsilon,\tau}-\frac1{\tau}\ln (1+\tau n_1^{\varepsilon,\tau})\kappa(c_1^{\varepsilon,\tau})-(u_1^{\varepsilon,\tau} \ast \rho^\varepsilon) \cdot \nabla c_1^{\varepsilon,\tau}\big\}(\cdot,s)\|_{L^{2}}ds\\
%&&\leq C(\kappa(R),\kappa^{'}(R),\kappa^{''}(R))T^{\frac12}\|(n_2^{\varepsilon,\tau},c_2^{\varepsilon,\tau},u_2^{\varepsilon,\tau})-(n_1^{\varepsilon,\tau},c_1^{\varepsilon,\tau},u_1^{\varepsilon,\tau})\|_{S}.
%\eeno
The estimate about $u^{\varepsilon,\tau}$ is similar, so we omit it.
Ultimately, one can obtain  
\beno
&&\|\Phi(n_{1}^{\varepsilon,\tau}, c_{1}^{\varepsilon,\tau}, u_{1}^{\varepsilon,\tau})(\cdot,t)-\Phi(n_{2}^{\varepsilon,\tau}, c_{2}^{\varepsilon,\tau}, u_{2}^{\varepsilon,\tau})(\cdot,t)\|_{S}\\&\leq~& C (T^{\frac18}+T^{\frac12}+T)\|(n_{1}^{\varepsilon,\tau}, c_{1}^{\varepsilon,\tau}, u_{1}^{\varepsilon,\tau})-(n_{2}^{\varepsilon,\tau}, c_{2}^{\varepsilon,\tau}, u_{2}^{\varepsilon,\tau})\|_{S},
\eeno
for all $ t\in (0,T), $ which shows that if $ T $ is chosen sufficiently small then $ \Phi $ acts as a contraction on $ S. $ Accordingly, Lemma \ref{lem:Banach fixed point} asserts 
%Using Lemma \ref{lem:3},
there exists $(n^{\varepsilon,\tau},c^{\varepsilon,\tau},u^{\varepsilon,\tau})\in S$, such that $\Phi(n^{\varepsilon,\tau},c^{\varepsilon,\tau},u^{\varepsilon,\tau})=(n^{\varepsilon,\tau},c^{\varepsilon,\tau},u^{\varepsilon,\tau})$.
%Therefore, we can prove that there exist a smooth function $P$ such that $(n^{\varepsilon,\tau},c^{\varepsilon,\tau},u^{\varepsilon,\tau},P)$ is classical solution of (\ref{eq:CNS}).\\

%{\bf The positivity of $n^\varepsilon$ and $c^\varepsilon$}\\

{\bf {Step III. The positivity of $n^{\varepsilon,\tau}$ and $c^{\varepsilon,\tau}$ in system $\eqref{eq:CNS}$.}}

We know that there exists a strong solution $(n^{\varepsilon,\tau},c^{\varepsilon,\tau},u^{\varepsilon,\tau})$ in $ (0,T) $ satisfies
\begin{equation*}
\begin{aligned}
(n^{\varepsilon,\tau}, c^{\varepsilon,\tau}, u^{\varepsilon,\tau})\in& L^{\infty}[(0,T);\{L^{\infty}(\mathbb{R}^{3})\cap W^{1,2}(\mathbb{R}^{3})\}\\&\times \{\dot{W}^{1,\infty}(\mathbb{R}^{3})\cap L^{2}(\mathbb{R}^{3})\cap \dot{W}^{2,2}(\mathbb{R}^{3})\}\times \{L^{2}(\mathbb{R}^{3})\cap \dot{W}^{2,2}(\mathbb{R}^{3})\}].
\end{aligned}
\end{equation*}
%Assume that the maximal interval time interval of the solution of $\eqref{eq:GKS}$ is $(0,T_\ast]$, then the strong solution $(n^{\varepsilon,\tau},c^{\varepsilon,\tau},u^{\varepsilon,\tau})$ satisfies
%\begin{equation}\label{rela:ncu}
%\begin{aligned}
%(n^{\varepsilon,\tau}, c^{\varepsilon,\tau}, u^{\varepsilon,\tau})&\in L^{\infty}((0,T_\ast];\{L^{\infty}(\mathbb{R}^{3})\cap W^{1,2}(\mathbb{R}^{3})\}\\&\times \{\dot{W}^{1,\infty}(\mathbb{R}^{3})\cap L^{2}(\mathbb{R}^{3})\cap \dot{W}^{2,2}(\mathbb{R}^{3})\}\times \{L^{2}(\mathbb{R}^{3})\cap \dot{W}^{2,2}(\mathbb{R}^{3})\}).
%\end{aligned}
%\end{equation}
Now we claim
\ben\label{ine:nc>0}
n^{\varepsilon,\tau}(x,t)\geq 0 \quad {\rm and} \quad c^{\varepsilon,\tau}(x,t) \geq 0 \quad {\rm in} \quad \mathbb{R}^3\times (0,T).
\een

Indeed, setting $n_{+}^{\varepsilon,\tau}= \max\{0,n^{\varepsilon,\tau}\}$, $-n_{-}^{\varepsilon,\tau}= \min\{0,n^{\varepsilon,\tau}\}$, then $n^{\varepsilon,\tau} = n_+^{\varepsilon,\tau} - n_-^{\varepsilon,\tau}$.  Multiply the first equation of (\ref{eq:CNS n+}) by $n_{-}^{\varepsilon,\tau}$, we arrive
\beno
\partial_t n_-^{\varepsilon,\tau} n_-^{\varepsilon,\tau} + (u^{\varepsilon,\tau} \ast \rho^\varepsilon) \cdot \nabla n_-^{\varepsilon,\tau}  n_-^{\varepsilon,\tau} - \Delta n_-^{\varepsilon,\tau}  n_-^{\varepsilon,\tau} =- \nabla \cdot \left(\frac{n_{+}^{\varepsilon,\tau}}{1+\tau n_{+}^{\varepsilon,\tau}} \nabla c^{\varepsilon,\tau} \chi(c^{\varepsilon,\tau})\right)  n_-^{\varepsilon,\tau}.
\eeno
%Since $\|n^{\varepsilon,\tau}\|_{L^\infty}$ is bounded, then we have $\|n_{-}^{\varepsilon,\tau}\|_{L^\infty}\leq C$, $ C $ is an absolute constant. We can choose sufficient small $\tau\in (0, \frac{1}{2\|n_{-}^{\varepsilon,\tau}\|_{L^{\infty}}}),$ 
%%such that $\|\varepsilon n_{-}^\varepsilon\|_{L^\infty} \leq \frac12$.
%which follows $\|\frac{1}{1-\tau n_{-}^{\varepsilon,\tau}}\|_{L^\infty}\leq C$.
Then integrating by parts on the above equation, we obtain
\begin{equation*}
\frac12 \frac d{dt} \|n_-^{\varepsilon,\tau}\|_{L^2}^2 + \|\nabla n_-^{\varepsilon,\tau}\|_{L^2}^2 =0,
\end{equation*}
which follows
%\beno
%\frac12 \frac d{dt} \|n_-^{\varepsilon,\tau}\|_{L^2}^2 + \|\nabla n_-^{\varepsilon,\tau}\|_{L^2}^2 &=& - \nabla \cdot \left(\frac{n^{\varepsilon,\tau}}{1+\tau n^{\varepsilon,\tau}} \nabla c^{\varepsilon,\tau} \chi(c^{\varepsilon,\tau})\right) n_-^{\varepsilon,\tau} \\
%&=&- \int_{\mathbb{R}^3} \frac{n_{-}^{\varepsilon,\tau}}{1-\tau n_{-}^{\varepsilon,\tau}} \nabla c^{\varepsilon,\tau} \chi(c^{\varepsilon,\tau}) \cdot \nabla n_-^{\varepsilon,\tau} \\
%&\leq& \frac12 \|\nabla n_-^{\varepsilon,\tau}\|_{L^2}^2 + C \|n_-^{\varepsilon,\tau}\|_{L^2}^2 \|\nabla c^{\varepsilon,\tau} \chi(c^{\varepsilon,\tau})\|_{L^\infty}^2 \\
%&\leq& \frac12 \|\nabla n_-^{\varepsilon,\tau}\|_{L^2}^2 + C \chi(\|c^{\varepsilon,\tau}\|_{L^\infty})^2 \|n_-^{\varepsilon,\tau}\|_{L^2}^2 \|\nabla c^{\varepsilon,\tau}\|_{L^{\infty}}^2.
%\eeno
%By Gronwall's inequality, for any $t \in (0,T_\ast]$, we have
%\beno
%\|n_-^{\varepsilon,\tau}(t)\|_{L^2}^2 + \int_0^t \|\nabla n_-^{\varepsilon,\tau}\|_{L^2}^2 \leq C \exp\left(\chi(\|c^{\varepsilon,\tau}\|_{L^\infty((0,T_{*}];L^{\infty})})^2 \int_0^t \| \nabla c^{\varepsilon,\tau}\|_{L^{\infty}}^2 \right) \|n_-^{\varepsilon,\tau}(\cdot,0)\|_{L^2}^2.
%\eeno
%By $\eqref{rela:ncu}$ and noting that $\|\nabla c^{\varepsilon,\tau}\|_{L^\infty((0,T_{*}];L^{\infty}(\mathbb{R}^{3}))}\leq C $, we have
\beno
\|n_-^{\varepsilon,\tau}(t)\|_{L^2}^2 + 2\int_0^t \|\nabla n_-^{\varepsilon,\tau}\|_{L^2}^2 =  \|n_-^{\varepsilon,\tau}(\cdot,0)\|_{L^2}^2=\|n_-^{\varepsilon}(\cdot,0)\|_{L^2}^2.
\eeno
Note that $ n_{-}^{\varepsilon}(\cdot,0)=0 $ due to $ n_{0}^{\varepsilon}=n_{0}\ast\rho^{\varepsilon}\geq 0, $ there holds 
%Since $n^{\varepsilon,\tau}_0 = n_0 \ast \rho^\varepsilon$ and $n_0 \geq 0$, we have
\beno
\|n_-^{\varepsilon,\tau}(t)\|_{L^2}^2 + 2\int_0^t \|\nabla n_-^{\varepsilon,\tau}\|_{L^2}^2 = 0,
\eeno
which mean $n_{-}^{\varepsilon,\tau} = 0$. Thus $ n^{\varepsilon,\tau}=n_{+}^{\varepsilon,\tau}, $ we have $ n^{\varepsilon,\tau}\geq 0. $ 

For the part of $c^{\varepsilon,\tau}$, setting $c_{+}^{\varepsilon,\tau}= \max\{0,c^{\varepsilon,\tau}\}$, $-c_{-}^{\varepsilon,\tau}= \min\{0,c^{\varepsilon,\tau}\}$, then $c^{\varepsilon,\tau} = c_+^{\varepsilon,\tau} - c_-^{\varepsilon,\tau}.$ For the second of (\ref{eq:CNS}), we multiply it by $c_{-}^{\varepsilon,\tau}$ and obtain
\beno
\partial_t c^{\varepsilon,\tau} \cdot c_{-}^{\varepsilon,\tau}+ u^{\varepsilon,\tau}  \cdot \nabla c^{\varepsilon,\tau} \cdot c_{-}^{\varepsilon,\tau}- \Delta c^{\varepsilon,\tau}\cdot c_{-}^{\varepsilon,\tau} + \frac1{\tau}\ln (1+\tau n^{\varepsilon,\tau})\kappa(c^{\varepsilon,\tau})\cdot c_{-}^{\varepsilon,\tau}=0.\\
%\kappa(c^{\varepsilon}) (n^\varepsilon \ast \rho^\varepsilon)
\eeno
Integrate by parts and note that $\kappa(0)=0$, $\kappa'(s) \geq 0$ and $\theta \in (0,1)$, we have
\beno
\frac12\frac d{dt}\| c_{-}^{\varepsilon,\tau}(t)\|_{L^2}^2+\|\nabla c_{-}^{\varepsilon,\tau}\|_{L^2}^2&=& \int_{\mathbb{R}^3}\frac1{\tau}\ln (1+\tau n^{\varepsilon,\tau})\kappa(c^{\varepsilon,\tau})\cdot c_{-}^{\varepsilon,\tau}\\
&=&- \frac1{\tau}\int_{\mathbb{R}^3}\frac{\kappa(c^{\varepsilon,\tau})-\kappa(0)}{c^{\varepsilon,\tau}}\cdot c_{-}^{\varepsilon,\tau}\ln (1+\tau n^{\varepsilon,\tau}) c_{-}^{\varepsilon,\tau}\\
&=&- \frac1{\tau}\int_{\mathbb{R}^3}\frac{\kappa^{'}(\theta c^{\varepsilon,\tau})c^{\varepsilon,\tau}}{c^{\varepsilon,\tau}} c_{-}^{\varepsilon,\tau}\ln (1+\tau n^{\varepsilon,\tau}) c_{-}^{\varepsilon,\tau}\\
&\leq&0,
\eeno
which means $c^{\varepsilon,\tau} \geq 0$. We prove the claim $\eqref{ine:nc>0}$.

{\bf{ Step IV. The mass.}}
We claim that
\ben\label{eq:n L^infty_t L^1}
||n^{\varepsilon,\tau}(x,t)||_{L^1} = ||n_0^{\varepsilon}(x)||_{L^1},
\een
and
\ben\label{ine:c L^1 L^infty'}
||c^{\varepsilon,\tau}(x,t)||_{L^1 \cap L^\infty} \leq ||c_0^{\varepsilon}(x)||_{L^1 \cap L^\infty},
\een
for  any $ (x,t)\in\mathbb{R}^{3}\times(0,T).$

{\bf Proof.}
% For the equation of (\ref{eq:n L^infty_t L^1}) and (\ref{ine:c L^1 L^infty'}),
Recall the first equation of $\eqref{eq:CNS}$ as follow:
\beno
\int_{\mathbb{R}^3}\partial_t n^{\varepsilon,\tau} + \int_{\mathbb{R}^3}(u^{\varepsilon,\tau} \ast \rho^\varepsilon) \cdot \nabla n^{\varepsilon,\tau} - \int_{\mathbb{R}^3}\Delta n^{\varepsilon,\tau} + 
\int_{\mathbb{R}^3}\nabla \cdot \left(\frac{n^{\varepsilon,\tau}}{1+\tau n^{\varepsilon,\tau}} \nabla c^{\varepsilon,\tau} \chi(c^{\varepsilon,\tau})\right)=0.
%\int_{\mathbb{R}^3} \nabla \cdot (n^\varepsilon (\nabla c^\varepsilon\chi(c^{\varepsilon})) \ast \rho^\varepsilon)=0.
\eeno
Integration by parts and  divergence theorem yield that
\beno
\frac d{dt}\int_{\mathbb{R}^3}n^{\varepsilon,\tau}=0,
\eeno
which means for any $t \in (0,T)$,
\beno
\|n^{\varepsilon,\tau}(x,t)\|_{L^1}-\|n^{\varepsilon}(x,0)\|_{L^1}=0.
\eeno
Similarly, we have
\ben\label{ine:c L 1}
\|c^{\varepsilon,\tau}(x,t)\|_{L^1}-\|c^{\varepsilon}(x,0)\|_{L^1}\leq 0.
\een
For the $L^\infty$ norm of $ c^{\varepsilon,\tau} $, since $\kappa(s) > 0$, by the maximum principle, we know that
\beno
\|c^{\varepsilon,\tau}(x,t)\|_{L^\infty} \leq \|c_0^{\varepsilon}(x)\|_{L^\infty}.
\eeno
The proof of (\ref{eq:n L^infty_t L^1}) and  (\ref{ine:c L^1 L^infty'}) is complete.

{\bf Step V.{ Higher regularity and global existence.}}
Let $ E^{\varepsilon,\tau}=(n^{\varepsilon,\tau},c^{\varepsilon,\tau},\nabla c^{\varepsilon,\tau},u^{\varepsilon,\tau}) $, $ E_{0}^{\varepsilon}=(n_{0}^{\varepsilon},c_{0}^{\varepsilon},\nabla c_{0}^{\varepsilon}, u_{0}^{\varepsilon}) $. Then the solution is global 
if
the following uniform estimates hold
\begin{equation}\label{uniform estimation}
\begin{aligned}
&\|E^{\varepsilon,\tau}\|_{L^{\infty}_{t}H^{2}_{x}}^{2}+\int_{0}^{t}\|E^{\varepsilon,\tau}\|_{H^{3}}^{2}
%+\|\nabla^{3}c^{\varepsilon,\tau}\|_{L^{\infty}_{t}L^{2}_{x}}^{2}+\int_{0}^{t}\|\nabla^{4}c^{\varepsilon,\tau}\|_{L^{2}}^{2}
\\\leq~&C(\varepsilon,\tau,\|\nabla\phi\|_{L^{\infty}},\|\chi\|_0, \|\kappa\|_0,\|n_{0}^{\varepsilon}\|_{L^{1}},\|E_{0}^{\varepsilon}\|_{H^{2}}  )(T+1)^{54},
\end{aligned}
\end{equation}
%where 
%\beno
%\|\chi\|_0=\|\chi\|_{L^\infty(0, \|c_0^{\varepsilon}\|_{L^{\infty}})}+\|\chi'\|_{L^\infty(0, \|c_0^{\varepsilon}\|_{L^{\infty}})}+\|\chi''\|_{L^\infty(0, \|c_0^{\varepsilon}\|_{L^{\infty}})},
%\eeno
%and
%\beno
%\|\kappa\|_0=\|\kappa\|_{L^\infty(0, \|c^{\varepsilon}_0\|_{L^{\infty}})}+\|\kappa'\|_{L^\infty(0, \|c_0^{\varepsilon}\|_{L^{\infty}})}
%+\|\kappa''\|_{L^\infty(0, \|c_0^{\varepsilon}\|_{L^{\infty}})},
%\eeno
%for any $ t\in(0,T)$ .

Next, we will estimate the above norms separately.

{\bf\underline{i. The energy norms of $ c^{\varepsilon,\tau} :$}}

\begin{equation}\label{c-energy}
\begin{aligned}
\|c^{\varepsilon,\tau}\|_{L^{\infty}_{t}L^{2}_{x}}^{2}+\int_0^t\|\nabla c^{\varepsilon,\tau}\|_{L^{2}}^{2} \leq  \|c_0^{\varepsilon}\|_{L^2}^2,~~{\rm {for}}~~{\rm {any}}~ t\in(0,T).
\end{aligned}
\end{equation}

For the second equation of $\eqref{eq:CNS}$, multiplying it by $ c^{\varepsilon,\tau} $ and integrating with respect to space variable $x$ over $ \mathbb{R}^{3} $, we obtain
\beno
\frac12\frac{d}{dt}\|c^{\varepsilon,\tau}\|_{L^{2}}^{2}+\|\nabla c^{\varepsilon,\tau}\|_{L^{2}}^{2}&=&-\int_{\R^{3}} u^{\varepsilon,\tau}
\cdot\nabla c^{\varepsilon,\tau} c^{\varepsilon,\tau}-\int_{\R^{3}}\frac{\ln(1+\tau n^{\varepsilon,\tau})}{\tau} \kappa(c^{\varepsilon,\tau})c^{\varepsilon,\tau}.
%&\leq& \|J_k \kappa(c^{k,\varepsilon})\|_{L^\infty} \|c^{k,\varepsilon}\|_{L^2} \|n^{k,\varepsilon}\|_{L^2}.
\eeno
Since $\int_{\mathbb{R}^{3}} u^{\varepsilon,\tau} \cdot\nabla c^{\varepsilon,\tau} c^{\varepsilon,\tau} = 0$, $\kappa(s) \geq 0$, by $\eqref{ine:nc>0}$ we have
\beno
\frac{d}{dt}\|c^{\varepsilon,\tau}\|_{L^{2}}^{2}+\|\nabla c^{\varepsilon,\tau}\|_{L^{2}}^{2} \leq 0,
\eeno
which means (\ref{c-energy}) holds.
%\ben\label{ine:c c_0}
%\|c^{\varepsilon,\tau}(t)\|_{L^{2}}^{2}+\int_0^t\|\nabla c^{\varepsilon,\tau}\|_{L^{2}}^{2} \leq C \|c_0^{\varepsilon,\tau}\|_{L^2}^2\leq C\|E_{0}^{\varepsilon,\tau}\|_{H^{3}}^{2}.
%\een

{\bf\underline{ii. The energy norms of $ n^{\varepsilon,\tau}: $}}
\begin{equation}\label{n-energy}
\begin{aligned}
&\|n^{\varepsilon,\tau}(t)\|_{L^{\infty}_{t}L^2_{x}}^2 + \int_0^t \|\nabla n^{\varepsilon,\tau}\|_{L^2}^2\\\leq~&C\left( \tau, \|\chi\|_{0}, \|n_{0}^{\varepsilon}\|_{L^{2}},\|c_{0}^{\varepsilon}\|_{L^{2}}  \right), ~~{\rm {for}}~~{\rm {any}}~ t\in(0,T).
\end{aligned}
\end{equation}
For the first equation of $\eqref{eq:CNS}$, multiplying it by $n^{\varepsilon,\tau}$ and integrating with respect to space variable $x$ over $\mathbb{R}^{3}$, we obtain
\beno
\frac12 \frac d{dt} \|n^{\varepsilon,\tau}\|_{L^2}^2 + \|\nabla n^{\varepsilon,\tau}\|_{L^2}^2
% &=&  -\int_{\mathbb{R}^3} \nabla \cdot \left(\frac{n^{\varepsilon,\tau}}{1+\tau n^{\varepsilon,\tau}} \nabla c^{\varepsilon,\tau} \chi(c^{\varepsilon,\tau})\right)n^{\varepsilon,\tau}\\
&=& \int_{\mathbb{R}^3} \frac{n^{\varepsilon,\tau}}{1+\tau n^{\varepsilon,\tau}} \nabla c^{\varepsilon,\tau}\chi(c^{\varepsilon,\tau}) \cdot \nabla n^{\varepsilon,\tau} \\
%&\leq& \frac12 \|\nabla n^{\varepsilon,\tau}\|_{L^2}^2 + C(\tau)\|\nabla c^{\varepsilon,\tau}\chi(c^{\varepsilon,\tau}) \|_{L^2}^2  \\
&\leq& \frac12 \|\nabla n^{\varepsilon,\tau}\|_{L^2}^2 + C(\tau)  \|\nabla c^{\varepsilon,\tau}\|_{L^2}^2 \|\chi(c^{\varepsilon,\tau})\|_{L^{\infty}}^{2},
\eeno
where we use $ \frac{n^{\varepsilon,\tau}}{1+\tau n^{\varepsilon,\tau}}<\tau^{-1}. $
Noting that $\chi\in C^{1}(\overline{\mathbb{R}^+})$, we know that $\|\chi(c^{\varepsilon,\tau})\|_{L^{\infty}} \leq \|\chi\|_{L^\infty(0, \|c_0^{\varepsilon}\|_{L^{\infty}})} $, which yields
\begin{equation*}
\frac d{dt} \|n^{\varepsilon,\tau}\|_{L^2}^2 + \|\nabla n^{\varepsilon,\tau}\|_{L^2}^2 \leq C\left(\tau, \|\chi\|_{0}  \right) \|\nabla c^{\varepsilon,\tau}\|_{L^2}^2.
\end{equation*}
Then (\ref{n-energy}) holds by (\ref{c-energy}).

{\bf\underline{iii. The energy norms of $ u^{\varepsilon,\tau}: $}}
\begin{equation}\label{u-energy}
\begin{aligned}
&\|u^{\varepsilon,\tau}(t)\|_{L^{\infty}_{t}L^2_{x}}^2 + \int_{0}^{t}\|\nabla u^{\varepsilon,\tau}\|_{L^2}^2
%	\\\leq~& \|u_0^{\varepsilon,\tau}\|_{L^2}^2 +C\left(T_{*}, \|\nabla\phi\|_{L^{\infty}}, \|E_{0}^{\varepsilon,\tau}\|_{H^{3}}   \right)+ C \|\nabla \phi\|_{L^\infty}^2 \int_0^t \|n^{\varepsilon,\tau}\|_{L^2}^2\\\leq~&C\left(T_{*}, \|\nabla\phi\|_{L^{\infty}}, \|E_{0}^{\varepsilon,\tau}\|_{H^{3}}   \right)+CT_{*}\|\nabla\phi\|_{L^{\infty}}^{2}\|n^{\varepsilon,\tau}\|_{L^{\infty}((0,T_{*});L^{2})}\\
\\\leq~& 
C\left(\tau,\|\nabla\phi\|_{L^{\infty}},\|\chi\|_{0}, \|n_{0}^{\varepsilon}\|_{L^{1}}, \|E_{0}^{\varepsilon}\|_{L^2}  \right)(T+1),~~{\rm {for}}~~{\rm {any}}~ t\in(0,T).
%	C \|u_0^{\varepsilon,\tau}\|_{L^2}^2 + C T_\ast \|\nabla \phi\|_{L^\infty}^2 \exp\left(C(\tau) \chi(\|c_0^{\varepsilon,\tau}\|_{L^\infty})^2 \|c_0^{\varepsilon,\tau}\|_{L^2}^2\right) \|n_0^{\varepsilon,\tau}\|_{L^2}^2.
\end{aligned}
\end{equation}

Multiplying the third equation of $\eqref{eq:CNS}$ by $ u^{\varepsilon,\tau} $ and integrating with respect to space variable $x$ over $ \mathbb{R}^{3} $, we obtain
\beno
\frac12 \frac d{dt} \|u^{\varepsilon,\tau}\|_{L^2}^2 + \|\nabla u^{\varepsilon,\tau}\|_{L^2}^2 &=&~- \int_{\mathbb{R}^3} (n^{\varepsilon,\tau} \nabla \phi) \ast \rho^\varepsilon \cdot u^{\varepsilon,\tau} \\
%&\leq~ \|(n^{\varepsilon,\tau} \nabla \phi) \ast \rho^\varepsilon\|_{L^2} \|u^{\varepsilon,\tau}\|_{L^2} \\
%&\leq~ C \|\nabla \phi\|_{L^\infty} \|n^{\varepsilon,\tau}\|_{L^2} \|u^{\varepsilon,\tau}\|_{L^2} \\
&\leq&~ C\left(\|\nabla\phi\|_{L^{\infty}} \right)\|u^{\varepsilon,\tau}\|_{L^6}\|n^{\varepsilon,\tau}\|_{L^\frac{6}{5}}\\
&\leq&~\frac12  \|\nabla u^{\varepsilon,\tau}\|_{L^2}^2+  C\left(\|\nabla\phi\|_{L^{\infty}} \right)\|n^{\varepsilon,\tau}\|_{L^\frac{6}{5}}^2\\
&\leq&~\frac12  \|\nabla u^{\varepsilon,\tau}\|_{L^2}^2+  C\left(\|\nabla\phi\|_{L^{\infty}} \right)\|n^{\varepsilon,\tau}\|_{L^1}^{\frac43}\|n^{\varepsilon,\tau}\|_{L^2}^{\frac23}
\eeno
%Multiplying the third equation of $\eqref{eq:CNS}$ by $ c^{\varepsilon} $ and integrating with respect to space variable $x$ over $ \mathbb{R}^{3} $, we obtain
%\beno
%\frac12 \frac d{dt} \|u^\varepsilon\|_{L^2}^2 + \|\nabla u^\varepsilon\|_{L^2}^2 &=& - \int_{\mathbb{R}^3} (n^\varepsilon \nabla \phi) \ast \rho^\varepsilon \cdot u^\varepsilon \\
%&\leq& \|(n^\varepsilon \nabla \phi) \ast \rho^\varepsilon\|_{L^2} \|u^\varepsilon\|_{L^2} \\
%&\leq& C \|\nabla \phi\|_{L^\infty} \|n^\varepsilon\|_{L^2} \|u^\varepsilon\|_{L^2} \\
%&\leq& \|u^\varepsilon\|_{L^2}^2 + C \|\nabla \phi\|_{L^\infty}^2 \|n^\varepsilon\|_{L^2}^2.
%\eeno
Integrating with respect to $t$, by (\ref{n-energy}) we get (\ref{u-energy}).

%\ben\label{ine:u u_0}
%\|u^{\varepsilon,\tau}(t)\|_{L^2}^2 + \|\nabla u^{\varepsilon,\tau}\|_{L^2}^2 &\leq& C \|u_0^{\varepsilon,\tau}\|_{L^2}^2 + C \|\nabla \phi\|_{L^\infty}^2 \int_0^t \|n^{\varepsilon,\tau}\|_{L^2}^2 \\ \nonumber
%&\leq& C \|u_0^{\varepsilon,\tau}\|_{L^2}^2 + C T_\ast \|\nabla \phi\|_{L^\infty}^2 \exp\left(C(\tau) \chi(\|c_0^{\varepsilon,\tau}\|_{L^\infty})^2 \|c_0^{\varepsilon,\tau}\|_{L^2}^2\right) \|n_0^{\varepsilon,\tau}\|_{L^2}^2.
%\een
%Collecting $\eqref{ine:c c_0}$, $\eqref{ine:n n_0}$ and $\eqref{ine:u u_0}$, for any $t \in (0,T_\ast)$, we have
%\ben\label{ine:E E_0}
%\|E^{\varepsilon,\tau}(t)\|_{L^2}^2 + \int_0^t \|\nabla E^{\varepsilon,\tau}\|_{L^2}^2 \leq C(\tau,T_\ast,\|E_0^{\varepsilon,\tau}\|_{L^{2}},\|\nabla \phi\|_{L^\infty}).
%\een

{\bf\underline{iv. the energy norms of $ \nabla u^{\varepsilon,\tau}: $}}
\begin{equation}\label{u'-energy}
\begin{aligned}
&\| u^{\varepsilon,\tau}\|_{L^{\infty}_{t}H^{1}_{x}}^{2}+\int_{0}^{t}\|\nabla u^{\varepsilon,\tau}\|_{H^{1}}^{2}\\\leq~&C\left(\varepsilon,\tau,\|\nabla\phi\|_{L^{\infty}},\|\chi\|_{0},\|n_{0}^{\varepsilon}\|_{L^{1}}, \|E_{0}^{\varepsilon}\|_{H^1} \right)(T+1)^2,~~{\rm {for}}~~{\rm {any}}~ t\in(0,T).
\end{aligned}
\end{equation}
%\|\nabla u_{0}^{\varepsilon,\tau}\|_{L^{2}} 
Rewrite the third equation of (\ref{eq:CNS}) as follow:
\begin{equation*}
\partial_t\partial_{i}u^{\varepsilon,\tau}+\partial_{i}\left(\mu(u^{\varepsilon,\tau}\ast\rho^{\varepsilon})\cdot\nabla u^{\varepsilon,\tau}  \right)-\partial_{i}\Delta u^{\varepsilon,\tau}+\partial_{i}\nabla P^{\varepsilon,\tau}=-\partial_{i}\left((n^{\varepsilon,\tau}\nabla\phi)\ast\rho^{\varepsilon}  \right),
\end{equation*}
with $ i=1,2,3. $ Multiplying it by $\partial_{i} u^{\varepsilon,\tau}$ and integrating with respect to space variable $x$ over $\mathbb{R}^{3}$, we obtain
\begin{equation*}
\begin{aligned}
\frac12\frac{d}{dt}\|\nabla u^{\varepsilon,\tau}\|_{L^{2}}^{2}+\|\nabla^{2}u^{\varepsilon,\tau}\|_{L^{2}}^{2}=~&-\int_{\mathbb{R}^{3}}\partial_{i}\left(\mu(u^{\varepsilon,\tau}\ast\rho^{\varepsilon})\cdot\nabla u^{\varepsilon,\tau}  \right)\cdot\partial_{i}u^{\varepsilon,\tau}\\&-~\int_{\mathbb{R}^{3}}\partial_{i}\left((n^{\varepsilon,\tau}\nabla\phi) \ast\rho^{\varepsilon} \right)\cdot\partial_{i}u^{\varepsilon,\tau}\\:=~&A_{1}+A_{2}.
\end{aligned}
\end{equation*}
%By the product derivative rule, we have
%\begin{equation*}
%	\begin{aligned}
%	A_{1}=~&-\int_{\mathbb{R}^{3}}\partial_{i}(u^{\varepsilon,\tau}\ast\rho^{\varepsilon})\cdot\nabla u^{\varepsilon,\tau}\cdot\partial_{i}u^{\varepsilon,\tau}-\int_{\mathbb{R}^{3}}(u^{\varepsilon,\tau}\ast\rho^{\varepsilon})\cdot\nabla\partial_{i}u^{\varepsilon,\tau}\cdot\partial_{i}u^{\varepsilon,\tau}.
%	\end{aligned}
%\end{equation*}
Integration by parts, H\"{o}lder's inequality and Young's inequality yield that
\begin{equation*}
\begin{aligned}
A_{1}=~&\mu\int_{\mathbb{R}^{3}}\left((u^{\varepsilon,\tau}\ast\rho^{\varepsilon})\cdot\nabla u^{\varepsilon,\tau} \right)\Delta u^{\varepsilon,\tau}\\
\leq~&\|u^{\varepsilon,\tau}\ast\rho^{\varepsilon}\|_{L^\infty}\|\nabla u^{\varepsilon,\tau}\|_{L^2}\|\Delta u^{\varepsilon,\tau}\|_{L^2}\\
%\leq~&\frac{1}{16}\|\nabla^{2}u^{\varepsilon,\tau}\|_{L^{2}}^{2}+C\|\nabla u^{\varepsilon,\tau}\|_{L^{2}}^{2}\|\rho^{\varepsilon}\|_{L^{2}}^{2}\|u^{\varepsilon,\tau}\|_{L^{2}}^{2}\\
\leq~&\frac{1}{16}\|\nabla^{2}u^{\varepsilon,\tau}\|_{L^{2}}^{2}+C(\varepsilon)\|\nabla u^{\varepsilon,\tau}\|_{L^{2}}^{2}\|u^{\varepsilon,\tau}\|_{L^{2}}^{2}.
\end{aligned}
\end{equation*}
%Note that $ \int_{\mathbb{R}^{3}}(u^{\varepsilon,\tau}\ast\rho^{\varepsilon})\cdot\nabla\partial_{i}u^{\varepsilon,\tau}\cdot\partial_{i}u^{\varepsilon,\tau}=0 $ due to $ \nabla\cdot(u^{\varepsilon,\tau}\ast\rho^{\varepsilon})=0. $ It follows that
%\begin{equation*}
%	\begin{aligned}
%	A_{1}=~&-\int_{\mathbb{R}^{3}}\partial_{i}(u_{j}^{\varepsilon,\tau}\ast\rho^{\varepsilon})\partial_{j} u_{k}^{\varepsilon,\tau}\partial_{i}u_{k}^{\varepsilon,\tau}\\=~&\int_{\mathbb{R}^{3}}\partial_{i}(u_{j}^{\varepsilon,\tau}\ast\rho^{\varepsilon})u_{k}^{\varepsilon,\tau}\partial_{ij}u_{k}^{\varepsilon,\tau}\\\leq~&\frac{1}{16}\|\nabla^{2}u^{\varepsilon,\tau}\|_{L^{2}}^{2}+C\|\partial_{i}(u^{\varepsilon,\tau}\ast\rho^{\varepsilon})\|_{L^{\infty}}^{2}\|u^{\varepsilon,\tau}\|_{L^{2}}^{2}\\\leq~&\frac{1}{16}\|\nabla^{2}u^{\varepsilon,\tau}\|_{L^{2}}^{2}+C\|\nabla u^{\varepsilon,\tau}\|_{L^{2}}^{2}\|\rho^{\varepsilon}\|_{L^{2}}^{2}\|u^{\varepsilon,\tau}\|_{L^{2}}^{2}\\\leq~&\frac{1}{16}\|\nabla^{2}u^{\varepsilon,\tau}\|_{L^{2}}^{2}+C(\varepsilon)\|\nabla u^{\varepsilon,\tau}\|_{L^{2}}^{2}\|u^{\varepsilon,\tau}\|_{L^{2}}^{2}.
%	\end{aligned}
%\end{equation*}
Similarly, for $ A_{2}, $ by integration by parts, H\"{o}lder's inequality and Young's inequality, we have
\begin{equation*}
\begin{aligned}
A_{2}%\leq~&\|\partial_{i}\left((n^{\varepsilon,\tau}\nabla\phi)\ast\rho^{\varepsilon}  \right)\|_{L^{2}}\|\nabla u^{\varepsilon,\tau}\|_{L^{2}}\\
%\leq~&\|n^{\varepsilon,\tau}\nabla\phi\|_{L^{2}}\|\rho^{\varepsilon}\|_{L^{1}}\|\nabla^2 u^{\varepsilon,\tau}\|_{L^{2}}\\\leq~&C\|\nabla\phi\|_{L^{\infty}}\|n^{\varepsilon,\tau}\|_{L^{2}}\|\nabla^2 u^{\varepsilon,\tau}\|_{L^{2}}\\
\leq~&\frac12\|\nabla^2 u^{\varepsilon,\tau}\|_{L^{2}}^{2}+C\left(\|\nabla\phi\|_{L^{\infty}}  \right)\|n^{\varepsilon,\tau}\|_{L^{2}}^{2}.
\end{aligned}
\end{equation*}
Combining $ A_{1} $ and $ A_{2}, $ there holds
\begin{equation*}
\begin{aligned}
&	\frac{d}{dt}\|\nabla u^{\varepsilon,\tau}\|_{L^{2}}^{2}+\|\nabla^{2}u^{\varepsilon,\tau}\|_{L^{2}}^{2}\\\leq~& C(\varepsilon) \|\nabla u^{\varepsilon,\tau}\|_{L^{2}}^{2}\|u^{\varepsilon,\tau}\|_{L^{2}}^{2}+C\left(\|\nabla\phi\|_{L^{\infty}} \right)\|n^{\varepsilon,\tau}\|_{L^{2}}^{2}.
\end{aligned}
\end{equation*}
%\begin{equation*}
%	\begin{aligned}
%&	\frac{d}{dt}\|\nabla u^{\varepsilon,\tau}\|_{L^{2}}^{2}+\|\nabla^{2}u^{\varepsilon,\tau}\|_{L^{2}}^{2}\\\leq~& C(\varepsilon)\left(1+\|u^{\varepsilon,\tau}\|_{L^{2}}^{2}  \right)\|\nabla u^{\varepsilon,\tau}\|_{L^{2}}^{2}+C\left(\varepsilon,\|\nabla\phi\|_{L^{\infty}} \right)\|n^{\varepsilon,\tau}\|_{L^{2}}^{2}.
%	\end{aligned}
%\end{equation*}
By  (\ref{n-energy}) and (\ref{u-energy}), we obtain (\ref{u'-energy}).

{\bf\underline{v. the energy norms of $ \nabla^{2}u^{\varepsilon,\tau}: $}}
\begin{equation}\label{u''-energy}
\begin{aligned}
&\| u^{\varepsilon,\tau}\|_{L^{\infty}_{t}H^{2}_{x}}^{2}+\int_{0}^{t}\|\nabla u^{\varepsilon,\tau}\|_{H^{2}}^{2}\\\leq~&C\left(\varepsilon,\tau,\|\nabla\phi\|_{L^{\infty}},\|\chi\|_{0},\|n_{0}^{\varepsilon}\|_{L^{1}}, \|E_{0}^{\varepsilon}\|_{H^{2}} \right)(T+1)^3,~~{\rm {for}}~~{\rm {any}}~ t\in(0,T).
\end{aligned}
\end{equation}
%\|\nabla u_{0}^{\varepsilon,\tau}\|_{L^2}
We write the third equation of $\eqref{eq:CNS}$ as follow:
\beno
\partial_t \partial_{ij} u^{\varepsilon,\tau}+ \mu\partial_{ij} ((u^{\varepsilon,\tau}\ast\rho^{\varepsilon})\cdot \nabla u^{\varepsilon,\tau}) +  \partial_{ij} \nabla P^{\varepsilon,\tau} = \partial_{ij} \Delta u^{\varepsilon,\tau} - \partial_{ij}((n^{\varepsilon,\tau}\nabla\phi)\ast\rho^{\varepsilon}).
\eeno
Multiplying above equation by $\partial_{ij} u^{\varepsilon,\tau} $, integrating with respect to space variable $ x $ over $ \mathbb{R}^{3} $ and noting that $\nabla \cdot u^{\varepsilon,\tau} = 0$, we obtain
\beno
\frac12\frac d{dt} \|\nabla^2 u^{\varepsilon,\tau}\|_{L^2}^2 + \| \nabla^3 u^{\varepsilon,\tau}\|_{L^2}^2 &=& - \mu\int_{\mathbb{R}^3} \partial_{ij} ((u^{\varepsilon,\tau}\ast\rho^{\varepsilon})\cdot\nabla u^{\varepsilon,\tau})\cdot \partial_{ij} u^{\varepsilon,\tau} \\
&&- \int_{\mathbb{R}^3} \partial_{ij}((n^{\varepsilon,\tau}\nabla\phi)\ast\rho^{\varepsilon})\cdot \partial_{ij} u^{\varepsilon,\tau}\\
&:=& B_1 + B_2.
\eeno
By the product derivative rule, we have
\beno
B_1 &=& -\mu \int_{\mathbb{R}^3}\partial_{ij}(u^{\varepsilon,\tau}\ast\rho^{\varepsilon})\cdot\nabla u^{\varepsilon,\tau} \cdot\partial_{ij}u^{\varepsilon,\tau}dx -\mu \int_{\mathbb{R}^3} (u^{\varepsilon,\tau}\ast\rho^{\varepsilon})\cdot\nabla\partial_{ij}u^{\varepsilon,\tau}\cdot\partial_{ij}u^{\varepsilon,\tau}dx \\
&&-\mu\int_{\mathbb{R}^3}\partial_{j}(u^{\varepsilon,\tau}\ast\rho^{\varepsilon})\cdot\nabla\partial_i u^{\varepsilon,\tau} \cdot\partial_{ij}u^{\varepsilon,\tau}dx-\mu \int_{\mathbb{R}^3}\partial_{i}(u^{\varepsilon,\tau}\ast\rho^{\varepsilon})\cdot\nabla\partial_ju^{\varepsilon,\tau} \cdot\partial_{ij}u^{\varepsilon,\tau}dx \\
&:=& B_{11} + B_{12} + B_{13} + B_{14}.
\eeno
For $B_{11}$, noting that
\beno
\|\partial_{ij}(u^{\varepsilon,\tau}\ast\rho^{\varepsilon})\|_{L^\infty} = \|(\partial_{ij}u^{\varepsilon,\tau})\ast\rho^{\varepsilon}\|_{L^\infty} \leq C \|\nabla^{2}u^{\varepsilon,\tau}\|_{L^2} \|\rho^{\varepsilon}\|_{L^2} \leq C(\varepsilon) \|\nabla^{2}u^{\varepsilon,\tau}\|_{L^2},
\eeno
integration by parts yields that
\beno
B_{11} &=& \mu\int_{\mathbb{R}^3}u^{\varepsilon,\tau} \partial_{ij}(u^{\varepsilon,\tau}\ast\rho^{\varepsilon})\cdot\nabla\partial_{ij}u^{\varepsilon,\tau}dx \\
&\leq& \frac{1}{16} \|\nabla^3u^{\varepsilon,\tau}\|_{L^{2}}^{2} + C \|u^{\varepsilon,\tau}\|_{L^2}^2 \|\partial_{ij}(u^{\varepsilon,\tau}\ast\rho^{\varepsilon})\|_{L^\infty}^2\\
&\leq& \frac{1}{16} \|\nabla^3u^{\varepsilon,\tau}\|_{L^{2}}^{2} + C(\varepsilon) \|u^{\varepsilon,\tau}\|_{L^2}^2 \|\nabla^{2}u^{\varepsilon,\tau}\|_{L^2}^{2}.
\eeno
For $B_{12}$, obviously
\[B_{12} = 0.\]
For $B_{13}$ and $B_{14}$, specially, we calculate $B_{13}$ in detail.
\beno
B_{13} &=\mu& \int_{\mathbb{R}^3} \partial_i u^{\varepsilon,\tau} \partial_{j}(u^{\varepsilon,\tau}\ast\rho^{\varepsilon}) \cdot\nabla\partial_{ij} u^{\varepsilon,\tau} dx \\
&\leq& \frac{1}{16} \|\nabla^3u^{\varepsilon,\tau}\|_{L^{2}}^{2} + C(\varepsilon) \|\nabla u^{\varepsilon,\tau}\|_{L^2}^2 \|u^{\varepsilon,\tau}\|_{L^2}^2.
\eeno
Collecting $B_{11}$ to $B_{14}$, we have
\beno
B_1 \leq \frac{3}{16} \|\nabla^3 u^{\varepsilon,\tau}\|_{L^{2}}^{2} + C(\varepsilon) \|u^{\varepsilon,\tau}\|_{L^2}^2 \|\nabla^{2}u^{\varepsilon,\tau}\|_{L^2}^{2}+ C(\varepsilon) \|\nabla u^{\varepsilon,\tau}\|_{L^2}^2\|u^{\varepsilon,\tau}\|_{L^2}^2.
\eeno
For the term $B_2$, using H\"{o}lder's inequality, we have
\begin{equation*}
\begin{aligned}
B_2 &\leq~\|\nabla^2 ((n^{\varepsilon,\tau}\nabla\phi)\ast\rho^{\varepsilon})\|_{L^2} \|\nabla^2 u^{\varepsilon,\tau}\|_{L^2}\\&\leq~C(\varepsilon)\|n^{\varepsilon,\tau}\nabla\phi\|_{L^2}  \|\nabla^2 u^{\varepsilon,\tau}\|_{L^2}\\&
\leq~C(\varepsilon)\|\nabla\phi\|_{L^{\infty}}^{2}\|\nabla^{2}u^{\varepsilon,\tau}\|_{L^{2}}^{2}+\frac12\|n^{\varepsilon,\tau}\|_{L^{2}}^{2}.
\end{aligned}
\end{equation*}
Collecting $B_1$ and $B_2$, we have
\begin{equation*}
\begin{aligned}
&\frac{d}{dt}\|\nabla^{2}u^{\varepsilon,\tau}\|_{L^{2}}^{2}+\|\nabla^{3}u^{\varepsilon,\tau}\|_{L^{2}}^{2}\\\leq~&C(\varepsilon)\left(\|\nabla\phi\|_{L^{\infty}}^{2}+\|u^{\varepsilon,\tau}\|_{L^{2}}^{2}  \right)\|\nabla^{2}u^{\varepsilon,\tau}\|_{L^{2}}^{2}+C(\varepsilon)\|\nabla u^{\varepsilon,\tau}\|_{L^{2}}^{2}\|u^{\varepsilon,\tau}\|_{L^{2}}^{2}+\frac12\|n^{\varepsilon,\tau}\|_{L^{2}}^{2}.
\end{aligned}
\end{equation*}
By  (\ref{n-energy}) and (\ref{u'-energy}), we obtain (\ref{u''-energy}).

{\bf\underline{vi. the energy norms of $ \nabla c^{\varepsilon,\tau}: $}}
\begin{equation}\label{c'-energy}
\begin{aligned}
&\| c^{\varepsilon,\tau}\|_{L^{\infty}_{t}H^{1}_{x}}^{2}+\int_{0}^{t}\|\nabla c^{\varepsilon,\tau}\|_{H^{1}}^{2}\\\leq~&C\left(\tau,\|\nabla\phi\|_{L^{\infty}},\|\chi\|_{0}, \|\kappa\|_{0}, \|n_{0}^{\varepsilon}\|_{L^{1}}, \|E_{0}^{\varepsilon}\|_{H^1}  \right)(T+1)^{\frac52},~~{\rm {for}}~~{\rm {any}}~ t\in(0,T).
\end{aligned}
\end{equation}

Rewrite the second equation of (\ref{eq:CNS}) as follow:
\begin{equation*}
\partial_t\partial_{i}c^{\varepsilon,\tau}+\partial_{i}\left(u^{\varepsilon,\tau}\cdot\nabla c^{\varepsilon,\tau}  \right)-\partial_{i}\Delta c^{\varepsilon,\tau}=-\partial_{i}\left( \frac{\ln(1+\tau n^{\varepsilon,\tau})}{\tau}\kappa(c^{\varepsilon,\tau}) \right)
\end{equation*}
with $ i=1,2,3. $ Multiplying it by $\partial_{i} c^{\varepsilon,\tau}$ and integrating with respect to space variable $x$ over $\mathbb{R}^{3}$, we obtain
\begin{equation*}
\begin{aligned}
\frac12\frac{d}{dt}\|\nabla c^{\varepsilon,\tau}\|_{L^{2}}^{2}+\|\nabla^{2}c^{\varepsilon,\tau}\|_{L^{2}}^{2}=~&-\int_{\mathbb{R}^{3}}\partial_{i}\left(u^{\varepsilon,\tau}\cdot\nabla c^{\varepsilon,\tau}  \right)\partial_{i} c^{\varepsilon,\tau}\\&-~\int_{\mathbb{R}^{3}}\partial_{i}\left(\frac{\ln(1+\tau n^{\varepsilon,\tau})}{\tau}\kappa(c^{\varepsilon,\tau})  \right)\partial_{i}c^{\varepsilon,\tau}\\:=~&C_{1}+C_{2}.
\end{aligned}
\end{equation*}
The term $ C_{1} $ is similar as $ A_{1}, $ we omit the proof and arrive
\begin{equation*}
\begin{aligned}
C_{1}\leq~\frac{1}{16}\|\nabla^{2}c^{\varepsilon,\tau}\|_{L^{2}}^{2}+C\|\nabla c^{\varepsilon,\tau}\|_{L^{2}}^{2}\|u^{\varepsilon,\tau}\|_{L^{\infty}}^{2}.
\end{aligned}
\end{equation*}
By the product derivative rule for $ C_{2}, $ we have
\begin{equation*}
\begin{aligned}
C_{2}=~&-\int_{\mathbb{R}^{3}}\frac{\partial_{i}n^{\varepsilon,\tau}}{1+\tau n^{\varepsilon,\tau}}\kappa(c^{\varepsilon,\tau})\partial_{i}c^{\varepsilon,\tau}-\int_{\mathbb{R}^{3}}\frac{\ln(1+\tau n^{\varepsilon,\tau})}{\tau}\kappa'(c^{\varepsilon,\tau})\partial_{i}c^{\varepsilon,\tau}\partial_{i}c^{\varepsilon,\tau}\\:=~&C_{21}+C_{22}.
\end{aligned}
\end{equation*}
For $ C_{21}, $ using H\"{o}lder's inequality and Young's inequality, we get
\begin{equation*}
\begin{aligned}
C_{21}\leq~\|\nabla n^{\varepsilon,\tau}\|_{L^{2}}\|\kappa(c^{\varepsilon,\tau})\|_{L^{\infty}}\|\nabla c^{\varepsilon,\tau}\|_{L^{2}}%\leq& \|\nabla n^{\varepsilon,\tau}\|_{L^{2}}\kappa(\|c_{0}^{\varepsilon}\|_{L^{\infty}})\|\nabla c^{\varepsilon,\tau}\|_{L^{2}}\\
%\leq~\frac12\|\nabla c^{\varepsilon,\tau}\|_{L^{2}}^{2}+C\left(\|\kappa\|_{0}
% \right)\|\nabla n^{\varepsilon,\tau}\|_{L^{2}}^{2}.
\end{aligned}
\end{equation*}
For $ C_{22}, $ note that
\begin{equation*}
\|\nabla c^{\varepsilon,\tau}\|_{L^{4}}\leq C\|c^{\varepsilon,\tau}\|_{L^{2}}^{\frac18}\|\nabla^{2}c^{\varepsilon,\tau}\|_{L^{2}}^{\frac78},
\end{equation*}
then
\begin{equation*}
\begin{aligned}
C_{22}\leq~&\|\frac{\ln(1+\tau n^{\varepsilon,\tau})}{\tau}\|_{L^{2}}\|\kappa'(c^{\varepsilon,\tau})\|_{L^{\infty}}\||\nabla c^{\varepsilon,\tau}|^{2}\|_{L^{2}}\\\leq~&C(\|\kappa\|_{0}) \|n^{\varepsilon,\tau}\|_{L^{2}}\|c^{\varepsilon,\tau}\|_{L^{2}}^{\frac14}\|\nabla^{2}c^{\varepsilon,\tau}\|_{L^{2}}^{\frac74}\\\leq~&\frac{1}{16}\|\nabla^{2}c^{\varepsilon,\tau}\|_{L^{2}}^{2}+C\left( \|\kappa\|_{0} \right)\|n^{\varepsilon,\tau}\|_{L^{2}}^{8}\|c^{\varepsilon,\tau}\|_{L^{2}}^{2}.
\end{aligned}
\end{equation*}
Collecting $ C_{21} $ and $ C_{22} $ yields that
\begin{equation*}
\begin{aligned}
C_{2}\leq~&\frac{1}{16}\|\nabla^{2}c^{\varepsilon,\tau}\|_{L^{2}}^{2}+\|\nabla n^{\varepsilon,\tau}\|_{L^{2}}\|\kappa(c^{\varepsilon,\tau})\|_{L^{\infty}}\|\nabla c^{\varepsilon,\tau}\|_{L^{2}}\\&+~C\left( \|\kappa\|_{0} \right)\|n^{\varepsilon,\tau}\|_{L^{2}}^{8}\|c^{\varepsilon,\tau}\|_{L^{2}}^{2}.
\end{aligned}
\end{equation*}
Then combining $ C_{1} $ and $ C_{2}, $ we have
\begin{equation*}
\begin{aligned}
&\frac{d}{dt}\|\nabla c^{\varepsilon,\tau}\|_{L^{2}}^{2}+\|\nabla^{2}c^{\varepsilon,\tau}\|_{L^{2}}^{2}\\\leq~&\|\nabla n^{\varepsilon,\tau}\|_{L^{2}}\|\kappa(c^{\varepsilon,\tau})\|_{L^{\infty}}\|\nabla c^{\varepsilon,\tau}\|_{L^{2}}+C\|\nabla c^{\varepsilon,\tau}\|_{L^{2}}^{2}\|u^{\varepsilon,\tau}\|_{L^{\infty}}^{2}\\&+C\left( \|\kappa\|_{0} \right)\|n^{\varepsilon,\tau}\|_{L^{2}}^{8}\|c^{\varepsilon,\tau}\|_{L^{2}}^{2}.
\end{aligned}
\end{equation*}
Note that $ \|u^{\varepsilon,\tau}\|_{L^{\infty}}\leq C\|u^{\varepsilon,\tau}\|_{L^{2}}^{\frac14}\|\nabla^{2}u^{\varepsilon,\tau}\|_{L^{2}}^{\frac34}, $ by  (\ref{n-energy}), (\ref{u-energy}),(\ref{u''-energy}) and (\ref{c-energy}), we obtain (\ref{c'-energy}).
%\begin{equation}\label{c'-energy}
%	\begin{aligned}
%	&\|\nabla c^{\varepsilon,\tau}\|_{L^{2}}^{2}+\int_{0}^{t}\|\nabla^{2}c^{\varepsilon,\tau}\|_{L^{2}}^{2}\\\leq~&C\left(\varepsilon,\tau,T_{*},\|\nabla\phi\|_{L^{\infty}},\chi(\|c_{0}^{\varepsilon,\tau}\|_{L^{\infty}}), \kappa(\|c_{0}^{\varepsilon,\tau}\|_{L^{\infty}}),\kappa'(\|c_{0}^{\varepsilon,\tau}\|_{L^{\infty}}),\|E_{0}^{\varepsilon,\tau}\|_{H^{3}}  \right). 
%	\end{aligned}
%\end{equation}

{\bf\underline{vii. the energy norms of $ |\nabla c^{\varepsilon,\tau}|^{\frac32}: $}}
\begin{equation}\label{c'' half energy}
\begin{aligned}
&\|\nabla c^{\varepsilon,\tau}\|_{L^{\infty}_{t}L^{3}_{x}}^{3}+\int_{0}^{t}\|\nabla(|\nabla c^{\varepsilon,\tau}|^{\frac32})\|_{L^{2}}^{2}\\\leq~&C\left(\tau,\|\nabla\phi\|_{L^{\infty}}, \|\chi\|_{0},\|\kappa\|_{0},\|n_{0}^{\varepsilon}\|_{L^{1}}, \|E_{0}^{\varepsilon}\|_{H^1}  \right)(T+1)^{\frac{21}{4}},
~~{\rm {for}}~~{\rm {any}}~ t\in(0,T).
\end{aligned}
\end{equation}

Rewrite the second equation of (\ref{eq:CNS}) as follow:
\begin{equation*}
\partial_t\partial_{i}c^{\varepsilon,\tau}+\partial_{i}\left(u^{\varepsilon,\tau}\cdot\nabla c^{\varepsilon,\tau}  \right)-\partial_{i}\Delta c^{\varepsilon,\tau}=-\partial_{i}\left( \frac{\ln(1+\tau n^{\varepsilon,\tau})}{\tau}\kappa(c^{\varepsilon,\tau}) \right)
\end{equation*}
with $ i=1,2,3. $ Multiplying it by $|\nabla c^{\varepsilon,\tau}|\partial_{i}c^{\varepsilon,\tau}$ and integrating with respect to space variable $x$ over $\mathbb{R}^{3}$, we obtain
\begin{equation*}
\begin{aligned}
\frac13\frac{d}{dt}\|\nabla c^{\varepsilon,\tau}\|_{L^{3}}^{3}+\frac89\|\nabla(|\nabla c^{\varepsilon,\tau}|^{\frac32})\|_{L^{2}}^{2}\leq~&\left|\int_{\mathbb{R}^{3}}\partial_{i}\left(u^{\varepsilon,\tau}\cdot\nabla c^{\varepsilon,\tau}  \right)|\nabla c^{\varepsilon,\tau}|\partial_{i}c^{\varepsilon,\tau}dx\right|\\&+~\left|\int_{\mathbb{R}^{3}}\partial_{i}\left(\frac{\ln(1+\tau n^{\varepsilon,\tau})}{\tau}\kappa(c^{\varepsilon,\tau})  \right)|\nabla c^{\varepsilon,\tau}|\partial_{i}c^{\varepsilon,\tau}dx\right|\\:=~&|D_{1}|+|D_{2}|.
\end{aligned}
\end{equation*}
For $ D_{1}, $ we have
\begin{equation*}
\begin{aligned}
D_{1}\leq~C\|\nabla u^{\varepsilon,\tau}\|_{L^6}\|\nabla c^{\varepsilon,\tau}\|_{L^{\frac{18}{5}}}^3\leq~C\|\nabla^{2}u^{\varepsilon,\tau}\|_{L^{2}}\|\nabla c^{\varepsilon,\tau}\|_{L^{\frac{18}{5}}}^3.
%~C(\varepsilon)\|u^{\varepsilon,\tau}\|_{L^2}\|\nabla c^{\varepsilon,\tau}\|_{L^3}^3.
\end{aligned}
\end{equation*}
By the product derivative rule for $ D_{2}, $ we obtain
\begin{equation*}
\begin{aligned}
D_{2}&=~\int_{\mathbb{R}^{3}}\partial_{i}\left(\frac{\ln(1+\tau n^{\varepsilon,\tau})}{\tau}\kappa(c^{\varepsilon,\tau})  \right)|\nabla c^{\varepsilon,\tau}|\partial_{i}c^{\varepsilon,\tau}\\&=~\int_{\mathbb{R}^{3}}\frac{\partial_{i}n^{\varepsilon,\tau}}{1+\tau n^{\varepsilon,\tau}}\kappa(c^{\varepsilon,\tau})|\nabla c^{\varepsilon,\tau}|\partial_{i}c^{\varepsilon,\tau}+\int_{\mathbb{R}^{3}}\frac{\ln(1+\tau n^{\varepsilon,\tau})}{\tau}\kappa'(c^{\varepsilon,\tau})\partial_{i}c^{\varepsilon,\tau}|\nabla c^{\varepsilon,\tau}|\partial_{i}c^{\varepsilon,\tau}\\:&=~D_{21}+D_{22}.
\end{aligned}
\end{equation*}
For $ D_{21}, $ by interpolation inequality $\|\nabla c^{\varepsilon,\tau}\|^{\frac32}_{L^4}\leq C\|\nabla c^{\varepsilon,\tau}\|^{\frac{15}{16}}_{L^3}\|\nabla (|\nabla c^{\varepsilon,\tau}|^{\frac32})\|^{\frac38}_{L^2}$ ,H\"{o}lder's inequality and Young's inequality, we get
\begin{equation*}
\begin{aligned}
D_{21}%&\leq~C\left(\kappa(\|c_{0}^{\varepsilon}\|_{L^{\infty}})\right)\int_{\mathbb{R}^{3}}|\partial_{i} n^{\varepsilon,\tau}||\partial_{i}c^{\varepsilon,\tau}|^{2}\\
&\leq~C\left(\|\kappa\|_{0}\right)\|\nabla n^{\varepsilon,\tau}\|_{L^{2}}\|\nabla c^{\varepsilon,\tau}\|_{L^{4}}^{2}\\&\leq~C\left(\|\kappa\|_{0}\right)\|\nabla n^{\varepsilon,\tau}\|_{L^{2}}\|\nabla c^{\varepsilon,\tau}\|^{\frac54}_{L^3}\|\nabla (|\nabla c^{\varepsilon,\tau}|^{\frac32})\|^{\frac12}_{L^2}
%\\&\leq~\frac{1}{16}\|\nabla (|\nabla c^{\varepsilon,\tau}|^{\frac32})\|^{2}_{L^2}+C\left(\|\kappa\|_{0}\right)\|\nabla n^{\varepsilon,\tau}\|_{L^{2}}^{\frac43}\|\nabla c^{\varepsilon,\tau}\|^{\frac53}_{L^3}
%\\&\leq~\frac{1}{16}\|\nabla |\nabla c^{\varepsilon,\tau}|^{\frac32}\|^{2}_{L^2}+C\left(\kappa(\|c_{0}^{\varepsilon}\|_{L^{\infty}})\right)\|\nabla n^{\varepsilon,\tau}\|_{L^{2}}^{\frac43}(\frac12\|\nabla c^{\varepsilon,\tau}\|^{3}_{L^3}+1^{\frac94})
\\&\leq~\frac{1}{16}\|\nabla (|\nabla c^{\varepsilon,\tau}|^{\frac32})\|^{2}_{L^2}+C\left(\|\kappa\|_{0}\right)\|\nabla n^{\varepsilon,\tau}\|_{L^{2}}^{\frac43}\|\nabla c^{\varepsilon,\tau}\|^{\frac53}_{L^3}
%\\&\leq~\frac{1}{16}\|\nabla |\nabla c^{\varepsilon,\tau}|^{\frac32}\|^{2}_{L^2}+C\left(\kappa(\|c_{0}^{\varepsilon}\|_{L^{\infty}})\right)\|\nabla n^{\varepsilon,\tau}\|_{L^{2}}^{\frac43}\|\nabla c^{\varepsilon,\tau}\|^{3}_{L^3}+C\left(\kappa(\|c_{0}^{\varepsilon}\|_{L^{\infty}})\right)\|\nabla n^{\varepsilon,\tau}\|_{L^{2}}^{\frac43}.
%\\&\leq~C\left(\varepsilon_1,\kappa(\|c_{0}^{\varepsilon,\tau}\|_{L^{\infty}})\right)\|\nabla n^{\varepsilon,\tau}\|_{L^{2}}^{2}+\varepsilon_1\|\nabla c^{\varepsilon,\tau}\|_{L^{4}}^{4}.
\end{aligned}
\end{equation*}
%\begin{equation*}
%	\begin{aligned}
%	D_{21}&\leq~C\left(\kappa(\|c_{0}^{\varepsilon,\tau}\|_{L^{\infty}})\right)\int_{\mathbb{R}^{3}}|\partial_{i} n^{\varepsilon,\tau}||\partial_{i}c^{\varepsilon,\tau}|^{2}\\&\leq~C\left(\kappa(\|c_{0}^{\varepsilon,\tau}\|_{L^{\infty}})\right)\|\nabla n^{\varepsilon,\tau}\|_{L^{2}}\|\nabla c^{\varepsilon,\tau}\|_{L^{4}}^{2} \\&\leq~C\left(\varepsilon_1,\kappa(\|c_{0}^{\varepsilon,\tau}\|_{L^{\infty}})\right)\|\nabla n^{\varepsilon,\tau}\|_{L^{2}}^{2}+\varepsilon_1\|\nabla c^{\varepsilon,\tau}\|_{L^{4}}^{4}.
%	\end{aligned}
%\end{equation*}
For $ D_{22}, $ similar to $ D_{21} $, by interpolation inequality $\|\nabla c^{\varepsilon,\tau}\|_{L^{6}}\leq C\|\nabla c^{\varepsilon,\tau}\|_{L^{3}}^{\frac14}\|\nabla (|\nabla c^{\varepsilon,\tau}|^{\frac32})\|^{\frac12}_{L^2}$, we obtain
\begin{equation*}
\begin{aligned}
D_{22}&\leq~ C\left(\|\kappa\|_{0}\right)\|n^{\varepsilon,\tau}\|_{L^{2}}\|\nabla c^{\varepsilon,\tau}\|_{L^{6}}^{3}\\&\leq~C\left(\|\kappa\|_{0}\right)\|n^{\varepsilon,\tau}\|_{L^{2}}\|\nabla c^{\varepsilon,\tau}\|_{L^{3}}^{\frac34}\|\nabla (|\nabla c^{\varepsilon,\tau}|^{\frac32})\|^{\frac32}_{L^2}\\&\leq~\frac{1}{16}\|\nabla (|\nabla c^{\varepsilon,\tau}|^{\frac32})\|^{2}_{L^2}+C\left(\|\kappa\|_{0}\right)\|n^{\varepsilon,\tau}\|^4_{L^{2}}\|\nabla c^{\varepsilon,\tau}\|_{L^{3}}^{3}.
\end{aligned}
\end{equation*}
%C\left(\varepsilon_2, \kappa'(\|c_{0}^{\varepsilon,\tau}\|_{L^{\infty}})\right)\|n^{\varepsilon,\tau}\|_{L^{2}}^{4}+\varepsilon_2\|\nabla c^{\varepsilon,\tau}\|_{L^{6}}^{4}.
Collecting $ D_{21} $ and $ D_{22}, $ we have
\begin{equation*}
\begin{aligned}
D_{2}\leq~&\frac18\|\nabla (|\nabla c^{\varepsilon,\tau}|^{\frac32})\|^{2}_{L^2}+C\left(\|\kappa\|_{0}\right)\|\nabla n^{\varepsilon,\tau}\|_{L^{2}}^{\frac43}\|\nabla c^{\varepsilon,\tau}\|^{\frac53}_{L^3}\\&+C\left(\|\kappa\|_{0}\right)\|n^{\varepsilon,\tau}\|^4_{L^{2}}\|\nabla c^{\varepsilon,\tau}\|_{L^{3}}^{3}.
%C\left(\varepsilon_2, \kappa'(\|c_{0}^{\varepsilon,\tau}\|_{L^{\infty}})\right)\|n^{\varepsilon,\tau}\|_{L^{2}}^{4}+\varepsilon_2\|\nabla c^{\varepsilon,\tau}\|_{L^{6}}^{4}.
\end{aligned}
\end{equation*}
%C\left(\varepsilon_1,\kappa(\|c_{0}^{\varepsilon,\tau}\|_{L^{\infty}})\right)\|\nabla n^{\varepsilon,\tau}\|_{L^{2}}^{2}+C\left(\varepsilon_2, \kappa'(\|c_{0}^{\varepsilon,\tau}\|_{L^{\infty}})\right)\|n^{\varepsilon,\tau}\|_{L^{2}}^{4}\\&+~\varepsilon_1\|\nabla c^{\varepsilon,\tau}\|_{L^{4}}^{4}+\varepsilon_2\|\nabla c^{\varepsilon,\tau}\|_{L^{6}}^{4}.  
Combining with $ D_{1} $ and $ D_{2}, $ there holds
\begin{equation*}
\begin{aligned}
&\frac{d}{dt}\|\nabla c^{\varepsilon,\tau}\|_{L^{3}}^{3}+\|\nabla(|\nabla c^{\varepsilon,\tau}|^{\frac32})\|_{L^{2}}^{2}\\
%\leq~&C(\varepsilon)\|u^{\varepsilon,\tau}\|_{L^2}\|\nabla c^{\varepsilon,\tau}\|_{L^3}^3+C\left(\kappa(\|c_{0}^{\varepsilon}\|_{L^{\infty}})\right)\|\nabla n^{\varepsilon,\tau}\|_{L^{2}}^{\frac43}\|\nabla c^{\varepsilon,\tau}\|^{3}_{L^3}\\&+~C\left(\kappa(\|c_{0}^{\varepsilon}\|_{L^{\infty}})\right)\|\nabla n^{\varepsilon,\tau}\|_{L^{2}}^{\frac43}+C\left(\kappa'(\|c_{0}^{\varepsilon}\|_{L^{\infty}})\right)\|n^{\varepsilon,\tau}\|^4_{L^{2}}\|\nabla c^{\varepsilon,\tau}\|_{L^{3}}^{3}\\
\leq~&C\|\nabla^{2}u^{\varepsilon,\tau}\|_{L^{2}}\|\nabla c^{\varepsilon,\tau}\|_{L^{\frac{18}{5}}}^3+C\left(\|\kappa\|_{0}\right)\|n^{\varepsilon,\tau}\|^4_{L^{2}}\|\nabla c^{\varepsilon,\tau}\|^{\frac32}_{L^2}\|\nabla^2 c^{\varepsilon,\tau}\|^{\frac32}_{L^2}
\\
&+C\left(\|\kappa\|_{0}\right)\|\nabla n^{\varepsilon,\tau}\|_{L^{2}}^{\frac43}\|\nabla c^{\varepsilon,\tau}\|^{\frac53}_{L^3}\\
\leq~&C\|\nabla^{2}u^{\varepsilon,\tau}\|_{L^{2}}\|\nabla c^{\varepsilon,\tau}\|_{L^{\frac{18}{5}}}^3+C\left(\|\kappa\|_{0}\right)\|n^{\varepsilon,\tau}\|^4_{L^{2}}\|\nabla c^{\varepsilon,\tau}\|^{\frac32}_{L^2}\|\nabla^2 c^{\varepsilon,\tau}\|^{\frac32}_{L^2}\\
&+C\left(\|\kappa\|_{0}\right)\|\nabla n^{\varepsilon,\tau}\|_{L^{2}}^{\frac43}\|\nabla c^{\varepsilon,\tau}\|^{\frac{20}{21}}_{L^2}\|\nabla(|\nabla c^{\varepsilon,\tau}|^{\frac32})\|_{L^{2}}^{\frac{10}{21}}.
\end{aligned}
\end{equation*}
By  (\ref{c-energy}), (\ref{n-energy}), (\ref{c'-energy}) and (\ref{u-energy}), we have (\ref{c'' half energy}).

{\bf\underline{viii. the energy norms of $ |\nabla c^{\varepsilon,\tau}|^{\frac52}: $}}
\begin{equation}\label{c'5 half energy}
\begin{aligned}
&\|\nabla c^{\varepsilon,\tau}\|_{L^{\infty}_{t}L^{5}_{x}}^{5}+\int_{0}^{t}\|\nabla(|\nabla c^{\varepsilon,\tau}|^{\frac52})\|_{L^{2}}^{2}\\\leq~&C\left(\tau,\|\nabla\phi\|_{L^{\infty}}, \|\chi\|_{0},\|\kappa\|_{0},\|n_{0}^{\varepsilon}\|_{L^{1}}, \|E_{0}^{\varepsilon}\|_{H^2}  \right)(T+1)^{\frac{41}{4}},
~~{\rm {for}}~~{\rm {any}}~ t\in(0,T).
\end{aligned}
\end{equation}
 
Rewrite the second equation of (\ref{eq:CNS}) as follow:
\begin{equation*}
\partial_t\partial_{i}c^{\varepsilon,\tau}+\partial_{i}\left(u^{\varepsilon,\tau}\cdot\nabla c^{\varepsilon,\tau}  \right)-\partial_{i}\Delta c^{\varepsilon,\tau}=-\partial_{i}\left( \frac{\ln(1+\tau n^{\varepsilon,\tau})}{\tau}\kappa(c^{\varepsilon,\tau}) \right)
\end{equation*}
with $ i=1,2,3. $ Multiplying it by $|\nabla c^{\varepsilon,\tau}|^3\partial_{i}c^{\varepsilon,\tau}$ and integrating with respect to space variable $x$ over $\mathbb{R}^{3}$, we obtain
\begin{equation*}
\begin{aligned}
\frac15\frac{d}{dt}\|\nabla c^{\varepsilon,\tau}\|_{L^{5}}^{5}+\frac{16}{25}\|\nabla(|\nabla c^{\varepsilon,\tau}|^{\frac52})\|_{L^{2}}^{2}\leq~&\left|\int_{\mathbb{R}^{3}}\partial_{i}\left(u^{\varepsilon,\tau}\cdot\nabla c^{\varepsilon,\tau}  \right)|\nabla c^{\varepsilon,\tau}|^3\partial_{i}c^{\varepsilon,\tau}dx\right|\\&+~\left|\int_{\mathbb{R}^{3}}\partial_{i}\left(\frac{\ln(1+\tau n^{\varepsilon,\tau})}{\tau}\kappa(c^{\varepsilon,\tau})  \right)|\nabla c^{\varepsilon,\tau}|^3\partial_{i}c^{\varepsilon,\tau}dx\right|\\:=~&|D_{1}'|+|D_{2}'|.
\end{aligned}
\end{equation*}
For $ D_{1}', $ note that $ \|\nabla c^{\varepsilon,\tau}\|_{L^{6}}^{5}\leq C\|\nabla c^{\varepsilon,\tau}\|_{L^{5}}^{\frac{15}{4}}\|\nabla|\nabla c^{\varepsilon,\tau}|^{\frac52}\|_{L^{2}}^{\frac12}, $ 
we have
\begin{equation*}
\begin{aligned}
D_{1}'\leq~C\|\nabla u^{\varepsilon,\tau}\|_{L^6}\|\nabla c^{\varepsilon,\tau}\|_{L^6}^5
\leq~C\|\nabla^{2} u^{\varepsilon,\tau}\|_{L^{2}}\|\nabla c^{\varepsilon,\tau}\|_{L^{5}}^{\frac{15}{4}}\|\nabla|\nabla c^{\varepsilon,\tau}|^{\frac52}\|_{L^{2}}^{\frac12}.
%~C(\varepsilon)\|u^{\varepsilon,\tau}\|_{L^2}\|\nabla c^{\varepsilon,\tau}\|_{L^5}^5.
\end{aligned}
\end{equation*}
By the product derivative rule for $ D_{2}', $ we obtain
\begin{equation*}
\begin{aligned}
D_{2}'&=~\int_{\mathbb{R}^{3}}\partial_{i}\left(\frac{\ln(1+\tau n^{\varepsilon,\tau})}{\tau}\kappa(c^{\varepsilon,\tau})  \right)|\nabla c^{\varepsilon,\tau}|^3\partial_{i}c^{\varepsilon,\tau}\\&=~\int_{\mathbb{R}^{3}}\frac{\partial_{i}n^{\varepsilon,\tau}}{1+\tau n^{\varepsilon,\tau}}\kappa(c^{\varepsilon,\tau})|\nabla c^{\varepsilon,\tau}|^3\partial_{i}c^{\varepsilon,\tau}+\int_{\mathbb{R}^{3}}\frac{\ln(1+\tau n^{\varepsilon,\tau})}{\tau}\kappa'(c^{\varepsilon,\tau})\partial_{i}c^{\varepsilon,\tau}|\nabla c^{\varepsilon,\tau}|^3\partial_{i}c^{\varepsilon,\tau}\\:&=~D_{21}'+D_{22}'.
\end{aligned}
\end{equation*}
For $ D_{21}', $  we get
\begin{equation*}
\begin{aligned}
D_{21}'%&\leq~C\left(\kappa(\|c_{0}^{\varepsilon}\|_{L^{\infty}})\right)\int_{\mathbb{R}^{3}}|\partial_{i} n^{\varepsilon,\tau}||\partial_{i}c^{\varepsilon,\tau}|^{2}\\
&\leq~C\left(\|\kappa\|_{0}\right)\|\nabla n^{\varepsilon,\tau}\|_{L^{2}}\|\nabla c^{\varepsilon,\tau}\|_{L^{8}}^{4}\\&\leq~C\left(\|\kappa\|_{0}\right)\|\nabla n^{\varepsilon,\tau}\|_{L^{2}}\|\nabla c^{\varepsilon,\tau}\|^{\frac74}_{L^5}\|\nabla (|\nabla c^{\varepsilon,\tau}|^{\frac52})\|^{\frac{9}{10}}_{L^2}
%\\&\leq~\frac{1}{16}\|\nabla (|\nabla c^{\varepsilon,\tau}|^{\frac32})\|^{2}_{L^2}+C\left(\|\kappa\|_{0}\right)\|\nabla n^{\varepsilon,\tau}\|_{L^{2}}^{\frac43}\|\nabla c^{\varepsilon,\tau}\|^{\frac53}_{L^3}
%\\&\leq~\frac{1}{16}\|\nabla |\nabla c^{\varepsilon,\tau}|^{\frac32}\|^{2}_{L^2}+C\left(\kappa(\|c_{0}^{\varepsilon}\|_{L^{\infty}})\right)\|\nabla n^{\varepsilon,\tau}\|_{L^{2}}^{\frac43}(\frac12\|\nabla c^{\varepsilon,\tau}\|^{3}_{L^3}+1^{\frac94})
\\&\leq~\frac{1}{16}\|\nabla (|\nabla c^{\varepsilon,\tau}|^{\frac52})\|^{2}_{L^2}+C\left(\|\kappa\|_{0}\right)\|\nabla n^{\varepsilon,\tau}\|_{L^{2}}^{\frac{20}{11}}\|\nabla c^{\varepsilon,\tau}\|^{\frac{35}{11}}_{L^5}
%\\&\leq~\frac{1}{16}\|\nabla |\nabla c^{\varepsilon,\tau}|^{\frac32}\|^{2}_{L^2}+C\left(\kappa(\|c_{0}^{\varepsilon}\|_{L^{\infty}})\right)\|\nabla n^{\varepsilon,\tau}\|_{L^{2}}^{\frac43}\|\nabla c^{\varepsilon,\tau}\|^{3}_{L^3}+C\left(\kappa(\|c_{0}^{\varepsilon}\|_{L^{\infty}})\right)\|\nabla n^{\varepsilon,\tau}\|_{L^{2}}^{\frac43}.
%\\&\leq~C\left(\varepsilon_1,\kappa(\|c_{0}^{\varepsilon,\tau}\|_{L^{\infty}})\right)\|\nabla n^{\varepsilon,\tau}\|_{L^{2}}^{2}+\varepsilon_1\|\nabla c^{\varepsilon,\tau}\|_{L^{4}}^{4}.
\end{aligned}
\end{equation*}
%\begin{equation*}
%	\begin{aligned}
%	D_{21}&\leq~C\left(\kappa(\|c_{0}^{\varepsilon,\tau}\|_{L^{\infty}})\right)\int_{\mathbb{R}^{3}}|\partial_{i} n^{\varepsilon,\tau}||\partial_{i}c^{\varepsilon,\tau}|^{2}\\&\leq~C\left(\kappa(\|c_{0}^{\varepsilon,\tau}\|_{L^{\infty}})\right)\|\nabla n^{\varepsilon,\tau}\|_{L^{2}}\|\nabla c^{\varepsilon,\tau}\|_{L^{4}}^{2} \\&\leq~C\left(\varepsilon_1,\kappa(\|c_{0}^{\varepsilon,\tau}\|_{L^{\infty}})\right)\|\nabla n^{\varepsilon,\tau}\|_{L^{2}}^{2}+\varepsilon_1\|\nabla c^{\varepsilon,\tau}\|_{L^{4}}^{4}.
%	\end{aligned}
%\end{equation*}
For $ D_{22}'$, we obtain
\begin{equation*}
\begin{aligned}
D_{22}'&\leq~ C\left(\|\kappa\|_{0}\right)\|n^{\varepsilon,\tau}\|_{L^{2}}\|\nabla c^{\varepsilon,\tau}\|_{L^{10}}^{5}\\&\leq~C\left(\|\kappa\|_{0}\right)\|n^{\varepsilon,\tau}\|_{L^{2}}\|\nabla c^{\varepsilon,\tau}\|_{L^{5}}^{\frac54}\|\nabla (|\nabla c^{\varepsilon,\tau}|^{\frac52})\|^{\frac32}_{L^2}\\&\leq~\frac{1}{16}\|\nabla (|\nabla c^{\varepsilon,\tau}|^{\frac52})\|^{2}_{L^2}+C\left(\|\kappa\|_{0}\right)\|n^{\varepsilon,\tau}\|^4_{L^{2}}\|\nabla c^{\varepsilon,\tau}\|_{L^{5}}^{5}.
\end{aligned}
\end{equation*}
%C\left(\varepsilon_2, \kappa'(\|c_{0}^{\varepsilon,\tau}\|_{L^{\infty}})\right)\|n^{\varepsilon,\tau}\|_{L^{2}}^{4}+\varepsilon_2\|\nabla c^{\varepsilon,\tau}\|_{L^{6}}^{4}.
Collecting $ D_{21}' $ and $ D_{22}', $ we have
\begin{equation*}
\begin{aligned}
D_{2}'\leq~&\frac18\|\nabla (|\nabla c^{\varepsilon,\tau}|^{\frac52})\|^{2}_{L^2}+C\left(\|\kappa\|_{0}\right)\|\nabla n^{\varepsilon,\tau}\|_{L^{2}}^{\frac{20}{11}}\|\nabla c^{\varepsilon,\tau}\|^{\frac{35}{11}}_{L^5}\\&+C\left(\|\kappa\|_{0}\right)\|n^{\varepsilon,\tau}\|^4_{L^{2}}\|\nabla c^{\varepsilon,\tau}\|_{L^{5}}^{5}.
%C\left(\varepsilon_2, \kappa'(\|c_{0}^{\varepsilon,\tau}\|_{L^{\infty}})\right)\|n^{\varepsilon,\tau}\|_{L^{2}}^{4}+\varepsilon_2\|\nabla c^{\varepsilon,\tau}\|_{L^{6}}^{4}.
\end{aligned}
\end{equation*}
%C\left(\varepsilon_1,\kappa(\|c_{0}^{\varepsilon,\tau}\|_{L^{\infty}})\right)\|\nabla n^{\varepsilon,\tau}\|_{L^{2}}^{2}+C\left(\varepsilon_2, \kappa'(\|c_{0}^{\varepsilon,\tau}\|_{L^{\infty}})\right)\|n^{\varepsilon,\tau}\|_{L^{2}}^{4}\\&+~\varepsilon_1\|\nabla c^{\varepsilon,\tau}\|_{L^{4}}^{4}+\varepsilon_2\|\nabla c^{\varepsilon,\tau}\|_{L^{6}}^{4}.  
Combining with $ D_{1}' $ and $ D_{2}', $ there holds
\begin{equation*}
\begin{aligned}
&\frac{d}{dt}\|\nabla c^{\varepsilon,\tau}\|_{L^{5}}^{5}+\|\nabla(|\nabla c^{\varepsilon,\tau}|^{\frac52})\|_{L^{2}}^{2}\\
%\leq~&C(\varepsilon)\|u^{\varepsilon,\tau}\|_{L^2}\|\nabla c^{\varepsilon,\tau}\|_{L^3}^3+C\left(\kappa(\|c_{0}^{\varepsilon}\|_{L^{\infty}})\right)\|\nabla n^{\varepsilon,\tau}\|_{L^{2}}^{\frac43}\|\nabla c^{\varepsilon,\tau}\|^{3}_{L^3}\\&+~C\left(\kappa(\|c_{0}^{\varepsilon}\|_{L^{\infty}})\right)\|\nabla n^{\varepsilon,\tau}\|_{L^{2}}^{\frac43}+C\left(\kappa'(\|c_{0}^{\varepsilon}\|_{L^{\infty}})\right)\|n^{\varepsilon,\tau}\|^4_{L^{2}}\|\nabla c^{\varepsilon,\tau}\|_{L^{3}}^{3}\\
\leq~&C\|\nabla^{2}u^{\varepsilon,\tau}\|_{L^{2}}\|\nabla c^{\varepsilon,\tau}\|_{L^{5}}^{\frac{15}{4}}\|\nabla|\nabla c^{\varepsilon,\tau}|^{\frac52}\|_{L^{2}}^{\frac12}
+C\left(\|\kappa\|_{0}\right)\|\nabla n^{\varepsilon,\tau}\|_{L^{2}}^{\frac{20}{11}}\|\nabla c^{\varepsilon,\tau}\|^{\frac{35}{11}}_{L^5}\\&+C\left(\|\kappa\|_{0}\right)\|n^{\varepsilon,\tau}\|^4_{L^{2}}\|\nabla c^{\varepsilon,\tau}\|_{L^{5}}^{5}.
\end{aligned}
\end{equation*}
Integrating with respect to $t$, we get
\beno
&&\sup_t\|\nabla c^{\varepsilon,\tau}\|_{L^{5}}^{5}+\int_0^t\|\nabla(|\nabla c^{\varepsilon,\tau}|^{\frac52})\|_{L^{2}}^{2}ds\\
&\leq~&C\sup_t\|\nabla^{2}u^{\varepsilon,\tau}\|_{L^2}\int_0^t\|\nabla c^{\varepsilon,\tau}\|_{L^5}^5ds
\\
&&+C\left(\|\kappa\|_{0}\right)\left(\int_0^t\|\nabla n^{\varepsilon,\tau}\|_{L^{2}}^2ds\right)^{\frac{10}{11}}T^{\frac{1}{11}}\sup_t\|\nabla c^{\varepsilon,\tau}\|^{\frac{35}{11}}_{L^5}\\&&+C\left(\|\kappa\|_{0}\right)\sup_t\|n^{\varepsilon,\tau}\|^4_{L^{2}}\int_0^t\|\nabla c^{\varepsilon,\tau}\|_{L^{5}}^{5}ds.
\eeno

By  (\ref{n-energy}), (\ref{u-energy}) and (\ref{c'' half energy}), the proof of (\ref{c'5 half energy}) is complete.

{\bf\underline{ix. the energy norms of $ \nabla n^{\varepsilon,\tau}: $}}
\begin{equation}\label{n'-energy}
\begin{aligned}
&\| n^{\varepsilon,\tau}\|_{L_{t}^{\infty}H^{1}_{x}}^{2}+\int_{0}^{t}\|\nabla n^{\varepsilon,\tau}\|_{H^{1}}^{2}\\\leq~&C\left(\varepsilon,\tau,\|\nabla\phi\|_{L^{\infty}}, \|\chi\|_{0},\|\kappa\|_{0},\|n_{0}^{\varepsilon}\|_{L^{1}}, \|E_{0}^{\varepsilon}\|_{H^2}  \right)(T+1)^{\frac{41}{4}},
%~~{\rm{for}}~~{\rm {any}}~t\in(0,T).
\end{aligned}
\end{equation}
for any $ t\in(0,T). $

Rewrite the first equation of (\ref{eq:CNS}) as follow:
\beno
\partial_t \partial_{i}n^{\varepsilon,\tau}+\partial_{i} ((u^{\varepsilon,\tau}\ast\rho^{\varepsilon})\cdot\nabla n^{\varepsilon,\tau})= \partial_{i} \Delta n^{\varepsilon,\tau}- \partial_{i}\nabla\cdot\left(
\frac{n^{\varepsilon,\tau}}{1+\tau n^{\varepsilon,\tau}}\nabla c^{\varepsilon,\tau}\chi(c^{\varepsilon,\tau})
\right),
\eeno
with $i = 1,2,3$. Multiplying it by $\partial_{i} n^{\varepsilon,\tau}$ and integrating with respect to space variable $x$ over $\mathbb{R}^{3}$, we obtain
\beno
\frac12 \frac d{dt} \|\nabla n^{\varepsilon,\tau}\|_{L^2}^2 + \|\nabla^2 n^{\varepsilon,\tau}\|_{L^2}^2 &=& -\int_{\mathbb{R}^3} \partial_{i} ((u^{\varepsilon,\tau}\ast\rho^{\varepsilon})\cdot\nabla n^{\varepsilon,\tau}) \partial_{i} n^{\varepsilon,\tau} \\
&&- \int_{\mathbb{R}^3} \partial_{i}\nabla\cdot\left(
\frac{n^{\varepsilon,\tau}}{1+\tau n^{\varepsilon,\tau}}\nabla c^{\varepsilon,\tau}\chi(c^{\varepsilon,\tau})
\right) \partial_{i} n^{\varepsilon,\tau}\\&:=& E_{1}+E_{2}.
\eeno
The term $ E_{1} $ is same as $ A_{1}, $ we arrive
\begin{equation*}
E_{1}\leq~\frac{1}{16}\|\nabla^{2}n^{\varepsilon,\tau}\|_{L^{2}}^{2}+C(\varepsilon)\| u^{\varepsilon,\tau}\|_{L^{2}}^{2}\|\nabla n^{\varepsilon,\tau}\|_{L^{2}}^{2}.
\end{equation*}
By the product derivative rule for $ E_{2}, $ we have
\begin{equation*}
\begin{aligned}
E_{2}%=~&\int_{\mathbb{R}^{3}}\partial_{i}\left( \frac{n^{\varepsilon,\tau}}{1+\tau n^{\varepsilon,\tau}}\nabla c^{\varepsilon,\tau}\chi(c^{\varepsilon,\tau}) \right)\cdot\nabla\partial_{i}n^{\varepsilon,\tau}\\
=~&\int_{\mathbb{R}^{3}}\frac{\partial_{i}n^{\varepsilon,\tau}}{(1+\tau n^{\varepsilon,\tau})^{2}}\nabla c^{\varepsilon,\tau}\chi(c^{\varepsilon,\tau})\cdot\nabla\partial_{i}n^{\varepsilon,\tau}\\+~&\int_{\mathbb{R}^{3}}\frac{n^{\varepsilon,\tau}}{1+\tau n^{\varepsilon,\tau}}\partial_{i}\nabla c^{\varepsilon,\tau}\chi(c^{\varepsilon,\tau})\cdot\nabla\partial_{i}n^{\varepsilon,\tau}\\+~&\int_{\mathbb{R}^{3}}\frac{n^{\varepsilon,\tau}}{1+\tau n^{\varepsilon,\tau}}\nabla c^{\varepsilon,\tau}\chi'(c^{\varepsilon,\tau})\partial_{i}c^{\varepsilon,\tau}\cdot\nabla\partial_{i}n^{\varepsilon,\tau}\\:=~&E_{21}+E_{22}+E_{23}.
\end{aligned}
\end{equation*}
For $ E_{21}, $ by $ \|\nabla n^{\varepsilon,\tau}\|_{L^{\frac{10}{3}}}\leq C\|\nabla n^{\varepsilon,\tau}\|_{L^2}^{\frac25} \|\nabla^2 n^{\varepsilon,\tau}\|_{L^2}^{\frac35}$,
H\"{o}lder's inequality and Young's inequality, we get
\begin{equation*}
\begin{aligned}
E_{21}\leq~& \|\nabla n^{\varepsilon,\tau}\|_{L^{\frac{10}{3}}}\|\nabla c^{\varepsilon,\tau}\|_{L^{5}}\|\chi(c^{\varepsilon,\tau})\|_{L^{\infty}}\|\nabla^{2}n^{\varepsilon,\tau}\|_{L^{2}}\\\leq~&\|\nabla n^{\varepsilon,\tau}\|_{L^2}^{\frac25} \|\nabla^2 n^{\varepsilon,\tau}\|_{L^2}^{\frac35}\|\nabla c^{\varepsilon,\tau}\|_{L^{5}}\|\chi(c^{\varepsilon,\tau})\|_{L^{\infty}}\|\nabla^{2}n^{\varepsilon,\tau}\|_{L^{2}}
\\\leq~&\frac{1}{16}\|\nabla^{2}n^{\varepsilon,\tau}\|_{L^{2}}^{2}+C\left(\|\chi\|_{0}  \right)\|\nabla n^{\varepsilon,\tau}\|_{L^2}^2\|\nabla c^{\varepsilon,\tau}\|_{L^{5}}^{5}.
\end{aligned}
\end{equation*}
For $ E_{22}, $ note that $ \frac{n^{\varepsilon,\tau}}{1+\tau n^{\varepsilon,\tau}}<\tau^{-1}, $ we obtain
\begin{equation*}
\begin{aligned}
E_{22}\leq~&\frac{1}{16}\|\nabla^{2}n^{\varepsilon,\tau}\|_{L^{2}}^{2}+C\left(\tau, \|\chi\|_{0} \right)\|\nabla^{2}c^{\varepsilon,\tau}\|_{L^{2}}^{2}.
\end{aligned}
\end{equation*}
Similarly, for $ E_{23} $, we get
\begin{equation*}
\begin{aligned}
E_{23}\leq~&C\left(\tau, \|\chi\|_{0}  \right)\|\nabla c^{\varepsilon,\tau}\|_{L^{4}}^{2}\|\nabla^{2}n^{\varepsilon,\tau}\|_{L^{2}}\\\leq~&\frac{1}{16}\|\nabla^{2}n^{\varepsilon,\tau}\|_{L^{2}}^{2}+C\left(\tau, \|\chi\|_{0}  \right)\|\nabla c^{\varepsilon,\tau}\|_{L^{4}}^{4}.
%\\\leq~&\frac{1}{16}\|\nabla^{2}n^{\varepsilon,\tau}\|_{L^{2}}^{2}+\frac12\|\nabla c^{\varepsilon,\tau}\|_{L^{4}}^{8}+C\left(\tau, \chi'(\|c_{0}^{\varepsilon}\|_{L^{\infty}})  \right).
\end{aligned}
\end{equation*}
Collecting $ E_{21}-E_{23}, $ we have
\begin{equation*}
\begin{aligned}
E_{2}\leq~&\frac{3}{16}\|\nabla^{2}n^{\varepsilon,\tau}\|_{L^{2}}^{2}+C\left(\|\chi\|_{0} \right)\|\nabla n^{\varepsilon,\tau}\|_{L^2}^2\|\nabla c^{\varepsilon,\tau}\|_{L^{5}}^{5}\\&+~C\left(\tau,\|\chi\|_{0} \right)\|\nabla^{2}c^{\varepsilon,\tau}\|_{L^{2}}^{2}+C\left(\tau,\|\chi\|_{0} \right)\|\nabla c^{\varepsilon,\tau}\|_{L^{4}}^{4}.
\end{aligned}
\end{equation*}
Combining $ E_{1} $ with $ E_{2} $, we can yield that
\begin{equation*}
\begin{aligned}
&\frac{d}{dt}\|\nabla n^{\varepsilon,\tau}\|_{L^{2}}^{2}+\|\nabla^{2}n^{\varepsilon,\tau}\|_{L^{2}}^{2}\\\leq~&
C(\varepsilon)\| u^{\varepsilon,\tau}\|_{L^{2}}^{2}\|\nabla n^{\varepsilon,\tau}\|_{L^{2}}^{2}
+C\left(\|\chi\|_{0} \right)\|\nabla n^{\varepsilon,\tau}\|_{L^2}^2\|\nabla c^{\varepsilon,\tau}\|_{L^{5}}^{5}\\+&~C\left(\tau,\|\chi\|_{0} \right)\|\nabla^{2}c^{\varepsilon,\tau}\|_{L^{2}}^{2}+C\left(\tau,\|\chi\|_{0} \right)\|\nabla c^{\varepsilon,\tau}\|_{L^{4}}^{4}.
\end{aligned}
\end{equation*} 
Noting $ \|\nabla c^{\varepsilon,\tau}\|_{L^{4}}\leq C\|\nabla c^{\varepsilon,\tau}\|_{L^2}^{\frac16} \|\nabla c^{\varepsilon,\tau}\|_{L^5}^{\frac56}$, by  (\ref{c'' half energy}), (\ref{c'-energy}), (\ref{u-energy}) and (\ref{n-energy}),
we can get (\ref{n'-energy}).

{\bf\underline{x. the energy norms of $ \nabla^{2}c^{\varepsilon,\tau}: $}}
\begin{equation}\label{c''-energy}
\begin{aligned}
&\|c^{\varepsilon,\tau}\|_{L^{\infty}_{t}H^{2}_{x}}^{2}+\int_{0}^{t}\|\nabla c^{\varepsilon,\tau}\|_{H^{2}}^{2}\\\leq~&C\left(\varepsilon,\tau,\|\nabla\phi\|_{L^{\infty}}, \|\chi\|_{0},\|\kappa\|_{0},\|n_{0}^{\varepsilon}\|_{L^{1}}, \|E_{0}^{\varepsilon}\|_{H^2}  \right)(T+1)^{\frac{11}{2}},
%~~{\rm{for}}~~{\rm {any}}~t\in(0,T).
\end{aligned}
\end{equation}
for any $ t\in(0,T). $

Rewrite the second equation of (\ref{eq:CNS}) as follow:
\begin{equation*}
\partial_t\partial_{ij}c^{\varepsilon,\tau}+\partial_{ij}\left(u^{\varepsilon,\tau}\cdot\nabla c^{\varepsilon,\tau}  \right)-\partial_{ij}\Delta c^{\varepsilon,\tau}=-\partial_{ij}\left( \frac{\ln(1+\tau n^{\varepsilon,\tau})}{\tau}\kappa(c^{\varepsilon,\tau}) \right)
\end{equation*}
with $ i,j =1,2,3. $ Multiplying it by $\partial_{ij} c^{\varepsilon,\tau}$ and integrating with respect to space variable $x$ over $\mathbb{R}^{3}$, we obtain
\begin{equation*}
\begin{aligned}
\frac12\frac{d}{dt}\|\nabla^{2} c^{\varepsilon,\tau}\|_{L^{2}}^{2}+\|\nabla^{3}c^{\varepsilon,\tau}\|_{L^{2}}^{2}=~&-\int_{\mathbb{R}^{3}}\partial_{ij}\left(u^{\varepsilon,\tau}\cdot\nabla c^{\varepsilon,\tau}  \right)\partial_{ij} c^{\varepsilon,\tau}\\&-~\int_{\mathbb{R}^{3}}\partial_{ij}\left(\frac{\ln(1+\tau n^{\varepsilon,\tau})}{\tau}\kappa(c^{\varepsilon,\tau})  \right)\partial_{ij}c^{\varepsilon,\tau}\\:=~&F_{1}+F_{2}.
\end{aligned}
\end{equation*}
The term $ F_{1} $ is similar as $ B_{1}, $ we arrive
\begin{equation*}
\begin{aligned}
F_{1}\leq~\frac{3}{16}\|\nabla^{3}c^{\varepsilon,\tau}\|_{L^{2}}^{2}+C\|u^{\varepsilon,\tau}\|_{L^{\infty}}^{2}\|\nabla^{2}c^{\varepsilon,\tau}\|_{L^{2}}^{2}+C\|\nabla c^{\varepsilon,\tau}\|_{L^{3}}^{2}\|\nabla^{2} u^{\varepsilon,\tau}\|_{L^{2}}^{2}.
\end{aligned}
\end{equation*}
For the term $ F_{2} $, integration by parts  yields that
\begin{equation*}
\begin{aligned}
F_{2}=~&\int_{\mathbb{R}^{3}}\frac{\ln(1+\tau n^{\varepsilon,\tau})}{\tau}\kappa^{'}(c^{\varepsilon,\tau})\partial_i c^{\varepsilon,\tau}\partial_{ijj}c^{\varepsilon,\tau}+\int_{\mathbb{R}^{3}}\frac{\partial_i n^{\varepsilon,\tau}}{1+\tau n^{\varepsilon,\tau}}\kappa(c^{\varepsilon,\tau})\partial_{ijj}c^{\varepsilon,\tau}\\
\leq~&\frac18\|\nabla^{3}c^{\varepsilon,\tau}\|_{L^2}^2+C(\|\kappa\|_{0})\|\nabla n^{\varepsilon,\tau}\|_{L^2}^2+C(\|\kappa\|_{0})\|\nabla c^{\varepsilon,\tau}\|_{L^5}^2\|n^{\varepsilon,\tau}\|_{L^\frac{10}3}^2.
\end{aligned}
\end{equation*}
Combining $ F_{1} $ with $ F_{2}$, there holds
\begin{equation*}
\begin{aligned}
&\frac12\frac{d}{dt}\|\nabla^{2} c^{\varepsilon,\tau}\|_{L^{2}}^{2}+\|\nabla^{3}c^{\varepsilon,\tau}\|_{L^{2}}^{2}\\\leq~&C\|u^{\varepsilon,\tau}\|_{L^{\infty}}^{2}\|\nabla^{2}c^{\varepsilon,\tau}\|_{L^{2}}^{2}+C\|\nabla c^{\varepsilon,\tau}\|_{L^{3}}^{2}\|\nabla^{2} u^{\varepsilon,\tau}\|_{L^{2}}^{2}.\\&+~C(\|\kappa\|_{0})\|\nabla n^{\varepsilon,\tau}\|_{L^2}^2+C(\|\kappa\|_{0})\|\nabla c^{\varepsilon,\tau}\|_{L^5}^2\|n^{\varepsilon,\tau}\|_{L^\frac{10}3}^2.
\end{aligned}
\end{equation*}
By (\ref{u-energy}), (\ref{c-energy}), (\ref{c'-energy}), (\ref{c'' half energy}), and (\ref{n-energy}), we have (\ref{c''-energy}).
%Note that  $ \int_{0}^{t}\|\nabla n^{\varepsilon,\tau}\|^2_{L^{\frac{10}{3}}}\leq C $ due to $ \nabla n^{\varepsilon,\tau}\in L^{\infty}_{t}L^{2}_{x}\cap L^{2}_{t}\dot{H}^{1}_{x}. $ Similarly,
%we know $ \int_{0}^{t}\|\nabla c^{\varepsilon,\tau}\|_{L^{5}}^{5}\leq C $ %and $ \int_{0}^{t}\|\nabla c^{\varepsilon,\tau}\|_{L^{6}}^{4}\leq C $ 
%due to $ |\nabla c^{\varepsilon,\tau}|^{\frac32}\in L^{\infty}_{t}L^{2}_{x}\cap L^{2}_{t}\dot{H}^{1}_{x}. $ 
%Moreover, by (\ref{u-energy}), (\ref{c'-energy}), (\ref{n'-energy}) and
%\begin{equation*}
%	\int_{0}^{t}\|\nabla^{2}c^{\varepsilon,\tau}\|_{L^{3}}^{4}\leq C\left(\|\nabla^{2}c^{\varepsilon,\tau}\|_{L^{\infty}_{t}L^{2}_{x}}^{4}+\int_{0}^{t}\|\nabla^{3}c^{\varepsilon,\tau}\|_{L^{2}}^{4}  \right),
%\end{equation*}
% gronwall's inequality, we have (\ref{c''-energy}).
% (\ref{c''})

{\bf\underline{xi. the energy norms of $ |\nabla n^{\varepsilon,\tau}|^{\frac32}: $}}
\begin{equation}\label{n'3-energy}
\begin{aligned}
&\|\nabla n^{\varepsilon,\tau}\|_{L_{t}^{\infty}L^{3}_{x}}^{3}+\int_{0}^{t}\|\nabla |\nabla n^{\varepsilon,\tau}|^{\frac32}\|_{L^{2}}^{2}\\\leq~&C\left(\varepsilon,\tau,\|\nabla\phi\|_{L^{\infty}}, \|\chi\|_{0},\|\kappa\|_{0},\|n_{0}^{\varepsilon}\|_{L^{1}}, \|E_{0}^{\varepsilon}\|_{H^2}  \right)(T+1)^{\frac{321}{8}},
%~~{\rm{for}}~~{\rm {any}}~t\in(0,T).
\end{aligned}
\end{equation}
for any $ t\in(0,T). $

Rewrite the first equation of (\ref{eq:CNS}) as follow:
\beno
\partial_t \partial_{i}n^{\varepsilon,\tau}+\partial_{i} ((u^{\varepsilon,\tau}\ast\rho^{\varepsilon})\cdot\nabla n^{\varepsilon,\tau})= \partial_{i} \Delta n^{\varepsilon,\tau}- \partial_{i}\nabla\cdot\left(
\frac{n^{\varepsilon,\tau}}{1+\tau n^{\varepsilon,\tau}}\nabla c^{\varepsilon,\tau}\chi(c^{\varepsilon,\tau})
\right),
\eeno
with $i = 1,2,3$. Multiplying it by $|\nabla  n^{\varepsilon,\tau}|\partial_{i} n^{\varepsilon,\tau}$ and integrating with respect to space variable $x$ over $\mathbb{R}^{3}$, we obtain
\beno
\frac13\frac{d}{dt}\|\nabla n^{\varepsilon,\tau}\|_{L^{3}}^{3}+\frac89\|\nabla(|\nabla n^{\varepsilon,\tau}|^{\frac32})\|_{L^{2}}^{2}&=& -\int_{\mathbb{R}^3} \partial_{i} ((u^{\varepsilon,\tau}\ast\rho^{\varepsilon})\cdot\nabla n^{\varepsilon,\tau}) |\nabla  n^{\varepsilon,\tau}|\partial_{i} n^{\varepsilon,\tau} \\
&&- \int_{\mathbb{R}^3} \partial_{i}\nabla\cdot\left(
\frac{n^{\varepsilon,\tau}}{1+\tau n^{\varepsilon,\tau}}\nabla c^{\varepsilon,\tau}\chi(c^{\varepsilon,\tau})
\right) |\nabla  n^{\varepsilon,\tau}|\partial_{i} n^{\varepsilon,\tau}\\&:=& S_{1}+S_{2}.
\eeno
For the term $ S_{1} $, H\"{o}lder inequality and interpolation inequality yield that
\begin{equation*}
S_{1}\leq~\|\partial_{i} (u^{\varepsilon,\tau}\ast\rho^{\varepsilon})\|_{L^\infty}\|\nabla n^{\varepsilon,\tau}\|_{L^3}^3
\leq~ C(\varepsilon)\|u^{\varepsilon,\tau}\|_{L^2}\|\nabla n^{\varepsilon,\tau}\|^{\frac32}_{L^2}\|\nabla^2 n^{\varepsilon,\tau}\|^{\frac32}_{L^2}.
\end{equation*}
By the product derivative rule for $ S_{2}, $ there holds
\beno
S_2&=&\int_{\mathbb{R}^3} \partial_{i}\left(
\frac{n^{\varepsilon,\tau}}{1+\tau n^{\varepsilon,\tau}}\right)\nabla c^{\varepsilon,\tau}\chi(c^{\varepsilon,\tau})
\nabla\left(|\nabla  n^{\varepsilon,\tau}|\partial_{i} n^{\varepsilon,\tau}\right)\\
&+&\int_{\mathbb{R}^3}
\frac{n^{\varepsilon,\tau}}{1+\tau n^{\varepsilon,\tau}}\partial_{i}\left(\nabla c^{\varepsilon,\tau}\right)\chi(c^{\varepsilon,\tau})
\nabla\left(|\nabla  n^{\varepsilon,\tau}|\partial_{i} n^{\varepsilon,\tau}\right)\\
&+&\int_{\mathbb{R}^3}
\frac{n^{\varepsilon,\tau}}{1+\tau n^{\varepsilon,\tau}}\nabla c^{\varepsilon,\tau}\partial_{i}\left(\chi(c^{\varepsilon,\tau})\right)
\nabla\left(|\nabla  n^{\varepsilon,\tau}|\partial_{i} n^{\varepsilon,\tau}\right)\\
&=:&S_{21}+S_{22}+S_{23}.
\eeno 
For $S_{21}$,
due to
\beno
\|\nabla n^{\varepsilon,\tau}\|_{L^{\frac92}}\leq \|\nabla n^{\varepsilon,\tau}\|_{L^{\frac{10}3}}^{\frac29}\|\nabla n^{\varepsilon,\tau}\|_{L^5}^{\frac79},
\eeno
there holds
\beno
S_{21}
&\leq&  C\int_{\mathbb{R}^3}|\nabla c^{\varepsilon,\tau}|\chi(c^{\varepsilon,\tau})|\nabla n^{\varepsilon,\tau}|^{\frac32}\left|\nabla\left(|\nabla n^{\varepsilon,\tau}|^{\frac32}\right)\right|\\
%&\leq& \frac18\int_{\mathbb{R}^3}\left|\nabla\left(|\nabla n^{\varepsilon,\tau}|^{\frac32}\right)\right|^2dx+C\int_{\mathbb{R}^3}|\nabla c^{\varepsilon,\tau}|^2 |\chi^2(c^{\varepsilon,\tau})| |\nabla n^{\varepsilon,\tau}|^{3}\\
&\leq& \frac18\int_{\mathbb{R}^3}\left|\nabla\left(|\nabla n^{\varepsilon,\tau}|^{\frac32}\right)\right|^2dx+C\|\chi(c^{\varepsilon,\tau})\|^2_{L^\infty}\|\nabla c^{\varepsilon,\tau}\|_{L^6}^2\|\nabla n^{\varepsilon,\tau}\|_{L^{\frac92}}^3\\
&\leq& \frac18\int_{\mathbb{R}^3}\left|\nabla\left(|\nabla n^{\varepsilon,\tau}|^{\frac32}\right)\right|^2dx+C\|\chi(c^{\varepsilon,\tau})\|^2_{L^\infty}\|\nabla c^{\varepsilon,\tau}\|_{L^6}^2\|\nabla n^{\varepsilon,\tau}\|_{L^{\frac{10}3}}^{\frac23}\|\nabla n^{\varepsilon,\tau}\|_{L^5}^{\frac73}.
\eeno
Similarly, for $S_{22}$ and $ S_{23}, $ we have
\beno
S_{22}+S_{23}&\leq& \frac18\int_{\mathbb{R}^3}\left|\nabla\left(|\nabla n^{\varepsilon,\tau}|^{\frac32}\right)\right|^2dx+  C(\tau)\int \left(|\nabla^2 c^{\varepsilon,\tau}|^2 +|\nabla c^{\varepsilon,\tau}|^4\right)|\nabla n^{\varepsilon,\tau}|.
\eeno
Then we obtain
\beno
&&\|\nabla n^{\varepsilon,\tau}\|_{L_{t}^{\infty}L^{3}_{x}}^{3}+\int_{0}^{t}\|\nabla |\nabla n^{\varepsilon,\tau}|^{\frac32}\|_{L^{2}}^{2}\\
&\leq&  C(\varepsilon)\|u^{\varepsilon,\tau}\|_{L^\infty L^2}\int_0^t\|\nabla n^{\varepsilon,\tau}\|^{\frac32}_{L^2}\|\nabla^2 n^{\varepsilon,\tau}\|^{\frac32}_{L^2}ds\\
&&+C\|\chi(c^{\varepsilon,\tau})\|^2_{L^\infty}\int_0^t\|\nabla c^{\varepsilon,\tau}\|_{L^6}^2\|\nabla n^{\varepsilon,\tau}\|_{L^{\frac{10}3}}^{\frac23}\|\nabla n^{\varepsilon,\tau}\|_{L^5}^{\frac73}ds\\
&&+ C(\tau)\int_0^t\int \left(|\nabla^2 c^{\varepsilon,\tau}|^2 +|\nabla c^{\varepsilon,\tau}|^4\right)|\nabla n^{\varepsilon,\tau}|dxds.
\eeno
By (\ref{u-energy}) and (\ref{n'-energy}), we know
\beno
\|u^{\varepsilon,\tau}\|_{L^\infty L^2}\int_0^t\|\nabla n^{\varepsilon,\tau}\|^{\frac32}_{L^2}\|\nabla^2 n^{\varepsilon,\tau}\|^{\frac32}_{L^2}ds\leq C T^{\frac{129}{8}}.\eeno
By
(\ref{c''-energy}), (\ref{n'-energy}) and the interpolation inequality
\ben\label{eq:interpolation}
\|f\|_{L^qL^p}^2\leq C \left( \|f\|_{L^\infty L^2}^2 +\|\nabla f\|_{L^2L^2}^2\right),\quad \frac2q+\frac3p=\frac32,\quad 2\leq p\leq 6,
\een
there holds
\beno
&&\int_0^t\|\nabla c^{\varepsilon,\tau}\|_{L^6}^2\|\nabla n^{\varepsilon,\tau}\|_{L^{\frac{10}3}}^{\frac23}\|\nabla n^{\varepsilon,\tau}\|_{L^5}^{\frac73}ds\\
&\leq & \sup_t\|\nabla^2 c^{\varepsilon,\tau}\|_{L^2}^2\|\nabla n^{\varepsilon,\tau}\|_{L^{\frac{10}{3}}_{t,x}}^{\frac23}\|\nabla n^{\varepsilon,\tau}\|_{L^{5}_{t,x}}^{\frac73}\\
&\leq & \frac14 \left(\|\nabla n^{\varepsilon,\tau}\|_{L_{t}^{\infty}L^{3}_{x}}^{3}+\int_{0}^{t}\|\nabla |\nabla n^{\varepsilon,\tau}|^{\frac32}\|_{L^{2}}^{2}\right)+
C\sup_t\|\nabla^2 c^{\varepsilon,\tau}\|_{L^2}^9\|\nabla n^{\varepsilon,\tau}\|_{L^{\frac{10}{3}}_{t,x}}^{3}\\
&\leq & \frac14 \left(\|\nabla n^{\varepsilon,\tau}\|_{L_{t}^{\infty}L^{3}_{x}}^{3}+\int_{0}^{t}\|\nabla |\nabla n^{\varepsilon,\tau}|^{\frac32}\|_{L^{2}}^{2}\right)+
CT^{\frac{99}{4}+\frac{123}{8}}.
\eeno
By
(\ref{c''-energy}), (\ref{n'-energy}) and the interpolation inequality
\beno
&&C(\tau)\int_0^t\int \left(|\nabla^2 c^{\varepsilon,\tau}|^2 +|\nabla c^{\varepsilon,\tau}|^4\right)|\nabla n^{\varepsilon,\tau}|dxds\\
&\leq & C(\tau) \|\nabla^2 c^{\varepsilon,\tau}\|_{L^{\frac{10}{3}}_{t,x}}^{2} \|\nabla n^{\varepsilon,\tau}\|_{L^{\frac{5}{2}}_{t,x}}+\sup_t\|\nabla^2 c^{\varepsilon,\tau}\|_{L^2}^4 \|\nabla n^{\varepsilon,\tau}\|_{L_{t}^{1}L^{3}_{x}}\\
&\leq& C T^{\frac{11}{2}+\frac{43}{8}}+CT^{11+\frac{47}{8}}.
\eeno
Thus the proof of  (\ref{n'3-energy}) is complete.

{\bf\underline{xii. the energy norms of $ \nabla^{3}c^{\varepsilon,\tau}: $}}
\begin{equation}\label{c'''-energy}
\begin{aligned}
&\|\nabla c^{\varepsilon,\tau}\|_{L^{\infty}_{t}H^{2}_{x}}^{2}+\int_{0}^{t}\|\nabla^{2}c^{\varepsilon,\tau}\|_{H^{2}}^{2}\\\leq~&C(\varepsilon,\tau,\|\nabla\phi\|_{L^{\infty}},\|\chi\|_{0},\|\kappa\|_{0},\|n_{0}^{\varepsilon}\|_{L^{1}},\|E_{0}^{\varepsilon}\|_{H^{2}}  )(T+1)^{54},
%~~{\rm{for}}~~{\rm {any}}~t\in(0,T).
\end{aligned}
\end{equation}
for any $ t\in(0,T). $

We write the second equation of $\eqref{eq:CNS}$ as follow:
\begin{equation*}
\partial_t \partial_{ijk}c^{\varepsilon,\tau}+\partial_{ijk}\left(u^{\varepsilon,\tau}\cdot\nabla c^{\varepsilon,\tau}  \right)-\partial_{ijk}\Delta c^{\varepsilon,\tau}=-\partial_{ijk}\left( \frac{1}{\tau}\ln(1+\tau n^{\varepsilon,\tau})\kappa(c^{\varepsilon,\tau}) \right),
\end{equation*}
with $ i,j,k=1,2,3. $
Multiplying above equation  by $\partial_{ijk}c^{\varepsilon,\tau} $ and integrating with respect to space variable $ x $ over $ \mathbb{R}^{3} $, we obtain
\beno
\frac12\frac d{dt} \|\nabla^3 c^{\varepsilon,\tau}\|_{L^2}^2 + \| \nabla^4 c^{\varepsilon,\tau}\|_{L^2}^2 &=& - \int_{\mathbb{R}^3} \partial_{ijk} (u^{\varepsilon,\tau}\cdot\nabla c^{\varepsilon,\tau}) \partial_{ijk}c^{\varepsilon,\tau}\\
&&- \int_{\mathbb{R}^3} \partial_{ijk}\left( \frac{1}{\tau}\ln(1+\tau n^{\varepsilon,\tau})\kappa(c^{\varepsilon,\tau}) \right)\partial_{ijk}c^{\varepsilon,\tau} \\
&:=& H_1 + H_2.
\eeno
For $ H_{1} $, by H\"{o}lder inequality, we have
\begin{equation*}
\begin{aligned}
H_{1}\leq~& C\big[\|\nabla^{2}u^{\varepsilon,\tau}\|_{L^{\frac{10}{3}}}\|\nabla c^{\varepsilon,\tau}\|_{L^5}\|\nabla^4 c^{\varepsilon,\tau}\|_{L^2}+\|\nabla u^{\varepsilon,\tau}\|_{L^3}\|\nabla^2 c^{\varepsilon,\tau}\|_{L^6}\|\nabla^4 c^{\varepsilon,\tau}\|_{L^2}\\
&+~\|u^{\varepsilon,\tau}\|_{L^\infty}\|\nabla^3 c^{\varepsilon,\tau}\|_{L^2}\|\nabla^4 c^{\varepsilon,\tau}\|_{L^2}\big]\\
\leq~&\frac{1}{16}\|\nabla^{4}c^{\varepsilon,\tau}\|_{L^{2}}^{2}+C(\|\nabla^{2}u^{\varepsilon,\tau}\|_{L^{\frac{10}{3}}}^{2}\|\nabla c^{\varepsilon,\tau}\|_{L^5}^{2}+\|\nabla u^{\varepsilon,\tau}\|_{L^{2}}\|\nabla^{2}u^{\varepsilon,\tau}\|_{L^{2}}\|\nabla^{3} c^{\varepsilon,\tau}\|_{L^{2}}^{2}\\&+~\|u^{\varepsilon,\tau}\|_{L^\infty}^{2}\|\nabla^3 c^{\varepsilon,\tau}\|_{L^2}^{2}).
\end{aligned}
\end{equation*}
Next, we deal with $ H_{2}. $ Similarly, we have 
\beno
H_{2}&\leq~ %C(\tau,\|\kappa\|_{0})\|\nabla^4 c^{\varepsilon,\tau}\|_{L^2}\big[\|\nabla^2 n^{\varepsilon,\tau}\|_{L^2}+\|\nabla n^{\varepsilon,\tau}\|^2_{L^4}\\&+~\| n^{\varepsilon,\tau}\|_{L^\infty}\left(\|\nabla^2 c^{\varepsilon,\tau}\|_{L^2}
%+\|\nabla c^{\varepsilon,\tau}\|^2_{L^4}\right)+\|\nabla n^{\varepsilon,\tau}\|_{L^4}\|\nabla c^{\varepsilon,\tau}\|_{L^4}\big]\\
& \frac{1}{16}\|\nabla^{4}c^{\varepsilon,\tau}\|_{L^{2}}^{2}+C(\tau,\|\kappa\|_{0})\big[\|\nabla^2 n^{\varepsilon,\tau}\|^2_{L^2}+\|\nabla n^{\varepsilon,\tau}\|^4_{L^4}\\&&+~ \left(\|n^{\varepsilon,\tau}\nabla^2 c^{\varepsilon,\tau}\|^2_{L^2}
+\|n^{\varepsilon,\tau}|\nabla c^{\varepsilon,\tau}|^2\|^2_{L^2}\right)+\|\nabla n^{\varepsilon,\tau}\|^2_{L^4}\|\nabla c^{\varepsilon,\tau}\|^2_{L^4}\big].
\eeno
Collecting $ H_{1}$ and $H_{2}, $ we have
\beno
&&\frac{d}{dt}\|\nabla^{3}c^{\varepsilon,\tau}\|_{L^{2}}^{2}+\|\nabla^{4}c^{\varepsilon,\tau}\|_{L^{2}}^{2}\\&\leq&
C(\|\nabla^{2}u^{\varepsilon,\tau}\|_{L^{\frac{10}{3}}}^{2}\|\nabla c^{\varepsilon,\tau}\|_{L^5}^{2}+\|\nabla u^{\varepsilon,\tau}\|_{L^{2}}\|\nabla^{2}u^{\varepsilon,\tau}\|_{L^{2}}\|\nabla^{3} c^{\varepsilon,\tau}\|_{L^{2}}^{2}+\|u^{\varepsilon,\tau}\|_{L^\infty}^{2}\|\nabla^3 c^{\varepsilon,\tau}\|_{L^2}^{2})\\
&&+C(\tau,\|\kappa\|_{0})\big[\|\nabla^2 n^{\varepsilon,\tau}\|^2_{L^2}+\|\nabla n^{\varepsilon,\tau}\|^4_{L^4}\\&&+~ \left(\|n^{\varepsilon,\tau}\nabla^2 c^{\varepsilon,\tau}\|^2_{L^2}
+\|n^{\varepsilon,\tau}|\nabla c^{\varepsilon,\tau}|^2\|^2_{L^2}\right)+\|\nabla c^{\varepsilon,\tau}\|^4_{L^4}\big].
\eeno
%By  (\ref{u-energy}), (\ref{c'-energy}) and (\ref{c''-energy}), we have
%\beno
%&&\int_0^t\left(\|\nabla^3c^{\varepsilon,\tau}\|_{L^{2}}^{2}+\|\nabla^{2}c^{\varepsilon,\tau}\|_{L^{2}}^{2}  +\|\nabla c^{\varepsilon,\tau}\|_{L^{2}}^{2}\right)\|u^{\varepsilon,\tau}\|_{L^{2}}^{2}ds\\
%&\leq &  C T^3.
%\eeno
(\ref{u'-energy}), (\ref{u''-energy}) and (\ref{c'5 half energy}) show that
\begin{equation*}
	\begin{aligned}
	&\int_{0}^{t}\left( \|\nabla^{2}u^{\varepsilon,\tau}\|_{L^{\frac{10}{3}}}^{2}\|\nabla c^{\varepsilon,\tau}\|_{L^5}^{2}+\|\nabla u^{\varepsilon,\tau}\|_{L^{2}}\|\nabla^{2}u^{\varepsilon,\tau}\|_{L^{2}}\|\nabla^{3} c^{\varepsilon,\tau}\|_{L^{2}}^{2}+\|u^{\varepsilon,\tau}\|_{L^\infty}^{2}\|\nabla^3 c^{\varepsilon,\tau}\|_{L^2}^{2} \right)ds\\\leq~&C\sup_t\|\nabla c^{\varepsilon,\tau}\|_{L^{5}}^{2}\int_{0}^{t}\|\nabla^{2}u^{\varepsilon,\tau}\|_{L^{\frac{10}{3}}}^{2}ds+C\sup_t\|\nabla u^{\varepsilon,\tau}\|_{L^{2}}\sup_t\|\nabla^{2} u^{\varepsilon,\tau}\|_{L^{2}}\int_{0}^{t}\|\nabla^{3}c^{\varepsilon,\tau}\|_{L^{2}}^{2}ds\\&+~C\sup_t\|u^{\varepsilon,\tau}\|_{L^{2}}^{\frac12}\sup_t\|\nabla^{2}u^{\varepsilon,\tau}\|_{L^{2}}^{\frac32}\int_{0}^{t}\|\nabla^{3}c^{\varepsilon,\tau}\|_{L^{2}}^{2}ds\\\leq~&CT^{8}.
	\end{aligned}
\end{equation*}
(\ref{n'-energy}),  (\ref{c'-energy}) and (\ref{eq:interpolation}) imply 
\beno
&&\int_0^t\|\nabla^2 n^{\varepsilon,\tau}\|^2_{L^2}+\|\nabla n^{\varepsilon,\tau}\|^4_{L^4}+\|\nabla c^{\varepsilon,\tau}\|^4_{L^4} ds\leq   C T^{54}.
\eeno
Moreover, (\ref{n'-energy}),  (\ref{c'-energy}) and (\ref{c''-energy})
\beno
&&\int_0^t\|n^{\varepsilon,\tau}\nabla^2 c^{\varepsilon,\tau}\|^2_{L^2}
+\|n^{\varepsilon,\tau}|\nabla c^{\varepsilon,\tau}|^2\|^2_{L^2}ds\\
& \leq& \|\nabla^2 c^{\varepsilon,\tau}\|^2_{L^{\frac{10}{3}}_{t,x}} T^{\frac25} \|n^{\varepsilon,\tau}\|^{\frac{2}{25}}_{L^1}\|\nabla n^{\varepsilon,\tau}\|^{\frac{48}{25}}_{L^2}+CT \|n^{\varepsilon,\tau}\|^{2}_{L^6} \|\nabla c^{\varepsilon,\tau}\|^{2}_{L^6}\\
&\leq & C T^{\frac{787}{50}}+CT^{\frac{55}{4}}.\eeno
The proof of (\ref{c'''-energy}) is complete.

{\bf\underline{xiii. the energy norms of $ \nabla^{2}n^{\varepsilon,\tau}: $}}
\begin{equation}\label{n''-energy}
\begin{aligned}
&\|n^{\varepsilon,\tau}\|_{L^{\infty}_{t}H^{2}_{x}}^{2}+\int_{0}^{t}\|\nabla n^{\varepsilon,\tau}\|_{H^{2}}^{2}\\\leq~&C\left(\varepsilon,\tau,\|\nabla\phi\|_{L^{\infty}}, \|\chi\|_{0},\|\kappa\|_{0},\|n_{0}^{\varepsilon}\|_{L^{1}}, \|E_{0}^{\varepsilon}\|_{H^2}  \right)(T+1)^{\frac{411}{8}},
%~~{\rm{for}}~~{\rm {any}}~t\in(0,T).
\end{aligned}
\end{equation}
for any $ t\in(0,T). $

Rewrite the first equation of (\ref{eq:CNS}) as follow:
\beno
\partial_t \partial_{ij}n^{\varepsilon,\tau}+\partial_{ij} ((u^{\varepsilon,\tau}\ast\rho^{\varepsilon})\cdot\nabla n^{\varepsilon,\tau})= \partial_{ij} \Delta n^{\varepsilon,\tau}- \partial_{ij}\nabla\cdot\left(
\frac{n^{\varepsilon,\tau}}{1+\tau n^{\varepsilon,\tau}}\nabla c^{\varepsilon,\tau}\chi(c^{\varepsilon,\tau})
\right),
\eeno
with $i,j = 1,2,3$. Multiplying it by $\partial_{ij} n^{\varepsilon,\tau}$ and integrating with respect to space variable $x$ over $\mathbb{R}^{3}$, we obtain
\beno
\frac12 \frac d{dt} \|\nabla^2 n^{\varepsilon,\tau}\|_{L^2}^2 + \|\nabla^3 n^{\varepsilon,\tau}\|_{L^2}^2 &=& -\int_{\mathbb{R}^3} \partial_{ij} ((u^{\varepsilon,\tau}\ast\rho^{\varepsilon})\cdot\nabla n^{\varepsilon,\tau}) \partial_{ij} n^{\varepsilon,\tau} \\
&&- \int_{\mathbb{R}^3} \partial_{ij}\nabla\cdot\left(
\frac{n^{\varepsilon,\tau}}{1+\tau n^{\varepsilon,\tau}}\nabla c^{\varepsilon,\tau}\chi(c^{\varepsilon,\tau})
\right) \partial_{ij} n^{\varepsilon,\tau}\\&:=& G_{1}+G_{2}.
\eeno
The term $ G_{1} $ is same as $ B_{1}, $ we arrive
\begin{equation*}
\begin{aligned}
G_{1}\leq \frac{3}{16}\|\nabla^{3}n^{\varepsilon,\tau}\|_{L^{2}}^{2}+C(\varepsilon)\|u^{\varepsilon,\tau}\|_{L^{2}}^{2}\|\nabla^{2}n^{\varepsilon,\tau}\|_{L^{2}}^{2}+C(\varepsilon)\|\nabla n^{\varepsilon,\tau}\|_{L^{2}}^{2}\|u^{\varepsilon,\tau}\|_{L^{2}}^{2}.
\end{aligned}
\end{equation*}
By the product derivative rule for $G_2$, we also have
\beno
G_2 &=& \int_{\mathbb{R}^{3}} \partial_{ij}\left(
\frac{n^{\varepsilon,\tau}}{1+\tau n^{\varepsilon,\tau}}\nabla c^{\varepsilon,\tau}\chi(c^{\varepsilon,\tau})
\right) \cdot\nabla\partial_{ij}n^{\varepsilon,\tau} dx \\
&=& \int_{\mathbb{R}^{3}} \partial_{ij}(\frac{n^{\varepsilon,\tau}}{1+\tau n^{\varepsilon,\tau}})(\nabla c^{\varepsilon,\tau}\chi(c^{\varepsilon,\tau})) \cdot\nabla\partial_{ij}n^{\varepsilon,\tau} dx \\&&+ \int_{\mathbb{R}^{3}} \partial_{i}(\frac{n^{\varepsilon,\tau}}{1+\tau n^{\varepsilon,\tau}})\partial_j(\nabla c^{\varepsilon,\tau}\chi(c^{\varepsilon,\tau})) \cdot\nabla\partial_{ij}n^{\varepsilon,\tau} dx \\
&&+ \int_{\mathbb{R}^{3}} \partial_{j}(\frac{n^{\varepsilon,\tau}}{1+\tau n^{\varepsilon,\tau}})\partial_i(\nabla c^{\varepsilon,\tau}\chi(c^{\varepsilon,\tau})) \cdot\nabla\partial_{ij}n^{\varepsilon,\tau} dx\\&& + \int_{\mathbb{R}^{3}}\frac{n^{\varepsilon,\tau}}{1+\tau n^{\varepsilon,\tau}}\partial_{ij}(\nabla c^{\varepsilon,\tau}\chi(c^{\varepsilon,\tau})) \cdot\nabla\partial_{ij}n^{\varepsilon,\tau} dx \\
&:=& G_{21} + G_{22} + G_{23} + G_{24}.
\eeno
For $ G_{21}, $ note that
\begin{equation*}
\begin{aligned}
\partial_{ij}\left(\frac{n^{\varepsilon,\tau}}{1+\tau n^{\varepsilon,\tau}} \right)&=~\partial_{i}\left( \frac{\partial_{j}n^{\varepsilon,\tau}}{(1+\tau n^{\varepsilon,\tau})^{2}} \right)\\&=~\frac{\partial_{ij}n^{\varepsilon,\tau}(1+\tau n^{\varepsilon,\tau})-2\tau\partial_{i}n^{\varepsilon,\tau}\partial_{j}n^{\varepsilon,\tau}}{(1+\tau n^{\varepsilon,\tau})^{3}},
\end{aligned}
\end{equation*}
we have
\begin{equation*}
\begin{aligned}
G_{21}=~&\int_{\mathbb{R}^{3}}\frac{\partial_{ij}n^{\varepsilon,\tau}}{(1+\tau n^{\varepsilon,\tau})^{2}}\nabla c^{\varepsilon,\tau}\chi(c^{\varepsilon,\tau})\cdot\nabla\partial_{ij}n^{\varepsilon,\tau}\\&-~\int_{\mathbb{R}^{3}}\frac{2\tau\partial_{i}n^{\varepsilon,\tau}\partial_{j}n^{\varepsilon,\tau}}{(1+\tau n^{\varepsilon,\tau})^{3}}\nabla c^{\varepsilon,\tau}\chi(c^{\varepsilon,\tau})\cdot\nabla\partial_{ij}n^{\varepsilon,\tau}\\\leq~&C\left(\|\chi\|_{0} \right)\|\nabla^{2}n^{\varepsilon,\tau}\|_{L^{2}}\|\nabla c^{\varepsilon,\tau}\|_{L^{\infty}}\|\nabla^{3}n^{\varepsilon,\tau}\|_{L^{2}}\\&+~C\left(\tau, \|\chi\|_{0} \right)\|\nabla n^{\varepsilon,\tau}\|_{L^{2}}\|\nabla c^{\varepsilon,\tau}\|_{L^{\infty}}\|\nabla^{3}n^{\varepsilon,\tau}\|_{L^{2}}\\\leq~&\frac{1}{8}\|\nabla^{3}n^{\varepsilon,\tau}\|_{L^{2}}^{2}+C\left(\|\chi\|_{0} \right)\|\nabla^{2}n^{\varepsilon,\tau}\|^2_{L^{2}}\|\nabla c^{\varepsilon,\tau}\|^2_{L^{\infty}}\\&+~C\left(\tau, \|\chi\|_{0} \right)\|\nabla n^{\varepsilon,\tau}\|^2_{L^{2}}\|\nabla c^{\varepsilon,\tau}\|^2_{L^{\infty}}.
%\\\leq~&\frac{1}{8}\|\nabla^{3}n^{\varepsilon,\tau}\|_{L^{2}}^{2}+\frac16\left(\|\nabla c^{\varepsilon,\tau}\|_{L^{4}}^{8}+1 \right)+C\left(\chi(\|c_{0}^{\varepsilon}\|_{L^{\infty}}) \right)\|\nabla^{2}n^{\varepsilon,\tau}\|_{L^{2}}^{2}\\&+~C\left(\tau, \chi(\|c_{0}^{\varepsilon}\|_{L^{\infty}})  \right)\|\nabla n^{\varepsilon,\tau}\|_{L^{6}}^{2}\|\nabla^{2}c^{\varepsilon,\tau}\|_{L^{2}}^{2}\|\nabla^{2}n^{\varepsilon,\tau}\|_{L^{2}}^{2}.	
\end{aligned}
\end{equation*}
For $ G_{22}, $ note that
\begin{equation*}
\begin{aligned}
\partial_{j}\left(\nabla c^{\varepsilon,\tau}\chi(c^{\varepsilon,\tau})  \right)=\partial_{j}\nabla c^{\varepsilon,\tau}\chi(c^{\varepsilon,\tau})+\nabla c^{\varepsilon,\tau}\chi'(c^{\varepsilon,\tau})\partial_{j}c^{\varepsilon,\tau},
\end{aligned}
\end{equation*}
we obtain
\begin{equation*}
\begin{aligned}
G_{22}=~&\int_{\mathbb{R}^{3}}\frac{\partial_{j}n^{\varepsilon,\tau}}{(1+\tau n^{\varepsilon,\tau})^{2}}\partial_{j}\nabla c^{\varepsilon,\tau}\chi(c^{\varepsilon,\tau})\cdot\nabla\partial_{ij}n^{\varepsilon,\tau}\\&+~\int_{\mathbb{R}^{3}}\frac{\partial_{j}n^{\varepsilon,\tau}}{(1+\tau n^{\varepsilon,\tau})^{2}}\nabla c^{\varepsilon,\tau}\chi'(c^{\varepsilon,\tau})\partial_{j}c^{\varepsilon,\tau}\cdot\nabla\partial_{ij}n^{\varepsilon,\tau}\\:=& G_{221}+G_{222}.
\end{aligned}
\end{equation*}
For $ G_{221}, $ 
\begin{equation*}
\begin{aligned}
G_{221}\leq~&C\left(\|\chi\|_{0} \right)\|\nabla n^{\varepsilon,\tau}\|_{L^{3}}\|\nabla^{2}c^{\varepsilon,\tau}\|_{L^{6}}\|\nabla^{3}n^{\varepsilon,\tau}\|_{L^{2}}\\
\leq~&\frac{1}{16}\|\nabla^{3}n^{\varepsilon,\tau}\|_{L^{2}}^{2}+\|\nabla n^{\varepsilon,\tau}\|^2_{L^{3}}\|\nabla^{3}c^{\varepsilon,\tau}\|^2_{L^{2}}.
%+C\left(\|\chi\|_{0} \right)\|\nabla n^{\varepsilon,\tau}\|_{L^{2}}^{\frac12}\|\nabla^{2}c^{\varepsilon,\tau}\|_{L^{2}}^{\frac12}\|\nabla^{2}n^{\varepsilon,\tau}\|_{L^{2}}^{\frac32}\|\nabla^{3}c^{\varepsilon,\tau}\|_{L^{2}}^{\frac32}\\\leq~&\frac{1}{16}\|\nabla^{3}n^{\varepsilon,\tau}\|_{L^{2}}^{2}+\frac14\|\nabla n^{\varepsilon,\tau}\|_{L^{2}}^{2}\|\nabla^{2}c^{\varepsilon,\tau}\|_{L^{2}}^{2}+C\left(\|\chi\|_{0} \right)\|\nabla^{3}c^{\varepsilon,\tau}\|_{L^{2}}^{2}\|\nabla^{2}n^{\varepsilon,\tau}\|_{L^{2}}^{2}.
\end{aligned}
\end{equation*}
For $ G_{222}, $
\begin{equation*}
\begin{aligned}
G_{222}\leq~&C\left(\|\chi\|_{0}  \right)\|\nabla n^{\varepsilon,\tau}\|_{L^{2}}\|\nabla c^{\varepsilon,\tau}\|_{L^{\infty}}\|\nabla^{3}n^{\varepsilon,\tau}\|_{L^{2}}\\\leq~&\frac{1}{16}\|\nabla^{3}n^{\varepsilon,\tau}\|_{L^{2}}^{2}+\|\nabla n^{\varepsilon,\tau}\|^2_{L^{2}}\|\nabla c^{\varepsilon,\tau}\|_{L^{\infty}}^{2}.
%C\left(\|\chi\|_{0} \right)\|\nabla c^{\varepsilon,\tau}\|_{L^{6}}^{4}\|\nabla^{2}n^{\varepsilon,\tau}\|_{L^{2}}^{2}.
\end{aligned}
\end{equation*}
Then
\begin{equation*}
\begin{aligned}
G_{22}\leq~&\frac18\|\nabla^3 n^{\varepsilon,\tau}\|_{L^{2}}^{2}+\|\nabla n^{\varepsilon,\tau}\|^2_{L^{3}}\|\nabla^{3}c^{\varepsilon,\tau}\|^2_{L^{2}}+\|\nabla n^{\varepsilon,\tau}\|^2_{L^{2}}\|\nabla c^{\varepsilon,\tau}\|_{L^{\infty}}^{2}.
%\frac14\|\nabla n^{\varepsilon,\tau}\|_{L^{2}}^{2}\|\nabla^{2}c^{\varepsilon,\tau}\|_{L^{2}}^{2}+C\left(\|\chi\|_{0} \right)\|\nabla^{3}c^{\varepsilon,\tau}\|_{L^{2}}^{2}\|\nabla^{2}n^{\varepsilon,\tau}\|_{L^{2}}^{2}\\&+~C\left(\|\chi\|_{0} \right)\|\nabla c^{\varepsilon,\tau}\|_{L^{6}}^{4}\|\nabla^{2}n^{\varepsilon,\tau}\|_{L^{2}}^{2}.
\end{aligned}
\end{equation*}
The estimation of term $ G_{23} $ is same as $ G_{22}, $ we omit it.
For $ G_{24}, $ note that
\begin{equation*}
\begin{aligned}
&\partial_{ij}\left(\nabla c^{\varepsilon,\tau}\chi(c^{\varepsilon,\tau})  \right)\\=~&\partial_{i}\left(\partial_{j}\nabla c^{\varepsilon,\tau}\chi(c^{\varepsilon,\tau})+\nabla c^{\varepsilon,\tau}\chi'(c^{\varepsilon,\tau})\partial_{j}c^{\varepsilon,\tau}   \right)\\=~&\partial_{ij}\nabla c^{\varepsilon,\tau}\chi(c^{\varepsilon,\tau})+\partial_{j}\nabla c^{\varepsilon,\tau}\chi'(c^{\varepsilon,\tau})\partial_{i}c^{\varepsilon,\tau}+\partial_{i}\nabla c^{\varepsilon,\tau}\chi'(c^{\varepsilon,\tau})\partial_{j}c^{\varepsilon,\tau}\\&+~\nabla c^{\varepsilon,\tau}\chi''(c^{\varepsilon,\tau})\partial_{i}c^{\varepsilon,\tau}\partial_{j}c^{\varepsilon,\tau}+\nabla c^{\varepsilon,\tau}\chi'(c^{\varepsilon,\tau})\partial_{ij}c^{\varepsilon,\tau}\\
\leq~& C\left(\|\chi\|_{0}  \right)\left(|\nabla^3 c^{\varepsilon,\tau}|+|\nabla^2 c^{\varepsilon,\tau}||\nabla c^{\varepsilon,\tau}|+|\nabla c^{\varepsilon,\tau}|^3\right),
\end{aligned}
\end{equation*}
we get
\beno
&&G_{24}\leq \int_{\mathbb{R}^{3}}\frac{n^{\varepsilon,\tau}}{1+\tau n^{\varepsilon,\tau}}\partial_{ij}(\nabla c^{\varepsilon,\tau}\chi(c^{\varepsilon,\tau})) \cdot\nabla\partial_{ij}n^{\varepsilon,\tau} dx \\
&\leq& \int_{\mathbb{R}^{3}}C\left(\tau,\|\chi\|_{0}  \right)\left(|\nabla^3 c^{\varepsilon,\tau}|+|\nabla^2 c^{\varepsilon,\tau}||\nabla c^{\varepsilon,\tau}|+|\nabla c^{\varepsilon,\tau}|^3\right)\nabla\partial_{ij}n^{\varepsilon,\tau} dx\\
&\leq& \|\nabla^{3} n^{\varepsilon,\tau}\|_{L^2}\left(\|\nabla^3 c^{\varepsilon,\tau}\|_{L^2}+\|\nabla^2 c^{\varepsilon,\tau}\|_{L^2}\|\nabla c^{\varepsilon,\tau}\|_{L^\infty}+\|\nabla c^{\varepsilon,\tau}\|_{L^6}^3\right)\\
&\leq&\frac18\|\nabla^{3} n^{\varepsilon,\tau}\|_{L^{2}}^{2}+C\left(\|\nabla^3 c^{\varepsilon,\tau}\|^2_{L^2}+\|\nabla^2 c^{\varepsilon,\tau}\|^2_{L^2}\|\nabla c^{\varepsilon,\tau}\|^2_{L^\infty}+\|\nabla c^{\varepsilon,\tau}\|_{L^6}^6\right).
\eeno
Collecting $ G_{21}-G_{24}, $ we have
\beno
G_{2}&\leq& \frac{1}{2}\|\nabla^{3}n^{\varepsilon,\tau}\|_{L^{2}}^{2}+C\left(\|\chi\|_{0} \right)\|\nabla^{2}n^{\varepsilon,\tau}\|^2_{L^{2}}\|\nabla c^{\varepsilon,\tau}\|^2_{L^{\infty}}\\&+&C\left(\tau, \|\chi\|_{0} \right)\|\nabla n^{\varepsilon,\tau}\|^2_{L^{2}}\|\nabla c^{\varepsilon,\tau}\|^2_{L^{\infty}}\\
&+&C\|\nabla n^{\varepsilon,\tau}\|^2_{L^{3}}\|\nabla^{3}c^{\varepsilon,\tau}\|^2_{L^{2}}+\|\nabla n^{\varepsilon,\tau}\|^2_{L^{2}}\|\nabla c^{\varepsilon,\tau}\|_{L^{\infty}}^{2}\\
&+&C\left(\|\nabla^3 c^{\varepsilon,\tau}\|^2_{L^2}+\|\nabla^2 c^{\varepsilon,\tau}\|^2_{L^2}\|\nabla c^{\varepsilon,\tau}\|^2_{L^\infty}+\|\nabla c^{\varepsilon,\tau}\|_{L^6}^6\right).
\eeno

%\begin{equation*}
%	\begin{aligned}
%	G_{2}\leq~&\frac{11}{16}\|\nabla^{3}n^{\varepsilon,\tau}\|_{L^{2}}^{2}+\frac16\|\nabla c^{\varepsilon,\tau}\|_{L^{4}}^{4}\\&+~C\left(\tau,\|\chi\|_{0} \right)\left(\|\nabla^{3}c^{\varepsilon,\tau}\|_{L^{2}}^{2}+\|\nabla c^{\varepsilon,\tau}\|_{L^{6}}^{4}+\|\nabla^{2}c^{\varepsilon,\tau}\|_{L^{2}}^{2}\|\nabla n^{\varepsilon,\tau}\|_{L^{6}}^{2}+1  \right)\|\nabla^{2}n^{\varepsilon,\tau}\|_{L^{2}}^{2}\\&+~\frac12\|\nabla n^{\varepsilon,\tau}\|_{L^{2}}^{2}\|\nabla^{2}c^{\varepsilon,\tau}\|_{L^{2}}^{2}+C\left(\tau,\|\chi\|_{0} \right)\|\nabla c^{\varepsilon,\tau}\|_{L^{2}}^{2}\\&+~C\left(\tau,\|\chi\|_{0} \right)\|\nabla^{2}c^{\varepsilon,\tau}\|_{L^{\frac{10}{3}}}^{\frac{10}{3}}+C\left(\tau,\|\chi\|_{0} \right)\|\nabla^{2}c^{\varepsilon,\tau}\|_{L^{2}}^{2}\|\nabla c^{\varepsilon,\tau}\|_{L^{6}}^{4}.
%	\end{aligned}
%\end{equation*}
Then combining $ G_{1} $ with $ G_{2}, $ we obtain
\beno
&&\frac{d}{dt}\|\nabla^{2}n^{\varepsilon,\tau}\|_{L^{2}}^{2}+\|\nabla^{3}n^{\varepsilon,\tau}\|_{L^{2}}^{2}\\
&\leq&\frac{11}{16}\|\nabla^{3}n^{\varepsilon,\tau}\|_{L^{2}}^{2}+C(\varepsilon)\|u^{\varepsilon,\tau}\|_{L^{2}}^{2}\|\nabla^{2}n^{\varepsilon,\tau}\|_{L^{2}}^{2}+C(\varepsilon)\|\nabla n^{\varepsilon,\tau}\|_{L^{2}}^{2}\|u^{\varepsilon,\tau}\|_{L^{2}}^{2}\\
&+&C\left(\|\chi\|_{0} \right)\|\nabla^{2}n^{\varepsilon,\tau}\|^2_{L^{2}}\|\nabla c^{\varepsilon,\tau}\|^2_{L^{\infty}}+C\left(\tau, \|\chi\|_{0} \right)\|\nabla n^{\varepsilon,\tau}\|^2_{L^{2}}\|\nabla c^{\varepsilon,\tau}\|^2_{L^{\infty}}\\
&+&C\|\nabla n^{\varepsilon,\tau}\|^2_{L^{3}}\|\nabla^{3}c^{\varepsilon,\tau}\|^2_{L^{2}}+\|\nabla n^{\varepsilon,\tau}\|^2_{L^{2}}\|\nabla c^{\varepsilon,\tau}\|_{L^{\infty}}^{2}\\
&+&C\left(\|\nabla^3 c^{\varepsilon,\tau}\|^2_{L^2}+\|\nabla^2 c^{\varepsilon,\tau}\|^2_{L^2}\|\nabla c^{\varepsilon,\tau}\|^2_{L^\infty}+\|\nabla c^{\varepsilon,\tau}\|_{L^6}^6\right).
\eeno
Due to
\beno
\|\nabla c^{\varepsilon,\tau} \|_{L^\infty}\leq C\|\nabla c^{\varepsilon,\tau} \|^{\frac14}_{L^2}\|\nabla^3 c^{\varepsilon,\tau} \|^{\frac34}_{L^2},
\eeno
there holds
\ben\label{nabla c infty}
\sup_t\|\nabla c^{\varepsilon,\tau} \|_{L^\infty}\leq CT^{\frac{329}{16}}+C.
\een
Integrating with respect to time variable $t$ from 0 to $T$,
%\begin{equation*}
%	\begin{aligned}
%&\frac{d}{dt}\|\nabla^{2}n^{\varepsilon,\tau}\|_{L^{2}}^{2}+\|\nabla^{3}n^{\varepsilon,\tau}\|_{L^{2}}^{2}\\\leq~&
%\frac16\|\nabla c^{\varepsilon,\tau}\|_{L^{4}}^{4}+C\left(\varepsilon,\tau,\|\chi\|_{0}  \right)\|\nabla^{2}n^{\varepsilon,\tau}\|_{L^{2}}^{2}\\&\times\left(\|\nabla^{3}c^{\varepsilon,\tau}\|_{L^{2}}^{2}+\|\nabla c^{\varepsilon,\tau}\|_{L^{6}}^{4}+\|\nabla^{2}c^{\varepsilon,\tau}\|_{L^{2}}^{2}\|\nabla n^{\varepsilon,\tau}\|_{L^{6}}^{2}+1+\|u^{\varepsilon,\tau}\|_{L^{2}}^{2} \right)\\&+~C\|\nabla n^{\varepsilon,\tau}\|_{L^{2}}^{2}\left(\|\nabla^{2}c^{\varepsilon,\tau}\|_{L^{2}}^{2}+\|u^{\varepsilon,\tau}\|_{L^{2}}^{2}\right)+C\left(\tau,\|\chi\|_{0} \right)\|\nabla c^{\varepsilon,\tau}\|_{L^{2}}^{2}\\&+~C\left(\tau,\|\chi\|_{0})  \right)\|\nabla^{2}c^{\varepsilon,\tau}\|_{L^{\frac{10}{3}}}^{\frac{10}{3}}+C\left(\tau,\|\chi\|_{0} \right)\|\nabla^{2}c^{\varepsilon,\tau}\|_{L^{2}}^{2}\|\nabla c^{\varepsilon,\tau}\|_{L^{6}}^{4}.
%%&+C(\varepsilon)\|u^{\varepsilon,\tau}\|_{L^{2}}^{2}\|\nabla^{2}n^{\varepsilon,\tau}\|_{L^{2}}^{2}+C(\varepsilon)\|\nabla n^{\varepsilon,\tau}\|_{L^{2}}^{2}\|u^{\varepsilon,\tau}\|_{L^{2}}^{2}.
%\end{aligned}
%\end{equation*}
by (\ref{u-energy}), (\ref{n'-energy}), (\ref{n'3-energy}), (\ref{c'''-energy}) and (\ref{nabla c infty}), we can get (\ref{n''-energy}).
To sum up the above estimates,  for any $ t\in(0,T), $ we achieve (\ref{uniform estimation}).

Then we complete the proof of Theorem \ref{lem:2}.

%By (\ref{u''-energy}), (\ref{n'-energy}), (\ref{c''-energy}) and (\ref{c'' half energy}),

\section{Uniform estimates for the regularized system}
For the problem (\ref{eq:CNS}), one can obtain the following uniform estimates independent of $\varepsilon$ and $\tau$.
\begin{proposition}\label{lem:same energy}
	Let $(n_0^{\varepsilon},c_0^{\varepsilon}, u_0^{\varepsilon})$ satisfying%$s > \frac 32$
	\ben\label{ine:initial condition assumption}
	\left\{
	\begin{array}{llll}
		\displaystyle n_0^{\varepsilon} \in L^1,\quad (n_0^{\varepsilon}+1)\ln (n_0^{\varepsilon}+1)\in L^1,\quad u_0^{\varepsilon} \in L^2_\sigma;\\
		\displaystyle \nabla \sqrt{c_0^{\varepsilon}} \in L^2, \quad c_0^{\varepsilon} \in L^1 \cap L^\infty ;\\
		\displaystyle n_0^{\varepsilon} \geq 0, \quad c_0^{\varepsilon} \geq 0.\\
	\end{array}
	\right.
	\een
	Assume that $(n^{\varepsilon,\tau},c^{\varepsilon,\tau},u^{\varepsilon,\tau})$ is a solution of (\ref{eq:CNS}) as in Theorem \ref{lem:2}, $\kappa$ and $\chi$ satisfy $\eqref{ine:chi}$, $\eqref{ine:kappa}$ and $ \eqref{equ:chi kappa}. $ 
	Then there exists a constant $C > 0$ independent of $\varepsilon$ and $\tau$ such that
%	\ben\label{eq:n L^infty_t L^1}
%	||n^{\varepsilon,\tau}(t)||_{L^1} = ||n_0||_{L^1};
%	\een
%	\ben\label{ine:c L^1 L^infty'}
%	||c^{\varepsilon,\tau}(t)||_{L^1 \cap L^\infty} \leq ||c_0||_{L^1 \cap L^\infty};
%	\een
	\ben\label{ine:  energy}
	&&\mathbb{U^{\varepsilon,\tau}}(t) + \int_0^t \mathbb{V^{\varepsilon,\tau}}(\tau) d\tau\\\nonumber& \leq& C(\|\nabla \phi\|_{L^\infty},\|n_0^{\varepsilon}\|_{L^1},\|c_{0}^{\varepsilon}\|_{L^{\infty}\cap L^{1}},\|u_{0}^{\varepsilon}\|_{L^{2}},\|(n_{0}^{\varepsilon}+1)\ln (n_{0}^{\varepsilon}+1)\|_{L^{1}},\|\nabla \sqrt{c_0^{\varepsilon}}\|_{L^{2}})(1+t),
	\een
	where
	\ben\label{eq:def U}
	\mathbb{U^{\varepsilon,\tau}}(t) = ||n^{\varepsilon,\tau}||_{L^1} + \|(n^{\varepsilon,\tau}+1)\ln (n^{\varepsilon,\tau}+1)\|_{L^1}+ \frac{2}{\Theta_0}||\nabla \sqrt{c^{\varepsilon,\tau}}||_{L^{2}}^{2} + ||u^{\varepsilon,\tau}||_{L^{2}}^{2},
	\een
	and
	\ben\label{eq:def V}
	\mathbb{V^{\varepsilon,\tau}}(t) \nonumber&=&||\nabla \sqrt{n^{\varepsilon,\tau}+1}||_{L^{2}}^{2} + \frac{4}{3\Theta_0}||\Delta \sqrt{c^{\varepsilon,\tau}}||_{L^{2}}^{2} + ||\nabla u^{\varepsilon,\tau}||_{L^{2}}^{2}
	\\&&+~\frac{1}{3\Theta_0}\int_{\mathbb{R}^3} (\sqrt{c^{\varepsilon,\tau}})^{-2} |\nabla \sqrt{c^{\varepsilon,\tau}}|^4.
	\een
\end{proposition}

{\bf Proof.} 
% For the equation of (\ref{eq:n L^infty_t L^1}) and (\ref{ine:c L^1 L^infty'}),
%recall the first equation of $\eqref{eq:CNS}$ as follow:
%\beno
%\int_{\mathbb{R}^3}\partial_t n^{\varepsilon,\tau} + \int_{\mathbb{R}^3}(u^{\varepsilon,\tau} \ast \rho^\varepsilon) \cdot \nabla n^{\varepsilon,\tau} - \int_{\mathbb{R}^3}\Delta n^{\varepsilon,\tau} + 
% \int_{\mathbb{R}^3}\nabla \cdot \left(\frac{n^{\varepsilon,\tau}}{1+\tau n^{\varepsilon,\tau}} \nabla c^{\varepsilon,\tau} \chi(c^{\varepsilon,\tau})\right)=0.
%%\int_{\mathbb{R}^3} \nabla \cdot (n^\varepsilon (\nabla c^\varepsilon\chi(c^{\varepsilon})) \ast \rho^\varepsilon)=0.
%\eeno
%Integration by parts and  divergence theorem yields that
%\beno
%\frac d{dt}\int_{\mathbb{R}^3}n^{\varepsilon,\tau}=0,
%\eeno
%which means for any $t \in (0,T_\ast]$,
%\beno
%\|n^{\varepsilon,\tau}(t)\|_{L^1}-\|n^{\varepsilon,\tau}(.,0)\|_{L^1}=0.
%\eeno
%Similarly, we have
%\ben\label{ine:c L 1}
%\|c^{\varepsilon,\tau}(t)\|_{L^1}-\|c^{\varepsilon,\tau}(.,0)\|_{L^1}\leq 0.
%\een
%For the $L^\infty-$norm, since $\kappa(s) > 0$, by the maximum principle, we know that
%\beno
%\|c^{\varepsilon,\tau}\|_{L^\infty} \leq \|c_0^{\varepsilon,\tau}\|_{L^\infty}.
%\eeno
%Next,
 First, we  rewrite the equation of $n^{\varepsilon,\tau}$ in $(\ref{eq:CNS})_1$ as
% \beno
% \partial_t n^{\varepsilon,\tau} + (u^{\varepsilon,\tau} \ast \rho^\varepsilon) \cdot \nabla n^{\varepsilon,\tau} - \Delta n^{\varepsilon,\tau} = - \nabla \cdot \left(\frac{n^{\varepsilon,\tau}}{1+\tau n^{\varepsilon,\tau}} \nabla c^{\varepsilon,\tau} \chi(c^{\varepsilon,\tau})\right)
% \eeno
\ben\label{n'}
&&\partial_t(n^{\varepsilon,\tau}+1)+(u^{\varepsilon,\tau}\ast\rho^\varepsilon)\cdot \nabla (n^{\varepsilon,\tau}+1)-\Delta (n^{\varepsilon,\tau}+1)\nonumber\\&&=- \nabla \cdot \left(\frac{n^{\varepsilon,\tau}+1}{1+\tau n^{\varepsilon,\tau}} \nabla c^{\varepsilon,\tau} \chi(c^{\varepsilon,\tau})\right)+ \nabla \cdot \left(\frac{1}{1+\tau n^{\varepsilon,\tau}} \nabla c^{\varepsilon,\tau} \chi(c^{\varepsilon,\tau})\right).
%-\nabla\cdot(\chi(c^{\varepsilon,\tau})(n^{\varepsilon,\tau}+1)\nabla(c^{\varepsilon,\tau}\ast\rho^\varepsilon))+\nabla\cdot(\chi(c^{\varepsilon,\tau})\nabla(c^\varepsilon\ast\rho^\varepsilon)).
\een
Multiplying $(1+\ln (n^{\varepsilon,\tau}+1))$ in $(\ref{n'})$, we have
\begin{equation*}
	\begin{aligned}
	&\int_{\mathbb{R}^3}\partial_t(n^{\varepsilon,\tau}+1) (1+\ln (n^{\varepsilon,\tau}+1)) +\int_{\mathbb{R}^3}(u^{\varepsilon,\tau}\ast\rho^\varepsilon)\cdot \nabla (n^{\varepsilon,\tau}+1) (1+\ln (n^{\varepsilon,\tau}+1)) \\
	-&\int_{\mathbb{R}^3}\Delta (n^{\varepsilon,\tau}+1)(1+\ln (n^{\varepsilon,\tau}+1)) +\int_{\mathbb{R}^3}\nabla \cdot \left(\frac{n^{\varepsilon,\tau}+1}{1+\tau n^{\varepsilon,\tau}} \nabla c^{\varepsilon,\tau} \chi(c^{\varepsilon,\tau})\right) (1+\ln (n^{\varepsilon,\tau}+1)) \\
	-&\int_{\mathbb{R}^3}\nabla \cdot \left(\frac{1}{1+\tau n^{\varepsilon,\tau}} \nabla c^{\varepsilon,\tau} \chi(c^{\varepsilon,\tau})\right)(1+\ln (n^{\varepsilon,\tau}+1)) =0.
	\end{aligned}
\end{equation*}
Integration by parts, we get
\ben\label{ine:energy estimate of n}
&& \frac{d}{dt}\int_{\mathbb{R}^3} (n^{\varepsilon,\tau}+1) \ln (n^{\varepsilon,\tau}+1)(\cdot,t)dx + \int_{\mathbb{R}^3}\frac1{n^{\varepsilon,\tau}+1} |\nabla (n^{\varepsilon,\tau}+1)|^2\\\nonumber
&=&\int_{\mathbb{R}^{3}}\frac1{1+\tau n^{\varepsilon,\tau}}\nabla c^{\varepsilon,\tau}\chi(c^{\varepsilon,\tau})\cdot\nabla n^{\varepsilon,\tau}
%&&=-\int_{\mathbb{R}^3}\nabla \cdot \left(\frac{n^{\varepsilon,\tau}+1}{1+\tau n^{\varepsilon,\tau}} \nabla c^{\varepsilon,\tau} \chi(c^{\varepsilon,\tau})\right) (1+\ln (n^{\varepsilon,\tau}+1)) dx\\\nonumber
+\int_{\mathbb{R}^3}\nabla \cdot \left(\frac{1}{1+\tau n^{\varepsilon,\tau}} \nabla c^{\varepsilon,\tau} \chi(c^{\varepsilon,\tau})\right)\ln (n^{\varepsilon,\tau}+1) .
\een
 %Firstly, multiplying $ (1+\ln n^{\varepsilon,\tau}) $ in $ (\ref{eq:CNS})_{1}, $   we have
%\beno
%&&\int_{\mathbb{R}^3}\partial_t n^{\varepsilon,\tau} (1+\ln n^{\varepsilon,\tau}) dx+\int_{\mathbb{R}^3}(u^{\varepsilon,\tau}\ast\rho^\varepsilon)\cdot \nabla n^{\varepsilon,\tau}\cdot (1+\ln n^{\varepsilon,\tau}) dx\\
%&&-\int_{\mathbb{R}^3}\Delta n^{\varepsilon,\tau}(1+\ln n^{\varepsilon,\tau}) dx+\int_{\mathbb{R}^3} \nabla \cdot \left(\frac{n^{\varepsilon,\tau}}{1+\tau n^{\varepsilon,\tau}} \nabla c^{\varepsilon,\tau} \chi(c^{\varepsilon,\tau})\right) (1+\ln n^{\varepsilon,\tau}) dx=0.\\
%\eeno
%Integration by parts, we have
%\ben\label{ine:energy estimate of n}
%\frac{d}{dt}\int_{\mathbb{R}^3} n^{\varepsilon,\tau} \ln n^{\varepsilon,\tau}(\cdot,t)dx+4\int_{\mathbb{R}^{3}}|\nabla \sqrt{n^{\varepsilon,\tau}}|^{2}=\int_{\mathbb{R}^{3}}\frac1{1+\tau n^{\varepsilon,\tau}}(\nabla c^{\varepsilon,\tau}\chi(c^{\varepsilon,\tau}))\cdot\nabla n^{\varepsilon,\tau}.
%\een

For the equation of $c^{\varepsilon,\tau}$ in (\ref{eq:CNS}), noting that
\beno
\Delta c^{\varepsilon,\tau} = 2 |\nabla \sqrt{c^{\varepsilon,\tau}}|^2 + 2 \sqrt{c^{\varepsilon,\tau}} \Delta \sqrt{c^{\varepsilon,\tau}},
\eeno
and dividing $2\sqrt{c^{\varepsilon,\tau}}$ on both sides, we obtain
\ben\label{eq:c tilde}
\partial_t\sqrt{c^{\varepsilon,\tau}}+u^{\varepsilon,\tau}\cdot \nabla\sqrt{c^{\varepsilon,\tau}}-\frac{|\nabla \sqrt{c^{\varepsilon,\tau}}|^2}{\sqrt{c^{\varepsilon,\tau}}}-\Delta\sqrt{c^{\varepsilon,\tau}}= - \frac{\kappa(c^{\varepsilon,\tau})}{2\tau\sqrt{c^{\varepsilon,\tau}}}~~\ln(1+\tau n^{\varepsilon,\tau}).
\een
Multiplying the above equation (\ref{eq:c tilde}) by $-\partial_i (\partial_i \sqrt{c^{\varepsilon,\tau}})$ , there holds
\begin{equation}\label{eq:c}
	\begin{aligned}
	\frac 12 \frac{d}{dt} ||\nabla \sqrt{c^{\varepsilon,\tau}(t)}||_{L^{2}}^{2} +||\nabla^2 \sqrt{c^{\varepsilon,\tau}(t)}||_{L^{2}}^{2} =& - \int_{\mathbb{R}^3} (\sqrt{c^{\varepsilon,\tau}})^{-1} |\nabla \sqrt{c^{\varepsilon,\tau}}|^2 \Delta \sqrt{c^{\varepsilon,\tau}} dx \\ &+ \int_{\mathbb{R}^3} u^{\varepsilon,\tau}  \cdot \nabla \sqrt{c^{\varepsilon,\tau}} \Delta \sqrt{c^{\varepsilon,\tau}} dx \\ &+\frac 12 \int_{\mathbb{R}^3}\frac{\kappa(c^{\varepsilon,\tau})}{\tau\sqrt{c^{\varepsilon,\tau}}}~~\ln(1+\tau n^{\varepsilon,\tau})\Delta \sqrt{c^{\varepsilon,\tau}} dx\\
	:=~& I_1 + I_2 + I_3.
	\end{aligned}
\end{equation}
The key is to estimate the term $I_1$. By integration by parts, we have
\begin{equation*}
	\begin{aligned}
	I_1 =& - \int_{\mathbb{R}^3} (\sqrt{c^{\varepsilon,\tau}})^{-1} (\partial_j \sqrt{c^{\varepsilon,\tau}})^2 \partial_{ii} \sqrt{c^{\varepsilon,\tau}}   \\
	=& - \int_{\mathbb{R}^3} (\sqrt{c^{\varepsilon,\tau}})^{-2} (\partial_j \sqrt{c^{\varepsilon,\tau}})^2 (\partial_i \sqrt{c^{\varepsilon,\tau}})^2  + 2 \int_{\mathbb{R}^3} (\sqrt{c^{\varepsilon,\tau}})^{-1} \partial_{ij} \sqrt{c^{\varepsilon,\tau}} \partial_i \sqrt{c^{\varepsilon,\tau}} \partial_j \sqrt{c^{\varepsilon,\tau}}   \\
	=& - \sum_{i,j} \int_{\mathbb{R}^3} (\sqrt{c^{\varepsilon,\tau}})^{-2} (\partial_j \sqrt{c^{\varepsilon,\tau}})^2 (\partial_i \sqrt{c^{\varepsilon,\tau}})^2   + 2 \sum_{i = j} \int_{\mathbb{R}^3} (\sqrt{c^{\varepsilon,\tau}})^{-1} \partial_{ij} \sqrt{c^{\varepsilon,\tau}} \partial_i \sqrt{c^{\varepsilon,\tau}} \partial_j \sqrt{c^{\varepsilon,\tau}}  \\
	&+~ 2 \sum_{i \neq j} \int_{\mathbb{R}^3} (\sqrt{c^{\varepsilon,\tau}})^{-1} \partial_{ij} \sqrt{c^{\varepsilon,\tau}} \partial_i \sqrt{c^{\varepsilon,\tau}} \partial_j \sqrt{c^{\varepsilon,\tau}}.
	\end{aligned}
\end{equation*}
Noting that
\beno
&&\sum_{i = j} \int_{\mathbb{R}^3} (\sqrt{c^{\varepsilon,\tau}})^{-1} \partial_{ij} \sqrt{c^{\varepsilon,\tau}} \partial_i \sqrt{c^{\varepsilon,\tau}} \partial_j \sqrt{c^{\varepsilon,\tau}} \\&&= - I_1 - \sum_{i \neq j} \int_{\mathbb{R}^3} (\sqrt{c^{\varepsilon,\tau}})^{-1} \partial_{ii} \sqrt{c^{\varepsilon,\tau}} \partial_j \sqrt{c^{\varepsilon,\tau}} \partial_j \sqrt{c^{\varepsilon,\tau}} ,
\eeno
we have
\beno
I_1 &=& - \frac 13 \sum_{i,j} \int_{\mathbb{R}^3} (\sqrt{c^{\varepsilon,\tau}})^{-2} (\partial_j \sqrt{c^{\varepsilon,\tau}})^2 (\partial_i \sqrt{c^{\varepsilon,\tau}})^2 \\&&- \frac 23 \sum_{i \neq j} \int_{\mathbb{R}^3} (\sqrt{c^{\varepsilon,\tau}})^{-1} \partial_{ii} \sqrt{c^{\varepsilon,\tau}} \partial_j \sqrt{c^{\varepsilon,\tau}} \partial_j \sqrt{c^{\varepsilon,\tau}}  \\
&&+ \frac 23 \sum_{i \neq j} \int_{\mathbb{R}^3} (\sqrt{c^{\varepsilon,\tau}})^{-1} \partial_{ij} \sqrt{c^{\varepsilon,\tau}} \partial_i \sqrt{c^{\varepsilon,\tau}} \partial_j \sqrt{c^{\varepsilon,\tau}}.
\eeno
Using Young inequality it follows that
%of $ab\leq \epsilon{a^2}+\frac{b^2}{4\epsilon}$, letting $\epsilon=\frac12$,
\beno
\frac 23 \sum_{i \neq j} \int_{\mathbb{R}^3} (\sqrt{c^{\varepsilon,\tau}})^{-1} \partial_{ij} \sqrt{c^{\varepsilon,\tau}} \partial_i \sqrt{c^{\varepsilon,\tau}} \partial_j \sqrt{c^{\varepsilon,\tau}}  &\leq& \frac 13 \sum_{i \neq j} \int_{\mathbb{R}^3} (\sqrt{c^{\varepsilon,\tau}})^{-2} |\partial_i \sqrt{c^{\varepsilon,\tau}}|^2 |\partial_j \sqrt{c^{\varepsilon,\tau}}|^2  \\
&&+ \frac 1{3} \sum_{i \neq j} \int_{\mathbb{R}^3} |\partial_{ij} \sqrt{c^{\varepsilon,\tau}}|^2 ,
\eeno
and
\beno
&&- \frac 23 \sum_{i \neq j} \int_{\mathbb{R}^3} (\sqrt{c^{\varepsilon,\tau}})^{-1} \partial_{ii} \sqrt{c^{\varepsilon,\tau}} \partial_j \sqrt{c^{\varepsilon,\tau}} \partial_j \sqrt{c^{\varepsilon,\tau}} \\
&=&- \frac 23 \sum_{i\neq j} \int_{\mathbb{R}^3} (\sqrt{c^{\varepsilon,\tau}})^{-2} (\partial_j \sqrt{c^{\varepsilon,\tau}})^2 (\partial_i \sqrt{c^{\varepsilon,\tau}})^2 +\frac 43 \sum_{i \neq j} \int_{\mathbb{R}^3} (\sqrt{c^{\varepsilon,\tau}})^{-1} \partial_{ij} \sqrt{c^{\varepsilon,\tau}} \partial_i \sqrt{c^{\varepsilon,\tau}} \partial_j \sqrt{c^{\varepsilon,\tau}} \\
&\leq& \frac 23  \sum_{i\neq j} \int_{\mathbb{R}^3} |\partial_{ij} \sqrt{c^{\varepsilon,\tau}}|^2.
\eeno
Then we have
\ben\label{I_1}
I_1 &\leq& - \frac 1{3} \sum_{i=j} \int_{\mathbb{R}^3} (\sqrt{c^{\varepsilon,\tau}})^{-2} (\partial_j \sqrt{c^{\varepsilon,\tau}})^2 (\partial_i \sqrt{c^{\varepsilon,\tau}})^2 + \sum_{i\neq j} \int_{\mathbb{R}^3} |\partial_{ij} \sqrt{c^{\varepsilon,\tau}}|^2.\nonumber\\
\een
Using $\kappa(0)=0$ and mean value theorem, for some $\eta \in (0,1)$, there holds
$$
\kappa(c^{\varepsilon,\tau}) = c^{\varepsilon,\tau} \frac{\kappa(c^{\varepsilon,\tau}) - \kappa(0)}{c^{\varepsilon,\tau}} = c^{\varepsilon,\tau} \kappa'(\eta c^{\varepsilon,\tau}).
$$
Noting that $\kappa''(s) \geq 0$ and $\kappa'(s) \geq 0$ from $\eqref{ine:kappa}$, there means $\kappa'(0) \geq 0$ and
$$
\kappa'(\eta c^{\varepsilon}) \leq 2\kappa'(c^{\varepsilon}),
$$
%By mean value theorem and the monotony of $\kappa'$,  using $$2\kappa'(c^{\varepsilon,\tau})-\kappa'(\eta c^{\varepsilon,\tau})>2\kappa'(\eta c^{\varepsilon,\tau})-\kappa'(\eta c^{\varepsilon,\tau})>\frac12\kappa'(\eta c^{\varepsilon,\tau})>\frac12\kappa'(0),$$ 
and we have
\begin{equation}\label{k''}
\begin{aligned}
&\int_{\mathbb{R}^3}\nabla\left(\frac{\kappa(c^{\varepsilon,\tau})}{\sqrt{c^{\varepsilon,\tau}}}\right)~~\ln(1+\tau n^{\varepsilon,\tau})\cdot \nabla\sqrt{c^{\varepsilon,\tau}}
\\=~&\frac1{\tau}\int_{\mathbb{R}^3}\ln(1+\tau n^{\varepsilon,\tau})\kappa'(c^{\varepsilon,\tau})|\nabla\sqrt{c^{\varepsilon,\tau}}|^{2}-\frac1{2\tau}\int_{\mathbb{R}^3}\ln(1+\tau n^{\varepsilon,\tau})\kappa(c^{\varepsilon,\tau})(\sqrt{c^{\varepsilon,\tau}})^{-2}|\nabla\sqrt{c^{\varepsilon,\tau}}|^{2}\\=~&\frac1{\tau}\int_{\mathbb{R}^3}\ln(1+\tau n^{\varepsilon,\tau})\kappa'(c^{\varepsilon,\tau})|\nabla\sqrt{c^{\varepsilon,\tau}}|^{2}-\frac1{2\tau}\int_{\mathbb{R}^3}\ln(1+\tau n^{\varepsilon,\tau})\kappa'(\eta c^{\varepsilon,\tau})|\nabla\sqrt{c^{\varepsilon,\tau}}|^{2}\\
\geq~& 0.
\end{aligned}
\end{equation}
Then for $I_3$, integration by parts yields that

\begin{equation}\label{I 3}
\begin{aligned}
I_3=~&\frac 12 \int_{\mathbb{R}^3}\frac{\kappa(c^{\varepsilon,\tau})}{\tau\sqrt{c^{\varepsilon,\tau}}}~~\ln(1+\tau n^{\varepsilon,\tau})\Delta \sqrt{c^{\varepsilon,\tau}} \\
=~&-\frac 12 \int_{\mathbb{R}^3}  \ln(1+\tau n^{\varepsilon,\tau})\nabla(\frac{\kappa( c^{\varepsilon,\tau})}{\tau\sqrt{c^{\varepsilon,\tau}}})\cdot\nabla\sqrt{c^{\varepsilon,\tau}} -  \frac 12 \int_{\mathbb{R}^3}\frac{\kappa(c^{\varepsilon,\tau})}{\tau \sqrt{c^{\varepsilon,\tau}}}\nabla(\ln(1+\tau n^{\varepsilon,\tau}))\cdot \nabla \sqrt{c^{\varepsilon,\tau}} \\\leq~&- \frac 1{4\tau} \int_{\mathbb{R}^3}\frac{\kappa(c^{\varepsilon,\tau})}{c^{\varepsilon,\tau}}\nabla(\ln(1+\tau n^{\varepsilon,\tau}))\cdot \nabla c^{\varepsilon,\tau}\\=~& -\frac 1{4} \int_{\mathbb{R}^3}\frac{\kappa(c^{\varepsilon,\tau})}{c^{\varepsilon,\tau}}\frac1{1+\tau n^{\varepsilon,\tau}}\nabla n^{\varepsilon,\tau}\cdot \nabla c^{\varepsilon,\tau}.
\end{aligned}
\end{equation}
Collecting $\eqref{eq:c}$, $\eqref{I_1}$ and $\eqref{I 3}$, we have
%and noting the relation $||\nabla^2 \sqrt{c^{\varepsilon,\tau}}||_2 = ||\Delta \sqrt{c^{\varepsilon,\tau}}||_2$,
\beno\label{ine:energy estimate of c}\nonumber
&&\frac 12 \frac{d}{dt} ||\nabla \sqrt{c^{\varepsilon,\tau}(t)}||_{L^{2}}^{2} + ||\nabla^2 \sqrt{c^{\varepsilon,\tau}(t)}||_{L^{2}}^{2}-\sum_{i\neq j} \int_{\mathbb{R}^3} |\partial_{ij} \sqrt{c^{\varepsilon,\tau}}|^2 \\ \nonumber
&&+\frac 1{3} \sum_{i=j} \int_{\mathbb{R}^3} (\sqrt{c^{\varepsilon,\tau}})^{-2} (\partial_j \sqrt{c^{\varepsilon,\tau}})^2 (\partial_i \sqrt{c^{\varepsilon,\tau}})^2 \\ \nonumber
&&\leq \int_{\mathbb{R}^3} u^{\varepsilon,\tau}  \cdot \nabla \sqrt{c^{\varepsilon,\tau}} \Delta \sqrt{c^{\varepsilon,\tau}} 
- \frac 1{4} \int_{\mathbb{R}^3}\frac{\kappa(c^{\varepsilon,\tau})}{c^{\varepsilon,\tau}}\frac1{1+\tau n^{\varepsilon,\tau}}\nabla n^{\varepsilon,\tau}\cdot \nabla c^{\varepsilon,\tau} .
\eeno
Since
\beno
\frac 13 \int_{\mathbb{R}^3} |\Delta \sqrt{c^{\varepsilon,\tau}}|^2 \leq \sum_{i=1}^3 \int_{\mathbb{R}^3} |\partial_{ii} \sqrt{c^{\varepsilon,\tau}}|^2 = ||\nabla^2 \sqrt{c^{\varepsilon,\tau}(t)}||_{L^{2}}^{2} - \sum_{i\neq j} \int_{\mathbb{R}^3} |\partial_{ij} \sqrt{c^{\varepsilon,\tau}}|^2,
\eeno
we have
\ben\label{ine:energy estimate of c}\nonumber
&&\frac 12 \frac{d}{dt} ||\nabla \sqrt{c^{\varepsilon,\tau}(t)}||_{L^{2}}^{2} + \frac13 ||\Delta \sqrt{c^{\varepsilon,\tau}(t)}||_{L^{2}}^{2}+\frac 1{3} \sum_{i=j} \int_{\mathbb{R}^3} (\sqrt{c^{\varepsilon,\tau}})^{-2} (\partial_j \sqrt{c^{\varepsilon,\tau}})^2 (\partial_i \sqrt{c^{\varepsilon,\tau}})^2\\
%+\frac14\kappa'(0)\int_{\mathbb{R}^3}(n^\varepsilon \ast \rho^\varepsilon)|\nabla\sqrt{c^{\varepsilon}}|^{2}\\
&&\leq \int_{\mathbb{R}^3} u^{\varepsilon,\tau}  \cdot \nabla \sqrt{c^{\varepsilon,\tau}} \Delta \sqrt{c^{\varepsilon,\tau}} 
- \frac 1{4} \int_{\mathbb{R}^3}\frac{\kappa(c^{\varepsilon,\tau})}{c^{\varepsilon,\tau}}\frac1{1+\tau n^{\varepsilon,\tau}}\nabla n^{\varepsilon,\tau}\cdot \nabla c^{\varepsilon,\tau} .
\een
Then $(\ref{ine:energy estimate of n})+\frac{4}{\Theta_0}\times(\ref{ine:energy estimate of c})$ yields that
\ben\label{n and c}\nonumber
&& \frac{d}{dt}\int_{\mathbb{R}^3} (n^{\varepsilon,\tau}+1) \ln (n^{\varepsilon,\tau}+1)(\cdot,t)dx + \int_{\mathbb{R}^3}\frac1{n^{\varepsilon,\tau}+1} |\nabla (n^{\varepsilon,\tau}+1)|^2\\\nonumber
&&+~ \frac2{\Theta_0}\frac{d}{dt} ||\nabla \sqrt{c^{\varepsilon,\tau}(t)}||_{L^{2}}^{2} 
+ \frac4{3\Theta_0}||\Delta \sqrt{c^{\varepsilon,\tau}(t)}||_{L^{2}}^{2}\\\nonumber
&& +~ \frac 4{3\Theta_0} \sum_{i=j} \int_{\mathbb{R}^3} (\sqrt{c^{\varepsilon,\tau}})^{-2} (\partial_j \sqrt{c^{\varepsilon,\tau}})^2 (\partial_i \sqrt{c^{\varepsilon,\tau}})^2 \\ \nonumber
&\leq&
%\int_{\mathbb{R}^{3}}\frac1{1+\tau n^{\varepsilon,\tau}}(\nabla c^{\varepsilon,\tau}\chi(c^{\varepsilon,\tau}))\cdot\nabla (n^{\varepsilon,\tau}+1)\\\nonumber
\int_{\mathbb{R}^3}\nabla \cdot \left(\frac{1}{1+\tau n^{\varepsilon,\tau}} \nabla c^{\varepsilon,\tau} \chi(c^{\varepsilon,\tau})\right)\ln (n^{\varepsilon,\tau}+1) dx\\
&&+~\frac4{\Theta_0}\int_{\mathbb{R}^3} u^{\varepsilon,\tau}  \cdot \nabla \sqrt{c^{\varepsilon,\tau}} \Delta \sqrt{c^{\varepsilon,\tau}}:= J_1+J_2.
%- \frac 1{C_0} \int_{\mathbb{R}^3}\frac{\kappa(c^{\varepsilon,\tau})}{c^{\varepsilon,\tau}}\frac1{1+\tau n^{\varepsilon,\tau}}\nabla n^{\varepsilon,\tau}\cdot \nabla c^{\varepsilon,\tau}
\een
%\ben\label{n and c}\nonumber
%&& \frac{d}{dt}\int_{\mathbb{R}^3} n^{\varepsilon,\tau} \ln n^{\varepsilon,\tau}(\cdot,t) + 4\int_{\mathbb{R}^3}|\nabla \sqrt{n^{\varepsilon,\tau}}|^{2}+ \frac2{C_0}\frac{d}{dt} ||\nabla \sqrt{c^{\varepsilon,\tau}(t)}||_2^2 \\ \nonumber
%&&+ \frac4{3C_0}||\Delta \sqrt{c^{\varepsilon,\tau}(t)}||_2^2 + \frac 4{3C_0} \sum_{i=j} \int_{\mathbb{R}^3} (\sqrt{c^{\varepsilon,\tau}})^{-2} (\partial_j \sqrt{c^{\varepsilon,\tau}})^2 (\partial_i \sqrt{c^{\varepsilon,\tau}})^2 \\ \nonumber
%%&&+\frac1{C_0}\kappa'(0)\int_{\mathbb{R}^3}(n^\varepsilon \ast \rho^\varepsilon)|\nabla\sqrt{c^{\varepsilon}}|^{2}\\
%&&\leq \frac4{C_0}\int_{\mathbb{R}^3} (u^{\varepsilon,\tau} \ast \rho^\varepsilon) \cdot \nabla \sqrt{c^{\varepsilon,\tau}} \Delta \sqrt{c^{\varepsilon,\tau}} dx:= J.
%%\int_{\mathbb{R}^{3}}\frac1{1+\varepsilon n^\varepsilon}(\nabla c^{\varepsilon}\chi(c^{\varepsilon}))\cdot\nabla n^{\varepsilon}
%\een
{\bf Estimate of $J_1.$} 
Noting that $ \frac{1}{1+\tau n^{\varepsilon,\tau}}\leq 1, $ and $ \|\chi(c^{\varepsilon,\tau})\|_{L^{\infty}}\leq \|\chi\|_{L^{\infty}(0,\|c^{\varepsilon}_{0}\|_{L^{\infty}})}\leq C,$ 
by integration by parts, H\"{o}lder inequality, Young inequality, (\ref{ine:nc>0}) and (\ref{ine:c L^1 L^infty'}), there holds
\ben\label{ine:J_1}\nonumber
&J_1&=\int_{\mathbb{R}^3}\nabla \cdot \left(\frac{1}{1+\tau n^{\varepsilon,\tau}} \nabla c^{\varepsilon,\tau} \chi(c^{\varepsilon,\tau})\right)\ln (n^{\varepsilon,\tau}+1) dx\\\nonumber
&&=-2\int_{\mathbb{R}^3}\frac1{1+\tau n^{\varepsilon,\tau}}(\sqrt{c^{\varepsilon,\tau}})^{\frac32}(\sqrt{c^{\varepsilon,\tau}})^{-\frac12}\nabla \sqrt{c^{\varepsilon,\tau}} \chi(c^{\varepsilon,\tau})\frac{\nabla n^{\varepsilon,\tau}}{1+n^{\varepsilon,\tau}}dx\\\nonumber
%&&\leq C\int_{\mathbb{R}^3}\left|(\sqrt{c^{\varepsilon,\tau}})^{\frac32}(\sqrt{c^{\varepsilon,\tau}})^{-\frac12}\nabla \sqrt{c^{\varepsilon,\tau}} \frac{\nabla n^{\varepsilon,\tau}}{(1+n^{\varepsilon,\tau})^2}\right|dx\\\nonumber
&&\leq C\|\sqrt{c^{\varepsilon,\tau}}^{\frac32}\|_{L^4}\left(\int_{\mathbb{R}^3}(\sqrt{c^{\varepsilon,\tau}})^{-2}|\nabla \sqrt{c^{\varepsilon,\tau}}|^4dx\right)^{\frac14}\left(\int_{\mathbb{R}^3}\left|\frac{\nabla n^{\varepsilon,\tau}}{(1+n^{\varepsilon,\tau})^2}\right|^2dx\right)^{\frac12}\\\nonumber
&&\leq C\|c^{\varepsilon,\tau}\|^{\frac14}_{L^1}\|c^{\varepsilon,\tau}\|^{\frac12}_{L^\infty}\left(\int_{\mathbb{R}^3}(\sqrt{c^{\varepsilon,\tau}})^{-2}|\nabla \sqrt{c^{\varepsilon,\tau}}|^4dx\right)^{\frac14}\left(\int_{\mathbb{R}^3}\left|\frac{\nabla n^{\varepsilon,\tau}}{(1+n^{\varepsilon,\tau})^2}\right|^2dx\right)^{\frac12}\\\nonumber
%&&\leq\epsilon_0\int_{\mathbb{R}^3}|\frac{\nabla n^{\varepsilon,\tau}}{(1+n^{\varepsilon,\tau})^2}|^2dx+C(\epsilon_0,\|c_{0}\|_{L^{\infty}},\|c_{0}\|_{L^{1}})\left(\int_{\mathbb{R}^3}(\sqrt{c^{\varepsilon,\tau}})^{-2}|\nabla \sqrt{c^{\varepsilon,\tau}}|^4dx\right)^{\frac12}\\\nonumber
%&&\leq\epsilon_0\int_{\mathbb{R}^3}|\frac{\nabla n^{\varepsilon,\tau}}{(1+n^{\varepsilon,\tau})^2}|^2dx+\epsilon_0^{'}\int_{\mathbb{R}^3}(\sqrt{c^{\varepsilon,\tau}})^{-2}|\nabla \sqrt{c^{\varepsilon,\tau}}|^4dx+C(\epsilon_0,\epsilon_{0}^{'},\|c_{0}\|_{L^{\infty}},\|c_{0}\|_{L^{1}})\\\nonumber
&&\leq\frac14\int_{\mathbb{R}^3}\frac1{1+n^{\varepsilon,\tau}} |\nabla (n^{\varepsilon,\tau}+1)|^2dx+\frac1{27\Theta_0 }\int_{\mathbb{R}^3}(\sqrt{c^{\varepsilon,\tau}})^{-2}|\nabla \sqrt{c^{\varepsilon,\tau}}|^4dx\\
&&~~+C(\Theta_0,\|c^{\varepsilon}_{0}\|_{L^{\infty}},\|c^{\varepsilon}_{0}\|_{L^{1}}).
\een
{\bf Estimate of $J_2.$} 
%Since $\|\Delta (u^{\varepsilon,\tau} \ast \rho^\varepsilon)\|_{L^2} \leq C \|\Delta u^{\varepsilon,\tau}\|_{L^2}$ and $\|\nabla (u^{\varepsilon,\tau} \ast \rho^\varepsilon)\|_{L^2} \leq C \|\nabla u^{\varepsilon,\tau}\|_{L^2}$ with that the constant $C$ is independent of $\tau$, 
For $J_{2}$, we derive
\begin{equation}\label{ine:J_2}
\begin{aligned}
J_{2}&=-\frac{4}{\Theta_0}\int_{\mathbb{R}^{3}}\partial_ju^{\varepsilon,\tau}\partial_{i} \sqrt{c^{\varepsilon,\tau}}\partial_{j} \sqrt{c^{\varepsilon,\tau}}-\frac{4}{\Theta_0}\int_{\mathbb{R}^{3}}u^{\varepsilon,\tau}\partial_{ij} \sqrt{c^{\varepsilon,\tau}}\partial_{j} \sqrt{c^{\varepsilon,\tau}}\\&=-\frac{4}{\Theta_0}\int_{\mathbb{R}^{3}}\partial_ju^{\varepsilon,\tau}\sqrt{c^{\varepsilon,\tau}}(\sqrt{c^{\varepsilon,\tau}})^{-1}\partial_{i} \sqrt{c^{\varepsilon,\tau}}\partial_{j} \sqrt{c^{\varepsilon,\tau}}\\&\leq \frac{8}{27 \Theta_0} \int_{\mathbb{R}^3} (\sqrt{c^{\varepsilon,\tau}})^{-2} |\nabla \sqrt{c^{\varepsilon,\tau}}|^4 dx + \frac{27}{2 \Theta_0} ||c^{\varepsilon}_0||_{L^\infty} \int_{\mathbb{R}^3} |\nabla u^{\varepsilon,\tau}|^2 dx.
\end{aligned}
\end{equation}
Noting that
\beno
\int_{\mathbb{R}^3}  (\sqrt{{c^{\varepsilon,\tau}}})^{-2} |\nabla \sqrt{{c^{\varepsilon,\tau}}}|^4 &=& \sum_{i,j} \int_{\mathbb{R}^3}  (\sqrt{{c^{\varepsilon,\tau}}})^{-2} (\partial_j \sqrt{{c^{\varepsilon,\tau}}})^2 (\partial_i \sqrt{{c^{\varepsilon,\tau}}})^2 \\
&\leq&  3\sum_{i=j} \int_{\mathbb{R}^3}  (\sqrt{{c^{\varepsilon,\tau}}})^{-2} (\partial_j \sqrt{{c^{\varepsilon,\tau}}})^2 (\partial_i \sqrt{{c^{\varepsilon,\tau}}})^2,
\eeno
%choosing $\epsilon_{0}=\frac14, 3 \epsilon_{0}^{'}=\frac1{9 \Theta_0 },3 \epsilon_{1} = \frac{8}{9 \Theta_0},$ 
%and $C(\epsilon_{1})||c_0||_{L^\infty} \leq \frac12$,
 combining $\eqref{ine:J_1}$ and $\eqref{ine:J_2}$, we achieve
%choose $ \epsilon_{1} $ is small enough, combine with the estimates of $ J $, we achieve
\ben\label{ine:energy estimate of c and n} \nonumber
&&\frac{d}{dt}\int_{\mathbb{R}^3} (n^{\varepsilon,\tau}+1) \ln (n^{\varepsilon,\tau}+1)(\cdot,t)dx +3 \int_{\mathbb{R}^3}|\nabla (n^{\varepsilon,\tau}+1)^{\frac12}|^2\\\nonumber 
&&+~ \frac2{\Theta_0}\frac{d}{dt} ||\nabla \sqrt{c^{\varepsilon,\tau}(t)}||_{L^{2}}^{2} + \frac4{3\Theta_0}||\Delta \sqrt{c^{\varepsilon,\tau}(t)}||_{L^{2}}^{2}+ \frac1{3\Theta_0} \int_{\mathbb{R}^3}(\sqrt{c^{\varepsilon,\tau}})^{-2} |\nabla \sqrt{c^{\varepsilon,\tau}}|^4  \\
&&\leq~ \frac{27}{2 \Theta_0}||c^{\varepsilon}_0||_{L^\infty} \int_{\mathbb{R}^3} |\nabla u^{\varepsilon,\tau}|^2 + C(\|c^{\varepsilon}_{0}\|_{L^{\infty}},\|c^{\varepsilon}_{0}\|_{L^{1}}).
\een

Finally, for the equation of $u^{\varepsilon,\tau}$, multiplying the equation $(\ref{eq:CNS})_3$ by $u^{\varepsilon,\tau}$ and integration by parts, we have
\beno
\frac 12 \frac d{dt} ||u^{\varepsilon,\tau}(t)||_{L^{2}}^{2} + ||\nabla u^{\varepsilon,\tau}(t)||_{L^{2}}^{2} = - \int_{\mathbb{R}^3} ((n^{\varepsilon,\tau} \nabla \phi) \ast \rho^\varepsilon) \cdot u^{\varepsilon,\tau} dx.
\eeno
Due to
\begin{equation*}
\begin{aligned}
- \int_{\mathbb{R}^3} ((n^{\varepsilon,\tau} \nabla \phi) \ast \rho^\varepsilon) \cdot u^{\varepsilon,\tau} dx = -\int_{\mathbb{R}^3}n^{\varepsilon,\tau}\nabla\phi\cdot (u^{\varepsilon,\tau}\ast\rho^{\varepsilon})dx,
\end{aligned}
\end{equation*}
it follows that
\begin{equation}\label{u}
	\frac 12 \frac d{dt} ||u^{\varepsilon,\tau}(t)||_{L^{2}}^{2} + ||\nabla u^{\varepsilon,\tau}(t)||_{L^{2}}^{2} =  - \int_{\mathbb{R}^3} n^{\varepsilon,\tau} \nabla \phi \cdot (u^{\varepsilon,\tau} \ast \rho^\varepsilon) dx.
\end{equation}
Note that
\begin{equation*}
	\begin{aligned}
	\|n^{\varepsilon,\tau}\|_{L^{3}}^{3}=~&\int_{\mathbb{R}^{3}}\left((\sqrt{n^{\varepsilon,\tau}+1}-1)^{2}+2\sqrt{n^{\varepsilon,\tau}+1}-2 \right)^{3}dx  \\\leq~&C\int_{\mathbb{R}^{3}}(\sqrt{n^{\varepsilon,\tau}+1}-1)^{6}+C\int_{\mathbb{R}^{3}}(\sqrt{n^{\varepsilon,\tau}+1}-1)^{3}\\
\leq~&C\left(\int_{\mathbb{R}^{3}}|\nabla\sqrt{n^{\varepsilon,\tau}+1}|^{2}\right)^3+C\int_{\{n^{\varepsilon,\tau}\leq 1\}}(\sqrt{n^{\varepsilon,\tau}+1}-1)^{3}\\&+~C\int_{\{n^{\varepsilon,\tau}> 1\}}(\sqrt{n^{\varepsilon,\tau}+1}-1)^{3}\\:=~&C\left(\int_{\mathbb{R}^{3}}|\nabla\sqrt{n^{\varepsilon,\tau}+1}|^{2}\right)^3+K_{1}+K_{2}.
	\end{aligned}
\end{equation*}
For $ K_{1}, $ due to $ \sqrt{n^{\varepsilon,\tau}+1}-1\leq 2n^{\varepsilon,\tau} $ for $ n^{\varepsilon,\tau}\leq 1, $ there holds
\begin{equation*}
	\begin{aligned}
	K_{1}&\leq~C\int_{\{n^{\varepsilon,\tau}\leq 1\}}(n^{\varepsilon,\tau})^{3}\leq C\int_{\{n^{\varepsilon,\tau}\leq 1\}}n^{\varepsilon,\tau}\leq C\|n^{\varepsilon}_{0}\|_{L^{1}}.
	\end{aligned}
\end{equation*} 
For $ K_{2}, $ clearly, $ \sqrt{n^{\varepsilon,\tau}+1}-1\leq\sqrt{n^{\varepsilon,\tau}+1} $ and $ \sqrt{n^{\varepsilon,\tau}+1}\leq 4(\sqrt{n^{\varepsilon,\tau}+1}-1) $ for $ n^{\varepsilon,\tau}>1, $ then we have
\begin{equation*}
	\begin{aligned}
	K_{2}&\leq~ C\int_{\{n^{\varepsilon,\tau}>1 \}}(\sqrt{n^{\varepsilon,\tau}+1})^{3}\\&\leq~C\|\sqrt{n^{\varepsilon,\tau}+1}\|_{L^{2}\{n^{\varepsilon,\tau}>1\}}^{\frac32}\|\sqrt{n^{\varepsilon,\tau}+1}\|_{L^{6}\{n^{\varepsilon,\tau}>1 \}}^{\frac32}\\&\leq~C\|\sqrt{2n^{\varepsilon,\tau}}\|_{L^{2}}^{\frac32}\|\sqrt{n^{\varepsilon,\tau}+1}-1\|_{L^{6}}^{\frac32}\\&\leq~ C\|n^{\varepsilon}_{0}\|_{L^{1}}^{\frac34}\|\nabla\sqrt{n^{\varepsilon,\tau}+1}\|_{L^{2}}^{\frac32}.
	\end{aligned}
\end{equation*}
Combining $ K_{1} $ and $ K_{2}, $ it follows that 
\begin{equation}\label{n L3}
	\begin{aligned}
	\|n^{\varepsilon,\tau}\|_{L^{3}}^{3}\leq C\|\nabla\sqrt{n^{\varepsilon,\tau}+1}\|_{L^{2}}^{6}+C\|n^{\varepsilon}_{0}\|_{L^{1}}^{\frac34}\|\nabla\sqrt{n^{\varepsilon,\tau}+1}\|_{L^{2}}^{\frac32}+C\|n^{\varepsilon}_{0}\|_{L^{1}}.
	\end{aligned}
\end{equation}
Using (\ref{n L3}), the right term of $ (\ref{u}) $ can be estimated as follow:
\begin{equation*}
\begin{aligned}
-&\int_{\mathbb{R}^3} n^{\varepsilon,\tau} \nabla \phi \cdot (u^{\varepsilon,\tau} \ast \rho^{\varepsilon}) dx\\\leq~& ||\nabla \phi||_{L^\infty} ||u^{\varepsilon,\tau}||_{L^6} ||n^{\varepsilon,\tau}||_{L^\frac 65}\\\leq~&C\|\nabla\phi\|_{L^{\infty}}\|\nabla u^{\varepsilon,\tau}\|_{L^{2}}\|n^{\varepsilon,\tau}\|_{L^{1}}^{\frac34}\|n^{\varepsilon,\tau}\|_{L^{3}}^{\frac14}\\\leq~&\frac12\|\nabla u^{\varepsilon,\tau}\|_{L^{2}}^{2}+C\|\nabla\phi\|_{L^{\infty}}^{2}\|n^{\varepsilon}_{0}\|_{L^{1}}^{\frac32}\|n^{\varepsilon,\tau}\|_{L^{3}}^{\frac12}\\\leq~&\frac12\|\nabla u^{\varepsilon,\tau}\|_{L^{2}}^{2}+C\|\nabla\phi\|_{L^{\infty}}^{2}\|n^{\varepsilon}_{0}\|_{L^{1}}^{\frac32}\left(\|\nabla\sqrt{n^{\varepsilon,\tau}+1}\|_{L^{2}}+\|n^{\varepsilon}_{0}\|_{L^{1}}^{\frac18}\|\nabla\sqrt{n^{\varepsilon,\tau}+1}\|_{L^{2}}^{\frac14}+\|n^{\varepsilon}_{0}\|_{L^{1}}^{\frac16}\right)\\\leq~&\frac12\|\nabla u^{\varepsilon,\tau}\|_{L^{2}}^{2}+\frac{1}{8(\frac{27}{2 \Theta_0}\|c^{\varepsilon}_0\|_{L^\infty} + 1)} \|\nabla \sqrt{n^{\varepsilon,\tau}+1}\|_{L^2}^2 \\+&~C\left(\|\nabla\phi\|_{L^{\infty}}^{4}\|n^{\varepsilon}_{0}\|_{L^{1}}^{3}+ \|\nabla\phi\|_{L^{\infty}}^{\frac{16}{7}}\|n^{\varepsilon}_{0}\|_{L^{1}}^{\frac{13}{7}}+\|\nabla\phi\|_{L^{\infty}}^{2}\|n^{\varepsilon}_{0}\|_{L^{1}}^{\frac53} \right).
\end{aligned}
\end{equation*}
Then from (\ref{u}) we can obtain that
\ben\label{K_1,K_2}
&&\frac{d}{dt}||u^{\varepsilon,\tau}(t)||_{L^{2}}^{2} + ||\nabla u^{\varepsilon,\tau}(t)||_{L^{2}}^{2} \\\nonumber &\leq& \frac{1}{8(\frac{27}{2 \Theta_0}\|c^{\varepsilon}_0\|_{L^\infty} + 1)} \|\nabla (n^{\varepsilon,\tau}+1)^\frac12\|_{L^2}^2 + C(\|\nabla \phi\|_{L^\infty},\|n_0^\varepsilon\|_{L^1}),
\een
where we use $ ||\nabla\phi||_{L^{\infty}}\leq C $ and
$ ||n^{\varepsilon}_{0}||_{L^{1}}\leq C. $ 

Therefore, (\ref{K_1,K_2})$\times 2(\frac{27}{2 \Theta_0}||c^{\varepsilon}_{0}||_{L^{\infty}}+1) $+(\ref{ine:energy estimate of c and n}) yields that
%The ri^{\varepsilon}ght term of the above equation can be estimated as follow:
%\begin{equation*}
%\begin{aligned}
%&-\int_{\mathbb{R}^3} n^{\varepsilon,\tau} \nabla \phi \cdot (u^{\varepsilon,\tau} \ast \rho^\varepsilon) dx\\\leq &~ C ||\nabla \phi||_{L^\infty} ||u^{\varepsilon,\tau}||_{L^6} ||n^{\varepsilon,\tau}||_{L^\frac 65} \\
%\leq& ~ C\|\nabla\phi\|_{L^{\infty}}\|\nabla u^{\varepsilon,\tau}\|_{L^{2}}\|\sqrt{n^{\varepsilon,\tau}}\|_{L^{\frac{12}{5}}}^{2}\\\leq&~ \epsilon_{2}\|\nabla u^{\varepsilon,\tau}\|_{L^{2}}^{2}+C(\epsilon_{2})\|\nabla\phi\|_{L^{\infty}}^{2}\|\sqrt{n^{\varepsilon,\tau}}\|_{L^{\frac{12}{5}}}^{4}\\\leq&~ \epsilon_{2}\|\nabla u^{\varepsilon,\tau}\|_{L^{2}}^{2}+C(\epsilon_{2})\|\nabla\phi\|_{L^{\infty}}^{2}\|\nabla\sqrt{n^{\varepsilon,\tau}}\|_{L^{2}}\|\sqrt{n^{\varepsilon,\tau}}\|_{L^{2}}^{3}\\\leq&~ \epsilon_{2}\|\nabla u^{\varepsilon,\tau}\|_{L^{2}}^{2}+\epsilon_{3}\|\nabla\sqrt{n^{\varepsilon,\tau}}\|_{L^{2}}^{2}+C(\epsilon_{2},\epsilon_{3})\|\nabla\phi\|_{L^{\infty}}^{4}\|n_{0}\|_{L^{1}}^{3}.
%\end{aligned}
%\end{equation*}
%Noting that $ ||\nabla\phi||_{L^{\infty}}\leq C $ and
%$ ||n_{0}||_{L^{1}}\leq C, $ there holds
%\ben\label{K_1,K_2}
%\frac{d}{dt}||u^{\varepsilon,\tau}||_2^2 + (2-2\epsilon_{2})||\nabla u^{\varepsilon,\tau}(t)||_2^2 \leq 2\epsilon_{3}\|\nabla\sqrt{n^{\varepsilon,\tau}}\|_{L^{2}}^{2} + C(\epsilon_{2},\epsilon_{3})\|\nabla\phi\|_{L^{\infty}}^{4}\|n_{0}\|_{L^{1}}^{3}.
%\een

\beno
&&\frac{d}{dt}||u^{\varepsilon,\tau}||_{L^{2}}^{2} + 2(\frac{27}{2 \Theta_0}||c^{\varepsilon}_{0}||_{L^{\infty}}+1)||\nabla u^{\varepsilon,\tau}(t)||_{L^{2}}^{2}\\
&&+\frac{d}{dt}\int_{\mathbb{R}^3} (n^{\varepsilon,\tau}+1) \ln (n^{\varepsilon,\tau}+1)(\cdot,t)dx + \int_{\mathbb{R}^3} |\nabla (n^{\varepsilon,\tau}+1)^{\frac12}|^2\\\nonumber 
&&+ \frac2{\Theta_0}\frac{d}{dt} ||\nabla \sqrt{c^{\varepsilon,\tau}(t)}||_{L^{2}}^{2}+ \frac4{3\Theta_0}||\Delta \sqrt{c^{\varepsilon,\tau}(t)}||_{L^{2}}^{2} + \frac1{3\Theta_0} \int_{\mathbb{R}^3}(\sqrt{c^{\varepsilon,\tau}})^{-2} |\nabla \sqrt{c^{\varepsilon,\tau}}|^4  \\
&\leq &\frac{27}{2 \Theta_0}||c^{\varepsilon}_0||_{L^\infty} \int_{\mathbb{R}^3} |\nabla u^{\varepsilon,\tau}|^2 + C\left(   \|\nabla\phi\|_{L^{\infty}},\|n^{\varepsilon}_{0}\|_{L^{1}},\|c^{\varepsilon}_{0}\|_{L^{\infty}},\|c^{\varepsilon}_{0}\|_{L^{1}}\right),
\eeno
which means
%Since $C(\epsilon_{1})||c_0||_{L^\infty}\leq \frac12$, we have
\ben\label{n,c,u}\nonumber 
&&\frac{d}{dt}||u^{\varepsilon,\tau}||_{L^{2}}^{2} + ||\nabla u^{\varepsilon,\tau}(t)||_{L^{2}}^{2}\\\nonumber 
&&+\frac{d}{dt}\int_{\mathbb{R}^3} (n^{\varepsilon,\tau}+1) \ln (n^{\varepsilon,\tau}+1)(\cdot,t)dx +\int_{\mathbb{R}^3} |\nabla (n^{\varepsilon,\tau}+1)^{\frac12}|^2\\\nonumber 
&&+ \frac2{\Theta_0}\frac{d}{dt} ||\nabla \sqrt{c^{\varepsilon,\tau}(t)}||_{L^{2}}^{2} + \frac4{3\Theta_0}||\Delta \sqrt{c^{\varepsilon,\tau}(t)}||_{L^{2}}^{2} + \frac1{3\Theta_0} \int_{\mathbb{R}^3}(\sqrt{c^{\varepsilon,\tau}})^{-2} |\nabla \sqrt{c^{\varepsilon,\tau}}|^4  \\\nonumber
&\leq & 
C(\|\nabla \phi\|_{L^\infty},\|n_0^{\varepsilon}\|_{L^1},\|c^{\varepsilon}_{0}\|_{L^{\infty}},\|c^{\varepsilon}_{0}\|_{L^{1}}).
\een
Integrating with respect to time from 0 to $t$, we get (\ref{ine:  energy}). 
%\begin{equation}\label{U V}
%	\begin{aligned}
%	&||u^{\varepsilon,\tau}||_{L^{2}}^{2} + \int_0^t||\nabla u^{\varepsilon,\tau}(t)||_{L^{2}}^{2}\\\nonumber 
%+~&\int_{\mathbb{R}^3} (n^{\varepsilon,\tau}+1) \ln (n^{\varepsilon,\tau}+1)(\cdot,t)dx + \int_0^t\int_{\mathbb{R}^3} |\nabla (n^{\varepsilon,\tau}+1)^{\frac12}|^2\\\nonumber 
%+~& \frac2{\Theta_0}||\nabla \sqrt{c^{\varepsilon,\tau}(t)}||_{L^{2}}^{2} +\frac4{3\Theta_0}\int_0^t||\Delta \sqrt{c^{\varepsilon,\tau}(t)}||_{L^{2}}^{2} + \frac1{3\Theta_0} \int_0^t \int_{\mathbb{R}^3}(\sqrt{c^{\varepsilon,\tau}})^{-2} |\nabla \sqrt{c^{\varepsilon,\tau}}|^4  \\
%\leq ~ &
%C(\|\nabla \phi\|_{L^\infty},\|n_0^{\varepsilon}\|_{L^1},\|c^{\varepsilon}_{0}\|_{L^{\infty}},\|c^{\varepsilon}_{0}\|_{L^{1}},\|u^{\varepsilon}_{0}\|_{L^{2}},\|n^{\varepsilon}_{0}\ln n^{\varepsilon}_{0}\|_{L^{1}},\|\nabla \sqrt{c_0^{\varepsilon}}\|_{L^{2}})t.
%	\end{aligned}
%\end{equation}
%Thus (\ref{ine:  energy}) is proved.

%\ben\label{nlnn varepsilon,tau }
%\int_{\Omega} (n^{\varepsilon,\tau}\ln n^{\varepsilon,\tau})(\cdot,t) dx + 2\alpha^{-1}e^{-1} C^{1-\alpha}\geq\int_{\Omega} (n^{\varepsilon,\tau} |\ln n^{\varepsilon,\tau}|)(\cdot,t) dx> 0.
%\een

{\bf Proof of (\ref{eq: good esimate of u}) when $\mu=0$.}

 Rewrite the equation of $ u^{\varepsilon,\tau} $ as follows:
\begin{equation*}
	\partial_{t}(-\Delta)^{\frac14}u^{\varepsilon,\tau}-\Delta(-\Delta)^{\frac14}u^{\varepsilon,\tau}+\nabla(-\Delta)^{\frac14}P^{\varepsilon,\tau}
=-(-\Delta)^{\frac14}(n^{\varepsilon,\tau}\nabla\phi)\ast \rho^{\varepsilon}.
\end{equation*} 
Multiplying $ (-\Delta)^{\frac14}u^{\varepsilon,\tau} $ in above equation, integration by parts yields that
\begin{equation*}
	\begin{aligned}
	&\|(-\Delta)^{\frac14}u^{\varepsilon,\tau}\|_{L^{2}}^{2}+2\int_{0}^{t}\|\nabla(-\Delta)^{\frac14}u^{\varepsilon,\tau}\|_{L^{2}}^{2}\\\leq~&
C\|\nabla\phi\|_{L^{\infty}}\int_{0}^{t}\|n^{\varepsilon,\tau}\|_{L^{\frac32}}\|(-\Delta)^{\frac12}u^{\varepsilon,\tau}\|_{L^{3}}ds+
\|(-\Delta)^{\frac14}u_{0}^{\varepsilon}\|_{L^{2}}^{2}\\\leq~&C\|\nabla\phi\|_{L^{\infty}}\left( \int_{0}^{t}\|n^{\varepsilon,\tau}\|_{L^{\frac32}}^{2} \right)^{\frac12}\left( \int_{0}^{t}\|(-\Delta)^{\frac12}u^{\varepsilon,\tau}\|_{L^{3}}^{2}\right)^{\frac12}+\|(-\Delta)^{\frac14}u_{0}^{\varepsilon}\|_{L^{2}}^{2}
\\\leq~&C\|\nabla\phi\|_{L^{\infty}}\|\sqrt{n^{\varepsilon,\tau}}\|_{L^{4}_{t}L^{3}_{x}}^{2}\|(-\Delta)^{\frac12}u^{\varepsilon,\tau}\|_{L^{2}_{t}L^{3}_{x}}+
\|(-\Delta)^{\frac14}u_{0}^{\varepsilon}\|_{L^{2}}^{2}\\\leq~&C\|\nabla\phi\|_{L^{\infty}}\|\sqrt{n^{\varepsilon,\tau}}\|_{L^{4}_{t}L^{3}_{x}}^{2}
\|\nabla(-\Delta)^{\frac14}u^{\varepsilon,\tau}\|_{L^{2}_{t}L^{2}_{x}}+\|(-\Delta)^{\frac14}u_{0}^{\varepsilon}\|_{L^{2}}^{2}
\\\leq~&\frac12\|\nabla(-\Delta)^{\frac14}u^{\varepsilon,\tau}\|_{L^{2}_{t}L^{2}_{x}}^{2}+C\|\nabla\phi\|_{L^{\infty}}^{2}
\|\sqrt{n^{\varepsilon,\tau}}\|_{L^{4}_{t}L^{3}_{x}}^{4}+\|(-\Delta)^{\frac14}u_{0}^{\varepsilon}\|_{L^{2}}^{2},
	\end{aligned}
\end{equation*}
which follows that
\begin{equation}\label{u better}
	\begin{aligned}
	&\|(-\Delta)^{\frac14}u^{\varepsilon,\tau}\|_{L^{2}}^{2}+\int_{0}^{t}\|\nabla(-\Delta)^{\frac14}u^{\varepsilon,\tau}\|_{L^{2}}^{2}\\\leq~&C\|\nabla\phi\|_{L^{\infty}}^{2}\|\sqrt{n^{\varepsilon,\tau}}\|_{L^{4}_{t}L^{3}_{x}}^{4}+\|(-\Delta)^{\frac14}u_{0}^{\varepsilon}\|_{L^{2}}^{2}.
	\end{aligned}
\end{equation}
Note that
\begin{equation*}
	\|n^{\varepsilon,\tau}\|_{L^{\frac32}}\leq \|n^{\varepsilon,\tau}\|_{L^{1}}^{\frac12}\|n^{\varepsilon,\tau}\|_{L^{3}}^{\frac12},
\end{equation*}
by (\ref{n L3}), we have
\begin{equation*}
\begin{aligned}
&\left(\int_{0}^{t}\|n^{\varepsilon,\tau}\|_{L^{\frac32}}^{2} \right)^{\frac12}\\\leq~&C\|n_{0}^{\varepsilon}\|_{L^{1}}^{\frac12}\left( \int_{0}^{t}\|n^{\varepsilon,\tau}\|_{L^{3}} \right)^{\frac12}\\\leq~&C\|n_{0}^{\varepsilon}\|_{L^{1}}^{\frac12}\left(\int_{0}^{t}\|\nabla\sqrt{n^{\varepsilon,\tau}+1}\|_{L^{2}}^{2}ds+ \|n_{0}^{\varepsilon}\|_{L^{1}}^{\frac14}\int_{0}^{t}\|\nabla\sqrt{n^{\varepsilon,\tau}+1}\|_{L^{2}}^{\frac12}ds +\|n_{0}^{\varepsilon}\|_{L^{1}}^{\frac13}t\right)^{\frac12}\\\leq~&C(\|\nabla \phi\|_{L^\infty},\|n_0^{\varepsilon}\|_{L^1},\|c_{0}^{\varepsilon}\|_{L^{\infty}},\|c_{0}^{\varepsilon}\|_{L^{1}},\|u_{0}^{\varepsilon}\|_{L^{2}},\|(n_{0}^{\varepsilon}+1)\ln (n_{0}^{\varepsilon}+1)\|_{L^{1}},\|\nabla \sqrt{c_0^{\varepsilon}}\|_{L^{2}})(t+1)^{\frac12}.
\end{aligned}
\end{equation*}
 Thus we have 
\ben\label{sqrt n 1}
\|\sqrt{n^{\varepsilon,\tau}}\|_{L^{4}_{t}L^{3}_{x}}^{4}\leq C(t+1).
\een
 By (\ref{sqrt n 1}),$ \nabla\phi\in L^{\infty} $ and the Sobolev inequality (also see Lemma \ref{Sobolev inequalities})
$$\|f\|_{L^{q}}\leq C\|(-\Delta)^{\frac{\alpha}{2}}f\|_{L^{p}},\quad{\rm for}\quad\frac1q=\frac1p-\frac{\alpha}{n} ,$$ (\ref{u better}) becomes
\begin{equation*}
	\begin{aligned}
&	\|(-\Delta)^{\frac14}u^{\varepsilon,\tau}\|_{L^{2}}^{2}+\int_{0}^{t}\|\nabla(-\Delta)^{\frac14}u^{\varepsilon,\tau}\|_{L^{2}}^{2}\\
\leq~&C(\|\nabla \phi\|_{L^\infty},\|n_0^{\varepsilon}\|_{L^1},\|c_{0}^{\varepsilon}\|_{L^{\infty}},\|c_{0}^{\varepsilon}\|_{L^{1}},\|u_{0}^{\varepsilon}\|_{L^{2}},\|(n_{0}^{\varepsilon}+1)\ln (n_{0}^{\varepsilon}+1)\|_{L^{1}},\|\nabla \sqrt{c_0^{\varepsilon}}\|_{L^{2}})(t+1).
	\end{aligned}
\end{equation*}
%where we use $$ \|\sqrt{n^{\varepsilon,\tau}}\|_{L^{4}_{t}L^{3}_{x}}^{4}\leq C\left(\|\sqrt{n^{\varepsilon,\tau}}\|_{L^{\infty}_{t}L^{2}_{x}}^{2}+\|\nabla\sqrt{n^{\varepsilon,\tau}+1}\|_{L^{2}_{t}L^{2}_{x}}^{2}
%\right)^{2}\leq C(T^2+1), $$ 

\section{Local energy inequality}

\begin{lemma}[local energy inequality]\label{local energy inequality}
	Let $ \Omega\subset\mathbb{R}^{3} $ is an arbitrary bounded domain and  $\psi$ is a cut-off function, which vanishes on the parabolic boundary of $\Omega\times(0,t) $.
	Assume $(n^{\varepsilon,\tau},c^{\varepsilon,\tau},u^{\varepsilon,\tau})$ is a regular solution of (\ref{eq:CNS}) in $\Omega\times(0,t)\subset\mathbb{R}^{3}\times\mathbb{R}^{+}$ with initial data
	$\eqref{ine:initial condition assumption}$. Assume that the function $\kappa$ and $\chi$ satisfy $\eqref{ine:chi}$ ,$\eqref{ine:kappa}$ and $ \eqref{equ:chi kappa}. $ Then the following inequality holds:

\ben\label{local energy}
\nonumber&& \int_{\Omega} (n^{\varepsilon,\tau} \ln (n^{\varepsilon,\tau}) \psi)(\cdot,t) + 4 \int_{(0,t)\times\Omega} |\nabla \sqrt{n^{\varepsilon,\tau}}|^2 \psi
	+ \frac{2}{\Theta_0}  \int_{\Omega} (|\nabla \sqrt{c^{\varepsilon,\tau}}|^2 \psi)(\cdot,t)\\ \nonumber&&
 +  \frac{4}{3\Theta_0}\int_{(0,t)\times\Omega} |\Delta \sqrt{c^{\varepsilon,\tau}}|^2 \psi 
+\frac{2}{3\Theta_0}  \int_{(0,t)\times\Omega} (\sqrt{c^{\varepsilon,\tau}})^{-2} |\nabla \sqrt{c^{\varepsilon,\tau}}|^4 \psi\\\nonumber
&&+\frac{18}{\Theta_0}\|{c^{\varepsilon}_{0}}\|_{L^{\infty}}\int_{\Omega} (|u^{\varepsilon,\tau}|^2)(\cdot,t) \psi + \frac{18}{\Theta_0}\|{c^{\varepsilon}_{0}}\|_{L^{\infty}}\int_{(0,t)\times\Omega} |\nabla u^{\varepsilon,\tau}|^2 \psi \\\nonumber
&&\leq \int_{(0,t)\times\Omega} n^{\varepsilon,\tau} \ln n^{\varepsilon,\tau} (\partial_t \psi + \Delta \psi) + \int_{(0,t)\times\Omega} n^{\varepsilon,\tau} \ln n^{\varepsilon,\tau} (u^{\varepsilon,\tau}\ast\rho^\varepsilon) \cdot \nabla \psi\\ \nonumber
&&+\int_{(0,t)\times\Omega}\left(\frac{n^{\varepsilon,\tau}}{1+\tau n^{\varepsilon,\tau}} \nabla c^{\varepsilon,\tau} \chi(c^{\varepsilon,\tau})\right)\ln n^{\varepsilon,\tau}\cdot\nabla\psi +\int_{(0,t)\times\Omega}\left(\frac{n^{\varepsilon,\tau}}{1+\tau n^{\varepsilon,\tau}} \nabla c^{\varepsilon,\tau} \chi(c^{\varepsilon,\tau})\right)\cdot\nabla\psi\\ \nonumber
%&&+\int_{(0,t)\times\Omega} \left(\frac{1}{1+\tau n^{\varepsilon,\tau}} \nabla c^{\varepsilon,\tau} \chi(c^{\varepsilon,\tau})\right) \cdot \nabla n^{\varepsilon,\tau} \psi-\int_{(0,t)\times\Omega}\frac{\chi(c^{\varepsilon,\tau})}{1+\tau n^{\varepsilon,\tau}}\nabla n^{\varepsilon,\tau}\cdot \nabla c^{\varepsilon,\tau}\psi\\ \nonumber
&&+ \frac{2}{\Theta_0}\int_{(0,t)\times\Omega} |\nabla \sqrt{c^{\varepsilon,\tau}}|^2 (\partial_t \psi + \Delta \psi)
+ \frac{2}{\Theta_0} \int_{(0,t)\times\Omega} |\nabla \sqrt{c^{\varepsilon,\tau}}|^2 {(u^{\varepsilon,\tau}\ast\rho^\varepsilon)} \cdot \nabla \psi\\ \nonumber %&&+\frac{18}{\Theta_0}\|{c^{\varepsilon}_{0}}\|_{L^{\infty}} \int_{(0,t)\times\Omega} |\nabla u^{\varepsilon,\tau}|^2 \psi.\\
&&+\frac{18}{\Theta_0}\|{c^{\varepsilon}_{0}}\|_{L^{\infty}}\int_{(0,t)\times\Omega} |u^{\varepsilon,\tau}|^2 \left(\partial_t \psi + \Delta \psi\right) + \frac{18\mu}{\Theta_0}\|{c^{\varepsilon}_{0}}\|_{L^{\infty}}\int_{(0,t)\times\Omega} |u^{\varepsilon,\tau}|^2(u^{\varepsilon,\tau}\ast\rho^\varepsilon)\cdot\nabla\psi\\&&+\frac{36}{\Theta_0}\|{c^{\varepsilon}_{0}}\|_{L^{\infty}}\int_{(0,t)\times\Omega} (P^{\varepsilon,\tau} - \bar{P}^{\varepsilon,\tau}) u^{\varepsilon,\tau} \cdot \nabla \psi
-\frac{36}{\Theta_0}\|{c^{\varepsilon}_{0}}\|_{L^{\infty}}\int_{(0,t)\times\Omega} (n^{\varepsilon,\tau}\nabla\phi)\ast\rho^\varepsilon \cdot u^{\varepsilon,\tau} \psi.
\een

\end{lemma}	
{\bf Proof.} 
Multiplying $(1+\ln n^{\varepsilon,\tau})\psi$ in $(\ref{eq:CNS})_1$, integration by parts yields that
\beno
&&\int_{0}^t\int_{\Omega}\partial_t n^{\varepsilon,\tau} (1+\ln n^{\varepsilon,\tau})\psi+\int_{0}^t\int_{\Omega} (u^{\varepsilon,\tau} \ast \rho^\varepsilon) \cdot \nabla n^{\varepsilon,\tau}(1+\ln n^{\varepsilon,\tau})\psi\\
&& -\int_{0}^t\int_{\Omega} \Delta n^{\varepsilon,\tau}(1+\ln n^{\varepsilon,\tau})\psi +\int_{0}^t\int_{\Omega}\nabla \cdot \left(\frac{n^{\varepsilon,\tau}}{1+\tau n^{\varepsilon,\tau}} \nabla c^{\varepsilon,\tau} \chi(c^{\varepsilon,\tau})\right)(1+\ln n^{\varepsilon,\tau})\psi\\
%&&+\int_{0}^t\int_{\Omega} \nabla \cdot (n^\varepsilon (\nabla c^\varepsilon\chi(c^{\varepsilon})) \ast \rho^\varepsilon)(1+\ln n^\varepsilon)\psi\\
&&\doteq T_1+\cdots+T_4=0,
\eeno
where
\beno
T_1&=&\int_{(0,t)\times\Omega} \partial_tn^{\varepsilon,\tau}{(1+\ln n^{\varepsilon,\tau})}\psi dxdt \\&=& \int_{\Omega}  (n^{\varepsilon,\tau}{\ln n^{\varepsilon,\tau}}\psi)(\cdot,t) dx-\int_{(0,t)\times\Omega}n^{\varepsilon,\tau}\ln n^{\varepsilon,\tau}\partial_t \psi dxdt;
\eeno
\beno
T_2&=&-\int_{(0,t)\times\Omega} (u^{\varepsilon,\tau} \ast \rho^\varepsilon)\cdot \nabla n^{\varepsilon,\tau} \psi dxdt-\int_{(0,t)\times\Omega}n^{\varepsilon,\tau}(1+\ln n^{\varepsilon,\tau})(u^{\varepsilon,\tau} \ast \rho^\varepsilon)\cdot\nabla\psi dxdt\\
&=&-\int_{(0,t)\times\Omega}n^{\varepsilon,\tau}\ln n^{\varepsilon,\tau} (u^{\varepsilon,\tau} \ast \rho^\varepsilon)\cdot\nabla\psi dxdt;
\eeno
\beno
T_3&=&\int_{(0,t)\times\Omega}\frac1{n^{\varepsilon,\tau}} \nabla n^{\varepsilon,\tau} \cdot \nabla n^{\varepsilon,\tau} \psi dxdt+\int_{(0,t)\times\Omega}\nabla n^{\varepsilon,\tau} \cdot \nabla \psi dxdt\\&&+\int_{(0,t)\times\Omega} \ln n^{\varepsilon,\tau} \nabla n^{\varepsilon,\tau} \cdot \nabla \psi dxdt\\
%&=&\int_{(0,t)\times\Omega}\frac1{n^{\varepsilon,\tau}} \nabla n^{\varepsilon,\tau} \cdot \nabla n^{\varepsilon,\tau} \psi dxdt+\int_{(0,t)\times\Omega}\nabla n^{\varepsilon,\tau} \cdot \nabla \psi dxdt\\
%&&-\int_{(0,t)\times\Omega} n^{\varepsilon,\tau} \frac1{n^{\varepsilon,\tau}} \nabla n^{\varepsilon,\tau} \cdot \nabla \psi dxdt-\int_{(0,t)\times\Omega} n^{\varepsilon,\tau} \ln n^{\varepsilon,\tau} \Delta\psi dxdt\\
&=&4\int_{(0,t)\times\Omega} |\nabla\sqrt{n^{\varepsilon,\tau}}|^2\psi dxdt-\int_{(0,t)\times\Omega}n^{\varepsilon,\tau}\ln n^{\varepsilon,\tau}\Delta\psi dxdt;
\eeno
and

\beno
T_4&=&-\int_{0}^t\int_{\Omega}\left(\frac{n^{\varepsilon,\tau}}{1+\tau n^{\varepsilon,\tau}} \nabla c^{\varepsilon,\tau} \chi(c^{\varepsilon,\tau})\right)(1+\ln n^{\varepsilon,\tau})\cdot\nabla\psi dxdt\\&&- \int_{0}^t\int_{\Omega} \left(\frac{1}{1+\tau n^{\varepsilon,\tau}} \nabla c^{\varepsilon,\tau} \chi(c^{\varepsilon,\tau})\right) \cdot \nabla n^{\varepsilon,\tau} \psi dxdt.\\
\eeno
%\beno
%T_4&=& - \int_{0}^t\int_{\Omega} (1+\ln n^\varepsilon) n^\varepsilon (\nabla c^\varepsilon\chi(c^{\varepsilon})) \ast \rho^\varepsilon \cdot \nabla \psi dxdt - \int_{0}^t\int_{\Omega} (\nabla c^\varepsilon\chi(c^{\varepsilon})) \ast \rho^\varepsilon \cdot \nabla n^\varepsilon \psi
%\eeno
To sum up, we obtain
\ben\label{n}
\nonumber&&\int_{\Omega}  (n^{\varepsilon,\tau}{\ln n^{\varepsilon,\tau}}\psi)(\cdot,t)+4\int_{(0,t)\times\Omega} |\nabla\sqrt{n^{\varepsilon,\tau}}|^2\psi\\ \nonumber
&=&\int_{(0,t)\times\Omega}n^{\varepsilon,\tau}\ln n^{\varepsilon,\tau}(\partial_t \psi+\Delta\psi)+\int_{(0,t)\times\Omega}n^{\varepsilon,\tau}\ln n^{\varepsilon,\tau} (u^{\varepsilon,\tau} \ast \rho^\varepsilon)\cdot\nabla\psi\\ \nonumber
&&+\int_{(0,t)\times\Omega}\left(\frac{n^{\varepsilon,\tau}}{1+\tau n^{\varepsilon,\tau}} \nabla c^{\varepsilon,\tau} \chi(c^{\varepsilon,\tau})\right)\ln n^{\varepsilon,\tau}\cdot\nabla\psi \\\nonumber&&+\int_{(0,t)\times\Omega}\left(\frac{n^{\varepsilon,\tau}}{1+\tau n^{\varepsilon,\tau}} \nabla c^{\varepsilon,\tau} \chi(c^{\varepsilon,\tau})\right)\cdot\nabla\psi\\
&&+\int_{(0,t)\times\Omega} \left(\frac{1}{1+\tau n^{\varepsilon,\tau}} \nabla c^{\varepsilon,\tau} \chi(c^{\varepsilon,\tau})\right) \cdot \nabla n^{\varepsilon,\tau} \psi.
\een
Note that
\beno
\Delta {c^{\varepsilon,\tau}} = 2 |\nabla \sqrt{{c^{\varepsilon,\tau}}}|^2 + 2 \sqrt{{c^{\varepsilon,\tau}}} \Delta \sqrt{{c^{\varepsilon,\tau}}},
\eeno
and dividing $2\sqrt{{c^{\varepsilon,\tau}}}$ on both sides, it follows from the equation $(\ref{eq:CNS})_2$ that
\ben\label{eq:c-tilde}
\partial_t\sqrt{c^{\varepsilon,\tau}}+u^{\varepsilon,\tau}\cdot \nabla\sqrt{c^{\varepsilon,\tau}}-\frac{|\nabla \sqrt{c^{\varepsilon,\tau}}|^2}{\sqrt{c^{\varepsilon,\tau}}}-\Delta\sqrt{c^{\varepsilon,\tau}}= - \frac{\kappa(c^{\varepsilon,\tau})}{2\tau\sqrt{c^{\varepsilon,\tau}}}~~\ln(1+\tau n^{\varepsilon,\tau}).
\een
Multiplying the above equation $(\ref{eq:CNS})_2$ by $-\partial_i (\partial_i \sqrt{{c^{\varepsilon,\tau}}} \psi)$ and integration by parts, there holds
\beno
&&-\int_{(0,t)\times\Omega}\partial_t\sqrt{c^{\varepsilon,\tau}}\partial_i (\partial_i \sqrt{c^{\varepsilon,\tau}} \psi)dxdt-\int_{(0,t)\times\Omega}u^{\varepsilon,\tau}\cdot \nabla\sqrt{c^{\varepsilon,\tau}}\partial_i (\partial_i \sqrt{c^{\varepsilon,\tau}} \psi)dxdt \\&&+\int_{(0,t)\times\Omega}\frac{|\nabla \sqrt{c^{\varepsilon,\tau}}|^2}{\sqrt{c^{\varepsilon,\tau}}}\partial_i (\partial_i \sqrt{c^{\varepsilon,\tau}} \psi)dxdt +\int_{(0,t)\times\Omega}\Delta\sqrt{c^{\varepsilon,\tau}}\partial_i (\partial_i \sqrt{c^{\varepsilon,\tau}} \psi)dxdt\\&&-\int_{(0,t)\times\Omega}\frac{\kappa(c^{\varepsilon,\tau})}{2\tau\sqrt{c^{\varepsilon,\tau}}}~~\ln(1+\tau n^{\varepsilon,\tau})\partial_i (\partial_i \sqrt{c^{\varepsilon,\tau}} \psi)dxdt \\
&& \doteq J_1+J_2+\cdots+J_5=0,
\eeno
where
\beno
J_1&=&\frac12\int_{(0,t)\times\Omega}\partial_t(\partial_i\sqrt{c^{\varepsilon,\tau}})^2\psi =\frac12\int_{\Omega}(|\nabla\sqrt{c^{\varepsilon,\tau}}|^2\psi)(\cdot,t) -\frac12\int_{(0,t)\times\Omega}|\nabla\sqrt{c^{\varepsilon,\tau}}|^2\partial_t \psi ;\\
J_2&=&\int_{(0,t)\times\Omega}\partial_iu^{\varepsilon,\tau}\partial_j\sqrt{c^{\varepsilon,\tau}}\partial_i\sqrt{c^{\varepsilon,\tau}}\psi+\int_{(0,t)\times\Omega}u^{\varepsilon,\tau}\partial_{ij}\sqrt{c^{\varepsilon,\tau}}\partial_i\sqrt{c^{\varepsilon,\tau}}\psi\\
&=&\int_{(0,t)\times\Omega} \nabla u^{\varepsilon,\tau} : (\nabla\sqrt{c^{\varepsilon,\tau}} \otimes \nabla\sqrt{c^{\varepsilon,\tau}})\psi -\frac12 \int_{(0,t)\times\Omega} u^{\varepsilon,\tau}\cdot\nabla\psi \nabla \sqrt{c^{\varepsilon,\tau}} \cdot \nabla \sqrt{c^{\varepsilon,\tau}};\\
J_3&=&-\int_{(0,t)\times\Omega}\partial_i\big((\sqrt{c^{\varepsilon,\tau}})^{-1}|\nabla\sqrt{c^{\varepsilon,\tau}}|^2\big)(\partial_i\sqrt{c^{\varepsilon,\tau}}\psi)\\
&=&\int_{(0,t)\times\Omega}\big((\sqrt{c^{\varepsilon,\tau}})^{-1}|\nabla\sqrt{c^{\varepsilon,\tau}}|^2\big)\Delta\sqrt{c^{\varepsilon,\tau}}\psi +\int_{(0,t)\times\Omega}\big((\sqrt{c^{\varepsilon,\tau}})^{-1}|\nabla\sqrt{c^{\varepsilon,\tau}}|^2\big)\nabla\sqrt{c^{\varepsilon,\tau}}\cdot\nabla\psi ;\\
J_4&=&-\int_{(0,t)\times\Omega} \partial_{j} \sqrt{c^{\varepsilon,\tau}} \partial_{ij} (\sqrt{c^{\varepsilon,\tau}} \psi) =\int_{(0,t)\times\Omega} |\nabla^2\sqrt{c^{\varepsilon,\tau}}|^2 \psi  + \int_{(0,t)\times\Omega} \partial_{ij} \sqrt{c^{\varepsilon,\tau}} \partial_i \sqrt{c^{\varepsilon,\tau}} \partial_j \psi \\
&=& \int_{(0,t)\times\Omega} |\nabla^2\sqrt{c^{\varepsilon,\tau}}|^2 \psi  - \frac12\int_{(0,t)\times\Omega}|\nabla\sqrt{c^{\varepsilon,\tau}}|^2\Delta \psi ;\\
\eeno
and
\beno
J_5%&=&-\int_{(0,t)\times\Omega}\frac{\kappa(c^{\varepsilon,\tau})}{2\tau\sqrt{c^{\varepsilon,\tau}}}~~\ln(1+\tau n^{\varepsilon,\tau})\partial_i (\partial_i \sqrt{c^{\varepsilon,\tau}} \psi)\\
%&=&\int_{(0,t)\times\Omega}\nabla\left(\frac{\kappa(c^{\varepsilon,\tau})}{2\tau\sqrt{c^{\varepsilon,\tau}}}\right)~~\ln(1+\tau n^{\varepsilon,\tau})\cdot \nabla\sqrt{c^{\varepsilon,\tau}}\psi \\&&+\int_{(0,t)\times\Omega}\frac{\kappa(c^{\varepsilon,\tau})}{2\tau\sqrt{c^{\varepsilon,\tau}}}~~\nabla(\ln(1+\tau n^{\varepsilon,\tau}))\cdot\nabla\sqrt{c^{\varepsilon,\tau}}\psi \\
&=&\int_{(0,t)\times\Omega}\nabla\left(\frac{\kappa(c^{\varepsilon,\tau})}{2\tau\sqrt{c^{\varepsilon,\tau}}}\right)~~\ln(1+\tau n^{\varepsilon,\tau})\cdot \nabla\sqrt{c^{\varepsilon,\tau}}\psi\\&&+\frac 1{4} \int_{\mathbb{R}^3}\frac{\kappa(c^{\varepsilon,\tau})}{c^{\varepsilon,\tau}}\frac1{1+\tau n^{\varepsilon,\tau}}\nabla n^{\varepsilon,\tau}\cdot \nabla c^{\varepsilon,\tau}\psi\\
&=&\frac1{\tau}\int_{(0,t)\times\Omega}\ln(1+\tau n^{\varepsilon,\tau})\kappa'(c^{\varepsilon,\tau})|\nabla\sqrt{c^{\varepsilon,\tau}}|^{2}\psi\\&&-\frac1{2\tau}\int_{(0,t)\times\Omega}\ln(1+\tau n^{\varepsilon,\tau})\kappa(c^{\varepsilon,\tau})(\sqrt{c^{\varepsilon,\tau}})^{-2}|\nabla\sqrt{c^{\varepsilon,\tau}}|^{2}\psi\\&&+ \frac 1{4} \int_{(0,t)\times\Omega}\frac{\kappa(c^{\varepsilon,\tau})}{c^{\varepsilon,\tau}}\frac1{1+\tau n^{\varepsilon,\tau}}\nabla n^{\varepsilon,\tau}\cdot \nabla c^{\varepsilon,\tau}\psi.
\eeno
Combining all of them, we have
\ben\label{ine:c 1} 
&&\frac12   \int_{\Omega} (|\nabla \sqrt{c^{\varepsilon,\tau}}|^2 \psi)(\cdot,t) + \int_{(0,t)\times\Omega} |\nabla^2 \sqrt{c^{\varepsilon,\tau}}|^2 \psi \nonumber\\
&& =\frac12 \int_{(0,t)\times\Omega} |\nabla \sqrt{c^{\varepsilon,\tau}}|^2 (\partial_t \psi + \Delta \psi) + \frac 12 \int_{(0,t)\times\Omega} |\nabla \sqrt{c^{\varepsilon,\tau}}|^2 u^{\varepsilon,\tau} \cdot \nabla \psi\nonumber \\
\nonumber
&&-\int_{(0,t)\times\Omega} \nabla u^{\varepsilon,\tau} : (\nabla\sqrt{c^{\varepsilon,\tau}} \otimes \nabla\sqrt{c^{\varepsilon,\tau}})\psi - \int_{(0,t)\times\Omega} (\sqrt{c^{\varepsilon,\tau}})^{-1} |\nabla \sqrt{c^{\varepsilon,\tau}}|^2 \Delta \sqrt{c^{\varepsilon,\tau}} \psi \\  \nonumber
&&- \int_{(0,t)\times\Omega} (\sqrt{c^{\varepsilon,\tau}})^{-1} |\nabla \sqrt{c^{\varepsilon,\tau}}|^2 \nabla \sqrt{c^{\varepsilon,\tau}} \cdot \nabla \psi-\frac1{\tau}\int_{(0,t)\times\Omega}\ln(1+\tau n^{\varepsilon,\tau})\kappa'(c^{\varepsilon,\tau})|\nabla\sqrt{c^{\varepsilon,\tau}}|^{2}\psi\\\nonumber&&+\frac1{2\tau}\int_{(0,t)\times\Omega}\ln(1+\tau n^{\varepsilon,\tau})\kappa(c^{\varepsilon,\tau})(\sqrt{c^{\varepsilon}})^{-2}|\nabla\sqrt{c^{\varepsilon,\tau}}|^{2}\psi\\&&- \frac 1{4} \int_{(0,t)\times\Omega}\frac{\kappa(c^{\varepsilon,\tau})}{c^{\varepsilon,\tau}}\frac1{1+\tau n^{\varepsilon,\tau}}\nabla n^{\varepsilon,\tau}\cdot \nabla c^{\varepsilon,\tau}\psi.
\een
We remark that the bad terms are those integrals without $\nabla\psi$.
One bad term of all the above terms is $I\doteq-\int_{(0,t)\times\Omega} (\sqrt{c^{\varepsilon,\tau}})^{-1} |\nabla \sqrt{c^{\varepsilon,\tau}}|^2 \Delta \sqrt{c^{\varepsilon,\tau}} \psi$. The estimation of $ I $ is similar to $ I_{1} $ from (\ref{eq:c}) (see also (2.6) in \cite{CLW 2022}).
Then we have
\ben\label{eq:I}
I &\leq& - \frac 1{3} \sum_{i=j} \int_{(0,t)\times\Omega} (\sqrt{c^{\varepsilon,\tau}})^{-2} (\partial_j \sqrt{c^{\varepsilon,\tau}})^2 (\partial_i \sqrt{c^{\varepsilon,\tau}})^2 \psi +  \sum_{i\neq j} \int_{(0,t)\times\Omega} |\partial_{ij} \sqrt{c^{\varepsilon,\tau}}|^2 \psi \nonumber\\
&+& \int_{(0,t)\times\Omega} (\sqrt{c^{\varepsilon,\tau}})^{-1} |\nabla \sqrt{c^{\varepsilon,\tau}}|^2 \nabla \sqrt{c^{\varepsilon,\tau}} \cdot \nabla \psi.
\een
Submitting it to (\ref{ine:c 1}), we get
\ben\label{ine:c 1'}  \nonumber
&&\frac12  \int_{\Omega} (|\nabla \sqrt{c^{\varepsilon,\tau}}|^2 \psi)(\cdot,t) + \frac13\int_{(0,t)\times\Omega} |\triangle \sqrt{c^{\varepsilon,\tau}}|^2 \psi \nonumber\\
&&+\int_{(0,t)\times\Omega}\nabla\left(\frac{\kappa(c^{\varepsilon,\tau})}{2\tau\sqrt{c^{\varepsilon,\tau}}}\right)~~\ln(1+\tau n^{\varepsilon,\tau})\cdot \nabla\sqrt{c^{\varepsilon,\tau}}\psi \\\nonumber&&+\frac 1{3} \sum_{i=j} \int_{(0,t)\times\Omega} (\sqrt{c^{\varepsilon,\tau}})^{-2} (\partial_j \sqrt{c^{\varepsilon,\tau}})^2 (\partial_i \sqrt{c^{\varepsilon,\tau}})^2 \psi\nonumber\\
&&\leq~ \frac12 \int_{(0,t)\times\Omega} |\nabla \sqrt{c^{\varepsilon,\tau}}|^2 (\partial_t \psi + \Delta \psi) + \frac 12 \int_{(0,t)\times\Omega} |\nabla \sqrt{c^{\varepsilon,\tau}}|^2 u^{\varepsilon,\tau} \cdot \nabla \psi\nonumber \\
\nonumber
&&- \int_{(0,t)\times\Omega} \nabla u^{\varepsilon,\tau} : (\nabla\sqrt{c^{\varepsilon,\tau}} \otimes \nabla\sqrt{c^{\varepsilon,\tau}})\psi  - \frac 1{4} \int_{(0,t)\times\Omega}\frac{\kappa(c^{\varepsilon,\tau})}{c^{\varepsilon,\tau}}\frac1{1+\tau n^{\varepsilon,\tau}}\nabla n^{\varepsilon,\tau}\cdot \nabla c^{\varepsilon,\tau}\psi.
\een
Note that $ \|c^{\varepsilon,\tau}\|_{L^{\infty}}\leq \|c^{\varepsilon}_{0}\|_{L^{\infty}}$ and
\beno
\int_{(0,t)\times\Omega} (\sqrt{{c}^{\varepsilon,\tau}})^{-2} |\nabla \sqrt{{c}^{\varepsilon,\tau}}|^4 \psi&=& \sum_{i,j} \int_{(0,t)\times\Omega} (\sqrt{{c}^{\varepsilon,\tau}})^{-2} (\partial_j \sqrt{{c}^{\varepsilon,\tau}})^2 (\partial_i \sqrt{{c}^{\varepsilon,\tau}})^2 \psi\\
&\leq&  3\sum_{i=j} \int_{(0,t)\times\Omega} (\sqrt{{c}^{\varepsilon,\tau}})^{-2} (\partial_j \sqrt{{c}^{\varepsilon,\tau}})^2 (\partial_i \sqrt{c^{\varepsilon,\tau}})^2 \psi,
\eeno
using Young's inequality, we obtain
\ben\label{ine:g}
&-& \int_{(0,t)\times\Omega} \nabla u^{\varepsilon,\tau} : (\nabla\sqrt{c^{\varepsilon,\tau}} \otimes \nabla\sqrt{c^{\varepsilon,\tau}})\psi \\&\leq&\frac1{18} \int_{(0,t)\times\Omega} (\sqrt{c^{\varepsilon,\tau}})^{-2} |\nabla \sqrt{c^{\varepsilon,\tau}}|^4\psi\nonumber
+\frac92 \int_{(0,t)\times\Omega} (\sqrt{c^{\varepsilon,\tau}})^2|\nabla u^{\varepsilon,\tau}|^2 \psi\nonumber
\\&\leq&\frac1{18} \int_{(0,t)\times\Omega} (\sqrt{c^{\varepsilon,\tau}})^{-2} |\nabla \sqrt{c^{\varepsilon,\tau}}|^4\psi \nonumber
+\frac92 \|{c^{\varepsilon}_{0}}\|_{L^{\infty}} \int_{(0,t)\times\Omega} |\nabla u^{\varepsilon,\tau}|^2 \psi \nonumber
\\\nonumber&\leq&\frac16\sum_{i=j} \int_{(0,t)\times\Omega} (\sqrt{c^{\varepsilon,\tau}})^{-2} (\partial_j \sqrt{c^{\varepsilon,\tau}})^2 (\partial_i \sqrt{c^{\varepsilon,\tau}})^2 \psi+\frac92 \|{c^{\varepsilon}_{0}}\|_{L^{\infty}} \int_{(0,t)\times\Omega} |\nabla u^{\varepsilon,\tau}|^2 \psi.
\een
Substitute these estimates (\ref{eq:I})-(\ref{ine:g}) to the inequality of (\ref{ine:c 1}),  we obtain\label{ine:c 2}
\begin{equation}\label{sqrt c varepsilon tau}
\begin{aligned}
&\frac12 \int_{\Omega} (|\nabla \sqrt{c^{\varepsilon,\tau}}|^{2} \psi)(\cdot,t) + \frac13\int_{(0,t)\times\Omega} |\Delta \sqrt{c^{\varepsilon,\tau}}|^2 \psi\\&+~ \int_{(0,t)\times\Omega}\nabla\left(\frac{\kappa(c^{\varepsilon,\tau})}{2\tau\sqrt{c^{\varepsilon,\tau}}}\right)~~\ln(1+\tau n^{\varepsilon,\tau})\cdot \nabla\sqrt{c^{\varepsilon,\tau}}\psi \\ 
&+\frac16\sum_{i=j} \int_{(0,t)\times\Omega} (\sqrt{c^{\varepsilon,\tau}})^{-2} (\partial_j \sqrt{c^{\varepsilon,\tau}})^2 (\partial_i \sqrt{c^{\varepsilon,\tau}})^2 \psi\\
\leq~ &\frac12 \int_{(0,t)\times\Omega} |\nabla \sqrt{c^{\varepsilon,\tau}}|^2 (\partial_t \psi + \Delta \psi) + \frac 12 \int_{(0,t)\times\Omega} |\nabla \sqrt{c^{\varepsilon,\tau}}|^2 (u^{\varepsilon,\tau}\ast\rho^\varepsilon) \cdot \nabla \psi \\&+~\frac92\|{c^{\varepsilon}_{0}}\|_{L^{\infty}} \int_{(0,t)\times\Omega} |\nabla u^{\varepsilon,\tau}|^2 \psi- \frac 1{4} \int_{(0,t)\times\Omega}\frac{\kappa(c^{\varepsilon,\tau})}{c^{\varepsilon,\tau}}\frac1{1+\tau n^{\varepsilon,\tau}}\nabla n^{\varepsilon,\tau}\cdot \nabla c^{\varepsilon,\tau}\psi.
\end{aligned}
\end{equation}
Since $ \kappa(s)=\Theta_0 s\chi(s), $ we have
\begin{equation}\label{k'}
\begin{aligned}
-\frac{1}{\Theta_0}\int_{(0,t)\times\Omega}\frac{\kappa(c^{\varepsilon,\tau})}{c^{\varepsilon,\tau}}\frac{\nabla n^{\varepsilon,\tau}}{1+\tau n^{\varepsilon,\tau}}\cdot \nabla c^{\varepsilon,\tau}\psi= -\int_{(0,t)\times\Omega}\frac{\chi(c^{\varepsilon,\tau})}{1+\tau n^{\varepsilon,\tau}}\nabla n^{\varepsilon,\tau}\cdot \nabla c^{\varepsilon,\tau}\psi.
\end{aligned}
\end{equation}
Recall the local estimate of $n^{\varepsilon,\tau}$ in (\ref{n}), by taking
$(\ref{n}) + \frac{4}{\Theta_0}\times(\ref{sqrt c varepsilon tau})$ and using by (\ref{k'}) and (\ref{k''}), we arrive at, 
\ben\label{ine:energy n1 c1} 
\nonumber&& \int_{\Omega} (n^{\varepsilon,\tau} \ln (n^{\varepsilon,\tau}) \psi)(\cdot,t) + 4 \int_{(0,t)\times\Omega} |\nabla \sqrt{n^{\varepsilon,\tau}}|^2 \psi
\\ \nonumber&&	+ \frac{2}{\Theta_0}  \int_{\Omega} (|\nabla \sqrt{c^{\varepsilon,\tau}}|^2 \psi)(\cdot,t) +  \frac{4}{3\Theta_0}\int_{(0,t)\times\Omega} |\Delta \sqrt{c^{\varepsilon,\tau}}|^2 \psi 
\\ \nonumber
&&+ \frac{2}{3\Theta_0}  \sum_{i=j} \int_{(0,t)\times\Omega} (\sqrt{c^{\varepsilon,\tau}})^{-2} (\partial_j \sqrt{c^{\varepsilon,\tau}})^2 (\partial_i \sqrt{c^{\varepsilon,\tau}})^2 \psi\\\nonumber
&&\leq \int_{(0,t)\times\Omega} n^{\varepsilon,\tau} \ln n^{\varepsilon,\tau} (\partial_t \psi + \Delta \psi) + \int_{(0,t)\times\Omega} n^{\varepsilon,\tau} \ln n^{\varepsilon,\tau} (u^{\varepsilon,\tau}\ast\rho^\varepsilon) \cdot \nabla \psi\\ \nonumber
&&+\int_{(0,t)\times\Omega}\left(\frac{n^{\varepsilon,\tau}}{1+\tau n^{\varepsilon,\tau}} \nabla c^{\varepsilon,\tau} \chi(c^{\varepsilon,\tau})\right)\ln n^{\varepsilon,\tau}\cdot\nabla\psi +\int_{(0,t)\times\Omega}\left(\frac{n^{\varepsilon,\tau}}{1+\tau n^{\varepsilon,\tau}} \nabla c^{\varepsilon,\tau} \chi(c^{\varepsilon,\tau})\right)\cdot\nabla\psi\\ \nonumber
%&&+\int_{(0,t)\times\Omega} \left(\frac{1}{1+\tau n^{\varepsilon,\tau}} \nabla c^{\varepsilon,\tau} \chi(c^{\varepsilon,\tau})\right) \cdot \nabla n^{\varepsilon,\tau} \psi-\int_{(0,t)\times\Omega}\frac{\chi(c^{\varepsilon,\tau})}{1+\tau n^{\varepsilon,\tau}}\nabla n^{\varepsilon,\tau}\cdot \nabla c^{\varepsilon,\tau}\psi\\ \nonumber
&&+ \frac{2}{\Theta_0}\int_{(0,t)\times\Omega} |\nabla \sqrt{c^{\varepsilon,\tau}}|^2 (\partial_t \psi + \Delta \psi)
+ \frac{2}{\Theta_0} \int_{(0,t)\times\Omega} |\nabla \sqrt{c^{\varepsilon,\tau}}|^2 {(u^{\varepsilon,\tau}\ast\rho^\varepsilon)} \cdot \nabla \psi\\ &&+\frac{18}{\Theta_0}\|{c^{\varepsilon}_{0}}\|_{L^{\infty}} \int_{(0,t)\times\Omega} |\nabla u^{\varepsilon,\tau}|^2 \psi.
\een

For the equation of $u^{\varepsilon,\tau}$, multiplying $(\ref{eq:CNS})_3$ by $2u^{\varepsilon,\tau} \psi$ and integration by parts, we have

\beno
&&\int_{(0,t)\times\Omega}\partial_t u^{\varepsilon,\tau}\cdot 2u^{\varepsilon,\tau}\psi +\mu\int_{(0,t)\times\Omega} (u^{\varepsilon,\tau} \ast \rho^\varepsilon) \cdot \nabla u^{\varepsilon,\tau} \cdot 2u^{\varepsilon,\tau}\psi-\int_{(0,t)\times\Omega} \Delta u^{\varepsilon,\tau}\cdot 2u^{\varepsilon,\tau}\psi\\
&&+\int_{(0,t)\times\Omega} \nabla P^{\varepsilon,\tau}\cdot 2u^{\varepsilon,\tau}\psi +\int_{(0,t)\times\Omega} (n^{\varepsilon,\tau} \nabla \phi) \ast \rho^{\varepsilon,\tau}\cdot 2u^{\varepsilon,\tau}\psi=0.
\eeno
Integration by parts achieves that
\begin{equation}\label{ine:energy u}
\begin{aligned}
&\int_{\Omega} (|u^{\varepsilon,\tau}|^2)(\cdot,t) \psi + 2\int_{(0,t)\times\Omega} |\nabla u^{\varepsilon,\tau}|^2 \psi \\=&\int_{(0,t)\times\Omega} |u^{\varepsilon,\tau}|^2 \left(\partial_t \psi + \Delta \psi\right) + \mu\int_{(0,t)\times\Omega} |u^{\varepsilon,\tau}|^2(u^{\varepsilon,\tau}\ast\rho^\varepsilon)\cdot\nabla\psi\\&+2\int_{(0,t)\times\Omega} (P^{\varepsilon,\tau} - \bar{P}^{\varepsilon,\tau}) u^{\varepsilon,\tau} \cdot \nabla \psi
-2\int_{(0,t)\times\Omega} (n^{\varepsilon,\tau}\nabla\phi)\ast\rho^\varepsilon \cdot u^{\varepsilon,\tau} \psi.
\end{aligned}
\end{equation}
Then $(\ref{ine:energy u})\times\frac{18}{\Theta_0}\|{c^{\varepsilon}_{0}}\|_{L^{\infty}}+(\ref{ine:energy n1 c1})$ yields (\ref{local energy}). The proof is complete.

\begin{remark} It follows from (\ref{local energy}) that 
\begin{equation}\label{sqrt n'}
\int_{(0,t)\times\Omega} |\nabla \sqrt{n^{\varepsilon,\tau}}|^2\leq C.
\end{equation}
Indeed, for the bounded domain $\Omega$, 
using $|\ln n| n^{\alpha}\leq \alpha^{-1}e^{-1}$ for $0<n<1$ and $0<\alpha<\frac{1}{20}$,  we have
\beno\label{nlnn'} \nonumber
&& \int_{\Omega} (n^{\varepsilon,\tau} |\ln n^{\varepsilon,\tau}|)(\cdot,t) dx\\
&\leq& \int_{\Omega} (n^{\varepsilon,\tau}\ln n^{\varepsilon,\tau} )(\cdot,t) dx -2\int_{\Omega \cap \{x;0<n^{\varepsilon,\tau}<1\}}(n^{\varepsilon,\tau}\ln n^{\varepsilon,\tau} )(\cdot,t) dx
\nonumber \\
&\leq &\int_{\Omega} (n^{\varepsilon,\tau}\ln n^{\varepsilon,\tau})(\cdot,t) dx + 2\alpha^{-1}e^{-1} \int_{\Omega}  (n^{\varepsilon,\tau(1-\alpha)})(\cdot,t) dx\nonumber \\
&\leq &\int_{\Omega} (n^{\varepsilon,\tau}\ln n^{\varepsilon,\tau})(\cdot,t) dx + 2\alpha^{-1}e^{-1} C^{1-\alpha}.\\
%&\leq &\int_{\Omega} (n^{\varepsilon,\tau}+1)\ln (n^{\varepsilon,\tau}+1)(\cdot,t) dx + 2\alpha^{-1}e^{-1} C^{1-\alpha}\\
%&\leq& C.
\eeno
thus adding $2\alpha^{-1}e^{-1} C^{1-\alpha}$ to both sides of the inequality of  (\ref{local energy}), the proof is complete.
\end{remark}

\section{The convergence of $u^{\varepsilon,\tau}$, $n^{\varepsilon,\tau}$ and $c^{\varepsilon,\tau}$}

In the following, we prove the strong convergence of $u^{\varepsilon,\tau}$, $ n^{\varepsilon,\tau} $ and $c^{\varepsilon,\tau}$ respectively by Aubin-Lions Lemma.

\begin{lemma}\label{lem:u}
	$u^{\varepsilon,\tau} \rightarrow u^{\varepsilon}$ in $L^q(0,t;L^q(\Omega)$) strongly for $q \in [2,\frac{10}3)$ as $ \tau\to 0 $.
\end{lemma}
{\bf Proof.} Recalling the equation of $u^{\varepsilon,\tau}$ of $(\ref{eq:GKS})_3$,
%\beno
%\partial_t u^{\varepsilon,\tau} + (u^{\varepsilon,\tau} \ast \rho^\varepsilon) \cdot \nabla u^{\varepsilon,\tau} - \Delta u^{\varepsilon,\tau} + \nabla P^{\varepsilon,\tau} = - (n^{\varepsilon,\tau }\nabla \phi) \ast \rho^\varepsilon.
%\eeno
for any smooth function $\zeta \in C_c^\infty(\Omega \times (0,t))$ with $\nabla \cdot \zeta = 0$, there holds
\beno
&&\int_{(0,t)\times\Omega} \partial_t u^{\varepsilon,\tau} \cdot \zeta \\
%&=& \int_{(0,t)\times\Omega} \Delta u^{\varepsilon,\tau} \cdot \zeta - \int_{(0,t)\times\Omega} (u^{\varepsilon,\tau} \ast \rho^\varepsilon) \cdot \nabla u^{\varepsilon,\tau} \cdot \zeta - \int_{(0,t)\times\Omega} (n^{\varepsilon,\tau} \nabla \phi) \ast \rho^\varepsilon \cdot \zeta \\
&\leq& C \|\nabla u^{\varepsilon,\tau}\|_{L^{2}(0,t;L^{2}(\Omega))} \|\nabla \zeta\|_{L^{2}(0,t;L^{2}(\Omega))} \\&&+~ C \|u^{\varepsilon,\tau} \ast \rho^\varepsilon\|_{L^{\frac83}(0,t;L^{4}(\Omega))} \|u^{\varepsilon,\tau}\|_{L^{\frac83}(0,t;L^{4}(\Omega))}  \|\nabla \zeta\|_{L^{4}(0,t;L^{2}(\Omega))} \\
&&+~ C \|(n^{\varepsilon,\tau} \nabla \phi) \ast \rho^\varepsilon\|_{L^{\frac53}(0,t;L^{\frac65}(\Omega))} \|\zeta\|_{L^{\frac52}(0,t;L^{6}(\Omega))} \\
&\leq& C \|\nabla u^{\varepsilon,\tau}\|_{L^{2}(0,t;L^{2}(\Omega))} \|\nabla \zeta\|_{L^{2}(0,t;L^{2}(\Omega))} + C \|u^{\varepsilon,\tau}\|_{L^{\frac83}(0,t;L^{4}(\Omega))}^2 \|\nabla \zeta\|_{L^{4}(0,t;L^{2}(\Omega))} \\
&&+ ~C \|\nabla \phi\|_{{L^\infty}(\Omega)} \|n^{\varepsilon,\tau}\|_{L^{\frac53}(0,t;L^{\frac65}(\Omega))} \|\nabla \zeta\|_{L^{\frac52}(0,t;L^{2}(\Omega))}.
\eeno
%Here, we use the strong $(p,p)-$type operator for the convolution operator and the boundary of the support set of $\zeta$ since $\zeta$ has a compact supported set.
By $ (\ref{ine:  energy}) $ and (\ref{sqrt n'}), we know  $n^{\varepsilon,\tau} \in L^\infty (0,t;L^1(\mathbb{R}^3))$ and $\nabla \sqrt{n^{\varepsilon,\tau}} \in L^2 (0,t;L^2(\Omega))$, note that
\beno
\|\sqrt{n^{\varepsilon,\tau}}\|_{L^q_tL^p_x} &\leq& C (\|\sqrt{n^{\varepsilon,\tau}}\|_{L^\infty_tL^2_x}+\|\nabla \sqrt{n^{\varepsilon,\tau}}\|_{L^2_tL^2_x}),\quad \frac3p+\frac2q=\frac32,\quad 2\leq p\leq 6,
\eeno
especially, we have
%\beno
%\|\sqrt{n^{\varepsilon,\tau}}\|_{L^\frac{10}3(0,t;L^\frac{10}3(\Omega))} \leq C,
%\eeno
%which means
\begin{equation}\label{ine:n}
\|n^{\varepsilon,\tau}\|_{L^\frac53(0,t; L^\frac53(\Omega))} \leq C.
\end{equation}
Since $\Omega$ is bounded, it follows that
\beno
\|n^{\varepsilon,\tau}\|_{L^{\frac53}(0,t;L^{\frac65}(\Omega))} \leq C(|\Omega|) \|n^{\varepsilon,\tau}\|_{L^{\frac53}(0,t;L^{\frac53}(\Omega))} \leq C.\eeno
Similarly, we have
\begin{equation}\label{ine:u1}
	\|u^{\varepsilon,\tau}\|_{L^{\frac83}(0,t;L^{4}(\Omega))}\leq C.
\end{equation}
Due to $u^{\varepsilon,\tau} \in L^\infty (0,t;L^2(\mathbb{R}^3))$ and $\nabla u^{\varepsilon,\tau} \in L^2 (0,t;L^2(\mathbb{R}^3))$ from (\ref{ine:  energy}), we have
\beno
\int_{(0,t)\times\Omega} \partial_t u^{\varepsilon,\tau} \cdot \zeta &\leq& C \left(\|\nabla \zeta\|_{L^{\frac52}(0,t;L^{2}(\Omega))} + \|\nabla \zeta\|_{L^{2}(0,t;L^{2}(\Omega))} + \|\nabla \zeta\|_{L^{4}(0,t;L^{2}(\Omega))}\right) \\&\leq& C(t,|\Omega|)\|\nabla \zeta\|_{L^4(0,t;L^2(\Omega))},
\eeno
which means that
\beno
\partial_t u^{\varepsilon,\tau} \in L^\frac43(0,t;(\dot{W}^{1,2}_\sigma(\Omega))').
\eeno
Noting that $u^{\varepsilon,\tau} \in L^2(0,t;\dot{W}^{1,2}_\sigma(\Omega))$, by Lemma \ref{Aubin-Lions lemma}, we have
\beno
u^{\varepsilon,\tau} \rightarrow u^{\varepsilon} \quad {\rm in} \quad L^2(0,t;L^2_\sigma(\Omega))
\eeno
as $ \tau\to 0. $ Furthermore, we can obtain
\beno
u^{\varepsilon,\tau} \rightarrow u^{\varepsilon} \quad {\rm in} \quad L^q(0,t;L^q(\Omega)) \quad {\rm for} \quad q\in [2,\frac{10}3)
\eeno
as $ \tau\to 0. $

\begin{lemma}\label{lem:c}
	$\nabla c^{\varepsilon,\tau} \rightarrow \nabla c^{\varepsilon}$ in $L^q(0,t;L^q(\Omega)$) strongly for $q \in [2,\frac{10}3)$ as $ \tau\to 0 $.
\end{lemma}

{\bf Proof.} Recall the equation of $c^{\varepsilon,\tau}$ of  $(\ref{eq:GKS})_2$.
Taking the derivative of the above equation, we have
\beno
\partial_t \nabla c^{\varepsilon,\tau} + \nabla (u^{\varepsilon,\tau}  \cdot \nabla c^{\varepsilon,\tau}) - \Delta \nabla c^{\varepsilon,\tau} = - \nabla \left(\frac1{\tau}\ln (1+\tau n^{\varepsilon,\tau})\kappa(c^{\varepsilon,\tau})\right).\\
\eeno
For any smooth function $\zeta \in C_c^\infty(\Omega \times (0,t))$, we have
\beno
&&\int_{(0,t)\times\Omega}\partial_t \nabla c^{\varepsilon,\tau} \cdot \zeta\\ 
%&=& \int_{(0,t)\times\Omega} - \nabla ((u^{\varepsilon,\tau} \ast \rho^\varepsilon) \cdot \nabla c^{\varepsilon,\tau}) \cdot \zeta + \int_{(0,t)\times\Omega}  \Delta \nabla c^{\varepsilon,\tau} \cdot \zeta \\&&-~\int_{(0,t)\times\Omega} \nabla\left(\frac1{\tau}\ln (1+\tau n^{\varepsilon,\tau})\kappa(c^{\varepsilon,\tau})\right) \cdot \zeta \\
&=& \int_{(0,t)\times\Omega} u^{\varepsilon,\tau}  \cdot \nabla c^{\varepsilon,\tau} \nabla \cdot \zeta - \int_{(0,t)\times\Omega}  \Delta c^{\varepsilon,\tau} \nabla\cdot \zeta \\&&+~ \int_{(0,t)\times\Omega} \left(\frac1{\tau}\ln (1+\tau n^{\varepsilon,\tau})\kappa(c^{\varepsilon,\tau})\right) \nabla \cdot \zeta \\
%&\leq& \|u^{\varepsilon,\tau} \ast \rho^\varepsilon\|_{L^{\frac83}(0,t;L^{4}(\Omega))} \|\nabla c^{\varepsilon,\tau}\|_{L^{\frac83}(0,t;L^{4}(\Omega))} \|\nabla \zeta\|_{L^{4}(0,t;L^{2}(\Omega))} \\&&+~ \|\Delta c^{\varepsilon,\tau}\|_{L^{2}(0,t;L^{2}(\Omega))} \|\nabla \zeta\|_{L^{2}(0,t;L^{2}(\Omega))} \\
%&&+ ~\|\kappa\|_{L^{\infty}(0,\|c_{0}^{\varepsilon}\|_{L^{\infty}})} \|\frac1{\tau}\ln (1+\tau n^{\varepsilon,\tau})\|_{L^{\frac43}(0,t;L^{2}(\Omega))} \|\nabla \zeta\|_{L^{4}(0,t;L^{2}(\Omega))} \\
&\leq& C \|u^{\varepsilon,\tau}\|_{L^{\frac83}(0,t;L^{4}(\Omega))} \|\nabla c^{\varepsilon,\tau}\|_{L^{\frac83}(0,t;L^{4}(\Omega))} \|\nabla \zeta\|_{L^{4}(0,t;L^{2}(\Omega))} \\&&+~ C \|\Delta c^{\varepsilon,\tau}\|_{L^{2}(0,t;L^{2}(\Omega))} \|\nabla \zeta\|_{L^{2}(0,t;L^{2}(\Omega))} \\
&&+~ C\|\kappa\|_{L^{\infty}(0,\|c_{0}^{\varepsilon}\|_{L^{\infty}})} \|n^{\varepsilon,\tau}\|_{L^{\frac43}(0,t;L^{2}(\Omega))} \|\nabla \zeta\|_{L^{4}(0,t;L^{2}(\Omega))}.
\eeno
Since $\nabla\sqrt{c^{\varepsilon,\tau}} \in L^\infty(0,t;L^2(\Omega) ) \cap L^2(0,t;\dot{H}^1(\Omega)), $ there holds
\begin{equation}\label{sqrt c}
\|\nabla \sqrt{c^{\varepsilon,\tau}}\|_{L^{\frac83}(0,t;L^{4}(\Omega))} \leq C.
\end{equation}
Similarly, we obtain
\begin{equation}\label{ine:n2}
\|n^{\varepsilon,\tau}\|_{L^{\frac43}(0,t;L^{2}(\Omega))}\leq C.
\end{equation}
Noting that
\begin{equation*}
\nabla c^{\varepsilon,\tau} = 2 \sqrt{c^{\varepsilon,\tau}} \nabla \sqrt{c^{\varepsilon,\tau}},
\end{equation*}
and $\|c^{\varepsilon,\tau}\|_{L^\infty(0,t;L^\infty(\Omega))} \leq \|c_0^{\varepsilon}\|_{L^\infty}\leq C,$
we have
\begin{equation}\label{ine:nable c}
\|\nabla c^{\varepsilon,\tau}\|_{L^{\frac83}(0,t;L^{4}(\Omega))} \leq 2 \|c^{\varepsilon,\tau}\|_{L^\infty(0,t;L^\infty(\Omega))}^\frac12 \|\nabla \sqrt{c^{\varepsilon,\tau}}\|_{L^{\frac83}(0,t;L^{4}(\Omega))}\leq C.
\end{equation}
%Claim that $\|\Delta c^{\varepsilon,\tau}\|_{L^{2}(0,t;L^{2}(\Omega))} \leq C(\|c_0\|_{L^\infty} )$.
%(\ref{ine:  energy}) implies that $\|u^\varepsilon\|_{L^{\frac83}(0,t;L^{4}(\Omega))} \leq C$ since $u^\varepsilon \in L^\infty(0,t;L^2(\Omega) ) \cap L^2(0,t;\dot{H}^1(\Omega)) $.
%Similarly, we have $\|\nabla \sqrt{c^\varepsilon}\|_{L^{\frac83}(0,t;L^{4}(\Omega))} \leq C$. Noting that
%\beno
%\nabla c^\varepsilon = 2 \sqrt{c^\varepsilon} \nabla \sqrt{c^\varepsilon},
%\eeno
%we have
%\beno
%\|\nabla c^\varepsilon\|_{L^{\frac83}(0,t;L^{4}(\Omega))} \leq 2 \|c^\varepsilon\|_{L^\infty(0,t;L^\infty(\Omega))}^\frac12 \|\nabla \sqrt{c^\varepsilon}\|_{L^{\frac83}(0,t;L^{4}(\Omega))}.
%\eeno
%Since $\|c^\varepsilon\|_{L^\infty(0,t;L^\infty(\Omega))} \leq \|c_0\|_{L^\infty}$, we have $\|\nabla c^\varepsilon\|_{L^{\frac83}(0,t;L^{4}(\Omega))} \leq C$.
%We also need the estimate $\|\Delta c^\varepsilon\|_{L^{2}(0,t;L^{2}(\Omega))} \leq C$. 
Noting that
\beno
\Delta c^{\varepsilon,\tau} = 2 \nabla \sqrt{c^{\varepsilon,\tau}} \otimes \nabla \sqrt{c^{\varepsilon,\tau}} + 2 \sqrt{c^{\varepsilon,\tau}} \Delta \sqrt{c^{\varepsilon,\tau}},
\eeno
%and $\|\nabla^2 \sqrt{c^\varepsilon}\|_{L^2(0,t;L^2(\Omega))} \leq C$ by Lemma \ref{lem:same energy},
we arrive
\ben\label{ine: delta c}
\|\Delta c^{\varepsilon,\tau}\|_{L^2(0,t;L^2(\Omega))}^2 &\leq& 2 \int_{(0,t)\times \Omega} |\nabla \sqrt{c^{\varepsilon,\tau}}|^4 + 2 \int_{(0,t)\times \Omega} c^{\varepsilon,\tau} |\Delta \sqrt{c^{\varepsilon,\tau}}|^2\\\nonumber
&\leq& 2 \|c^{\varepsilon,\tau}\|_{L^\infty(0,t;L^\infty(\Omega))} \int_{(0,t)\times \Omega} (\sqrt{c^{\varepsilon,\tau}})^{-2} |\nabla \sqrt{c^{\varepsilon,\tau}}|^4 \\\nonumber
&&+~ 2 \|c^{\varepsilon,\tau}\|_{L^\infty(0,t;L^\infty(\Omega))} \int_{(0,t)\times \Omega} |\Delta \sqrt{c^{\varepsilon,\tau}}|^2.
\een
By (\ref{ine:  energy}) of Lemma \ref{lem:same energy}, we have $\|\Delta c^{\varepsilon,\tau}\|_{L^2(0,t;L^2(\Omega))} \leq C$,
where use $ \Delta\sqrt{c^{\varepsilon,\tau}}\in L^{2}(0,t;L^{2}(\Omega)) $ and $ \int_{(0,t)\times \Omega}(\sqrt{c^{\varepsilon,\tau}})^{-2}|\nabla\sqrt{c^{\varepsilon,\tau}}|^{4}\leq C $. Combining (\ref{ine:u1}) and (\ref{ine:n2})-(\ref{ine: delta c}), we have

%Moreover, due to $ \kappa(s)\in C(\mathbb{R})$ and $\|c^\varepsilon\|_{L^\infty(0,t;L^\infty(\Omega))} \leq \|c_0\|_{L^\infty}$, we have $ \|\kappa(c^{\varepsilon})\|_{L^{\infty}(0,t;L^{\infty}(\Omega))}\leq C. $

%By the priori estimates of $n^\varepsilon$, $c^\varepsilon$ and $u^\varepsilon$, we have
\beno
\int_0^t \int_\Omega \partial_t \nabla c^{\varepsilon,\tau} \cdot \zeta &\leq& C \left(\|\nabla \zeta\|_{L^{2}(0,t;L^{2}(\Omega))} + \|\nabla \zeta\|_{L^{4}(0,t;L^{2}(\Omega))}\right) \\&\leq& C(t)\|\nabla \zeta\|_{L^{4}(0,t;L^{2}(\Omega))},
\eeno
which means that
\beno
\partial_t \nabla c^{\varepsilon,\tau} \in  L^\frac43(0,t;(\dot{W}^{1,2}(\Omega))').
\eeno
Noting that $\nabla c^{\varepsilon,\tau} \in L^2(0,t;\dot{W}^{1,2}(\Omega))$, by Lemma \ref{Aubin-Lions lemma}, we have
\beno
\nabla c^{\varepsilon,\tau} \rightarrow \nabla c^{\varepsilon} \quad {\rm in} \quad  L^2(0,t;L^2(\Omega))
\eeno
as $ \tau\to 0. $
Similar to Lemma \ref{lem:u}, we obtain
\beno
\nabla c^{\varepsilon,\tau} \rightarrow \nabla c^{\varepsilon} \quad {\rm in} \quad L^q(0,t;L^q(\Omega)) \quad {\rm for} \quad q\in [2,\frac{10}3)
\eeno
as $ \tau\to 0. $

\begin{lemma}\label{lem:cc}
	$\nabla \sqrt{c^{\varepsilon,\tau}} \rightarrow \nabla \sqrt{c^{\varepsilon}}$ in $L^q(0,t;L^q(\Omega)$) strongly for $q \in [2,\frac{10}3)$ as $ \tau\to 0 $.
	%$u^\varepsilon \rightarrow u$ in $L^q(0,T;L^q(\mathbb{R}^3))$ strongly for $q \in [2,\frac{10}3)$;
	%$n^\varepsilon \rightarrow n$ in $L^{2m}(0,T;L^2(\mathbb{R}^3))$ strongly.
\end{lemma}
{\bf Proof.} Consider the equation of $\sqrt{c^{\varepsilon,\tau}}$ as follows:

\beno
\partial_t\sqrt{c^{\varepsilon,\tau}}+u^{\varepsilon,\tau}\cdot \nabla\sqrt{c^{\varepsilon,\tau}}-\frac{|\nabla \sqrt{c^{\varepsilon,\tau}}|^2}{\sqrt{c^{\varepsilon,\tau}}}-\Delta\sqrt{c^{\varepsilon,\tau}}=- \frac{\kappa(c^{\varepsilon,\tau})}{2\tau\sqrt{c^{\varepsilon,\tau}}}~~\ln(1+\tau n^{\varepsilon,\tau}). 
\eeno
Taking the derivative of the above equation, we have
\beno
&&\partial_t\nabla\sqrt{c^{\varepsilon,\tau}}+\nabla \left(u^{\varepsilon,\tau}\cdot \nabla\sqrt{c^{\varepsilon,\tau}}\right)-\nabla\frac{|\nabla \sqrt{c^{\varepsilon,\tau}}|^2}{\sqrt{c^{\varepsilon,\tau}}}-\Delta\nabla\sqrt{c^{\varepsilon,\tau}}\\&=& - \frac 1{2\tau} \nabla \left( \frac{\kappa(c^{\varepsilon,\tau})}{\sqrt{c^{\varepsilon,\tau}}}~~\ln(1+\tau n^{\varepsilon,\tau})\right).
\eeno
Now for any smooth function $\zeta \in C_c^\infty(\Omega \times (0,t))$, we know
\beno
&&\int_{(0,t)\times\Omega} \partial_t \nabla \sqrt{c^{\varepsilon,\tau}} \cdot \zeta \\
%&=& \int_{(0,t)\times\Omega} - \nabla ((u^{\varepsilon,\tau} \ast \rho^\varepsilon) \cdot \nabla \sqrt{c^{\varepsilon,\tau}}) \cdot \zeta + \int_{(0,t)\times\Omega}  \Delta \nabla \sqrt{c^{\varepsilon,\tau}} \cdot \zeta \\
%&&-~ \frac1{2\tau} \int_{(0,t)\times\Omega} \nabla \left( \frac{\kappa(c^{\varepsilon,\tau})}{\sqrt{c^{\varepsilon,\tau}}}~~\ln(1+\tau n^{\varepsilon,\tau})\right) \cdot \zeta +  \int_{(0,t)\times\Omega}\nabla\frac{|\nabla \sqrt{c^{\varepsilon,\tau}}|^2}{\sqrt{c^{\varepsilon,\tau}}}\cdot \zeta\\
&=&  \int_{(0,t)\times\Omega} u^{\varepsilon,\tau} \cdot \nabla \sqrt{c^{\varepsilon,\tau}} \nabla \cdot \zeta -  \int_{(0,t)\times\Omega}  \nabla^2 \sqrt{c^{\varepsilon,\tau}} : \nabla \zeta \\
&&+~ \frac 1{2\tau}  \int_{(0,t)\times\Omega} \frac{\kappa(c^{\varepsilon,\tau})-\kappa(0)}{c^{\varepsilon,\tau}} \sqrt{c^{\varepsilon,\tau}} \ln(1+\tau n^{\varepsilon,\tau}) \nabla \cdot \zeta
- \int_{(0,t)\times\Omega}\frac{|\nabla \sqrt{c^{\varepsilon,\tau}}|^2}{\sqrt{c^{\varepsilon,\tau}}} \nabla\cdot\zeta\\
&\leq&  \big|\int_{(0,t)\times\Omega} u^{\varepsilon,\tau}  \cdot \nabla \sqrt{c^{\varepsilon,\tau}} \nabla \cdot \zeta\big|+  \big|\int_{(0,t)\times\Omega}  \nabla^2 \sqrt{c^{\varepsilon,\tau}} : \nabla \zeta\big| \\
&&+~ \frac 1{2\tau}  \big|\int_{(0,t)\times\Omega} \frac{\kappa(c^{\varepsilon,\tau})-\kappa(0)}{c^{\varepsilon,\tau}} \sqrt{c^{\varepsilon,\tau}}(\tau n^{\varepsilon,\tau}) \nabla \cdot \zeta\big|
+\big|\int_{(0,t)\times\Omega}\frac{|\nabla \sqrt{c^{\varepsilon,\tau}}|^2}{\sqrt{c^{\varepsilon,\tau}}} \nabla\cdot\zeta\big|\\
%&\leq& \|u^{\varepsilon,\tau} \ast \rho^\varepsilon\|_{L^{\frac83}(0,t;L^{4}(\Omega))} \|\nabla \sqrt{c^{\varepsilon,\tau}}\|_{L^{\frac83}(0,t;L^{4}(\Omega))} \|\nabla \zeta\|_{L^{4}(0,t;L^{2}(\Omega))}\\&& +~ \|\nabla^2 \sqrt{c^{\varepsilon,\tau}}\|_{L^{2}(0,t;L^{2}(\Omega))} \|\nabla \zeta\|_{L^{2}(0,t;L^{2}(\Omega))} \\
%&&+~ \|\kappa'\|_{L^{\infty}(0,\|c_{0}\|_{L^{\infty}})}\|c^{\varepsilon,\tau}\|_{L^{\infty}(0,t;L^{\infty}(\Omega))}^{\frac12} \|n^{\varepsilon,\tau}\|_{L^{\frac43}(0,t;L^{2}(\Omega))} \|\nabla \zeta\|_{L^{4}(0,t;L^{2}(\Omega))} \\
%&&+~\|\frac{|\nabla \sqrt{c^{\varepsilon,\tau}}|^2}{\sqrt{c^{\varepsilon,\tau}}}\|_{L^{2}(0,t;L^{2}(\Omega))}\|\nabla \zeta\|_{L^{2}(0,t;L^{2}(\Omega))}\\
&\leq& C \|u^{\varepsilon,\tau} \|_{L^{\frac83}(0,t;L^{4}(\Omega))} \|\nabla \sqrt{c^{\varepsilon,\tau}}\|_{L^{\frac83}(0,t;L^{4}(\Omega))} \|\nabla \zeta\|_{L^{4}(0,t;L^{2}(\Omega))}\\
&& +~ \|\nabla^2 \sqrt{c^{\varepsilon,\tau}}\|_{L^{2}(0,t;L^{2}(\Omega))} \|\nabla \zeta\|_{L^{2}(0,t;L^{2}(\Omega))}  \\
&&+~ C\|\kappa'\|_{L^{\infty}(0,\|c^{\varepsilon}_{0}\|_{L^{\infty}})}\|c^{\varepsilon,\tau}\|_{L^{\infty}(0,t;L^{\infty}(\Omega))}^{\frac12} \|n^{\varepsilon,\tau} \|_{L^{\frac43}(0,t;L^{2}(\Omega))} \|\nabla \zeta\|_{L^{4}(0,t;L^{2}(\Omega))} \\
&&+~\|\frac{|\nabla \sqrt{c^{\varepsilon,\tau}}|^2}{\sqrt{c^{\varepsilon,\tau}}}\|_{L^{2}(0,t;L^{2}(\Omega))}\|\nabla \zeta\|_{L^{2}(0,t;L^{2}(\Omega))}.
\eeno
%Now we deal with the term $\|n^{\varepsilon,\tau} \|_{L^{\frac43}(0,t;L^{2}(\Omega))}.$ 
Using (\ref{ine:  energy}) and interpolation inequality, we have
\begin{equation*}
\begin{aligned}
\|n^{\varepsilon,\tau}\|_{L^{\frac43}(0,t;L^{2}(\Omega))}=\|\sqrt{n^{\varepsilon,\tau}}\|_{L^{\frac83}(0,t;L^{4}(\Omega))}\leq C.
\end{aligned}
\end{equation*}
By (\ref{ine:c L^1 L^infty'}), (\ref{ine:u1}), (\ref{sqrt c}) -(\ref{ine:nable c}) and
$\kappa(s) \in C^2(\overline{\mathbb{R^{+}}})$, we have
\beno
\int_{(0,t)\times\Omega} \partial_t \nabla \sqrt{c^{\varepsilon,\tau}} \cdot \zeta &\leq& C \left(\|\nabla \zeta\|_{L^2(0,t;L^2(\Omega))}+ \|\nabla \zeta\|_{L^4(0,t;L^2(\Omega))}\right)\\& \leq &C \|\nabla \zeta\|_{L^4(0,t;L^2(\Omega))},
\eeno
which means that
\beno
\partial_t \nabla \sqrt{c^{\varepsilon,\tau}} \in L^\frac43(0,t;(\dot{W}^{1,2}(\Omega))').
\eeno
Noting that $\nabla \sqrt{c^{\varepsilon,\tau}} \in L^2(0,t;\dot{W}^{1,2}(\Omega))$, by Lemma \ref{Aubin-Lions lemma}, we have
\beno
\nabla \sqrt{c^{\varepsilon,\tau}} \rightarrow \nabla \sqrt{c^{\varepsilon}} \quad {\rm in} \quad L^2(0,t;L^2(\Omega))
\eeno
as $ \tau\to 0. $
Similarly, we have
\beno
\nabla \sqrt{c^{\varepsilon,\tau}} \rightarrow \nabla\sqrt{c^{\varepsilon}} \quad {\rm in} \quad L^q(0,t;L^q(\Omega)) \quad {\rm for} \quad q\in [2,\frac{10}3)
\eeno
as $ \tau\to 0. $

\begin{lemma}\label{lem:n}
	$ n^{\varepsilon,\tau} \rightarrow  n^{\varepsilon}$ strongly in $L^q(0,t;L^q(\Omega))$ for $q \in [1,\frac53)$ as $ \tau\to 0 $.
\end{lemma}
{\bf Proof.} For any smooth function $\zeta \in C_c^\infty(\Omega \times (0,t))$, we know
\begin{equation*}
\begin{aligned}
&\int_{(0,t)\times\Omega}\partial_t  n^{\varepsilon,\tau}\cdot \zeta\\=~&- \int_{\mathbb{R}^3} \left((u^{\varepsilon,\tau} \ast \rho^\varepsilon) \cdot \nabla n^{\varepsilon,\tau} \right)\cdot \zeta+ \int_{(0,t)\times\Omega}\Delta n^{\varepsilon,\tau}\cdot \zeta\\&-\int_{(0,t)\times\Omega} \nabla \cdot \left(\frac{n^{\varepsilon,\tau}}{1+\tau n^{\varepsilon,\tau}} \nabla c^{\varepsilon,\tau} \chi(c^{\varepsilon,\tau})\right)\cdot \zeta\\\leq~&\|u^{\varepsilon,\tau}\|_{L^{\frac{10}{3}}(0,t;L^{\frac{10}{3}}(\Omega))}\|n^{\varepsilon,\tau}\|_{L^{\frac{5}{3}}(0,t;L^{\frac{5}{3}}(\Omega))}\|\nabla\zeta\|_{L^{10}(0,t;L^{10}(\Omega))}\\&+C\|\sqrt{n^{\varepsilon,\tau}}\|_{L^{\frac{10}{3}}(0,t;L^{\frac{10}{3}}(\Omega))}\|\nabla\sqrt{n^{\varepsilon,\tau}}\|_{L^{2}(0,t;L^{2}(\Omega))}\|\nabla\zeta\|_{L^{5}(0,t;L^{5}(\Omega))}\\&+C\|c_{0}^{\varepsilon}\|_{L^{\infty}(\Omega)}^{\frac12}\|\chi\|_{L^{\infty}(0,\|c_{0}^{\varepsilon}\|_{L^{\infty}})}\|\nabla\sqrt{c^{\varepsilon,\tau}}\|_{L^{\frac{10}{3}}(0,t;L^{\frac{10}{3}}(\Omega))}\|n^{\varepsilon,\tau}\|_{L^{\frac{5}{3}}(0,t;L^{\frac{5}{3}}(\Omega))}\|\nabla\zeta\|_{L^{10}(0,t;L^{10}(\Omega))}
\end{aligned}
\end{equation*}
By uniform estimate (\ref{ine:  energy}), we have the boundedness of $\|u^{\varepsilon,\tau}\|_{L^{\frac{10}{3}}(0,t;L^{\frac{10}{3}}(\Omega))}$,
$\|\nabla\sqrt{c^{\varepsilon,\tau}}\|_{L^{\frac{10}{3}}(0,t;L^{\frac{10}{3}}(\Omega))}$ and $\|n^{\varepsilon,\tau}\|_{L^{\frac{5}{3}}(0,t;L^{\frac{5}{3}}(\Omega))}$, thus
\beno
&&\int_{(0,t)\times\Omega} \partial_t n^{\varepsilon,\tau}\cdot \zeta \leq C\|\nabla\zeta\|_{L^{10}(0,t;L^{10}(\Omega))}.
\eeno
It follows that $\partial_t n^{\varepsilon,\tau}\in L^\frac{10}{9}(0,t;(\dot{W}^{1,10}(\Omega))').$
Using H\"{o}lder inequality, we get
\begin{equation*}
\begin{aligned}
\int_{(0,t)\times\Omega}|\nabla n^{\varepsilon,\tau}|^{\frac54}&=~\int_{(0,t)\times\Omega}|\sqrt{n^{\varepsilon,\tau}}\nabla \sqrt{n^{\varepsilon,\tau}}|^{\frac54}\\&\leq~ \left(\int_{(0,t)\times\Omega}|\nabla\sqrt{n^{\varepsilon,\tau}}|^2\right)^{\frac 58}\left(\int_{(0,t)\times\Omega}|n^{\varepsilon,\tau}|^{\frac53}\right)^{\frac 38}\leq C.
\end{aligned}
\end{equation*}
We choose $X=W^{1,\frac54}(\Omega)$, $Y=L^2{(\Omega)}$ and $Z=W^{-1,\frac{10}9}(\Omega)$, $p=\frac54$, $r=\frac{10}{9}$.
By Lemma \ref{Aubin-Lions lemma}, we have 
\beno
n^{\varepsilon,\tau}\rightarrow  n^{\varepsilon} \quad L^{\frac54}(0,t;L^2(\Omega)).
\eeno
Similarly, we have
\beno
n^{\varepsilon,\tau}\rightarrow  n^{\varepsilon} \quad {\rm in} \quad L^q(0,t;L^q(\Omega)) \quad {\rm for} \quad q\in [1,\frac{5}3)
\eeno
as $ \tau\to 0. $
\endproof

\begin{lemma}\label{lem:weak convergence of n}
	$(n^{\varepsilon,\tau},c^{\varepsilon,\tau},u^{\varepsilon,\tau})$ is a weak solution to the Cauchy problem (\ref{eq:CNS}), there holds
	\beno
	n^{\varepsilon,\tau} \ln n^{\varepsilon,\tau} \in L^{\frac32}(0,t;L^{\frac32}(\Omega)).
	\eeno
\end{lemma}

{\bf Proof.} Set $Q = (0,t) \times \Omega$, $ \Omega\subset\mathbb{R}^{3} $ is an arbitrary bounded domain, we have
\beno
\int_{0}^{t} \int_\Omega |n^{\varepsilon,\tau} \ln n^{\varepsilon,\tau}|^\frac32 &=& \int_Q |n^{\varepsilon,\tau} \ln n^{\varepsilon,\tau}|^\frac32 \\&=& \int_{Q \cap \{n^{\varepsilon,\tau} > e\}} |n^{\varepsilon,\tau} \ln n^{\varepsilon,\tau}|^\frac32 + \int_{Q \cap \{n^{\varepsilon,\tau} \leq e\}} |n^{\varepsilon,\tau} \ln n^{\varepsilon,\tau}|^\frac32.
\eeno
Since
\beno
\lim_{n^{\varepsilon,\tau} \rightarrow \infty} {(n^{\varepsilon,\tau})}^{-\frac16}|\ln n^{\varepsilon,\tau}|^{\frac32} = 0,
\eeno
we have
\beno
\int_{Q \cap \{n^{\varepsilon,\tau} > e\}} |n^{\varepsilon,\tau} \ln n^{\varepsilon,\tau}|^\frac32 \leq C \int_Q |n^{\varepsilon,\tau}|^\frac53.
\eeno
Noting that $n^{\varepsilon,\tau} \in L^\infty (0,t;L^1(\Omega))$ and $\nabla \sqrt{n^{\varepsilon,\tau}} \in L^2 (0,t;L^2(\Omega))$, we have
%\beno
%\|\sqrt{n^{\varepsilon,\tau}}\|_{L^\frac{10}3(0,t; L^\frac{10}3(\Omega))} \leq C.
%%&\leq& C \|\sqrt{n^{\varepsilon,\tau}}\|_{L^\infty(0,t;L^2(\Omega))}^\frac25 \|\nabla \sqrt{n^{\varepsilon,\tau}}\|_{L^2(0,t;L^2(\Omega))}^\frac35 \\
%%&&+~ C \|\sqrt{n^{\varepsilon,\tau}}\|_{L^\infty(0,t;L^2(\Omega))}
%\eeno
%Due to
%\beno
%&&\|\sqrt{n^\varepsilon+1}\|_{L^\infty(0,t; L^2(\Omega))}^2 = \int_{0}^{t} \int_\Omega (n^\varepsilon+1) dxdt= \int_{0}^{t} \int_\Omega n^\varepsilon dxdt + |(0,t) \times \Omega| \\
%&&\leq C(\|n_0\|_{L^1(0,t;L^1(\mathbb{R}^3))},|\Omega|,t).
%\eeno

\beno
\|\sqrt{n^{\varepsilon,\tau}}\|_{L^\frac{10}3(0,t;L^\frac{10}3(\Omega))} \leq C,
\eeno
which means
\beno
\|n^{\varepsilon,\tau}\|_{L^\frac53(0,t; L^\frac53(\Omega))} \leq C.
\eeno
%since $\sqrt n^\varepsilon \leq \sqrt{n^\varepsilon+1}$.
Noting that
\beno
\lim_{n^{\varepsilon,\tau} \rightarrow 0} (n^{\varepsilon,\tau})^\frac16 \ln n^{\varepsilon,\tau} = 0,
\eeno
we have
\beno
\int_{Q \cap \{n^{\varepsilon,\tau} \leq e\}} |n^{\varepsilon,\tau} \ln n^{\varepsilon,\tau}|^\frac32 \leq C \int_Q |n^{\varepsilon,\tau}|^\frac43 \leq C |Q|^\frac15 \left(\int_Q |n^{\varepsilon,\tau}|^\frac53\right)^\frac45 \leq C.
\eeno
It follows that
\beno
n^{\varepsilon,\tau} \ln n^{\varepsilon,\tau} \in L^{\frac32}(0,t;L^{\frac32}(\Omega)).
\eeno
\endproof

\begin{lemma}\label{n lnn convergence}
	$n^{\varepsilon,\tau} \ln n^{\varepsilon,\tau} \rightarrow n^{\varepsilon} \ln n^{\varepsilon}$ strongly in $L^q(0,t;L^q(\Omega))$ for $q \in [1,\frac32)$ as $ \tau\to 0 $.
\end{lemma}

{\bf Proof.} Firstly, we divide the following inequality into four parts as follows. We denote $Q=(0,t)\times\Omega$ for simplicity.
\beno
&&\int_{Q}|n^{\varepsilon,\tau} \ln n^{\varepsilon,\tau}-n^{\varepsilon} \ln n^{\varepsilon}|dxdt\\
&\leq& \int_{\{n^{\varepsilon,\tau}<1\}\bigcap\{n^{\varepsilon,\tau}>n^{\varepsilon}\}\bigcap Q}|n^{\varepsilon,\tau} \ln n^{\varepsilon,\tau}-n^{\varepsilon} \ln n^{\varepsilon}|\\
&&+\int_{\{n^{\varepsilon,\tau}<1\}\bigcap\{n^{\varepsilon,\tau}<n^{\varepsilon}\}\bigcap Q}|n^{\varepsilon,\tau} \ln n^{\varepsilon,\tau}-n^{\varepsilon} \ln n^{\varepsilon}|\\
&&+\int_{\{n^{\varepsilon,\tau}>1\}\bigcap Q}|(n^{\varepsilon,\tau}-n^{\varepsilon}) \ln n^{\varepsilon,\tau}|
+\int_{\{n^{\varepsilon,\tau}>1\}\bigcap Q}|n^{\varepsilon}(\ln n^{\varepsilon,\tau}-\ln n^{\varepsilon})|\\
&:=&N_1+N_2+N_3+N_4.
\eeno
For the term of $N_1$, we arrive
\beno
N_1&\leq& C\int_{\{n^{\varepsilon}<n^{\varepsilon,\tau}<1\}\bigcap Q}|(n^{\varepsilon,\tau}-n^{\varepsilon}) \ln n^{\varepsilon,\tau}|
+C \int_{\{n^{\varepsilon}<n^{\varepsilon,\tau}<1\}\bigcap Q }|n^{\varepsilon} \ln(\frac{ n^{\varepsilon,\tau}}{n^{\varepsilon}}) |\\
&:=&N_{11}+N_{12}.
\eeno
For $N_{11}$, note that 
\ben\label{n_0}
\lim_{n \rightarrow 0} n^\frac{1}{9} \ln n = 0,
\een
then by H\"{o}lder inequality, (\ref{n_0}) and  Lemma \ref{lem:n}, we get
\beno
N_{11}&&\leq C\int_{\{n^{\varepsilon}<n^{\varepsilon,\tau}<1\}\bigcap Q}|(n^{\varepsilon,\tau}-n^{\varepsilon}) (n^{\varepsilon,\tau})^{-\frac{1}{9}}|\\\nonumber&&\leq C\int_{\{n^{\varepsilon}<n^{\varepsilon,\tau}<1\}\bigcap Q}|(n^{\varepsilon,\tau}-n^{\varepsilon})(n^{\varepsilon,\tau}-n^{\varepsilon})^{-\frac19}|\\\nonumber
&&\leq C\int_{\{n^{\varepsilon}<n^{\varepsilon,\tau}<1\}\bigcap Q}|(n^{\varepsilon,\tau}-n^{\varepsilon})^{\frac89}|\\\nonumber
%&&\leq C\int_{\{n^{\varepsilon}<n^{\varepsilon,\tau}<1\}\bigcap\{Q\}}|n^{\varepsilon,\tau}-n^{\varepsilon}|^{\frac43}\\
&&\leq C |Q|^{\frac19}\|n^{\varepsilon,\tau}-n^{\varepsilon}\|_{L^1(Q)}\rightarrow 0 \quad {\rm as} \quad\tau\rightarrow0.
\eeno
As for $N_{12}$, by inequality $\ln(1+f)\leq f$ when $f>-1$ and  Lemma \ref{lem:n}, we have
\beno
N_{12}&&\leq C\int_{\{n^{\varepsilon}<n^{\varepsilon,\tau}<1\}\bigcap Q}|n^{\varepsilon,\tau}-n^{\varepsilon}|\rightarrow 0 \quad {\rm as} \quad\tau\rightarrow0.
%&&\leq C|Q|^{\frac14}\|n^{\varepsilon,\tau}-n^{\varepsilon}\|_{L^2(Q)}^{\frac32}
\eeno
For the term of $ N_{2}, $ we have
\beno
N_2&\leq& C\int_{\{n^{\varepsilon,\tau}<1\}\bigcap\{n^{\varepsilon,\tau}<n^{\varepsilon}\}\bigcap Q}|(n^{\varepsilon}-n^{\varepsilon,\tau}) \ln n^{\varepsilon}|
\\&&+C \int_{\{n^{\varepsilon,\tau}<1\}\bigcap\{n^{\varepsilon,\tau}<n^{\varepsilon}\}\bigcap Q}|n^{\varepsilon,\tau} \ln(\frac{n^{\varepsilon}}{ n^{\varepsilon,\tau}}) |\\
&\leq& C\int_{\{n^{\varepsilon,\tau}<1\}\bigcap\{n^{\varepsilon,\tau}<n^{\varepsilon}\}\bigcap Q\bigcap \{n^{\varepsilon}<1\}}|(n^{\varepsilon}-n^{\varepsilon,\tau}) \ln n^{\varepsilon}|\\
&&+C\int_{\{n^{\varepsilon,\tau}<1\}\bigcap\{n^{\varepsilon,\tau}<n^{\varepsilon}\}\bigcap Q\bigcap \{n^{\varepsilon}\geq1\}}|(n^{\varepsilon}-n^{\varepsilon,\tau}) \ln n^{\varepsilon}|\\
&&+C \int_{\{n^{\varepsilon,\tau}<1\}\bigcap\{n^{\varepsilon,\tau}<n^{\varepsilon}\}\bigcap Q}|n^{\varepsilon,\tau} \ln(\frac{n^{\varepsilon}}{ n^{\varepsilon,\tau}}) |\\
&:=&N_{21}+N_{22}+N_{23}.
\eeno
For $N_{21}$, by (\ref{n_0}) and Lemma \ref{lem:n}, we obtain
\beno
N_{21}&&\leq C\int_{\{n^{\varepsilon,\tau}<n^{\varepsilon}\}\bigcap Q}|(n^{\varepsilon}-n^{\varepsilon,\tau})^{\frac89} |\\
&&\leq C |Q|^{\frac19}\|n^{\varepsilon,\tau}-n^{\varepsilon}\|_{L^1(Q)}\rightarrow 0 \quad {\rm as} \quad\tau\rightarrow0.
%&&\leq C|Q|^{\frac59}\|n^{\varepsilon,\tau}-n^{\varepsilon}\|_{L^3(Q)}^{\frac43}\rightarrow 0 \quad {\rm as} \quad\tau\rightarrow0,
\eeno
Note that 
\ben\label{n_1}
\lim_{n \rightarrow \infty} n^{-\frac16}|\ln n|^{\frac32} = 0,
\een
 there holds
\beno
N_{22}&&\leq C\int_{\{n^{\varepsilon,\tau}<1\}\bigcap Q}|n^{\varepsilon}-n^{\varepsilon,\tau}||n^{\varepsilon}|^{\frac19} \\
&&\leq  C\|n^{\varepsilon,\tau}-n^{\varepsilon}\|_{L^{\frac98}(Q)}\|n^{\varepsilon}\|_{L^{1}(Q)}^{\frac19}\rightarrow 0 \quad {\rm as} \quad\tau\rightarrow0.
\eeno
Similarly, for the term of $N_{23}$, there holds
\beno
N_{23}\leq \int_{\{n^{\varepsilon,\tau}<1\}\bigcap\{n^{\varepsilon,\tau}<n^{\varepsilon}\}\bigcap Q}|n^{\varepsilon}-n^{\varepsilon,\tau}|\rightarrow 0 \quad {\rm as} \quad\tau\rightarrow0.
\eeno
%\beno
%N_{23}\leq C |Q|^{\frac14} \left(\int_{\{n^{\varepsilon,\tau}<1\}\bigcap\{n^{\varepsilon,\tau}<n^{\varepsilon}\}\bigcap\{Q\}}|n^{\varepsilon}-n^{\varepsilon,\tau}|^2\right)^{\frac34}\rightarrow 0 \quad {\rm as} \quad\tau\rightarrow0.
%\eeno
For the term of $N_3$, 
 by H\"{o}lder inequality, (\ref{ine:  energy}) , (\ref{n_1}) and  Lemma \ref{lem:n}, we arrive
\beno
N_3&&\leq \int_{\{n^{\varepsilon,\tau}>1\}\bigcap Q}|(n^{\varepsilon,\tau}-n^{\varepsilon}) \ln n^{\varepsilon,\tau}|\\
&&\leq \int_{\{n^{\varepsilon,\tau}>1\}\bigcap Q}|n^{\varepsilon,\tau}-n^{\varepsilon}| |n^{\varepsilon,\tau}|^{\frac19}\\
&&\leq C \|n^{\varepsilon,\tau}-n^{\varepsilon}\|_{L^{\frac98}(Q)}\|n^{\varepsilon,\tau}\|_{L^{1}(Q)}^{\frac19}
%&&\leq  C \|n^{\varepsilon,\tau}-n^{\varepsilon}\|_{L^{\frac{30}{11}}(Q)}^{\frac32}\|n^{\varepsilon,\tau}\|_{L^{\frac{10}3}(Q)}^{\frac32}
\rightarrow 0 \quad {\rm as} \quad\tau\rightarrow0.
\eeno
Moreover, we get
\beno
N_4&&=\int_{\{n^{\varepsilon,\tau}>1\}\bigcap Q}|n^{\varepsilon}\ln(\frac{n^{\varepsilon,\tau}-n^{\varepsilon}}{n^{\varepsilon}}+1)\\&&\leq \int_{\{n^{\varepsilon,\tau}>1\}\bigcap Q}|n^{\varepsilon,\tau}-n^{\varepsilon}|\rightarrow 0 \quad {\rm as} \quad\tau\rightarrow0.
%\leq C |Q|^{\frac14}\|n^{\varepsilon,\tau}-n^{\varepsilon}\|_{{L^2}(Q)}^{\frac32}
\eeno
Collecting $N_1-N_4$, which follows that $n^{\varepsilon,\tau} \ln n^{\varepsilon,\tau} \rightarrow n^{\varepsilon} \ln n^{\varepsilon}~~ {\rm in} ~~ L^{1}(0,t;L^{1}(\Omega))$ as $ \tau\to 0. $

Furthermore, we can obtain 
\beno
n^{\varepsilon,\tau}\ln n^{\varepsilon,\tau} \rightarrow  n^{\varepsilon}\ln n^{\varepsilon} \quad {\rm in} \quad L^q(0,t;L^q(\Omega)) \quad {\rm for} \quad q\in [1,\frac{3}2)
\eeno
as $ \tau\to 0. $

\section{The convergence of local energy inequality}

We will show that $(u^{\varepsilon,\tau}, c^{\varepsilon,\tau}, n^{\varepsilon,\tau})$ has a limit $(u^{\varepsilon},c^{\varepsilon},n^{\varepsilon})$ satisfying local energy inequality
(\ref{local energy}).

For the left part of (\ref{local energy}), using properties of weak convergence, we see that
%Fatou Lemma
\begin{equation*}
\begin{aligned}
\int_{(0,t)\times\Omega}|\nabla\sqrt{n^{\varepsilon}}|^{2}\psi\leq \liminf_{\tau \to 0} \left(\int_{(0,t)\times\Omega} |\nabla\sqrt{n^{\varepsilon,\tau}}|^{2} \psi \right);
\end{aligned}
\end{equation*}

\beno
 \int_{\Omega} (|\nabla \sqrt{{c^{\varepsilon}}}|^2 \psi)(\cdot,t) \leq    \liminf_{\tau\to0}   \int_{\Omega} (|\nabla \sqrt{{c}^{\varepsilon,\tau}}|^2 \psi)(\cdot,t) ;
\eeno

\beno
\int_{(0,t)\times\Omega} |\Delta \sqrt{{c^{\varepsilon}}}|^2 \psi \leq  \liminf_{\tau\to0} \int_{(0,t)\times\Omega} |\Delta \sqrt{{c}^{\varepsilon,\tau}}|^2 \psi;
\eeno
%\beno
%\int_{(0,t)\times\Omega} (\sqrt{\tilde{c}})^{-2} (\partial_j \sqrt{\tilde{c}})^2 (\partial_i \sqrt{\tilde{c}})^2 \psi \leq \liminf_{\varepsilon\to0} \int_{(0,t)\times\Omega} (\sqrt{\tilde{c}^\varepsilon})^{-2} (\partial_j \sqrt{\tilde{c}^\varepsilon})^2 (\partial_i \sqrt{\tilde{c}^\varepsilon})^2 \psi,
%\eeno
\beno
\int_{\Omega}|u^{\varepsilon}|^2 \psi \leq \liminf_{\tau\to0}  \int_{\Omega}|u^{\varepsilon,\tau}|^2 \psi;
\eeno

%For the term $\int_{(0,t)\times\Omega} |\nabla \sqrt{c^\varepsilon}|^2 (n^\varepsilon\ast\rho^{\varepsilon}) \psi$, since
%\beno
%\lim_{\varepsilon\rightarrow 0} |\nabla \sqrt{c^\varepsilon\ast\rho^{\varepsilon}}|^2 n^\varepsilon = |\nabla \sqrt{c}|^2 n \quad {\rm in} \quad L^q((0,t); L^q(\Omega)) \quad {\rm for ~~~ all~~} q\in(1,+\infty),
%\eeno
%by Fatou's Lemma, we have
%\beno
%\int_0^t \int_{\Omega} |\nabla \sqrt{c}|^2 n \psi \leq \liminf_{\varepsilon\to0} \int_0^t \int_{\Omega} |\nabla \sqrt{c^\varepsilon}|^2 (n^\varepsilon\ast\rho^{\varepsilon}) \psi.
%\eeno
\begin{equation*}
	\begin{aligned}
	\int_{(0,t)\times\Omega} |\nabla u^{\varepsilon}|^2 \psi \leq \liminf_{\tau\to0} \int_{(0,t)\times\Omega} |\nabla u^{\varepsilon,\tau}|^2 \psi.
	\end{aligned}
\end{equation*}
and
\begin{equation*}
	\begin{aligned}
\int_{(0,t)\times\Omega} (\sqrt{c^{\varepsilon}})^{-2} |\nabla \sqrt{c^{\varepsilon}}|^4 \psi\leq \liminf_{\tau\to0}  \int_{(0,t)\times\Omega} (\sqrt{c^{\varepsilon,\tau}})^{-2} |\nabla \sqrt{c^{\varepsilon,\tau}}|^4 \psi. 
 \end{aligned}
\end{equation*}
 Noting that $n^{\varepsilon,\tau} \ln n^{\varepsilon,\tau} > -e^{-1}$ and the domain $\Omega$ is bounded, by Fatou's Lemma we have
\begin{equation*}
\begin{aligned}
\int_{\Omega}(n^{\varepsilon}\ln n^{\varepsilon}\psi+2e^{-1}\psi)&\leq \liminf_{\tau \to 0}\int_{\Omega}(n^{\varepsilon,\tau}\ln n^{\varepsilon,\tau}\psi+2e^{-1}\psi)dx,
%\\&=\liminf_{\tau\to 0}\int_{\Omega}n^{\varepsilon,\tau}\ln n^{\varepsilon,\tau}\psi+\int_{\Omega}2e^{-1}\psi,
\end{aligned}
\end{equation*}
that is
\beno
\int_\Omega n^{\varepsilon} \ln n^{\varepsilon}\psi  \leq \liminf_{\tau \to 0} \int_\Omega n^{\varepsilon,\tau} \ln n^{\varepsilon,\tau}\psi,
\eeno
where we used $\|n^{\varepsilon,\tau} \ln n^{\varepsilon,\tau}\|_{L^1}\leq C(t)$, which means that uniformly bounded sequences have spot-convergent subcolumns.
%\beno
%\int_\Omega (n \ln n + 2e^{-1}) \leq \liminf_{\varepsilon \to 0} \int_\Omega (n^\varepsilon \ln n^\varepsilon + 2e^{-1})=\liminf_{\varepsilon \to 0} \int_\Omega n^\varepsilon \ln n^\varepsilon + \int_\Omega 2e^{-1},
%\eeno
%Thus
%\beno
%\int_\Omega n \ln n  \leq \liminf_{\varepsilon \to 0} \int_\Omega n^\varepsilon \ln n^\varepsilon.
%\eeno
%where we used $\|n^{\varepsilon,\tau} \ln n^{\varepsilon,\tau}\|_{L^1}\leq C(t)$, which means that uniformly bounded sequences have spot-convergent subcolumns.

For the right part of (\ref{local energy}), by a priori estimate, Lemma \ref{lem:u}-Lemma 5.13, we have
\beno
&&u^{\varepsilon,\tau}\rightarrow u^{\varepsilon} \quad {\rm in}~~ L^p((0,t);L^p(\Omega)),~~ p\in[2,\frac{10}3)~~{\rm as}~~\tau\rightarrow 0;\\
&&\nabla{c^{\varepsilon,\tau}}\rightarrow \nabla{c^{\varepsilon}} \quad {\rm in} ~~ L^p((0,t);L^p(\Omega)),~~ p\in[2,\frac{10}3)~~{\rm as}~~\tau\rightarrow 0;\\
&&\nabla\sqrt{c^{\varepsilon,\tau}}\rightarrow \nabla\sqrt{c^{\varepsilon}} \quad {\rm in} ~~ L^p((0,t);L^p(\Omega)),~~ p\in[2,\frac{10}3)~~{\rm as}~~\tau\rightarrow 0;\\
%&&n^{\varepsilon,\tau}\rightharpoonup n^{\varepsilon} \quad {\rm in}~~ L^{\frac53}((0,t);L^{\frac53}(\Omega)),~~{\rm as}~~\tau\rightarrow 0;\\
%&&(u^\varepsilon\ast\rho^\varepsilon)u\rightarrow u^2 \quad  {\rm in}~~ L^{\frac32}((0,t);L^{\frac32}(\Omega)),~~{\rm as}~~\varepsilon\rightarrow\infty;\\
%&&n^{\varepsilon,\tau} \ln n^{\varepsilon,\tau} \rightharpoonup n^{\varepsilon} \ln n^{\varepsilon} \quad  {\rm in}~~ L^{1}((0,t);L^1(\Omega)),~~{\rm as}~~\tau\rightarrow 0.\\
&& n^{\varepsilon,\tau} \rightarrow  n^{\varepsilon} \quad {\rm in} ~~ L^p((0,t);L^p(\Omega)),~~ p\in[2,\frac{10}3)~~{\rm as}~~\tau\rightarrow 0;\\
&& n^{\varepsilon,\tau} \ln n^{\varepsilon,\tau} \rightarrow n^{\varepsilon} \ln n^{\varepsilon} \quad {\rm in} ~~ L^p((0,t);L^p(\Omega)),~~ p\in[1,\frac{3}2)~~{\rm as}~~\tau\rightarrow 0.
%&& P^\varepsilon - \bar{P^\varepsilon} \rightharpoonup p  \quad {\rm in}~~ L^{\frac3}((0,t);L^{\frac3}(\Omega)),~~{\rm as}~~\varepsilon\rightarrow\infty.\\
\eeno
%\beno
%\rightharpoonup~{\rm in}~~ L^{\frac32}((0,T)\times \Omega),~~{\rm as}~~\varepsilon\rightarrow\infty;
%\een
We also need the weak convergence of $P^{\varepsilon,\tau}$. By the equation of $u^{\varepsilon,\tau}$, we have
\beno
-\Delta (P^{\varepsilon,\tau}-\bar {P}^{\varepsilon,\tau}) = \partial_i\partial_j((u^{\varepsilon,\tau}_i\ast\rho^\varepsilon)u^{\varepsilon,\tau}_j)+\nabla\cdot ((n^{\varepsilon,\tau}\nabla\phi) \ast \rho^\varepsilon).
\eeno
By Calder\'{o}n-Zygmund theorem, we have
\ben\label{ine:P varepsilon}
\int_{\mathbb{R}^3} |P^{\varepsilon,\tau}-\bar {P}^{\varepsilon,\tau}|^\frac53 \leq C \int_{\mathbb{R}^3} |(u^{\varepsilon,\tau} \ast \rho^\varepsilon) \otimes u^{\varepsilon,\tau}|^\frac53 + C \left(\int_{\mathbb{R}^3} |(n^{\varepsilon,\tau}\nabla\phi) \ast \rho^\varepsilon|^\frac{15}{14} \right)^\frac{14}{15}.
\een
Using H\"{o}lder's inequality and Young's inequality, and noting that the strong convergence of convolution, we have
\beno
\int_{\mathbb{R}^3} |(u^{\varepsilon,\tau} \ast \rho^\varepsilon) \otimes u^{\varepsilon,\tau}|^\frac53 \leq \int_{\mathbb{R}^3} |u^{\varepsilon,\tau}|^\frac{10}3,
\eeno
and
\beno
\left(\int_{\mathbb{R}^3} |(n^{\varepsilon,\tau}\nabla\phi) \ast \rho^\varepsilon|^\frac{15}{14} \right)^\frac{14}{15} \leq \|\nabla \phi\|_{L^\infty}  \left(\int_{\mathbb{R}^3} |n^{\varepsilon,\tau}|^\frac{15}{14} \right)^\frac{14}{15}.
\eeno
Integrating both sides of (\ref{ine:P varepsilon}) with respect of $t$ from $0$ to $t$, we have
\ben\label{ine:P 53}
\int_0^t \int_{\mathbb{R}^3} |P^{\varepsilon,\tau}-\bar {P}^{\varepsilon,\tau}|^\frac53 &\leq& C \int_0^t \int_{\mathbb{R}^3} |u^{\varepsilon,\tau}|^\frac{10}3 + C \|\nabla \phi\|_{L^\infty} \int_0^t \left(\int_{\mathbb{R}^3} |n^{\varepsilon,\tau}|^\frac{15}{14} \right)^\frac{14}{15}\\\nonumber&:=&I+II.
\een
For $ I, $ by interpolation inequality, we have
\begin{equation*}
I\leq C\left(\|u^{\varepsilon,\tau}\|_{L^{\infty}(0,t;L^{2}(\mathbb{R}^{3}))}+\|\nabla u^{\varepsilon,\tau}\|_{L^{2}(0,t;L^{2}(\mathbb{R}^{3}))}\right)^{\frac{10}{3}}.
\end{equation*}
For $ II, $ noting that $n^{\varepsilon,\tau} \geq 0$, by H\"{o}lder's and interpolation inequalities, we have
\begin{equation}\label{ine:n1}
	\begin{aligned}
	&\int_0^t\left(\int_{\mathbb{R}^3} |n^{\varepsilon,\tau}|^\frac{15}{14} dx \right)^\frac{14}{15} dt\\ =~&\int_0^t \left(\int_{\mathbb{R}^3} \left(\sqrt{n^{\varepsilon,\tau}}\right)^\frac{15}{7} dx \right)^\frac{14}{15} dt \\
	\leq~& t^\frac9{10} \left(\int_0^t \left(\int_{\mathbb{R}^3} \left(\sqrt{n^{\varepsilon,\tau}}\right)^\frac{15}{7} dx\right)^\frac{28}{3} dt\right)^\frac1{10}\\
	\leq~&C t^{\frac9{10}}\left( \|\sqrt{n^{\varepsilon,\tau}} \|_{L^\infty(0,t;L^2(\mathbb{R}^3))}^2 + \|\nabla \sqrt{n^{\varepsilon,\tau}}\|_{L^2(0,t;L^2(\mathbb{R}^3))}^2\right).
	\end{aligned}
\end{equation}
Then from (\ref{ine:  energy}),
we have
\beno
\int_0^t \int_{\mathbb{R}^3} |P^{\varepsilon,\tau}-\bar {P}^{\varepsilon,\tau}|^\frac53 \leq C(t,\|n_0^{\varepsilon}\|_{L^1},\|\nabla \phi\|_{L^\infty}),
\eeno
which means
\beno
P^{\varepsilon,\tau}-\bar {P}^{\varepsilon,\tau} \rightharpoonup P^{\varepsilon}-\bar {P}^{\varepsilon} ~~{\rm in}~~ L^{\frac32}(0,t;L^{\frac32}(\Omega)),~~{\rm as}~~\tau\rightarrow 0.
\eeno
%\Beno
%P^\Varepsilon\Rightharpoonup P ~~{\Rm In}~~ L^{\Frac32},~~{\Rm As}~~\Varepsilon\Rightarrow\Infty;
%\Eeno

Next, we prove the convegence of the right terms of the local energy inequality (\ref{local energy}). Let us call the right terms of the local energy inequality (\ref{local energy}) as $M_1$, $M_2$, $\cdots$, $M_{10}$ term by term.
%\begin{equation}\label{local energy 1}
%\begin{aligned}\triangle\Delta
%&\int_{(0,t)\times\Omega} n^\varepsilon \ln n^\varepsilon (\partial_t \psi + \Delta \psi) + \int_{(0,t)\times\Omega} n^\varepsilon \ln n^\varepsilon (u^\varepsilon\ast\rho^\varepsilon) \cdot \nabla \psi \\&+M_{0}\int_{(0,t)\times\Omega} n^\varepsilon \ln n^\varepsilon \nabla (c^\varepsilon\ast\rho^\varepsilon) \cdot \nabla \psi
%+M_{0}\int_{(0,t)\times\Omega}(n^\varepsilon) \nabla (c^\varepsilon\ast\rho^\varepsilon) \cdot \nabla \psi\\&+ 2\int_{(0,t)\times\Omega} |\nabla \sqrt{c^\varepsilon}|^2 (\partial_t \psi + \Delta \psi) + 2 \int_{(0,t)\times\Omega} |\nabla \sqrt{c^\varepsilon}|^2 {(u^\varepsilon\ast\rho^\varepsilon)} \cdot \nabla \psi\\&+\frac{20}3 \int_{(0,t)\times\Omega} (\sqrt{c^\varepsilon})^{-1} |\nabla \sqrt{c^\varepsilon}|^2 |\nabla \sqrt{c^\varepsilon} \cdot \nabla \psi|
%+108||c^{\varepsilon}||_{\infty}\int_{(0,t)\times\Omega} |u^\varepsilon|^2 \left(\partial_t \psi +\Delta \psi\right) \\&+ 108||c^{\varepsilon}||_{\infty}\int_{(0,t)\times\Omega} |u^{\varepsilon}|^2(u^\varepsilon\ast\rho^\varepsilon)\cdot\nabla\psi
%+216||c^{\varepsilon}||_{\infty}\int_{(0,t)\times\Omega} (P^{\varepsilon} - \bar{P^{\varepsilon}}) u^\varepsilon \cdot \nabla \psi\\& - 216||c^{\varepsilon}||_{\infty}\int_{(0,t)\times\Omega} (n^\varepsilon\nabla\phi)\ast\rho^\varepsilon \cdot u^{\varepsilon} \psi\\&:=M_1+\cdots+M_{11}.
%\end{aligned}
%\end{equation}

{\bf \underline{For $M_1$:}}
Since $n^{\varepsilon,\tau} \ln n^{\varepsilon,\tau} \rightarrow n^{\varepsilon} \ln n^{\varepsilon} ~~ {\rm in} ~~ L^{\frac{10}{7}}((0,t);L^{\frac{10}{7}}(\Omega))~~{\rm as}~~\tau\rightarrow 0$,
%By (\ref{ine:  energy}), we get $||n^{\varepsilon,\tau}\ln n^{\varepsilon,\tau}||_{L^1 (0,t;L^1(\Omega))} \leq  C(t)$ with the constant $C$ is independent of $\tau$, due to 
%uniformly bounded sequences consist of weakly convergent subcolumns, there exists a subsequence such that
%\beno
%n^{\varepsilon,\tau} \ln n^{\varepsilon,\tau} \rightharpoonup n^{\varepsilon} \ln n^{\varepsilon} \quad {\rm in} \quad L^1(0,t;L^1(\Omega))~~{\rm as}~~\tau\rightarrow0.
%\eeno
and $\psi \in C_c^\infty$, we have $\partial_t \psi + \Delta \psi \in C_c^\infty$. Therefore,
\beno
\int_{(0,t)\times\Omega} n^{\varepsilon,\tau} \ln n^{\varepsilon,\tau} (\partial_t \psi + \Delta \psi) \rightarrow \int_{(0,t)\times\Omega} n^{\varepsilon} \ln n^{\varepsilon} (\partial_t \psi + \Delta \psi)~~{\rm as}~~\tau\rightarrow0.
\eeno

%\beno
%&&\int_{(0,t)\times\Omega} n^\varepsilon \ln n^\varepsilon (\partial_t \psi + \Delta \psi)-\int_{(0,t)\times\Omega} n \ln n (\partial_t \psi + \Delta \psi)\\
%&\leq& C\int_{(0,t)\times\Omega} n^\varepsilon \ln n^\varepsilon-n \ln n\\
%&\rightarrow& 0
%\eeno
%as $\varepsilon\rightarrow\infty$.
%By
%\beno
%||n\ln n||_{L^1 L^1}\leq C
%\eeno
%where C is independent of $\varepsilon$.

{\bf \underline{For $M_2$:}}
\beno
&&\int_{(0,t)\times\Omega} n^{\varepsilon,\tau} \ln n^{\varepsilon,\tau} (u^{\varepsilon,\tau}\ast\rho^\varepsilon) \cdot \nabla \psi - \int_{(0,t)\times\Omega} n^{\varepsilon} \ln n^{\varepsilon} (u^{\varepsilon}\ast\rho^{\varepsilon}) \cdot \nabla \psi\\
&\leq& \int_{(0,t)\times \Omega}n^{\varepsilon,\tau}\ln n^{\varepsilon,\tau}[(u^{\varepsilon,\tau}-u^{\varepsilon})\ast\rho^{\varepsilon}]\cdot\nabla\psi\\&&+\int_{(0,t)\times \Omega}(n^{\varepsilon,\tau}\ln n^{\varepsilon,\tau}-n^{\varepsilon}\ln n^{\varepsilon})(u^{\varepsilon}\ast\rho^{\varepsilon})\cdot\nabla\psi\\
&\leq& \|n^{\varepsilon,\tau}\ln n^{\varepsilon,\tau}\|_{L^{\frac32}(0,t;L^{\frac32}(\Omega))}\|u^{\varepsilon,\tau}- u^{\varepsilon}\|_{L^3(0,t;L^3(\Omega_1))} \|\nabla \psi\|_{L^\infty(0,t;L^\infty(\Omega))}\\
&&+\|n^{\varepsilon,\tau}\ln n^{\varepsilon,\tau}-n^{\varepsilon}\ln n^{\varepsilon}\|_{L^{\frac{10}7}(0,t;L^{\frac{10}7}(\Omega))}\|u^{\varepsilon}\ast\rho^{\varepsilon}\|_{L^{\frac{10}3}(0,t;L^{\frac{10}3}(\Omega_1))}\|\nabla \psi\|_{L^\infty(0,t;L^\infty(\Omega))}\\
&:=&M_{21}+M_{22},
\eeno
where $\Omega_1=\bigcup_{x\in\Omega}B_1(x)$.
Since $u^{\varepsilon,\tau} \rightarrow u^{\varepsilon}$ in $L^3(0,t;L^3(\Omega_1))$, we have $M_{21}\to 0$ ~~{\rm as}~~$\tau\rightarrow0$. 
%Noting that
%\beno
%\|(u^{\varepsilon}\ast\rho^{\varepsilon}) \cdot \nabla \psi\|_{L^{3}(0,t;L^{3}(\Omega))} \leq  \|\nabla \psi\|_{L^{\infty}(0,t;L^{\infty}(\Omega))} \|u^{\varepsilon}\|_{L^{3}(0,t;L^{3}(\Omega_{1}))}\leq C,
%\eeno
By strong convergence of $n^{\varepsilon,\tau} \ln n^{\varepsilon,\tau}$, we have $M_{22}\to 0 ~~{\rm as}~~\tau\rightarrow0$.
%Using the strong convergence of convolution, we have $M_{23}\to 0 ~~{\rm as}~~\tau\rightarrow0$.

%Then
%\beno
%M_2 \to \int_{(0,t)\times\Omega} n^{\varepsilon} \ln n^{\varepsilon} (u^{\varepsilon}\ast\rho^{\varepsilon}) \cdot \nabla \psi~~{\rm as}~~\tau\rightarrow0.
%\eeno

%%By the same method of \cite{CLW}, we have
%\beno
%\|n^\varepsilon \ln n^\varepsilon\|_{L^\frac32(0,T;L^\frac32(\Omega))} \leq C,
%\eeno
%and
%\beno
%\|u \cdot \nabla \phi\|_{L^3(0,T;L^3(\Omega))} \leq C.
%\eeno
%
%\beno
%\|u \cdot \nabla \phi\|_{L^1(0,T;L^1(\Omega))} \leq C.
%\eeno
%%\beno
%%&&\int_{(0,t)\times\Omega} n^\varepsilon \ln n^\varepsilon (u^\varepsilon\ast\rho^\varepsilon) \cdot \nabla \psi-\int_{(0,t)\times\Omega} n \ln n u \cdot \nabla \psi\\
%%&=&\int_{(0,t)\times\Omega} n^\varepsilon \ln n^\varepsilon (u^\varepsilon\ast\rho^\varepsilon) \cdot \nabla \psi-\int_{(0,t)\times\Omega} n^\varepsilon \ln n^\varepsilon u \cdot \nabla \psi+\int_{(0,t)\times\Omega} n^\varepsilon \ln n^\varepsilon u \cdot \nabla \psi\\
%%&-&\int_{(0,t)\times\Omega} n \ln n u \cdot \nabla \psi\\
%%&\leq&\| n^\varepsilon \ln n^\varepsilon\|_{L^{\frac32}}\|u^\varepsilon\ast\rho^\varepsilon-u\|_{L^3}||\nabla \psi||_{L^\infty}+\| n^\varepsilon \ln n^\varepsilon-n \ln n\|_{L^{\frac32}}||u||_{L^3}||\nabla \psi||_{L^\infty}\\
%%&\rightarrow& 0.
%%\eeno
%where we used the strong convergence of $u$ and the weak convergence of $n\ln n$ respectively.

{\bf \underline{For $M_3$:}}
For the term 
\beno
M_3 = \int_{(0,t)\times\Omega}\left(\frac{n^{\varepsilon,\tau}}{1+\tau n^{\varepsilon,\tau}} \nabla c^{\varepsilon,\tau} \chi(c^{\varepsilon,\tau})\right)\ln n^{\varepsilon,\tau}\cdot\nabla\psi,
\eeno
we rewrite $M_3$ as follow:
\beno
&&M_3-\int_{(0,t)\times\Omega}n^{\varepsilon}\ln n^{\varepsilon}\nabla c^{\varepsilon} \chi(c^{\varepsilon})\cdot\nabla\psi\\ &=& \int_{(0,t)\times\Omega}\left(\frac{1}{1+\tau n^{\varepsilon,\tau}}-1\right)n^{\varepsilon,\tau}\ln n^{\varepsilon,\tau}\nabla c^{\varepsilon,\tau}\chi(c^{\varepsilon,\tau})\cdot \nabla\psi\\
&&+~\int_{(0,t)\times\Omega}\left(n^{\varepsilon,\tau}\ln n^{\varepsilon,\tau}
-n^{\varepsilon}\ln n^{\varepsilon}\right)\nabla c^{\varepsilon,\tau}\chi(c^{\varepsilon,\tau})\cdot \nabla\psi\\&&+~\int_{(0,t)\times\Omega}n^{\varepsilon}\ln n^{\varepsilon}(\nabla c^{\varepsilon,\tau}\chi(c^{\varepsilon,\tau})-\nabla c^{\varepsilon}\chi(c^{\varepsilon}))\cdot\nabla\psi\\
&=:& M_{31}+M_{32}+M_{33}.
\eeno
For $ M_{31}, $ note that
\begin{equation*}
	\begin{aligned}
	&\int_{(0,t)\times \Omega}\left|n^{\varepsilon,\tau}\ln n^{\varepsilon,\tau}\nabla c^{\varepsilon,\tau}\chi(c^{\varepsilon,\tau})\cdot\nabla\psi\right|\\\leq~&\|n^{\varepsilon,\tau}\ln n^{\varepsilon,\tau}\|_{L^\frac32(0,t;L^\frac32(\Omega))}\|\nabla c^{\varepsilon,\tau}\|_{L^3(0,t;L^3(\Omega))}\|\chi\|_{L^{\infty}(0,\|c_{0}^{\varepsilon}\|_{L^{\infty}})}\|\nabla\psi\|_{L^\infty(0,t;L^\infty(\Omega))}\leq~C,
	\end{aligned}
\end{equation*}
and
\begin{equation*}
	\frac{1}{1+\tau n^{\varepsilon,\tau}}-1\to 0~~{\rm as}~~\tau\to 0,
\end{equation*}
which implies $ M_{31}\to 0 $ as $\tau\to 0.$
%For $M_{31}$, noting that $\frac{n^{\varepsilon,\tau}}{1+\tau n^{\varepsilon,\tau}}\rightarrow n^{\varepsilon,\tau}$ pointwise as $\tau\rightarrow 0$, we know that
%\begin{equation*}
%	\begin{aligned}
%	M_{31}\leq& \|\frac{n^{\varepsilon,\tau}}{1+\tau n^{\varepsilon,\tau}}\ln n^{\varepsilon,\tau}-n^{\varepsilon,\tau}\ln n^{\varepsilon,\tau}\|_{L^\frac32(0,t;L^\frac32(\Omega))} \|\nabla c^{\varepsilon,\tau}\|_{L^3(0,t;L^3(\Omega))}\\&\cdot \|\chi(c^{\varepsilon,\tau})\|_{L^\infty(0,t;L^\infty(\Omega))}\|\nabla \psi\|_{L^\infty(0,t;L^\infty(\Omega))}
%\rightarrow 0,
%	\end{aligned}
%\end{equation*}
%where we used the boundedness of $\|\nabla c^{\varepsilon,\tau}\|_{L^3(0,t;L^3(\Omega))}$, $\|\chi(c^{\varepsilon,\tau})\|_{L^\infty(0,t;L^\infty(\Omega))}$ and $\|\nabla \psi\|_{L^\infty(0,t;L^\infty(\Omega))}$.

For $M_{32}$, by the strong convergence of $n^{\varepsilon,\tau} \ln n^{\varepsilon,\tau}$, 
\beno
M_{32}&\leq &\|(n^{\varepsilon,\tau}\ln n^{\varepsilon,\tau}-n^{\varepsilon}\ln n^{\varepsilon})\|_{L^{\frac{10}7}(0,t;L^{\frac{10}7}(\Omega))}\|\nabla c^{\varepsilon,\tau}\|_{L^{\frac{10}3}(0,t;L^{\frac{10}3}(\Omega))}\\
&&\cdot\|\chi\|_{L^{\infty}(0,\|c_{0}^{\varepsilon}\|_{L^{\infty}})}\|\nabla \psi\|_{L^\infty(0,t;L^\infty(\Omega))}
\\&\rightarrow& 0~~{\rm as}~~\tau\to 0.
\eeno
%\beno
%n^{\varepsilon,\tau}\ln n^{\varepsilon,\tau} \rightharpoonup n^{\varepsilon} \ln n^{\varepsilon} \quad {\rm in} \quad L^\frac32(0,t;L^\frac32(\Omega))~~{\rm as}~~\tau\rightarrow0,
%\eeno
%we have
%\beno
%M_{32} \to 0 \quad {\rm as} \quad \tau \to 0.
%\eeno
For $M_{33}$, we estimate it as follow:
\begin{equation*}
	\begin{aligned}
	M_{33} 
	\leq~& \|n^{\varepsilon} \ln n^{\varepsilon}\|_{L^\frac32(0,t;L^\frac32(\Omega))} \|\nabla c^{\varepsilon,\tau} \chi(c^{\varepsilon,\tau})-\nabla c^{\varepsilon} \chi(c^{\varepsilon})\|_{L^3(0,t;L^3(\Omega))} \|\nabla \psi\|_{L^\infty(0,t;L^\infty(\Omega))}\\
	\leq~ &\|n^{\varepsilon} \ln n^{\varepsilon}\|_{L^\frac32(0,t;L^\frac32(\Omega))}\|(\nabla c^{\varepsilon,\tau}-\nabla c^{\varepsilon}) \chi(c^{\varepsilon,\tau})\|_{L^3(0,t;L^3(\Omega))} \|\nabla \psi\|_{L^\infty(0,t;L^\infty(\Omega))}\\
	&+~\|n^{\varepsilon} \ln n^{\varepsilon}\|_{L^\frac32(0,t;L^\frac32(\Omega))}\|\nabla c^{\varepsilon}(\chi(c^{\varepsilon,\tau})-\chi (c^{\varepsilon}))\|_{L^3(0,t;L^3(\Omega))} \|\nabla \psi\|_{L^\infty(0,t;L^\infty(\Omega))}\\
	=:~&M_{33}^{'}+M_{33}^{''}.
	\end{aligned}
\end{equation*}
For the first term of $M_{33}^{'}$, by the continuity of function $\chi$ , H\"{o}lder inequality and the strong convergence of $\nabla c^{\varepsilon,\tau}$, we arrive
\begin{equation*}
\begin{aligned}
M_{33}^{'}\leq~& \|n^{\varepsilon} \ln n^{\varepsilon}\|_{L^\frac32(0,t;L^\frac32(\Omega))}\|\chi(c^{\varepsilon,\tau})\|_{L^\infty(0,t;L^\infty(\Omega))} \|\nabla c^{\varepsilon,\tau}-\nabla c^{\varepsilon}\|_{L^3(0,t;L^3(\Omega))} \|\nabla \psi\|_{L^\infty(0,t;L^\infty(\Omega))}\\
\leq~& \|n^{\varepsilon} \ln n^{\varepsilon}\|_{L^\frac32(0,t;L^\frac32(\Omega))}\|\chi\|_{L^{\infty}(0,\|c^{\varepsilon}_{0}\|_{L^{\infty}(\Omega)})} \|\nabla c^{\varepsilon,\tau}-\nabla c^{\varepsilon}\|_{L^3(0,t;L^3(\Omega))} \|\nabla \psi\|_{L^\infty(0,t;L^\infty(\Omega))}\\
&\rightarrow 0~~{\rm as}~~\tau\to 0.
\end{aligned}
\end{equation*}
For $M_{33}^{''}$, by H\"{o}lder inequality and mean value theorem, for $\theta\in(0,1)$, we achieve 
\begin{equation}\label{M_{31}}
	\begin{aligned}
	M_{33}^{''}\leq~&\|n^{\varepsilon} \ln n^{\varepsilon}\|_{L^\frac32(0,t;L^\frac32(\Omega))}\|\nabla c^{\varepsilon}(\chi(c^{\varepsilon,\tau})-\chi (c^{\varepsilon}))\|_{L^3(0,t;L^3(\Omega))} \|\nabla \psi\|_{L^\infty(0,t;L^\infty(\Omega))}\\
%\leq~&\|n^{\varepsilon} \ln n^{\varepsilon}\|_{L^\frac32(0,t;L^\frac32(\Omega))}\|\nabla c\|_{L^{\frac{10}{3}}(0,t;L^{\frac{10}{3}}(\Omega))}\|\chi(c^{\varepsilon,\tau})-\chi (c^{\varepsilon})\|_{L^{30}(0,t;L^{30}(\Omega))}\|\nabla \psi\|_{L^\infty(0,t;L^\infty(\Omega))}\\
	\leq~&\|n^{\varepsilon} \ln n^{\varepsilon}\|_{L^\frac32(0,t;L^\frac32(\Omega))}\|\nabla c^{\varepsilon}\|_{L^{\frac{10}{3}}(0,t;L^{\frac{10}{3}}(\Omega))}\\&\cdot\|\chi^{'}(\theta c^{\varepsilon,\tau}+(1-\theta)c^{\varepsilon})(c^{\varepsilon,\tau}-c^{\varepsilon})\|_{L^{30}(0,t;L^{30}(\Omega))}\|\nabla \psi\|_{L^\infty(0,t;L^\infty(\Omega))}\\
	\leq~&\|n^{\varepsilon} \ln n^{\varepsilon}\|_{L^\frac32(0,t;L^\frac32(\Omega))}\|\nabla c^{\varepsilon}\|_{L^{\frac{10}{3}}(0,t;L^{\frac{10}{3}}(\Omega))}\|\chi^{'}(\theta c^{\varepsilon,\tau}+(1-\theta)c^{\varepsilon})\|_{L^\infty(0,t;L^\infty(\Omega))} \\&\cdot\|c^{\varepsilon,\tau}-c^{\varepsilon}\|_{L^{30}(0,t;L^{30}(\Omega))} \|\nabla \psi\|_{L^\infty(0,t;L^\infty(\Omega))}\\
	\leq~&\|n^{\varepsilon} \ln n^{\varepsilon}\|_{L^\frac32(0,t;L^\frac32(\Omega))}\|\nabla c^{\varepsilon}\|_{L^{\frac{10}{3}}(0,t;L^{\frac{10}{3}}(\Omega))}\|\chi\|_{L^{\infty}(0,2\|c^{\varepsilon}_{0}\|_{L^{\infty}(\Omega)})} \\&\cdot\|c^{\varepsilon,\tau}-c^{\varepsilon}\|_{L^{30}(0,t;L^{30}(\Omega))} \|\nabla \psi\|_{L^\infty(0,t;L^\infty(\Omega))}.
	\end{aligned}
\end{equation}
Using embedding inequality
\begin{equation}\label{c 30}
	\begin{aligned}
	\|c^{\varepsilon,\tau}-c^{\varepsilon}\|_{L^{30}(0,t;L^{30}(\Omega))}\nonumber=~&\max |c^{\varepsilon,\tau}-c^{\varepsilon}|^{\frac{4}{5}}\left(\int_{(0,t)\times\Omega}|c^{\varepsilon,\tau}-c^{\varepsilon}|^6\right)^{\frac1{30}}\\
	%\leq~& \max |c^{\varepsilon,\tau}-c^{\varepsilon}|^{\frac{4}{5}}\left(\int_0^t(\int_{\Omega}|\nabla c^{\varepsilon,\tau}-\nabla c^{\varepsilon}|^2dx)^{3}dt\right)^{\frac1{30}}\\
	\leq~& \|c^{\varepsilon,\tau}-c^{\varepsilon}\|_{L^\infty}^{\frac{4}{5}}\|\nabla(c^{\varepsilon,\tau}- c^{\varepsilon})\|_{L^{6}(0,t;L^{2}(\Omega))}^{\frac15}\\
	\leq~& \|2 c^{\varepsilon}_{0}\|_{L^{\infty}}^{\frac{4}{5}}\|\nabla(c^{\varepsilon,\tau}- c^{\varepsilon})\|_{L^{\infty}(0,t;L^{2}(\Omega))}^{\frac{2}{15}}\|\nabla(c^{\varepsilon,\tau}- c^{\varepsilon})\|_{L^{2}(0,t;L^{2}(\Omega))}^{\frac{1}{15}}.
	\end{aligned}
\end{equation}
Combining (\ref{M_{31}}) and (\ref{c 30}), by the strong convergence of $\nabla c^{\varepsilon,\tau}$,
we have $M_{33}^{''}\to 0 ~~{\rm as}~~\tau\rightarrow0$.
To sum up, $M_{33}\to 0 ~~{\rm as}~~\tau\rightarrow0$.
Collecting $M_{31}$ to $M_{33}$, we have
\beno
M_{3} \to \int_{(0,t)\times\Omega}n^{\varepsilon}\ln n^{\varepsilon}\nabla c^{\varepsilon} \chi(c^{\varepsilon})\cdot\nabla\psi \quad {\rm as} \quad \tau \to 0.
\eeno

{\bf \underline{For $M_4$:}} The analysis of $ M_{4} $ is similar to $ M_{3}, $ we omit it. Ultimately, one can obtain
\beno
M_{4} \to \int_{(0,t)\times \Omega} n^{\varepsilon}\nabla c^{\varepsilon}\chi(c^{\varepsilon})\cdot\nabla\psi \quad {\rm as} \quad \tau \to 0.
\eeno

{\bf \underline{For $M_5$:}}
Since $\psi \in C_c^\infty$, and
\beno
&&\nabla\sqrt{{c}^{\varepsilon,\tau}}\rightarrow \nabla\sqrt{{c^{\varepsilon}}} \quad {\rm in} ~~ L^p(0,t;L^p(\Omega)),~~ p\in[2,\frac{10}3)~~{\rm as}~~\tau\rightarrow0.
\eeno
Then
\beno
\int_{(0,t)\times\Omega} |\nabla \sqrt{{c}^{\varepsilon,\tau}}|^2 (\partial_t \psi + \Delta \psi)-\int_{(0,t)\times\Omega} |\nabla \sqrt{{c^{\varepsilon}}}|^2 (\partial_t \psi + \Delta \psi)\rightarrow 0~~{\rm as}~~\tau\rightarrow0.
\eeno

%\leq C||\nabla \sqrt{c^\varepsilon}-\nabla \sqrt{c}||^2_{L^2}
{\bf \underline{For $M_6$:}}
\beno
&&\int_{(0,t)\times\Omega} |\nabla \sqrt{{c}^{\varepsilon,\tau}}|^2 {(u^{\varepsilon,\tau}\ast\rho^\varepsilon)} \cdot \nabla \psi-\int_{(0,t)\times\Omega} |\nabla \sqrt{{c^{\varepsilon}}}|^2 ({u^{\varepsilon}}\ast\rho^{\varepsilon}) \cdot \nabla \psi\\
&\leq& \int_{(0,t)\times \Omega}|\nabla \sqrt{{c}^{\varepsilon,\tau}}|^2[(u^{\varepsilon,\tau}-u^{\varepsilon})\ast\rho^{\varepsilon}]\cdot\nabla\psi+\int_{(0,t)\times \Omega}(|\nabla \sqrt{{c}^{\varepsilon,\tau}}|^2-|\nabla \sqrt{{c^{\varepsilon}}}|^2 )(u^{\varepsilon}\ast\rho^{\varepsilon})\cdot\nabla\psi\\
&\leq& \|\nabla \sqrt{{c^{\varepsilon,\tau}}}\|^2_{L^{3}(0,t;L^{3}(\Omega))}\|u^{\varepsilon,\tau}- u^{\varepsilon}\|_{L^3(0,t;L^3(\Omega))} \|\nabla \psi\|_{L^\infty(0,t;L^\infty(\Omega))}\\
&&+~\int_{(0,t)\times \Omega}(|\nabla \sqrt{{c}^{\varepsilon,\tau}}|^2-|\nabla \sqrt{{c^{\varepsilon}}}|^2 )(u^{\varepsilon}\ast\rho^{\varepsilon})\cdot\nabla\psi\\
&:=&M_{61}+M_{62}.
\eeno
Similar to $M_2$, by  $ u^{\varepsilon,\tau}\to u^{\varepsilon}  $ in $ L^{3}(0,t;L^{3}(\Omega)) $ and the strong convergence of $\nabla \sqrt{{c}^{\varepsilon,\tau}}$, we get 
\beno
\int_{(0,t)\times\Omega} |\nabla \sqrt{{c}^{\varepsilon,\tau}}|^2 {(u^{\varepsilon,\tau}\ast\rho^\varepsilon)} \cdot \nabla \psi\rightarrow\int_{(0,t)\times\Omega} |\nabla \sqrt{{c^{\varepsilon}}}|^2 ({u^{\varepsilon}}\ast\rho^{\varepsilon}) \cdot \nabla \psi~~{\rm as}~~\tau\rightarrow0.
\eeno

{\bf \underline{For $M_{7}$:}}
Since $\psi \in C_c^\infty$, we have $\partial_t \psi + \Delta \psi \in C_c^\infty$ and
\beno
u^{\varepsilon,\tau} \to u^{\varepsilon} \quad {\rm in} \quad L^3(0,t;L^3(\Omega))~~{\rm as}~~\tau\rightarrow0.
\eeno
By the definition of strong convergence, we have

\beno
&&\int_{(0,t)\times\Omega} |u^{\varepsilon,\tau}|^2 \left(\partial_t \psi + \Delta \psi\right)-\int_{(0,t)\times\Omega} |u^{\varepsilon}|^2 \left(\partial_t \psi + \Delta \psi\right)\\&\leq &C(t,|\Omega|) ||u^{\varepsilon,\tau}-u^{\varepsilon}||^2_{L^3(0,t;L^3(\Omega))}\rightarrow0~~{\rm as}~~\tau\rightarrow0.
\eeno

{\bf \underline{For $M_{8}$:}}

\begin{equation*}
\begin{aligned}
&\int_{(0,t)\times \Omega}|u^{\varepsilon,\tau}|^{2}(u^{\varepsilon,\tau}\ast\rho^{\varepsilon})\cdot\nabla\phi-\int_{(0,t)\times \Omega}|u^{\varepsilon}|^{2}(u^{\varepsilon}\ast\rho^{\varepsilon})\cdot\nabla\phi\\=~&\int_{(0,t)\times \Omega}|u^{\varepsilon,\tau}|^{2}[(u^{\varepsilon,\tau}-u^{\varepsilon})\ast\rho^{\varepsilon}]\cdot\nabla\phi+\int_{(0,t)\times \Omega}(|u^{\varepsilon,\tau}|^{2}-|u^{\varepsilon}|^{2})(u^{\varepsilon}\ast\rho^{\varepsilon})\cdot\nabla\phi\\\leq~&\|u^{\varepsilon,\tau}\|_{L^3(0,t;L^3(\Omega))}^2
||(u^{\varepsilon,\tau}-u^{\varepsilon})\ast\rho^{\varepsilon}||_{L^3(0,t;L^3(\Omega_1))}\|\nabla \phi\|_{L^\infty(0,t;L^{\infty}(\Omega))}\\&+~\int_{(0,t)\times \Omega}(|u^{\varepsilon,\tau}|^{2}-|u^{\varepsilon}|^{2})(u^{\varepsilon}\ast\rho^{\varepsilon})\cdot\nabla\phi.
\end{aligned}
\end{equation*}
By the convergence property of convolution, $ u^{\varepsilon,\tau}\to u^{\varepsilon}  $ in $ L^{3}(0,t;L^{3}(\Omega)), $  and the strong convergence of $|u^{\varepsilon,\tau}|^{2}$, we can achieve
\begin{equation*}
\begin{aligned}
\int_{(0,t)\times \Omega}|u^{\varepsilon,\tau}|^{2}(u^{\varepsilon,\tau}\ast\rho^{\varepsilon})\cdot\nabla\phi\rightarrow\int_{(0,t)\times \Omega}|u^{\varepsilon}|^{2}(u^{\varepsilon}\ast\rho^{\varepsilon})\cdot\nabla\phi
\end{aligned}
\end{equation*}
as $\tau\to 0.$ 
%Thus $M_{10}$ convergence 0 as $\varepsilon\to 0.$

{\bf \underline{For $M_{9}$:}}
\beno
&&\int_{(0,t)\times\Omega} (P^{\varepsilon,\tau} - \bar{P}^{\varepsilon,\tau}) u^{\varepsilon,\tau} \cdot \nabla \psi-\int_{(0,t)\times\Omega} (P^{\varepsilon} - \bar{P}^{\varepsilon}) u^{\varepsilon} \cdot \nabla \psi\\
&=&\int_{(0,t)\times\Omega} (P^{\varepsilon,\tau} - \bar{P}^{\varepsilon,\tau}) (u^{\varepsilon,\tau}-u^{\varepsilon}) \cdot \nabla \psi+
\int_{(0,t)\times\Omega} [(P^{\varepsilon,\tau} - \bar{P}^{\varepsilon,\tau})-(P^{\varepsilon} - \bar{P}^{\varepsilon})] u^{\varepsilon} \cdot \nabla \psi\\
&\leq& ||P^{\varepsilon,\tau} - \bar{P^{\varepsilon,\tau}}||_{L^\frac32(0,t;L^\frac32(\Omega))}||u^{\varepsilon,\tau}-u^{\varepsilon}||_{L^3(0,t;L^3(\Omega))}\|\nabla \psi\|_{L^\infty(0,t;L^{\infty}(\Omega))}\\&&
+~
\int_{(0,t)\times\Omega} [(P^{\varepsilon,\tau} - \bar{P}^{\varepsilon,\tau})-(P^{\varepsilon} - \bar{P}^{\varepsilon})] u^{\varepsilon} \cdot \nabla \psi.
\eeno
Since
\beno
u^{\varepsilon,\tau} \to u^{\varepsilon} \quad {\rm in} \quad L^3(0,t;L^3(\Omega))~~{\rm as}~~\tau\rightarrow0.
\eeno
Then
\beno
||P^{\varepsilon,\tau} - \bar{P}^{\varepsilon,\tau}||_{L^\frac32(0,t;L^\frac32(\Omega))}||u^{\varepsilon,\tau}-u^{\varepsilon}||_{L^3(0,t;L^3(\Omega))}\|\nabla \psi\|_{L^\infty}\rightarrow0~~{\rm as}~~\tau\rightarrow0.
\eeno
Since $|u^{\varepsilon} \cdot \nabla \psi| \leq C |\nabla \psi| |u^{\varepsilon}| \in L^3(0,t;L^3(\Omega))$, by the the definition of weak convergence, we have
\beno
\int_{(0,t)\times\Omega} [(P^{\varepsilon,\tau} - \bar{P}^{\varepsilon,\tau})-(P^{\varepsilon} - \bar{P}^{\varepsilon})] u^{\varepsilon} \cdot \nabla \psi\rightarrow0~~{\rm as}~~\tau\rightarrow0.
\eeno
{\bf \underline{For $M_{10}$:}}

\begin{equation*}
\begin{aligned}
&\int_{(0,t)\times\Omega} (n^{\varepsilon,\tau}\nabla\phi)\ast\rho^\varepsilon \cdot u^{\varepsilon,\tau}  \psi-\int_{(0,t)\times\Omega} (n^{\varepsilon}\nabla\phi)\ast\rho^{\varepsilon} \cdot u^{\varepsilon}  \psi\\=~&\int_{(0,t)\times \Omega}[(n^{\varepsilon,\tau}\nabla\phi-n^{\varepsilon}\nabla\phi)\ast\rho^{\varepsilon}]\cdot u^{\varepsilon,\tau}\psi+\int_{(0,t)\times \Omega}((n^{\varepsilon}\nabla\phi)\ast\rho^{\varepsilon})\cdot(u^{\varepsilon,\tau}-u^{\varepsilon})\psi\\\leq~ &\int_{(0,t)\times \Omega}(n^{\varepsilon,\tau}-n^{\varepsilon})\nabla\phi\cdot(u^{\varepsilon,\tau}\psi)\ast\rho^{\varepsilon}\\&+~\|n^{\varepsilon}\|_{L^\frac32(0,t;L^\frac32(\Omega))}||u^{\varepsilon,\tau}-u^{\varepsilon}||_{L^3(0,t;L^3(\Omega))}||\psi||_{L^\infty(0,t;L^{\infty}(\Omega))}\|\nabla\phi\|_{L^\infty(0,t;L^{\infty}(\Omega))}.
\end{aligned}
\end{equation*}
Since
\beno
&&n^{\varepsilon,\tau}\rightarrow n^{\varepsilon} \quad {\rm in}~~ L^{\frac32}(0,t;L^{\frac32}(\Omega))~~{\rm as}~~\tau\rightarrow 0,\\
\eeno
and

$$|\nabla\phi \cdot (u^{\varepsilon,\tau}  \psi)\ast\rho^{\varepsilon}|  \in L^{3}(0,t;L^{3}(\Omega_1)).$$
Using strong convergence of $n^{\varepsilon,\tau}$, we get
\beno
\int_{(0,t)\times\Omega} (n^{\varepsilon,\tau}\nabla\phi-n^{\varepsilon}\nabla\phi)\ast\rho^\varepsilon \cdot u^{\varepsilon,\tau} \psi\rightarrow0~~{\rm as}~~\tau\rightarrow0.
\eeno
Combining with $
u^{\varepsilon,\tau} \to u^{\varepsilon}$ in $ L^3(0,t;L^3(\Omega))$ as $\tau\rightarrow0 $, we have
$$ M_{10}\to \int_{(0,t)\times\Omega} (n^{\varepsilon}\nabla\phi)\ast\rho^{\varepsilon} \cdot u^{\varepsilon}  \psi ~~{\rm as}~~\tau\rightarrow0 .  $$

To sum up, we have

\begin{equation}\label{energy inequality}
\begin{aligned}
&\int_{\Omega} (n^{\varepsilon} \ln (n^{\varepsilon}) \psi)(\cdot,t) + 4 \int_{(0,t)\times\Omega} |\nabla \sqrt{n^{\varepsilon}}|^2 \psi+\frac{2}{\Theta_0}  \int_{\Omega} (|\nabla \sqrt{c^{\varepsilon}}|^2 \psi)(\cdot,t)\\&+\frac{4}{3\Theta_0}\int_{(0,t)\times\Omega} |\Delta \sqrt{c^{\varepsilon}}|^2 \psi+\frac{2}{3\Theta_0}  \int_{(0,t)\times\Omega} (\sqrt{c^{\varepsilon}})^{-2} |\nabla \sqrt{c^{\varepsilon}}|^4 \psi\\
&+\frac{18}{\Theta_0}\|{c^{\varepsilon}_{0}}\|_{L^{\infty}}\int_{\Omega}(|u^{\varepsilon}|^2)(\cdot,t) \psi + \frac{18}{\Theta_0}\|{c^{\varepsilon}_{0}}\|_{L^{\infty}}\int_{(0,t)\times\Omega} |\nabla u^{\varepsilon}|^2 \psi
\\\leq&\int_{(0,t)\times\Omega} n^{\varepsilon} \ln n^{\varepsilon} (\partial_t \psi + \Delta \psi) + \int_{(0,t)\times\Omega} n^{\varepsilon} \ln n^{\varepsilon} (u^{\varepsilon}\ast\rho^{\varepsilon}) \cdot \nabla \psi \\&+\int_{(0,t)\times\Omega}n^{\varepsilon}\chi(c^{\varepsilon})\nabla c^{\varepsilon}\cdot \nabla \psi +\int_{(0,t)\times \Omega}n^{\varepsilon}\ln n^{\varepsilon} \chi(c^{\varepsilon})\nabla c^{\varepsilon}\cdot\nabla\psi\\&+ \frac{2}{\Theta_0}\int_{(0,t)\times\Omega} |\nabla \sqrt{c^{\varepsilon}}|^2 (\partial_t \psi + \Delta \psi) + \frac{2}{\Theta_0} \int_{(0,t)\times\Omega} |\nabla \sqrt{c^{\varepsilon}}|^2 (u^{\varepsilon}\ast\rho^{\varepsilon}) \cdot \nabla \psi\\&
+\frac{18}{\Theta_0}||c^{\varepsilon}_{0}||_{L^{\infty}}\int_{(0,t)\times\Omega} |u^{\varepsilon}|^2 \left(\partial_t \psi +\Delta \psi\right) + \frac{18\mu}{\Theta_0}||c^{\varepsilon}_{0}||_{L^{\infty}}\int_{(0,t)\times\Omega} |u^{\varepsilon}|^2 (u^{\varepsilon}\ast\rho^{\varepsilon})\cdot\nabla\psi
\\&+\frac{36}{\Theta_0}||c^{\varepsilon}_{0}||_{L^{\infty}}\int_{(0,t)\times\Omega} (P^{\varepsilon} - \bar{P}^{\varepsilon}) u^{\varepsilon} \cdot \nabla \psi - \frac{36}{\Theta_0}||c^{\varepsilon}_{0}||_{L^{\infty}}\int_{(0,t)\times\Omega} (n^{\varepsilon}\nabla\phi)\ast\rho^{\varepsilon} \cdot u^{\varepsilon} \psi.
\end{aligned}
\end{equation}

\section{proof of Theorem \ref{suitable weak solution}}

{\bf Proof.}
By convergence, we know that $ (u^{\varepsilon}, c^{\varepsilon}, n^{\varepsilon}) $ satisfies the equation as follows:
\begin{eqnarray}\label{eq:CNS1}
\left\{
\begin{array}{llll}
\displaystyle \partial_t n^{\varepsilon} + (u^{\varepsilon} \ast \rho^\varepsilon) \cdot \nabla n^{\varepsilon} - \Delta n^{\varepsilon} = - \nabla \cdot \left(n^{\varepsilon} \nabla c^{\varepsilon} \chi(c^{\varepsilon})\right),\\
\displaystyle \partial_t c^{\varepsilon} + u^{\varepsilon}  \cdot \nabla c^{\varepsilon} - \Delta c^{\varepsilon} = - n^{\varepsilon}\kappa(c^{\varepsilon}),\\
\displaystyle \partial_t u^{\varepsilon} + \mu(u^{\varepsilon} \ast \rho^\varepsilon) \cdot \nabla u^{\varepsilon} - \Delta u^{\varepsilon} + \nabla P^{\varepsilon} = - (n^{\varepsilon} \nabla \phi) \ast \rho^\varepsilon,\\
\displaystyle \nabla \cdot u^{\varepsilon} = 0,\\
\displaystyle n^{\varepsilon}(x,0) = n^{\varepsilon}_0(x), \quad c^{\varepsilon}(x,0) = c^{\varepsilon}_0(x), \quad u^{\varepsilon}(x,0) = u^{\varepsilon}_0(x).\\
\end{array}
\right.
\end{eqnarray}
%where $\rho^\varepsilon$ is a standard mollifier. Here, the initial value $n_0^{\varepsilon}= n_0 \ast \rho^\varepsilon$ satisfies
%\ben\label{ine:n L log L1}
%n_0^{\varepsilon} \to n_0 \quad {\rm in} \quad L \log L \quad {\rm as} \quad \varepsilon \to 0.
%\een
%$c_0^{\varepsilon}$ and $u_0^{\varepsilon}$ can be defined by $c^{\varepsilon}_0(x) = \left(\sqrt{c_0} \ast \rho^\varepsilon\right)^2$ and $u^{\varepsilon}_0(x) = u_0 \ast \rho^\varepsilon$.
By completely similar arguments as the convergence of $ (u^{\varepsilon,\tau}, c^{\varepsilon,\tau}, n^{\varepsilon,\tau}) $  in the subsection, we can get the convergence of $ (u^{\varepsilon}, c^{\varepsilon},n^{\varepsilon}) .$ We omitted it. Thereby, $(u, c, n)$ satisfies the following local energy inequality:
\begin{equation}\label{energy inequality 1}
\begin{aligned}
&\int_{\Omega} (n \ln (n) \psi)(\cdot,t) + 4 \int_{(0,t)\times\Omega} |\nabla \sqrt{n}|^2 \psi+\frac{2}{\Theta_0}  \int_{\Omega} (|\nabla \sqrt{c}|^2 \psi)(\cdot,t)\\&+\frac{4}{3\Theta_0}\int_{(0,t)\times\Omega} |\Delta \sqrt{c}|^2 \psi
+\frac{18}{\Theta_0}\|{c_{0}}\|_{L^{\infty}}\int_{\Omega}(|u|^2)(\cdot,t) \psi \\
&+ \frac{18}{\Theta_0}\|{c_{0}}\|_{L^{\infty}}\int_{(0,t)\times\Omega} |\nabla u|^2 \psi+\frac{2}{3\Theta_0}  \int_{(0,t)\times\Omega} (\sqrt{c})^{-2} |\nabla \sqrt{c}|^4 \psi
\\\leq~&\int_{(0,t)\times\Omega} n \ln (n) (\partial_t \psi + \Delta \psi) + \int_{(0,t)\times\Omega} n \ln (n) u \cdot \nabla \psi \\&+\int_{(0,t)\times\Omega}n\chi(c)\nabla c\cdot \nabla \psi +\int_{(0,t)\times \Omega}n\ln n \chi(c)\nabla c\cdot\nabla\psi\\&+ \frac{2}{\Theta_0}\int_{(0,t)\times\Omega} |\nabla \sqrt{c}|^2 (\partial_t \psi + \Delta \psi) + \frac{2}{\Theta_0} \int_{(0,t)\times\Omega} |\nabla \sqrt{c}|^2 u \cdot \nabla \psi\\&
+\frac{18}{\Theta_0}||c_{0}||_{L^{\infty}}\int_{(0,t)\times\Omega} |u|^2 \left(\partial_t \psi +\Delta \psi\right) + \frac{18\mu}{\Theta_0}||c_{0}||_{L^{\infty}}\int_{(0,t)\times\Omega} |u|^2 u\cdot\nabla\psi
\\&+\frac{36}{\Theta_0}||c_{0}||_{L^{\infty}}\int_{(0,t)\times\Omega} (P - \bar{P}) u \cdot \nabla \psi - \frac{36}{\Theta_0}||c_{0}||_{L^{\infty}}\int_{(0,t)\times\Omega} n\nabla\phi \cdot u \psi.
\end{aligned}
\end{equation}
Thus the third condition of Definition \ref{sws} about suitable weak solution is proved.

To sum up, the proof of Theorem {\ref{suitable weak solution}} is complete.
%It is obvious that $(n,c,u)$ satisfied the second condition. Next we show the first condition of suitable weak solution also holds.

%Lemma \ref{lem:2} shows that $ (n^\varepsilon,c^\varepsilon,u^\varepsilon) \in \left(L^\infty(0,T;H^2(\mathbb{R}^3)) \cap L^2(0,T;H^3(\mathbb{R}^3))\right)^3 $ is a unique global solution of system \ref{eq:CNS}. Noting that
%$$ u^{\varepsilon}\to u~~{\rm as}~~\varepsilon\rightarrow 0;$$ $$c^{\varepsilon}\to c~~{\rm as}~~\varepsilon\rightarrow 0;$$ and $$ n^{\varepsilon}\to n~~{\rm as}~~\varepsilon\rightarrow 0$$ in $ (n^\varepsilon,c^\varepsilon,u^\varepsilon) \in \left(L^\infty(0,T;H^2(\mathbb{R}^3)) \cap L^2(0,T;H^3(\mathbb{R}^3))\right)^3. $

%Moreover, the space $ (n^\varepsilon,c^\varepsilon,u^\varepsilon) \in \left(L^\infty(0,T;H^2(\mathbb{R}^3)) \cap L^2(0,T;H^3(\mathbb{R}^3))\right)^3 $
%is complete. Thus $(n,c,u) \in \left(L^\infty(0,T;H^2(\mathbb{R}^3)) \cap L^2(0,T;H^3(\mathbb{R}^3))\right)^3.$
%In conclusion, the theorem \ref{suitable weak solution} is being proved.
%\end{proof}

\section{Appendix}

\begin{lemma}[See Theorem A.1 in \cite{Vicol}]\label{lem:Banach fixed point}
	Let $ X $ be a Banach space, let $ B\subset X $ be a closed bounded subset, and let $ \Phi: B\to B $ be Lipschitz continuous in the norm topology with Lipschitz constant $ L<1; $ i.e., $ \|\Phi(x)-\Phi(y)\|\leq L\|x-y\| $ and $ L<1. $ Then $ \Phi $ has a unique fixed point in $ B. $
\end{lemma}
%\begin{lemma}[\cite{Galdi}, Lemma \uppercase\expandafter{\romannumeral2}.3.3]\label{GN}
%	Let $ \Omega_{0}\subset\mathbb{R}^{n}(n\geq 2) $ be a bounded smooth domain. Assume that $ 1\leq q,r\leq\infty, $ and $ j, m $ are arbitrary integers satisfying $ 0\leq j<m. $ If $ v\in W^{m,r}(\Omega_{0})\cap L^{q}(\Omega_{0}), $ then we have 
%	\begin{equation*}
%	\|D^{j}v\|_{L^{p}}\leq C\|v\|_{L^{q}}^{1-a}\|v\|_{W^{m,r}}^{a},
%	\end{equation*}
%	where
%	\begin{equation*}
%	-j+\frac{n}{p}=(1-a)\frac{n}{q}+a(-m+\frac{n}{r}),
%	\end{equation*}
%	and
%	\begin{equation*}
%	a\in\left\{
%	\begin{aligned}
%	&[\frac{j}{m},1), ~~if~~(m-j-\frac{n}{r})~~is~an~nonnegative~integer,\\& [\frac{j}{m},1],~~otherwise,
%	\end{aligned}
%	\right.
%	\end{equation*}
%	the constant $ C $ depends only on $ m,j,q,r,a $ and $\Omega_{0}$.
%\end{lemma}
%\begin{lemma}[See \cite{1972}]\label{BS}
%	Let $ \{T_{n}\} $ be a bounded linear operator from Banach space $ X $ to normed linear space $ Y. $ If for any $ x\in X, $  $ \lim_{n\to\infty}T_{n}x $ exists, then there %exists $ T\in B(X,Y) $ such that
%	\begin{equation*}
%		Tx=\lim_{n\to\infty}T_{n}x,\quad \forall x\in X,
%	\end{equation*}
%	and
%	$$ \|T\|\leq \underline{\lim\limits}_{n\to\infty}\|T_{n}\|. $$
%\end{lemma}
\begin{lemma}[See \cite{Vicol} Proposition A.16.]\label{e Delta f }
$\nu$ is a nonzero constant, $n$ is the dimension of space, the operator $e^{\nu t\Delta}$ satisfies the following a priori estimate for $ 1\leq q\leq p\leq\infty $ and $ 0\leq s\leq r<\infty $
	\begin{equation*}
		\|e^{\nu t\Delta}f\|_{\dot{W}^{r,p}(\mathbb{R}^{n})}\leq C(r,s)(\nu t)^{\frac{n}{2}(\frac{1}{p}-\frac{1}{q})-\frac{(r-s)}{2}}\|f\|_{\dot{W}^{s,q}(\mathbb{R}^{n})},
	\end{equation*}
	as soon as $ t>0. $

\begin{lemma}[See \cite{Tsai 2018} Lemma 3.7]\label{e Delta f }\label{Aubin-Lions lemma}
Let $X\subset\subset Y\subset Z$ be reflexive Banach spaces, if $0<T<\infty$, $1\leq p, 1<r$, $f\in L^p(0,T;X)$ satisfy
\beno
\int_0^T\|f(t)\|_{X}^pdt\leq C_1,\quad \int_0^T\|\frac{d}{dt}f(t)\|_{Z}^rdt\leq C_2.
\eeno
for some constants $C_1, C_2<\infty$, then ${f}_j, j\in\mathbb{N}$  is relatively compact in $L^p(0,T;Y)$.
\end{lemma}

\begin{lemma}[See \cite{LR} Theorem 2.4]\label{Sobolev inequalities}
Let $\alpha\in(0,d)$, $d$ is dimension, then $\frac1{(\sqrt{-\Delta})^\alpha}$ is bounded:
\begin{itemize}
\item[(1)]  from $L^1(\mathbb{R}^d)$ to $L^{\frac d{d-\alpha},\infty}(\mathbb{R}^d)$;
\item[(2)] from $L^{\frac d{\alpha},1}(\mathbb{R}^d)$ to $L^{\infty}(\mathbb{R}^d)$;
\item[(3)] from $L^{p,q}(\mathbb{R}^d)$ to $L^{p_1,q}(\mathbb{R}^d)$ for $1<p<\frac d{\alpha}$ and $\frac1{p_1}=\frac1p-\frac{\alpha}d$.
\end{itemize}

\end{lemma}
%	\label{e1,2,3,4}
%Let $(e^{t\Delta})_{t\geq0}$ be the Neumann heat semigroup in $\Omega$, $\Omega$ is smooth bounded domain, let $\lambda_1>0$ denote the first nonzero eigenvalue of $-\Delta$ in $\Omega$ under Neumann boundary conditions, then there have $C_1$, $C_2$, $C_3$, $C_4$ depending on $\Omega$ only,which has the following properties:
%
%\begin{itemize}
%\item if $1\leq q\leq p\leq \infty$, $f\in L^q(\Omega)$ satisfying $\int_{\Omega}f=0$, then
%\beno
%\|e^{t\Delta}f\|_{L^p(\Omega)}\leq C_1(1+t^{-\frac n2(\frac1q-\frac1p)})e^{-\lambda_1t}\|f\|_{L^q(\Omega)}\quad {\rm for~~~ all} \quad t>0; 
%\eeno 
%\item if $1\leq q\leq p\leq \infty$, then for each $f\in L^q(\Omega)$,
%\beno\label{e2}
%\|\nabla e^{t\Delta}f\|_{L^p(\Omega)}\leq C_2(1+t^{-\frac12-\frac n2(\frac1q-\frac1p)})e^{-\lambda_1t}\|f\|_{L^q(\Omega)}\quad {\rm for~~~ all} \quad t>0; 
%\eeno
%\item if $2\leq p<\infty$, for $f\in W^{1,p}(\Omega)$, then 
%\beno\label{e3}
%\|\nabla e^{t\Delta}f\|_{L^p(\Omega)}\leq C_3e^{-\lambda_1t}\|\nabla f\|_{L^p(\Omega)}\quad {\rm for~~~ all} \quad t>0; 
%\eeno
%\item if $1< q\leq p<\infty$, then for each $f\in C_0^{\infty}(\Omega)$, there holds
%\beno\label{e4}
%\|e^{t\Delta}\nabla \cdot f\|_{L^p(\Omega)}\leq C_4(1+t^{-\frac12-\frac n2(\frac1q-\frac1p)})e^{-\lambda_1t}\|f\|_{L^q(\Omega)}\quad {\rm for~~~ all} \quad t>0; 
%\eeno
%\end{itemize}
 \end{lemma}

 \noindent {\bf Acknowledgments.}
The author would like to thank Professors Sining Zheng, Zhaoyin Xiang and Wei Wang for some helpful communications. W. Wang was supported by NSFC under grant 12071054,  National Support Program for Young Top-Notch Talents and by Dalian High-level Talent Innovation Project (Grant 2020RD09).

 \end{document}